\documentclass[reqno]{amsart}
\usepackage{lmodern}

\usepackage{amsmath,amsfonts,latexsym,amssymb, mathtools}
\usepackage{colonequals}%for :=
\usepackage[inline]{enumitem}
\usepackage{comment,url,soul,microtype,color}
\usepackage{xcolor}
\usepackage{thm-restate}
\usepackage[backref=page, bookmarks=false, 
colorlinks=true,
linkcolor=blue,
citecolor=blue,
filecolor=blue,
menucolor=blue,
urlcolor=blue]{hyperref}
\usepackage[nameinlink]{cleveref}
\usepackage{graphicx}
\graphicspath{{./figures/figures_pdf/}}
\usepackage{bm}%for boldface math letters
\usepackage{subcaption}

\calclayout

\let\oldtocsection=\tocsection
\let\oldtocsubsection=\tocsubsection
\let\oldtocsubsubsection=\tocsubsubsection

\renewcommand{\tocsection}[2]{\hspace{0em}\oldtocsection{#1}{#2}}
\renewcommand{\tocsubsection}[2]{\hspace{3em}\oldtocsubsection{#1}{#2}}
\renewcommand{\tocsubsubsection}[2]{\hspace{6em}\oldtocsubsubsection{#1}{#2}}

\newtheorem{lemma}{Lemma}[section]
\newtheorem*{lemma*}{Lemma}
\newtheorem{proposition}[lemma]{Proposition}
\newtheorem*{proposition*}{Proposition}
\newtheorem{definition/proposition}[lemma]{Definition/Proposition}
\newtheorem{corollary}[lemma]{Corollary}
\newtheorem{theorem}[lemma]{Theorem}
\newtheorem*{theorem*}{Theorem}

\theoremstyle{definition}
\newtheorem{definition}[lemma]{Definition}
\newtheorem{remark}[lemma]{Remark}

\newtheorem{construction}[lemma]{Construction}
\newtheorem{definition/lemma}[lemma]{Definition/Lemma}

\newtheorem{question}{Question}[section]

\newtheorem{claim}[lemma]{Claim}

\newtheorem{theorema}{Theorem}

\crefname{question}{question}{questions}
\Crefname{question}{Question}{Questions}
\crefname{theorem}{theorem}{theorems}
\Crefname{theorem}{Theorem}{Theorems}
\crefname{theorema}{theorem}{theorems}
\Crefname{Theorema}{Theorem}{Theorems}
\crefname{proposition}{proposition}{propositions}
\Crefname{proposition}{Proposition}{Propositions}

  \makeatletter
  \newcommand{\linkdest}[1]{\Hy@raisedlink{\hypertarget{#1}{}}}
  \makeatother

\makeatletter
\renewcommand{\subsection}{%
  \@startsection{subsection}{2}{\z@}%
    {-1.5ex plus -.5ex minus -.2ex}%    vertical space before subsection heading
    {1ex plus .2ex}%             vertical space after subsection heading
    {\normalfont\bfseries}%
}
\makeatother

  \newcommand{\linkref}[2]{\hyperlink{#1}{#2}}
  
\newcommand{\RomanCaps}[1]
    {\MakeUppercase{\romannumeral #1}}

\NewDocumentCommand{\gvec}{me{_^}}{\makegvec{#1}{#2}{#3}}

\makeatletter
\@for\thm@name:=lemma,proposition,definition/proposition,corollary,theorem,%
  definition,remark,observation,construction,definition/lemma,example,%
  conjecture,hyp,claim\do{%
  \expandafter\def\csname theH\thm@name\endcsname{\thesection.\arabic{lemma}}}
\makeatother

\ExplSyntaxOn
\NewDocumentCommand{\makegvec}{mmm}
 {
  \vec{#1}
  \IfValueT{#2}{\sb{#2}}
  \IfValueT{#3}{
  \prop_get:NnNTF \g_gvec_kern_prop {#1} \l_tmpa_tl
    { \sp{\mspace{\l_tmpa_tl mu}#3} }
    { \sp{\mspace{4mu}#3} }
}
}
\prop_new:N \g_gvec_kern_prop
\prop_gset_from_keyval:Nn \g_gvec_kern_prop
 {
  \imath=5,
  f=3,
  n=1,
  o=3,
 }
\ExplSyntaxOff

\begin{document}
\title{Counting Elements and Conjugacy Classes of Fixed Genus in Hyperbolic Groups}
\author{Ritwik Chakraborty}
\address{School of Mathematics, Tata Institute of Fundamental Research, 1, Homi Bhabha Road,
Mumbai 400005, India}
\email{r4physics@gmail.com, ritwik@math.tifr.res.in}

\author{Inhyeok Choi}
\address{School of Mathematics, KIAS, Hoegi-ro 85, Dongdaemun-gu, Seoul 02455,
South Korea}
\email{inhyeokchoi48@gmail.com}

\author{Biswajit Nag}
\email{biswajit.nag.math@gmail.com}
\begin{abstract}
     Let $G$ be a non-elementary hyperbolic group equipped with a word metric. A group element in $G$ has orientable genus $n$ if it can be written as the product of $n$ commutators and no fewer. Similarly, a group element in $G$ has non-orientable genus $m$ if it can be written as the product of $m$ squares and no fewer. We obtain the asymptotic order of growth for both the number of group elements and the number of conjugacy classes in $G$, with length at most $L$, and fixed genus (orientable or non-orientable) $n$. For free groups, our methods give precise asymptotics for each of the counts.  
\end{abstract}

\maketitle
\tableofcontents
\section{Introduction}\label{sec: Intro}
Counting geodesics in negatively curved spaces is a classical problem with a rich history. The asymptotic distribution of closed geodesics in closed hyperbolic surfaces was derived by Huber \cite{Huber1959}, and in general closed negatively curved manifolds by Margulis \cite{Margulis2004}. Long closed geodesics can look very complex, so it is natural to find suitable measures of complexity and understand how closed geodesics of a given complexity are asymptotically distributed. For example, the asymptotic distribution of simple closed geodesics, and more generally, geodesics of a given type in hyperbolic surfaces was derived by Mirzakhani \cite{Mirzakhani2008}, \cite{Mirzakhani2016}. The following theorem counts closed geodesics in closed hyperbolic manifolds with prescribed homology: 
\begin{theorem}[\cite{katsuda1988homology},\cite{PhillipsSarnak1987}]\label{thm: homologicalcounting}
Let $M$ be a closed hyperbolic $n$-manifold. Let $\alpha\in H_1(M,\mathbb{Z})$ and let $\psi:\pi_1(M)\rightarrow H_1(M,\mathbb{Z})$ be the Hurewicz homomorphism. Then:
        \begin{align*}
            \left | \begin{Bmatrix}
            \gamma\subset M, \text{ primitive closed geodesic}\\l_M(\gamma)\leq L
            \text{ and } \psi([\gamma])=\alpha
            \end{Bmatrix}\right | \sim c(M)\frac{e^{(n-1)L}}{mL^{r/2+1}}
        \end{align*}
        where $r=b_1(M)$, $m=|\mathrm{tors} ~H_1(M)|$, and $c(M)$ is a constant depending on $M$.
\end{theorem}
Note that the asymptotics here do not depend on the choice of $\alpha$. Homologically trivial curves come with a natural order of complexity. Any homologically trivial curve can be represented as a product of commutators in $\pi_1(M)$, up to free homotopy. The \textit{commutator length} of $\gamma$, denoted by $\text{cl}(\gamma)$, is defined to be the minimum number of commutators required to represent $\gamma$. Let us also note that the integer $\text{cl}(\gamma)$ equals the genus of the smallest filling of $\gamma$, that is, the smallest connected oriented surface $S$ with one boundary component $\partial S$ for which there is a map $S\rightarrow \Sigma$, mapping $\partial S$ to $\gamma$ up to free homotopy. Thus the commutator length of $\gamma$ can also be thought of as the \textit{genus} of $\gamma$. This paper is directly inspired by the following result of Erlandsson-Souto:
\begin{theorem}[\cite{erlandsson2023counting}]\label{thm: erlandssonsouto}
    Let $\Sigma$ be a closed, connected, oriented hyperbolic surface. Then for all $n\geq 1$:
    \[
    \left\{  
    \begin{matrix}
        \gamma\subset\Sigma\textnormal{ closed geodesic with } \\
        \textnormal{cl}(\gamma)=n, \textnormal{ and } l_{\Sigma}(\gamma)\leq L
    \end{matrix}
    \right\}
    \sim \frac{2}{12^n\cdot n!\cdot (3n-2)!\cdot \textnormal{vol}(T^1\Sigma)^{2n-1}}L^{6n-4}e^{\frac{L}{2}}\]
\end{theorem}
It is natural to ask for a group-theoretic analogue of the above theorem. We formulate this precisely below. Let $G$ be a countable group and let $X$ be a proper geodesic metric space with a proper, non-elementary action of $G$ by isometries. Fixing some basepoint $x_0\in X$, the \textit{growth rate} or \textit{volume entropy} is defined as:
\[h=h(G,X)\coloneq \limsup_{L\to\infty}\frac{\log |\{g\in G : d(x_0,g\cdot x_0)\leq L\}| }{L}\]
We will assume from now that the growth rate is finite. The asymptotic order of growth of the number of primitive conjugacy classes in non-elementary hyperbolic groups was obtained by Coornaert and Knieper \cite{CoornaertKnieper2002}, \cite{CoornaertKnieper2004} and in the significantly generalized case of statistically convex-cocompact actions with contracting elements by Gekhtman-Yang \cite{GekhtmanYang2022}. Precise asymptotics for the number of conjugacy classes have been obtained for geometrically finite groups of isometries of $\textnormal{CAT}(-1)$ spaces (with finite BMS measure) by Roblin \cite{Roblin2003}, mapping class groups $\textnormal{Mod}(\Sigma_g)$ of surfaces of genus at least $2$ acting on the Teichm\"{u}ller space $\mathcal{T}(\Sigma_g)$ \cite{EskinMirzakhani2011}.

The \textit{orientable genus} of a group element $g$ in $[G,G]$, $\text{genus}^+(g)$, is defined to be $n$ if it can be written as the product of $n$ commutators and no fewer. Similarly, the \textit{non-orientable genus} of a group element $g$ in $G^2$, $\text{genus}^-(g)$, is defined to be $m$ if it can be written as the product of $m$ squares and no fewer. The (orientable and non-orientable) genus of a group element is a conjugacy invariant and therefore well-defined for conjugacy classes. Moreover, for any path-connected space $Y$ with $\pi_1(Y)=G$, a conjugacy class $[g]$ has orientable genus at most $n$ if and only if there is map $f:S_{n,1}\rightarrow Y$ such that $f(\partial S_{n,1})$ is freely homotopic to a loop representing $[g]$. Similarly, a conjugacy class $[g]$ has non-orientable genus at most $k$ if and only if there is map $f:N_{k,1}\rightarrow Y$ such that $f(\partial N_{k,1})$ is freely homotopic to a loop representing $[g]$, where $N_{k,1}$ denotes the non-orientable surface of genus $k$ and connected boundary.

Let us note that non-orientable genus of a conjugacy class or the square length of a closed geodesic is also a measure of complexity for homologically trivial curves, except that one is using $\mathbb{Z}/2\mathbb{Z}$-coefficients.

We define:
\begin{align*}
    B^{\pm}_n(L)&\coloneq \{g\in G\mid \text{genus}^{\pm}(g)=n, d(x_0,g\cdot x_0)\leq L\}\\
    \mathcal{C}^{\pm}_n(L)&\coloneq \{[g]\in G\mid \text{genus}^{\pm}([g])=n, ~l([g])\leq L\}
\end{align*}
where $l([g])$ is the minimum displacement of $x_0$ caused by an element in $[g]$. Inspired by \Cref{thm: erlandssonsouto}, we have the following problems, asking for asymptotic counting with fixed genus:

\begin{question}\label[question]{ques: element counting by genus}
    Given $n\geq 1$, are there integers $\alpha_{\pm}$ and constants $c_{\pm}>0$, such that:
    \begin{align*}
        |\mathcal{C}^{\pm}_n(L)|\sim c_{\pm}L^{\alpha_{\pm}}e^{\frac{hL}{2}}\quad?
    \end{align*}
\end{question}

\begin{question}\label[question]{ques: conj class counting by genus}
    Given $n\geq 1$, are there integers $\beta_{\pm}$ and constants $c'_{\pm}>0$, such that:
    \begin{align*}
        |B^{\pm}_n(L)|\sim c'_{\pm}L^{\beta_{\pm}}e^{\frac{hL}{2}}\quad?
    \end{align*}
\end{question}

\subsection{Main results}\label{subsec: Main results} In this paper, we investigate \Cref{ques: element counting by genus,ques: conj class counting by genus} for the case when $G$ is a non-elementary hyperbolic group and $X$ is its Cayley graph, and for the rest of this section $B^{\pm}_n(L)$ and $\mathcal{C}^{\pm}_n(L)$ will be understood accordingly, that is, elements and conjugacy classes are counted with respect to a word metric. 

\begin{theorema}\label{intro thm: orientable element count}
For any non-elementary hyperbolic group $G$, there are constants $c_1, c_2,L_0>0$ depending only on $\delta, n, G$ such that for all $n\geq 1$:

\[
    \text{for all } L\geq L_0,\quad  c_1 L^{6n-2}e^{\frac{hL}{2}} \leq |B^{+}_n(L)| \leq  c_2 L^{6n-2} e^{\frac{hL}{2}}
\]

where $h$ is the volume entropy of $G$.
\end{theorema}

\begin{theorema}\label{intro thm: orientable conjugacy count}
    For any non-elementary hyperbolic group $G$, there are constants $c_1', c_2',L'_0 > 0$ depending only on $\delta,n,G$ such that for all $n\geq 1$:
    \[
        \text{for all } L\geq L'_0,\quad c_1'L^{6n-4}e^{\frac{hL}{2}} \leq |\mathcal{C}^{+}_n(L)| \leq  c_2' L^{6n-4} e^{\frac{hL}{2}}
    \]
    where $h$ is the volume entropy of $G$.
\end{theorema}

\begin{theorema}\label{intro thm: non-orientable element count}

For any non-elementary hyperbolic group $G$, there are constants $c_3, c_4, L_0>0$ depending only on $\delta, n, G$ such that for all $n\geq 1$:

\[
    \text{for all } L\geq L_0,\quad  c_3 L^{3n-2}e^{\frac{hL}{2}} \leq |B^{-}_n(L)| \leq  c_4 L^{3n-2} e^{\frac{hL}{2}}
\]

where $h$ is the volume entropy of $G$.
\end{theorema}

\begin{theorema}\label{intro thm: non-orientable conjugacy count}
    For any non-elementary hyperbolic group $G$, there are constants $c_3', c_4',L'_0 > 0$ depending only on $\delta,n,G$ such that for all $n\geq 1$:
    \[
        \text{for all } L\geq L'_0,\quad c_3'L^{3n-4}e^{\frac{hL}{2}} \leq |\mathcal{C}^{-}_n(L)| \leq  c_4' L^{3n-4} e^{\frac{hL}{2}}
    \]
    where $h$ is the volume entropy of $G$.
\end{theorema}

    While \Crefrange{intro thm: orientable element count}{intro thm: non-orientable conjugacy count} are stated for hyperbolic groups equipped with word metrics, our methods also go through for non-elementary geometric actions of groups on $\delta$-hyperbolic geodesic metric spaces by isometries. Therefore, in particular, \Cref{intro thm: orientable conjugacy count} (resp. \Cref{intro thm: non-orientable conjugacy count}) gives the asymptotic order of growth of the number of closed geodesics of fixed commutator length (resp. square length) $n$ and length at most $L$ in closed hyperbolic manifolds of any dimension, as $L$ grows to infinity.

    When \(G\) is a finitely generated non-abelian free group, our methods strengthen to yield precise asymptotic estimates of each of the counts involved above. In order to state these results, we need the notions of \emph{orientable} and \emph{non-orientable} \emph{pre-Wicks} and \emph{Wicks forms} of a given \emph{genus}. Wicks introduced genus $1$ Wicks forms to prescribe an algorithm to decide whether a given element in a free group is a commutator (cf. \cite{wicks1962commutators}). Culler \cite{culler1981using} used higher genus Wicks forms to classify homotopy classes of maps from surfaces with connected boundary into roses and to parametrize solution sets of quadratic equations over free groups (also see \cite{ComerfordEdmunds1994}). We refer the reader to the work of Vdovina and collaborators for results on enumeration and applications of Wicks forms (\cite{vdovina1995constructing}, \cite{vdovina1994products}, \cite{Vdovina2008},\cite{StoimenowVdovina2005}). Other incarnations of  Wicks forms include $1$-vertex triangulations of surfaces, ribbon or fat graphs, and unicellular maps.
    
    Let $\mathcal{X}$ be an alphabet, whose elements we call \emph{unknowns}. A \emph{pre-Wicks form} $W$ is a reduced word in $F(\mathcal{X})$, the free group generated by $\mathcal{X}$ which is: (a) \emph{quadratic} - every unknown $X$ appearing in $W$  appears exactly twice, either as $X$ and/or $X^{-1}$, (b) \emph{irredundant} - if $X$ and $Y$ are two unknowns and $\varepsilon_X,\varepsilon_Y\in\{\pm 1\}$, and $(X^{\varepsilon_X}Y^{\varepsilon_Y})$ is a subword of $W$, then $(X^{\varepsilon_X}Y^{\varepsilon_Y})^{\pm 1}$ does not appear as another subword of $W$. A \emph{Wicks form} $W$ is a reduced cyclic word which is quadratic and cyclically irredundant. For further definitions, we refer the reader to \Cref{sec: Prelims}.

\begin{theorema}\label{intro thm: element counts free group}
    Let $n$ and $r$ be positive integers, with $r\geq 2$. Let $F$ be a free group of rank $r$ equipped with a word metric associated with the standard symmetric generating set. Set $h_r=\log (2r-1)$. Then, for all $n\geq 1$:
    \begin{align*}
        \lim_{\substack{L\to\infty\\ L\textnormal{ even}}}\frac{| B^{+}_n(L)|}{L^{6n-2}e^{\frac{h_rL}{2}}}&= \frac{\tau^+_n(r) \mathcal{N}_{\text{pre}}^{+}(n)}{2^{6n-2}(6n-2)!(2r)^{6n-1}(2r-1)^{6n-2}(2r-2)}\\
        \lim_{\substack{L\to\infty\\ L\textnormal{ even}}}\frac{| B^{-}_n(L)|}{L^{3n-2}e^{\frac{h_rL}{2}}}&= \frac{\tau^-_n(r) \mathcal{N}_{\text{pre}}^{-}(n)}{2^{3n-2}(3n-2)!(2r)^{3n-1}(2r-1)^{3n-2}(2r-2)}
    \end{align*}
    where: 
    \begin{align*}
        \tau^+_n&=\tau^+_n(r)\coloneq 2r(2r(2r-1)(2r-2))^{4n-1}\\ \tau^-_n&=\tau^-_n(r)\coloneq 2r(2r(2r-1)(2r-2))^{2n-1} 
    \end{align*}
    and, $\mathcal{N}_{\text{pre}}^{+}(n)$ denotes the number of maximal pre-Wicks forms of orientable genus $n$, while $\mathcal{N}_{\text{pre}}^{-}(n)$ denotes the number of maximal pre-Wicks forms of non-orientable genus $n$.
\end{theorema}

Let $n$ be a positive integer. Define the \emph{mass} of orientable (and resp. non-orientable) Wicks forms of genus $n$ as:
    \[\mathfrak{M}_n^+\coloneq \sum_{\substack{\text{genus}^+(W)=n\\\text{maximal}\\\text{orientable}}}\frac{1}{|\text{Aut}(W)|},\quad \mathfrak{M}_n^-\coloneq \sum_{\substack{\text{genus}^-(W)=n\\\text{maximal}\\\text{non-orientable}}}\frac{1}{|\text{Aut}(W)|}\]

\begin{theorema}\label{intro thm: conjugacy class counts free group}
    Let $n$ and $r$ be positive integers, with $r\geq 2$. Let $F_r$ be the free group of rank $r$ equipped with the word metric associated with the standard symmetric generating set. Set $h_r=\log (2r-1)$. Then:
    \begin{align*}
        \text{for all }n\geq 1,\lim_{\substack{L\to\infty\\ L\textnormal{ even}}}\frac{|\mathcal{C}^{+}_n(L)|}{L^{6n-4}e^{\frac{h_rL}{2}}}&= \frac{\kappa^+_n(r)\mathfrak{M}^+_n}{2^{6n-4}(6n-4)!(2r)^{6n-3}(2r-1)^{6n-4}(2r-2)}\\
        \text{for all }n> 1,\lim_{\substack{L\to\infty\\ L\textnormal{ even}}}\frac{| \mathcal{C}^{-}_n(L)|}{L^{3n-4}e^{\frac{h_rL}{2}}}&= \frac{\kappa^-_n(r)\mathfrak{M}^-_n}{2^{3n-4}(3n-4)!(2r)^{3n-3}(2r-1)^{3n-4}(2r-2)}\\
        \lim_{\substack{L\to\infty\\ L\textnormal{ even}}}\frac{|\mathcal{C}^-_1(L)|}{L^{-1}e^{\frac{h_rL}{2}}}&=\frac{2r-1}{r-1}
    \end{align*}
    where:
    \begin{align*}
        \kappa^+_n&=\kappa^+_n(r)\coloneq (2r(2r-1)(2r-2))^{4n-2}\\ \kappa^-_n&=\kappa^-_n(r)\coloneq (2r(2r-1)(2r-2))^{2n-2} 
    \end{align*}
\end{theorema}

\begin{remark}
    Precise asymptotics for the number of conjugacy classes of fixed orientable genus $1$ and length $L$, that is, $\partial\mathcal{C}^+_1(L)$ in free groups were obtained by \cite{park2019conjugacy}. The limits in \Cref{intro thm: element counts free group,intro thm: conjugacy class counts free group} are taken over even values of $L$ because elements and conjugacy classes of a given genus have even length.
    
    $\mathfrak{M}^+_n$ and $\mathfrak{M}^-_n$ have closed expressions. By \cite[Theorem 1.1]{BACHER200213}, for all $n\geq 1$:
    \[\mathfrak{M}_n^+=\frac{2}{12^n}\frac{(6n-5)!}{n!(3n-3)!}\]
    On the other hand, using \cite[Corollary 8 and 9]{BernardiChapuy}, we have for all $j,k\geq 1$:
    \begin{align*}
        \mathfrak{M}^{-}_{2j}=\frac{c_j(3j-1)!\binom{6j-2}{3j-1}}{6^{j+1}(2j-1)} \quad \mathfrak{M}^{-}_{2k+1}=\frac{2^{6k-2}(3k)!}{3^{k+1}k\cdot k!}
        \quad \text{where }c_j=3\cdot 2^{3j-2}\frac{j!}{(2j)!}\sum_{l=0}^{j-1}\binom{2l}{l}16^{-l}
    \end{align*}
    Note that $\bigcup_{L}\bigsqcup_{n\geq 1}B^+_n(L)=[G,G]-{1_{G}}$. The growth rate of the commutator subgroup $[G,G]$ of a non-elementary hyperbolic group $G$ is the same as that of $G$ (cf. \cite{CoulonDalBoSambusetti2018}). For free groups, precise asymptotics for the number of reduced words and  the number of cyclically reduced words of length $L$ in the commutator subgroup were obtained by Sharp \cite{Sharp2001} and Rivin \cite{Rivin2010} respectively. For non-elementary hyperbolic groups, a recent work of Cantrell-Sharp \cite{CantrellSharp2021} obtains growth rate asymptotics for the number of group elements of length $L$ in a normal subgroup $N$ of $G$ with $G/N\cong \mathbb{Z}^{d}$.
\end{remark}

\subsection{Interpreting asymptotics using the moduli space of metric fat graphs}\label{subsec: Interpreting asymptotics}

We will briefly sketch a proof of the orientable genus part of \Cref{intro thm: conjugacy class counts free group} with the aim to explain the structure of the asymptotic expression. With that in mind, we fix positive integers $r\geq 2$, $n\geq 1$, and consider the free group of rank $r$, $F_r$, equipped with the standard word metric.

We begin with the observation that Wicks forms can also be thought as \emph{fat/ribbon graphs}. For a detailed account on fat graphs, the reader is referred to \cite{norbury2010counting}, \cite{andersen2020kontsevich}. A fat graph of genus $n$ is a graph $\Gamma$ with vertices of degree at least $3$ together with an embedding into an oriented surface $\Sigma_n$ of genus $n$ such that $\Sigma_n-\Gamma$ is a disk. An automorphism between two fat graphs $i_1:\Gamma_1\rightarrow\Sigma_n$ and $i_2:\Gamma_2\rightarrow\Sigma_n$ is a graph automorphism $\varphi: \Gamma_1\rightarrow \Gamma_2$ together with an orientation-preserving homeomorphism $h$ of $\Sigma_n$ such that $h\circ i_1=i_2\circ \varphi$. It is worth noting here that a fat graph is therefore a graph equipped with a cyclic ordering on the half edges incident at every vertex, and a fat graph automorphism is a graph automorphism that preserves the cyclic orderings. A metric on a fat graph $\Gamma$ is an assignment of non-negative real numbers $l(e)$ to its edges $e$. Equipped with a metric, the length of its boundary is defined as $l(\partial \Gamma)\coloneq 2\sum_{e\in E(\Gamma)}l(e)$. For every fat graph $\Gamma$, consider the polytope $V_\Gamma(L)\coloneq \{\vec{l}\in \mathbb{R}_{+}^{E(\Gamma)}\mid l(\partial \Gamma)=L \}$. The \textit{combinatorial moduli space of metric fat graphs} of genus $n$ with boundary length $L$, $\mathcal{M}^{\text{comb}}_n(L)$ is an orbicell-complex defined as follows:
\[\mathcal{M}^{\text{comb}}_n(L)\coloneq\left(\bigsqcup_{\Gamma}\frac{V_{\Gamma}(L)}{\text{Aut}(\Gamma)}\right)/\sim \]
where the cells glue together naturally via edge degeneration, that is, when an edge length contracts, the edge contracts to give a metric fat graph with fewer edges. There is a homeomorphism between $\mathcal{M}^{\text{comb}}_n=\bigcup \mathcal{M}^{\text{comb}}_n (L)$ and $\mathcal{M}_{g,1}\times \mathbb{R}_{+}$, the decorated moduli space of genus $n$ Riemann surfaces with one marked point, using Jenkins-Strebel quadratic differentials due to Mumford, Thurston, Harer, Penner (cf. \cite{Harer1985modulispace}). The set of integer points in $\mathcal{M}^{\text{comb}}_n(L)$ will be denoted by $\mathcal{M}^{\text{comb},\mathbb{Z}}_n(L)$. 

Given a Wicks form $W$ of genus $n$, consider a polygon $P_W$ with edges labeled by the unknowns appearing in $W$ such that $\partial P_W$, when treated as a labeled edge path, corresponds to $W$. Identifying sides labeled by the same unknown gives a compact orientable surface of genus $n$ and the image of $\partial P_W$ in the surface is a fat graph, $\Gamma_W$. This establishes a one-to-one correspondence between Wicks forms of genus $n$ and fat graphs of genus $n$. The top dimensional cells in $\mathcal{M}^{\text{comb}}_n(L)$ have dimension $6n-4$ and they correspond to maximal Wicks forms.

The following is a theorem due to Goldstein-Turner, Culler, and Comerford-Edmunds (also see \Cref{lem: canonical forms for conjugacy classes}).
\vspace{-0.5cm}
\begin{theorem*}[\cite{GoldsteinTurner},\cite{culler1981using},\cite{ComerfordEdmunds1994}]
    Let $[g]$ be a non-trivial conjugacy class in $F_r$ having orientable genus $n$. Then there exists an orientable  Wicks form $W$ of genus $n$ and a non-canceling substitution $\mathfrak{g}$ for $W$ such that $[g]=[\mathfrak{g}(W)]$, where $\mathfrak{g}(W)$ denotes the group element obtained by evaluating the substitution $\mathfrak{g}$ on $W$.
\end{theorem*}
\noindent  Therefore, given a conjugacy class $[g]$ of orientable genus $n$ and length $L$, we get a triple $(W,\vec{l},\mathfrak{g})$ where $W$ is a Wicks form of genus $n$, $\mathfrak{g}$ is a non-canceling substitution for $W$ such that $[\mathfrak{g}(W)]=[g]$, and $\vec{l}$ is an assignment of a positive integer $l_X$ to each unknown $X$ appearing in $W$, recording the length of the group element $g_X$ assigned by the substitution.  For a Wicks form $W$, let $\Lambda(W,L)$ be the set of triples $(W,\vec{l},\mathfrak{g})$ where $\vec{l}$ is a length type for $W$ satisfying $2\sum_{X}l_X=L$, and $\mathfrak{g}$ is a non-canceling substitution for $W$ such that $l(g_X)=l_X$ for every unknown $X$ appearing in $W$. One can also think of $\Lambda(W,L)$ as a decorated $V_{\Gamma_W}(L)$, that is, metric ribbon graphs decorated with group elements comprising a non-canceling substitution. Let $\Lambda_n(L)$ be the union of $\Lambda(W,L)$ for all Wicks forms $W$ of genus $n$. Consider the map $\Psi: \Lambda_n(L)\rightarrow \partial\mathcal{C}_{\leq n}^{+}(L)$ defined by $\Psi(W,\vec{l},\mathfrak{g})=[\mathfrak{g}(W)]$, where $\partial\mathcal{C}_{\leq n}^{+}(L)$ denotes the set of conjugacy classes of length $L$ and orientable genus at most $n$. Note that for any given $W$, $\Psi|_{\Lambda(W,L)}$ is $\text{Aut}(W)$-invariant. Thus we have an induced map $\overline{\Psi}:\overline{\Lambda}_n(L)\rightarrow \partial\mathcal{C}_{\leq n}^{+}(L)$, where $\overline{\Lambda}_n(L)\coloneq \bigsqcup_W \Lambda(W,L)/\text{Aut}(W) $. Most of the work in proving \Cref{intro thm: conjugacy class counts free group} is in showing that this map $\overline{\Psi}$ is essentially bijective onto $\partial \mathcal{C}^+_n(L)$ (cf. \Cref{lem: key lemma for conj class lower bounds}) using the results established in \Cref{subsec: combinatorialreduction} and \ref{subsec: Combreduction revisited}. Proceeding ahead, it suffices to make asymptotic estimates for $|\overline{\Lambda}_n(L)|$. Consider the projection map that forgets the substitution, from $\overline{\Lambda}_n(L)$ onto $\mathcal{M}^{\text{comb}}_n(L)$. As $L\to\infty$, most of the contribution to $|\overline{\Lambda}_n(L)|$ comes from the fibers above points in the top dimensional cells of  $\mathcal{M}^{\text{comb}}_n(L)$, and the size of such a fiber over a generically chosen point $[(\Gamma_W,\vec{l})]$ equals the number of non-canceling substitutions for $W$ of length type $\vec{l}$, which asymptotically equals:
\[\frac{\kappa^+_n(r)e^{\frac{h_rL}{2}}}{(2r)^{6n-3}(2r-1)^{6n-3}}\]
Therefore, as $L\to\infty$, $|\partial\mathcal{C}_{n}^{+}(L)|$ asymptotically equals:
\[|\mathcal{M}^{\text{comb},\mathbb{Z}}_n(L)|\times\frac{\kappa^+_n(r)e^{\frac{h_rL}{2}}}{(2r)^{6n-3}(2r-1)^{6n-3}}\]
and $|\mathcal{M}^{\text{comb},\mathbb{Z}}_n(L)|$ asymptotically equals: $\sum_{\Gamma}\frac{1}{|\text{Aut}(\Gamma)|}\frac{L^{6n-4}}{2^{6n-4}(6n-4)!}$, where the sum runs over all trivalent fat graphs of genus $n$, or equivalently, maximal Wicks forms of genus $n$. In the non-orientable genus case, one can similarly interpret the conjugacy class count asymptotics using metric M\"{o}bius graphs and their combinatorial moduli space (cf. \cite{GouldenHarerJackson2001}, \cite{chidambaram2026refinedlatticepointcounting}).

\subsection{Overview of proofs}

\subsubsection{Organization of the paper}

The proofs of \Crefrange{intro thm: orientable element count}{intro thm: conjugacy class counts free group} are broadly divided into upper bounds and lower bounds which can be read independently. All upper bound proofs only require the results proved in \nameref{sec: wicks form moves}. Element count lower bounds in the case of free groups require the results proved in \nameref{subsec: combinatorialreduction}. For the case of general hyperbolic groups, \Cref{subsec: relating combinatorial reduction to degeneracy in substitutions - hyperbolic} is needed in addition. Conjugacy count lower bounds need the results proved in \nameref{subsec: Combreduction revisited}.
\begin{itemize}
    \item \Cref{intro thm: orientable element count} follows from \Cref{thm: element count upper bound orientable hyperbolic} and \Cref{thm: lowerboundhyp}.
    \item \Cref{intro thm: non-orientable element count} follows from \Cref{thm: element count upper bound non-orientable hyperbolic} and \Cref{thm: lowerboundhyp nonorientable}.
    \item \Cref{intro thm: orientable conjugacy count} follows from \Cref{thm: conjugacy count upper bound orientable hyperbolic} and \Cref{thm: lowerboundhyp orientable conj class}.
    \item \Cref{intro thm: non-orientable conjugacy count} follows from \Cref{thm: conjugacy count upper bound non-orientable hyperbolic} and \Cref{thm: lowerboundhyp nonorientable conj class}.
\end{itemize}

\noindent Let $\partial B^{\pm}_n(L)$ (resp. $\partial \mathcal{C}^{\pm}_n(L)$) denote the set of group elements $g$ (resp. conjugacy classes $[g]$) with $\text{genus}^{\pm}(g)=n$ (resp. $\text{genus}^{\pm}([g])=n$) and length equal to $L$. For free groups, we first obtain precise asymptotics for the `shell' counts: $|\partial B^{\pm}_n(L)|$ and $|\partial \mathcal{C}^{\pm}_n(L)|$. The precise asymptotics for $|B^{\pm}_n(L)|$ and $| \mathcal{C}^{\pm}_n(L)|$ follow from the shell asymptotics and the fact that : $\sum_{\substack{l\leq L\\ l\text{ even}}}l^{\alpha}e^{\frac{h_rl}{2}}\sim \frac{2r-1}{2r-2}L^{\alpha}e^{\frac{h_rL}{2}}$. Keeping this in mind, we have:

\begin{itemize}
    \item \Cref{intro thm: element counts free group} follows from \Cref{thm: free grp shell upper bound asymptote,thm: free grp shell upper bound asymptote non-orientable}, and \Cref{thm: free grp shell lower bound asymptote,thm: free grp shell lower bound asymptote non-orientable}.
    \item \Cref{intro thm: conjugacy class counts free group} follows from \Cref{thm: free grp shell conj class upper bound asymptote,thm: free grp shell conj class upper bound asymptote non-orientable}, and \Cref{thm: free grp shell conj class lower bound asymptote,thm: free grp shell conj class lower bound asymptote non-or}.
\end{itemize}
All the constants appearing in this article have local scope, that is, they are meant to have scope within a section, unless otherwise specified.
\subsubsection{Upper bounds} 
Our approach relies on constructing quasi-geodesic representatives for genus $n$ group elements. Let $G$ be a hyperbolic group equipped with a word metric. Given any group element $g$ of orientable genus $n$, we construct a \emph{tracking pair}, which is a tuple $(W_t,\mathfrak{g}_t)$, where $W_t$ is an orientable  pre-Wicks form of genus $n$, and $\mathfrak{g}_t$ is a substitution evaluating out to $g$ such that the associated piece-wise geodesic $\gamma(\mathfrak{g}_t)$ \textit{$(c,D,l_0)$-tracks} any geodesic segment $\gamma_g$ joining the identity to $g$ (cf. \Cref{lem: longedgeparallelfellow-travelling}) for some uniform constants $c,D,l_0>0$ depending only on $n$ and the choice of word metric. To elaborate, let $W_t$ be a word of length $2m$ in the unknowns $X_1,\dots,X_m$. Then for each $1\leq j\leq m$, whenever $X_j$ is \emph{long}, i.e., $l(\gamma_{X_j})\geq D$, there is a decomposition of $\gamma_{X_j}$ into \emph{initial}, \emph{middle}, and \emph{final} pieces, that is, $\gamma_{X_j}=\gamma_{X_j}^{(i)}\gamma_{X_j}^{(m)}\gamma_{X_j}^{(f)}$. Moreover, there are subsegments: $\gamma'_{X_j}$, $\gamma'_{X^{-1}_j}$, of $\gamma_g$ such that: 
        \begin{enumerate}[label=(\roman*)]
            \item $l(\gamma_{X_j}^{(i)}),l(\gamma_{X_j}^{(f)})\leq c$
            \item $\gamma'_{X_j}$ and $\gamma_{X_j}^{(m)}$ $l_0$-fellow travel, $\gamma'_{X_{\overline{j}}}$ and $\gamma_{X_{\overline{j}}}^{(m)}$ $l_0$-fellow travel. 
            \item The collection of interiors of $\gamma'_{X_j}$ and $\gamma'_{X_{\overline{j}}}$ for $1\leq j\leq m$, is mutually disjoint.
        \end{enumerate}
This tracking property implies that $W_t$ admits a decomposition as: $W_t=S_1L_1S_2\dots S_{k}L_{k}S_{k+1}$, where $L_1,\dots,L_k$ are \emph{long} letters of $W_t$ and $S_1,\dots,S_{k+1}$ are products of consecutive \emph{short} letters (allowed to be the empty word) of $W_t$. This decomposition together with the fellow-traveling property and the fact that the length of an orientable (resp. non-orientable) pre-Wicks form of genus $n$ is at most $12n-2$ (resp. $6n-2$), enables us to obtain upper bounds for $B^{\pm}_n(L)$ and $\mathcal{C}^{\pm}_n(L)$.

Tracking pairs are constructed using three kinds of \textit{moves}, which take as input a tuple $(W,\mathfrak{g})$ and whenever applicable, they output another tuple $(W',\mathfrak{g}')$ such that:
\begin{enumerate}[label=(\alph*)]
    \item $\mathfrak{g}(W)=\mathfrak{g}'(W')$
    \item The total length of the substitution strictly decreases: $\sum_{X'}l(g_{X'})<\sum_X l(g_X)$.
    \item If $W$ is an orientable pre-Wicks form with $\text{genus}^+(W)=\text{genus}^+(\mathfrak{g}(W))=n$, then $W'$ is also an orientable pre-Wicks form with $\text{genus}^+(W')=n$.
    \item If $W$ is a non-orientable pre-Wicks form with $\text{genus}^-(W)=\text{genus}^-(\mathfrak{g}(W))=n$, then $W'$ is either a non-orientable pre-Wicks form with $\text{genus}^-(W')=n$ or an orientable pre-Wicks form with $\text{genus}^+(W')=(n-1)/2$.
\end{enumerate}

The construction of the tracking pair $(W_t,\mathfrak{g}_t)$ can be thought of as an extension of Comerford-Edmunds \cite[Theorem 1]{ComerfordEdmunds1994} results to the case of general hyperbolic groups. A much more detailed study of the existence of Wicks form representatives in hyperbolic groups has been undertaken by Duncan-Fulthorp  \cite{duncan2015orientablenonorientablegenusn} and much of the observations made in \Cref{subsec: properties of the tracking pair} have been inspired by \cite[Section 4]{duncan2015orientablenonorientablegenusn}.

\subsubsection{Lower Bounds}
 The proofs for lower bounds use the method of combinatorial reduction. For the sake of clarity and simplicity, we describe our methods in the simplest possible case: establishing lower bounds for the element count in free groups in the orientable genus case (\Cref{intro thm: element counts free group}). For the rest of this subsection, we fix some $r\geq 2$ and let $S$ denote a standard symmetric generating set for the free group of rank $r$.

To begin with, we note the following fact: For any non-trivial group element $g\in [F_r,F_r]$, there is an orientable pre-Wicks form $W$ and a non-canceling substitution $\mathfrak{g}$ such that $\mathfrak{g}(W)=g$ (cf. \cite[Theorem 1]{ComerfordEdmunds1994}, and \Cref{lem: canonical forms for group elements}). In order to use this observation to count group elements of a fixed orientable genus, we need to understand how often a group element arises from two different tuples $(W,\mathfrak{g})$ and $(W',\mathfrak{g}')$. Given a pre-Wicks form $W$, a length type $\vec{l}$ for $W$ is an assignment of positive integers to the unknowns appearing in $W$, and $||\vec{l}||\coloneq 2\sum_{X}l_X$ where the sum is over the unknowns $X$ appearing in $W$. We consider the set:
\[G_W(\vec{l})\coloneq \{\mathfrak{g}(W) \mid \mathfrak{g}\textnormal{ non-canceling substitution of length type }\vec{l}\}\]
We prove that given any two maximal pre-Wicks forms $W$ and $W'$, and two generically chosen respective length types $\vec{l}$ and $\vec{l'}$, with $(W,\vec{l})\neq (W',\vec{l'})$ and $||\vec{l}||=||\vec{l'}||=L$, the proportion of elements that arise from a substitution of type $(W',\vec{l'})$ among those that arise from a substitution of type $(W,\vec{l})$ goes to zero as $L\to\infty$. More precisely: 
\begin{lemma*}[\Cref{lem: non-cancellingcount2}]
    Let $W$ be a maximal orientable pre-Wicks form of genus $n$. Then, for all $L\geq 1$, there exists a collection of length types $\mathcal{L}=\mathcal{L}(W,n,r,L)$ and there are constants $L_0,C\geq 1$ depending only on $r$ and $n$ satisfying, for all $L\geq L_0$:
    \begin{enumerate}[font=\normalfont]
        \item For all $\vec{l}\in \mathcal{L}(W,n,r,L)$, $||\vec{l}||=L$ and \[\lim_{\substack{L\to\infty\\L\textnormal{ even}}}\frac{|\mathcal{L}(W,n,r,L)|}{L^{6n-2}}=\frac{1}{2^{6n-2}(6n-2)!}\]
    \end{enumerate}
    Let $W'$ be another orientable pre-Wicks form of genus $n$, and let $\vec{l}\in \mathcal{L}(W,n,r,L)$, $\vec{l}'\in \mathcal{L}(W',n,r,L)$.
    \begin{enumerate}[font=\normalfont, resume]
        \item If $(W,\vec{l})\neq(W',\vec{l}')$, then:
        \[|G_W(\vec{l})\cap G_{W'}(\vec{l}')|\leq C\frac{e^{\frac{h_rL}{2}}}{L^{6n}}\]
    \end{enumerate}
\end{lemma*}

In order to obtain \Cref{lem: non-cancellingcount2}, we define the notion of \emph{combinatorial reduction}, which we introduce briefly below.

Given two non-canceling substitutions \(\mathfrak{g}_1\) and \(\mathfrak{g}_2\) of length types \(\vec{l}\) and \(\vec{l}'\), the condition \(g = \mathfrak{g}_1(W) = \mathfrak{g}_2(W')\) forces $\mathfrak{g}_1$ and $\mathfrak{g}_2$ to have \textit{degeneracies}: there are subwords in the substitution that either repeat several times or completely determine other parts of the substitution. Establishing statistics of degeneracies in non-canceling substitutions for group elements that arise in different ways is key to obtaining the estimate in \Cref{lem: non-cancellingcount2} stated above. While it is possible to fully investigate such degeneracies in specific low genus pre-Wicks forms, the main question here is: 

\begin{quote}
    Is there a conceptual way of understanding how degeneracies arise from the condition $\mathfrak{g}_1(W) = \mathfrak{g}_2(W')$ for arbitrary pre-Wicks forms $W$ and $W'$?
\end{quote}  

Given a pre-Wicks form $W=X_1\cdots X_{2m}$ and a length type $\vec{l}$, there is a canonical partition induced on the set of \textit{positions} $\{1,\dots, L\}$, where $L=||\vec{l}||$: namely, a position $i$ is \textit{contained in} a letter $X_j$ if: $s_{X_{j}}\coloneq \sum_{k<j}l_{X_k}+1\leq i\leq t_{X_j}\coloneq \sum_{k\leq j}l_{X_k}$. Moreover, for any such position $i$, there is a \textit{matching position} defined by $\underline{i}\coloneq t_{X_{\overline{j}}}+s_{X_{j}}-i$, where $X_{\overline{j}}$ is the unique letter in $W$ \textit{paired} with $X_j$, that is, $X_{\overline{j}}=X_j^{-1}$. Given a triple $(W,\mathfrak{g},\vec{l})$ where $\mathfrak{g}$ is a non-canceling substitution and $||\vec{l}||=L$, we identify the set of positions $\{1,\dots, L\}$ with the letters (elements of $S$, the generating set) labeling the geodesic joining the identity element to $\mathfrak{g}(W)$. 

Given two triples $(W,\mathfrak{g},\vec{l})$ and $(W',\mathfrak{g}',\vec{l}')$ with $\mathfrak{g}(W)=\mathfrak{g}'(W')$ and $||\vec{l}||=||\vec{l'}||=L$, any position $i \in \{1, \dots, L\}$ defines matching positions \(\underline{i}\) and \(\underline{i}'\) relative to the length types \(\vec{l}\) and \(\vec{l}'\). The matching positions are defined such that if $x_i, x_{\underline{i}}$ and $x_{\underline{i}'}$ are the letters (elements of $S$, the generating set) of $g$ at the positions $i$, $\underline{i}$, and $\underline{i}'$ respectively, then $x_i = x_{\underline{i}}^{-1} = x_{\underline{i}'}^{-1}$. We therefore have the following equation/constraint for \(\mathfrak{g}\):
\[ 
    x_{\underline{i}} = x_{\underline{i}'} 
\]
When \(\underline{i}\neq \underline{i}'\), the above constraint reduces the degree of freedom we have in choosing \(\mathfrak{g}\) by one letter. Thus, if there are $N$ positions $i_1,\dots, i_N$ so that $\underline{i_1}\neq \underline{i_1}'$, $\dots$, $\underline{i_N}\neq \underline{i_N}'$, we say that there is a \emph{reduction} (in degrees of freedom) \emph{by $N$}. The notion of reduction only depends on the data $(W,\vec{l})$ and $(W',\vec{l}')$. \Cref{lem: non-cancellingcount2} follows from the following key proposition (cf. \Cref{thm: the main reduction theorem,cor: the main reduction theorem for free groups}):
\begin{proposition*}
    Let $N\geq 1$ be a positive integer. Let $W$ and $W'$ be maximal pre-Wicks forms of lengths $2m$ and $2m'$ respectively. Let $\vec{l}$ and $ \vec{l}'$ be length types for $W$ and $W'$ respectively, satisfying:
    \par\medskip
    \noindent
    \begin{enumerate*}[font=\normalfont, label=(\alph*), itemjoin=\qquad]
        \item $||\vec{l}||=||\vec{l}'||$
        \item 
         For all $1\leq i\leq 2m$, $1\leq j\leq 2m'$, $l_i,l'_j>6N+4$.
    \end{enumerate*}
    \par\medskip
    \noindent Then, either:
    \begin{enumerate*}[font=\normalfont, label=(\arabic*), itemjoin=\qquad]
        \item There is reduction by $N$, or
        \item $W=W'$, therefore \(m = m'\), and moreover $\vec{l} = \vec{l'}$.
    \end{enumerate*}
    
\end{proposition*}

The absence of reduction turns out to be a highly restrictive condition, as we describe in \Cref{subsec: combinatorialreduction}. As a first step, we show that the absence of reduction by $N$ between $(W,\vec{l})$ and $(W',\vec{l'})$ forces $W=W'$ and this does not require $W$ or $W'$ to be maximal. Moreover, $(W,\vec{l})$, $(W,\vec{l'})$ are forced to satisfy a \textit{shift equation}. The shift of a letter $X_j$ in $W$ is the difference between their ending positions: $\textnormal{shift}(j)\coloneq t_{X_j}-t'_{X_j}$ and  the absence of reduction implies: $\textnormal{shift}(j)=-\textnormal{shift}(\overline{j}-1)$, for paired letters $X_j$ and $X_{\overline{j}}$. To see how shift equations are useful, suppose $W$ has the following form:
\[\cdots XY\cdots Y^{-1}Z \cdots Z^{-1}X^{-1}\cdots\]
Then, simply applying the shift equation thrice gives:
\[\textnormal{shift}(X)=-\textnormal{shift}(Z^{-1})=(-1)^2\textnormal{shift}(Y^{-1})=(-1)^3\textnormal{shift}(X)\]
That is, $\textnormal{shift}(X)=0$. Therefore, if we further assume that $W$ is a maximal Wicks form, so that the associated fat graph is trivalent, the shift equations imply that all shifts must vanish and $\vec{l}=\vec{l'}$. So far we had assumed that the pre-Wicks forms involved are all orientable. In the non-orientable case, the analysis is more involved: a few more shift equations are obtained, and applying the shift equations proceeds according to the type of $1$-neighborhood of a vertex in the associated fat graph (cf. \Cref{fig: circuit in 1-neighbourhood}).

In the setup of general non-elementary hyperbolic groups, we do not have any non-canceling substitutions, so we prepare a large class of substitutions with controlled cancellation called \emph{chain substitutions}. Let $W$ be any pre-Wicks form of orientable/non-orientable genus $n$ and length $2m$, in the unknowns $X_1,\dots,X_m$. There exist constants $K,R_0>0$ depending only on $G$ and the choice of word metric, and a Schottky set $\mathcal{S}$, such that for any substitution $\mathfrak{g}$ with $l(g_{X_i})\geq R\geq R_0$ , 
     there exist $l_i,r_i\in\mathcal{S}$ and so that the substitution $\mathfrak{g}'$, defined by: $\text{for all }1\leq i\leq m, \quad g'_{X_i}=l_ig_{i}r_i$,
      is a $(K+2\delta,R)$-chain substitution (cf. \Cref{def: chainsubstitution}, \Cref{def: alignedchain}), which in particular, implies that the associated piecewise geodesic is a local quasigeodesic, so it fellow-travels with any geodesic segment $\gamma_g$ joining the identity to $g=\mathfrak{g}(W)$. Therefore one can encode the algebraic data of the Wicks form and its length type data into $\gamma_g$ (cf. \Cref{construction: geodesicencoding}). Upgrading the above strategy to the case of general hyperbolic groups requires generalizing the notion of combinatorial reduction appropriately, and this is the generality at which it is defined in \Cref{subsec: combinatorialreduction}. With this generalized combinatorial reduction, the free group lower bound strategy can be implemented using the encoding.

% {\color{gray}(Is it appropriate to include a table of notation here or at some other appropriate places?)}

\subsection{Acknowledgments}
Initial discussions on this work began during the program Geometry in Groups (code: ICTS/GiG2024/07) at the International Centre for Theoretical Sciences (ICTS), Bengaluru, and the authors are grateful to ICTS for its hospitality, to the organizers, and especially the participants in the second week of the program. RC would like to thank the Korea Institute for Advanced Study (KIAS) for its hospitality and support during his visit, where part of this work was carried out. RC is also grateful to Radhika Gupta for taking the time to read earlier drafts of this paper. BN thanks Indranil Bhattacharyya for several helpful discussions on the upper bound arguments. RC and BN were supported by the Department of Atomic Energy, Government of India (project 12-R\&D-TFR-5.01-0500). IC was supported by the Mid-Career Researcher Program (RS-2023-00278510) through the National Research Foundation funded by the government of Korea, and by the KIAS individual grant (MG091901) at KIAS. 

\subsection{AI Disclosure}

The entire paper was first written without aid of AI. After completing the writeup, it was fed into a publicly available LLM for proofreading, grammar check and typo corrections. The authors then fixed some typos based on the LLM's output.

\section{Preliminaries}\label{sec: Prelims}
Let $G$ be a group and $S$ be a finite generating set. We denote by $d=d_S$ the word metric on $G$. The length of a group element $g\in G$ is denoted by $l(g):=d(1,g)$.  

\begin{definition}
    For each conjugacy class $c$ in $G$, we define its \emph{(algebraic) length} as:
    \begin{equation*}
        l(c):=\min\{l(g)| g\in c\}
    \end{equation*}
\end{definition}

    Consider the Cayley graph $\Gamma=\Gamma(G,S)$ of $G$ with respect to the generating set $S$. Then $G$ acts on $\Gamma$ from the left, by isometries. It's easy to see that for a conjugacy class $c$ in $G$, $l(c)$ records the minimum displacement of a vertex in the Cayley graph by an element in the conjugacy class $c$, i.e. for any $g\in c$:
    \begin{equation*}
        l(c)=\min\{d(x,g\cdot x)| x\in V(\Gamma)\}
    \end{equation*}

\begin{definition}
\label{def: genus}
An element $g\in [G,G]$ is said to have \emph{orientable genus }$n$ if it can be expressed as a product of $n$ (and no fewer) commutators. That is, $n$ is the least positive integer such that the following equation has a solution:
    \begin{equation*}
        [a_1,b_1]\cdots [a_n,b_n]=g
    \end{equation*}
where $a_1,b_1,\dots,a_n,b_n$ are unknowns taking values in $G\backslash \{1_{G}\}$.

$g\in G^2$ is said to have \emph{non-orientable genus }$m$ if it can be expressed as a product of $m$ (and no fewer) squares. That is, $m$ is the least positive integer such that the following equation has a solution:
    \begin{equation*}
        e_1^2\cdots e_{m}^2=g
    \end{equation*}
where $e_1,\dots, e_{m}$ are unknowns taking values in $G\backslash \{1_{G}\}$. The orientable (resp. non-orientable) genus of a group element $g$ will be denoted as $\textnormal{genus}^+(g)$ (resp. $\textnormal{genus}^-(g)$). Moreover, it is invariant under conjugation and is therefore well-defined for conjugacy classes. 
\end{definition}

In any group $G$, $[G,G]$ is a subgroup of $G^2$ because any non-trivial commutator can be written as a product of $3$ squares: $[x,y]=x^2(x^{-1}y)^2{(y^{-1}})^2$. A classical result of Lyndon-Newman \cite{LyndonNewman1973} says that in the free group generated by $x$ and $y$, $\textnormal{genus}^-([x,y])=3$ (this can also be seen purely topologically, using the fact that there are no odd degree maps from the Klein bottle to the torus). More generally, we have the following fact, which is essentially a consequence of the fact that $\mathbb{RP}^2\#T^2$ is homeomorphic to $\mathbb{RP}^2\#\mathbb{RP}^2\#\mathbb{RP}^2$ (also see \cite{Edmunds1975} for an explicit sequence of Nielsen automorphisms transforming $x^2[y,z]$ to $x^2y^2z^2$ in $F(X,Y,Z)$).

\begin{proposition}\label{prop: square length upper bound}
For any $g\in [G,G]$, $\textnormal{genus}^{-}(g)\leq 2\textnormal{genus}^{+}(g)+1$.    
\end{proposition}

\begin{definition}
    \label{def: word} Let $\mathcal{X}$ be a given set of unknowns. By a word with unknowns in $\mathcal{X}$, we will mean a reduced word $W$ in the free group $F(\mathcal{X})$ generated by $\mathcal{X}$. A cyclic word refers to a cyclically reduced word, and all its cyclic permutations. 
    
\end{definition}

\begin{definition}\label{def: quadratic sign orientable non-orientable}
    A word $W$ in a set of unknowns $\mathcal{X}$ is said to be \emph{quadratic} if every unknown $X\in X$ which appears in $W$, appears exactly twice, as either $X$ or $X^{-1}$ each time.
    
    Given a quadratic word $W$, we define the \emph{sign} of an unknown $X$ appearing in $W$ to be the product of the exponents in its two occurrences in $W$. We denote this by $\text{sign}(X)$.
    
    A quadratic word $W$ is said to be \emph{non-orientable} if there is at least one unknown $X$ appearing in $W$ with $\text{sign}(X)=1$. Otherwise, it is said to be \emph{orientable}.
    
    Let $W=X_1\dots X_n$ be a quadratic word. For every letter $X_i$ in $W$, there exists a unique letter $X_{\bar{i}}$ such that $X_{\bar{i}}=X_i^{-1}$ or $X_i$, according to the sign of $X_i$. We call $X_{\bar{i}}$, \emph{the letter paired with} $X_i$. Note that for any such pair, we have $\text{sign}(X_i) = \text{sign}(X_{\overline{i}})$.
\end{definition}

\begin{definition}\label{def: irredundant/irreducible}
    A quadratic word $W$ in a set of unknowns $\mathcal{X}$ is said to \emph{contain no redundancies} or \emph{irredundant} if \(W\) does not contain two subwords of the form \(X_1^{\epsilon_1} X_2^{\epsilon_2}\) and \((X_1^{\epsilon_1} X_2^{\epsilon_2})^{\pm 1}\), where \(X_1, X_2\) are unknowns and \(\epsilon_1, \epsilon_2 \in \{\pm 1\}\). Equivalently, no substitution of the form \(X_1^{\epsilon_1} X_2^{\epsilon_2}\mapsto X\) is possible in \(W\). We'll say that $W$ is \emph{cyclically irredundant} if every cyclic permutation of $W$ is irredundant.
\end{definition}

\begin{definition}\label{def: Wicks form}
    A word $W$ in a set of unknowns $\mathcal{X}$ is said to be a \emph{Wicks form} if it is quadratic, cyclically reduced, and cyclically irredundant. A Wicks form can either be orientable or non-orientable. 
\end{definition}

\begin{definition}\label{def: wicks form genus}
    Let \(W = X_1 \cdots X_{2m}\) be a Wicks form of length \(2m\) with unknowns in \(\mathcal{X}\), and with \(X_i \in \mathcal{X} \cup \mathcal{X}^{-1}\). Consider a \(2m\)-sided polygon $P_W$ whose edges are oriented and labelled by the letters on \(W\) in the order which they appear. The \emph{topological genus} of the Wicks form \(W\) is the genus of the closed surface $\Sigma_{W}$ obtained by gluing the edges of the above \(2m\)-gon according to their labelling. The image of $\partial P_{W}$ in $\Sigma_{W}$ is a oriented, labelled graph $\Gamma_{W}$ embedded in $\Sigma_{W}$.
    Note that $\Sigma_W$ is orientable iff $W$ is orientable.
    Now suppose $W$ is an orientable Wicks form. The \emph{algebraic genus} of \(W\) is the smallest number \(n\) such that \(W\) can be written as a product of \(n\) commutators in \(F(\mathcal{X})\). If $W$ is a non-orientable Wicks form, its algebraic genus is defined to be the smallest number $m$ such that $W$ can be written as a product of $m$ squares in \(F(\mathcal{X})\).

\end{definition}

The next observation follows promptly from the definition.
\begin{lemma}\cite[1.1]{culler1981using}\label{lem: topoologicalandalgebraicgenus}
      Let $W$ be a Wicks form. 
      \leavevmode
      \begin{enumerate}[font=\normalfont]
          \item If $W$ is orientable, then the topological genus of $W$  coincides with its algebraic genus.
          \item If $W$ is non-orientable, then the topological genus of $W$  is one-half of its algebraic genus.
      \end{enumerate}
\end{lemma}

In the sequel, when we say genus of an orientable (resp. non-orientable) Wicks form $W$, we refer to its algebraic genus and denote it by $\text{genus}^{+}(W)$ (resp. $\text{genus}^{-}(W)$ ).

\begin{definition}\label{def: equivalenceofWicksforms}
    We define two notions of equivalence for Wicks forms, one for linear words, and another for cyclic words. The former will be used in \Cref{sec: Element count: Upper bounds} and \Cref{sec: Element count: Lower bounds} for element counts, while the latter will be used in \Cref{sec: Conjugacy class count: Upper bounds} and \Cref{sec: Conjugacy class count: Lower bounds}, for conjugacy class counts.
    
    Let $\mathcal{X}$ and $\mathcal{X}'$ be two sets of unknowns. Let $W$ and $W'$ be (linear) Wicks forms with unknowns in $\mathcal{X}$ and $\mathcal{X}'$ respectively. $W'$ is said to be \emph{linearly equivalent} to $W$ if there is a bijection $\varphi:\mathcal{X}\cup\mathcal{X}^{-1}\rightarrow \mathcal{X}'\cup(\mathcal{X}')^{-1}$ such that for every unknown $X\in\mathcal{X}$, $\varphi(X^{-1})=\varphi(X)^{-1}$, and $W'=\varphi(W)$.
    
    $W$ and $W'$ are said to be \emph{cyclically equivalent} if there is a cyclic permutation $\sigma$ of $W'$ such that $W$ and $\sigma(W')$ are linearly equivalent. Thus, if $W$ and $W'$ are considered as cyclic words, they are said to be equivalent if we can find linear words representing them that are equivalent.
\end{definition}

\begin{lemma}\label{lem: wicks form length}
    Let \(W\) be a Wicks form of genus \(n\). If $W$ is orientable, then $4n\leq l(W)\leq 12n-6$. If $W$ is non-orientable and $n>1$, then $2n\leq l(W)\leq 6n-6$. Moreover, the only non-orientable Wicks form of genus $1$ is $X^2$.
\end{lemma}

\begin{proof}
Consider a Wicks form of genus $n$ and length $2m$. Suppose $W$ is orientable. We look at the oriented, labelled $2m$-sided polygon $P_{W}$ associated with $W$. Identifying pairs of sides in the $2m$-gon according to the labelling gives us a closed orientable surface $\Sigma_W$ of genus $n$. The image of $\partial P_{W}$ in $\Sigma_{W}$ is a graph $\Gamma_{W}$ embedded in $\Sigma_{W}$.  Thus we have:
\begin{align*}
    \chi(\Sigma_{W})=2-2n=v-e+1
\end{align*}
Now $v\geq1$ gives us the lower bound $l(W)  = 2e\geq 4n$ with equality iff $\Gamma_{W}$ has exactly one vertex with degree $4n$. Because $W$ is cyclically reduced and irredundant, $\Gamma_W$ does not have any vertices of degree $1$ and $2$ respectively. Recall the handshake lemma, which says that $2e=\sum_{i\in V(\Gamma_W)}d(i)$ where $d(i)$ denotes the degree of vertex $i$. This gives us the upper bound $l(W)\leq 6(2n-1)$ with equality iff $\Gamma_W$ has $2(2n-1)$ vertices of degree $3$. We can argue completely analogously when $W$ is non-orientable, having genus greater than $1$, with the key difference appearing in the expression for the Euler characteristic of $\Sigma_W$, because the genus of $\Sigma_W$ in this case is $\frac{n}{2}$ (\Cref{lem: topoologicalandalgebraicgenus}). In the non-orientable genus $1$ case, the $\chi(\Sigma_W)=1$, so $v=e$ and the degree sum formula implies that all vertices of $\Gamma_W$ must have degree exactly $2$. This implies that $v=e=2$, and $W=X^2$ for some unknown $X$.
\end{proof}

Due to \Cref{lem: wicks form length}, it makes sense to call a Wicks form maximal and minimal whenever it achieves the upper and lower bounds that appear in the lemma. It is clear that minimal Wicks forms exist for every genus $n\geq 1$, one such minimal Wicks form being the product of $n$ commutators in the orientable case and the product of $n$ squares in the non-orientable case. It is not too hard to see that maximal Wicks forms also exist for every genus. We record the latter as a theorem below:
\begin{theorem}\cite[Section 2]{comerford1991powers}\label{thm: existence of maximal Wicks form}
    For every $n\geq 1$, there are maximal Wicks forms of genus $n$.
\end{theorem}

In fact, any Wicks form can be obtained from a maximal Wicks form by substituting some of its unknowns by the identity element (cf. \cite{ol1989homomorphism},\cite{comerford1991powers}). In \cite{vdovina1995constructing} and \cite{vdovina1994products}, Vdovina describes an inductive procedure for building all possible maximal Wicks forms of a given genus from those of lower genus. She also derives asymptotic counts for the number of maximal Wicks forms of a given genus. 

\begin{theorem}\cite[Theorem 2]{vdovina1995constructing}\cite[Theorem 1.1]{vdovina1994products}
    \label{thm: countWicksforms}
    Let $\mathcal{N}^+(n)$ denote the number of non-equivalent orientable maximal Wicks forms of genus $n\geq 1$. Let $\mathcal{N}^-(n)$ denote the number of non-equivalent non-orientable maximal Wicks forms of genus $n\geq 1$. Then:
    \begin{align*}
        \lim_{n\to\infty}\frac{\log \mathcal{N}^+(n)}{2n \log n} = 1\quad\text{and}\quad  \lim_{n\to\infty}\frac{\log \mathcal{N}^-(n)}{n \log n} = 1 
    \end{align*}
\end{theorem}

Let us note that in \Cref{thm: countWicksforms}, $\mathcal{N}^{\pm}(n)$ counts Wicks forms as linear words. While in \cite{vdovina1994products} and \cite{vdovina1995constructing}, Wicks forms are counted as cyclic words. Nonetheless, the asymptotics match because each cyclic equivalence class of Wicks forms contains at most $12n-6$ linearly non-equivalent linear Wicks forms in the orientable case, and at most $6n-6$ many of them in the non-orientable case.

\begin{definition}\label{def: pre-Wicks form}
    A word $W$ in a set of unknowns $X$ is said to be a \emph{pre-Wicks form} if it is quadratic, reduced, and irredundant. A pre-Wicks form is either orientable or non-orientable. 
\end{definition}

% {\color{gray}(It might be convenient to define \(G\)-valued substitutions somewhere around here and say that given a quadratic word \(W\) and a substitution \(\mathfrak{g}\), there is a natural way of turning it into a tuple \((W', \mathfrak{g}')\) such that \(W'\) is a pre-Wicks form with \(\mathfrak{g}'(W') = \mathfrak{g}(W)\). This is, of course, done by reducing and removing redundancies from \(W\))}

\begin{proposition}\label{prop: what does a preWicks form look like}
    Let $W$ be a pre-Wicks form. Then $W$ has one among the following forms necessarily:
    \begin{enumerate}[font=\normalfont]
        \item\label{item: preWicks form 1 - Wicks} $W$ is a Wicks form.
        \item\label{item: preWicks form 2 - tail reduction} $W=DUD^{-1}$ where $U$ is a Wicks form.
        \item\label{item: preWicks form 3 - tail redundancy} $W=AU_1A^{-1}B^{-1}U_2B$ and the word $U_1X^{-1}U_2X$ is a Wicks form. Or, $W=AU_1BAU_2B$ and the word $U_1XU_2X$ is a Wicks form. In both cases, \(U_1\), \(U_2\) are subwords of \(W\) with \(A, B\) and $X$ being an unknown that doesn't appear in $U_1$ or $U_2$. 
        \item\label{item: preWicks form 4 - tail redundancy nested in tail reduction} $W=DAU_1A^{-1}B^{-1}U_2BD^{-1}$ and the word $U_1X^{-1}U_2X$ is a Wicks form. Or, $W=DAU_1BAU_2BD^{-1}$ and the word $U_1XU_2X$ is a Wicks form. As before, \(U_1, U_2\) are subwords of \(W\) with \(A, B, D\) and $X$ being unknowns that don't appear in $U_1$ or $U_2$.
    \end{enumerate}
\end{proposition}

\begin{proof}
    By definition, a Wicks form is cyclically reduced and cyclically irredundant. A pre-Wicks form is only reduced and irredundant, so it could have a redundancy and/or a reduction at its ends. Accordingly $(2)$ is when it only has a tail reduction. $(3)$ is when it only has a tail redundancy, and $(4)$ describes the case when there is a tail redundancy nested with a tail reduction.
\end{proof}

By \Cref{prop: what does a preWicks form look like}, the notions of genus and equivalence of Wicks forms, together with the length bounds (\Cref{lem: wicks form length}) and growth rates (\Cref{thm: countWicksforms}) extend to pre-Wicks forms readily. 

\begin{definition}\label{def: genus of a preWicks form}
    The \emph{genus} of a pre-Wicks form $W$ is defined to be the genus of the Wicks form obtained from $W$ by removing any tail reduction and/or redundancy, as described in \Cref{prop: what does a preWicks form look like}.
\end{definition}

\begin{definition}\label{def: equivalence of preWicks forms}

    Let $\mathcal{X}$ and $\mathcal{X}'$ be two sets of unknowns. Let $W$ and $W'$ be pre-Wicks forms with unknowns in $\mathcal{X}$ and $\mathcal{X}'$ respectively. $W'$ is said to be \emph{linearly equivalent} to $W$ if there is a bijection $\varphi:\mathcal{X}\cup\mathcal{X}^{-1}\rightarrow \mathcal{X}'\cup(\mathcal{X}')^{-1}$ such that for every unknown $X\in\mathcal{X}$, $\varphi(X^{-1})=\varphi(X)^{-1}$, and $W'=\varphi(W)$.
    
\end{definition}

\begin{lemma}\label{lem: length of a preWicks form}
    Let $W$ be a pre-Wicks form of genus $n$. If $W$ is orientable, then $4n\leq l(W)\leq 12n-2$. If $W$ is non-orientable, then $2n\leq l(W)\leq 6n-2$ 
\end{lemma}

\begin{corollary}\label{cor: existence of maximal pre-Wicks form}
    For every $n\geq 1$, there are maximal pre-Wicks forms of genus $n$.
\end{corollary}

\begin{corollary}\label{thm: countpreWicksforms}
    Let $\mathcal{N}_{\text{pre}}^+(n)$ denote the number of non-equivalent maximal orientable pre-Wicks forms of genus $n\geq 1$. Let $\mathcal{N}_{\text{pre}}^-(n)$ denote the number of non-equivalent maximal non-orientable pre-Wicks forms of genus $n\geq 1$. Then:
    \begin{align*}
        \lim_{n\to\infty}\frac{\log \mathcal{N}_{\text{pre}}^+(n)}{2n \log n} = 1\quad\text{and}\quad  \lim_{n\to\infty}\frac{\log \mathcal{N}_{\text{pre}}^-(n)}{n \log n} = 1 
    \end{align*}
\end{corollary}

\subsection{Some background on hyperbolicity}\label{subsec: hyperbolicprelims}

For the rest of this section, $(Y,d)$ will denote a \(\delta\)-hyperbolic metric space.

\begin{definition}[Fellow-travelling geodesics]\label{def: fellow-travelling}
Two geodesic segments $\gamma_1, \gamma_2$, length parametrized, are said to $l_0$-fellow-travel each other for $l_0 \geq 0$ if $d(\gamma_1(t), \gamma_2(t)) \leq l_0$ for all $0 \leq t \leq l(\gamma_1) = l(\gamma_2).$
\end{definition}

It will be useful to keep the following fact in mind.

\begin{proposition}\label{prop: ftbound}
    Let $\varepsilon_1,\varepsilon_2\geq 0$. Let $\gamma_1$ and $\gamma_2$ be two length parametrized geodesic segments. If $d(i(\gamma_1),i(\gamma_2))\leq \varepsilon_1$ and $d(\tau(\gamma_1),\tau(\gamma_2))\leq \varepsilon_2$, then for all $t$, we have:
    \begin{equation*}
        d(\gamma(t),\gamma'(t))\leq 2\max(\varepsilon_1,\varepsilon_2)+4\delta 
    \end{equation*}
\end{proposition}

We will require the following direct consequence of Gromov's tree approximation for hyperbolic polygons (see \cite[6.1]{gromov1987hyperbolic} and \cite[Chapter 2, Theorem 12]{ghys2013groupes}).

\begin{figure}[h]
    \centering
    \includegraphics[width=0.5\columnwidth]{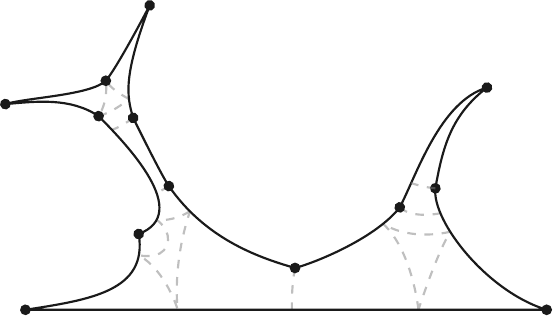}
    \caption{Pairing segments of an $(n+1)$-sided polygon: paired segments \(2 \delta \lceil\log_2 n\rceil\)-fellowtravel each other, have equal length, and there's no ``crisscrossing" between two separate pairs.}
    \label{fig: polygon pairing}
\end{figure}

\begin{lemma}\label{lem: polygonpairing}
    Every geodesic $(n+1)$-gon in $Y$ is $\delta \lceil \log_2n \rceil$-thin.
    
    More precisely, let \(\gamma_0\cdot\gamma_1 \cdots \gamma_n\) be a geodesic polygon in \(Y\). There is a partitioning of each \(\gamma_i = \gamma_{i, 1}\cdots \gamma_{i, m_i}, m_i \geq 1\) and a pairing (that is, a partitioning into subsets of size 2) on the set \(\cup_i \{\gamma_{i,1}, \cdots , \gamma_{i, m_i}\}\) such that the following holds (see \Cref{fig: polygon pairing}):
    \begin{enumerate}[font=\normalfont]
      \item Any such pair \(\gamma_{i_1, m_{j_1}}, \gamma_{i_2, m_{j_2}}\) \(2 \delta \lceil\log_2 n\rceil\)-fellow-travel each other and have the same length.
      \item Let \(\{\gamma_{i_1, m_{j_1}}, \gamma_{i_2, m_{j_2}}\}\) and \(\{\gamma_{i_3, m_{j_3}}, \gamma_{i_4, m_{j_4}}\}\) two such distinct pairs. Then, if we remove $\gamma_{i_1, m_{j_1}}$ and $\gamma_{i_2, m_{j_2}}$ from the polygon \(\gamma_1 \cdots \gamma_n\), the interiors of \(\gamma_{i_3, m_{j_3}}\) and \(\gamma_{i_4, m_{j_4}}\) lie in the same component of the remainder.
    \end{enumerate}
\end{lemma}

Next we define a slightly quantitative version of geodesic tracking which will play a key role in our counting arguments in Section~\ref{sec: Element count: Upper bounds}. 
% Also, see \Cref{fig: tracking}
\begin{definition}\label{def: tracking}
    Let $c,D,l_0>0$. In \(Y\), consider a geodesic $\gamma$ and a piece-wise geodesic $\alpha = \gamma_1'\cdots\gamma_n'$, between two points in \(Y\). We will say that $\alpha$ \emph{$(c,D,l_0)$-tracks} $\gamma$ if the following holds:
    \begin{enumerate}
        \item Every geodesic segment $\gamma'_j$ that is at least $D$-long, i.e. $l(\gamma'_j)\geq D$, \(\gamma'_j\) can be written as a concatenation of three collinear geodesic segments, $\gamma'_j=(\gamma'_j)^{i}(\gamma'_j)^{m}(\gamma'_j)^{f}$ with $l((\gamma'_j)^{i}),l((\gamma'_j)^{f})\leq c$. Moreover, there is a subsegment $\gamma_j$ of $\gamma$ such that, $\gamma_j$ and $(\gamma'_j)^{m}$ $l_0$-fellow-travel.
        \item The subsegments $\{\gamma_j\}$ corresponding to the long geodesic segments $\gamma'_j$ as above, are mutually disjoint. 
    \end{enumerate}
\end{definition}
    % \begin{figure}[h]
    % \centering
    % \includegraphics[width=0.7\columnwidth]{tracking}
    % \caption{$(c,D,l_0)$-tracking}
    % \label{fig: tracking}
    % \end{figure}

\begin{lemma}\cite[Lemma 10]{Lysenok}\label{lemma: conjugationbound}
    Let $G$ be a $\delta$-hyperbolic group equipped with the word metric induced by some finite generating set $S$. If $g,h\in G$ are conjugate, then there exists an element $w\in G$ such that $g=whw^{-1}$, and:
    \begin{equation}
        l(w)\leq \frac{1}{2}(l(g)+l(h))+M
    \end{equation}
    Here $M$ denotes the number of elements in $G$ represented by words in $F(S)$ of length at most $4\delta$.
\end{lemma}

The following definition and lemma appear in \cite{gromov1987hyperbolic}. We cite the version in \cite{gouezel2023exponential} for convenience.

\begin{definition}\label{def: alignedchain}
    Let $C,D\geq 0$. A sequence of points $x_0,x_1,\dots, x_n$ will be called a $(C,D)$-chain if $(x_{i-1},x_{i+1})_{x_i}\leq C$ for all $1\leq i<n$, and $d(x_i,x_{i+1})\geq D$ for all $0\leq i<n$.  
\end{definition}

\begin{lemma}\label{lem: alignmentandprogress}
    Let $x_0,x_1,\dots x_n$ be a $(C,D)$-chain with $D\geq 2C+2\delta+1$, then:
    \begin{align*}
        d(x_0,x_n)\geq \sum_{i=0}^{n-1}(d(x_i,x_{i+1})-(2C+2\delta))
    \end{align*}
    Moreover, if we have $D\geq 2C+4\delta+1$, then for all $0<i<n$, $(x_0,x_n)_{x_i}\leq C+2\delta$.
\end{lemma}

\begin{remark}
    We will often say that the points $x,y,z$ are $C$-aligned if $(x,z)_y\leq C$.
\end{remark}

\begin{proposition}\label{prop: chainfellowtravelling}
    Let $C>0$ and $m\geq 2$ be a positive integer. Let $x_0,\dots,x_{m}$ be a $(C,D)$-chain with $D\geq 5(C+30\delta)$. Let $\gamma_j$ denote a geodesic segment $[x_{j-1},x_{j}]$ for all $1\leq j\leq m$ and let $\gamma$ denote a geodesic segment $[x_0,x_{m}]$. Then:
    \begin{enumerate}[font=\normalfont]
        \item $\cup_{i=1}^{m}\gamma_i\subset N_{C+30\delta}(\gamma)$ and $\gamma\subset N_{C+30\delta}(\cup_{i=1}^{m}\gamma_i)$.
    \end{enumerate}
    There are decompositions of $\gamma$ and $\gamma_j$ $(1\leq j\leq m)$: $\gamma=\gamma'_1\beta'_1\gamma'_2\dots\beta'_{m-1}\gamma'_{m}$, $\gamma_1=\alpha_1\beta_1$, $\gamma_m=\beta_m\alpha_m$, and $\gamma_j=\beta_j^{(0)}\alpha_j\beta_j^{(1)}$ $(2\leq j \leq m-1)$ such that for all $1\leq j\leq m$,
    \begin{enumerate}[font=\normalfont, resume]
        \item  $l(\alpha_j)=l(\gamma'_j)$.
        \item $\alpha_j$ and $\gamma'_j$ $4(C+30\delta)$-fellow travel.
        \item Lengths of $\beta_1,\beta_m$, and  $\beta_j^{(0)},\beta_j^{(1)}$ are at most $2(C+30\delta)$. 
        \item $l(\beta'_j)=2(C+30\delta)$.
    \end{enumerate}
\end{proposition}
\begin{proof}
    Proof of $(1)$ essentially follows along the same lines as that of \cite[H, Theorem 1.13]{bridson2013metric}. $(2),(3),(4)$ follow from $(1)$ and \Cref{prop: ftbound}.
\end{proof}

\begin{definition}\label{def: Schottkyset}
    Suppose $0<\eta<1$, $K,D>0$. $S\subset G$ is said to be a $(\eta,K,D)$-Schottky set if:
    \begin{enumerate}
        \item For all $g,h\in G$, $\#\{s\in S | (g,sh)\leq K\}\geq (1-\eta)\#S$.
        \item For all $g,h\in G$, $\#\{s\in S | (g,s^{-1}h)\leq K\}\geq (1-\eta)\#S$.
        \item For all $s\in S$, $l(s)\geq D$.
    \end{enumerate}
\end{definition}

\begin{proposition}\cite[Proposition 3.12]{gouezel2023exponential}\label{prop: existenceofSchottkyset}
    Let $G$ be a hyperbolic group, equipped with the word metric induced by a finite generating set $S$. Let $g_1$ and $g_2$ be two independent loxodromic elements of $G$, that is, their respective sets of fixed points at infinity are disjoint. For all $\eta>0$, there exists a $K>0$ such that for all $D>0$, there exists a positive integer $n$ and a $(\eta,K,D)$-Schottky set in $\{w_1\dots w_n|w_i\in\{g_1,g_2\}, 1\leq i\leq n\}$.
\end{proposition}
\begin{remark}\label{rmk: rmktoexistenceofSchottkyset}
    Gou\"{e}zel's construction of Schottky sets in \cite{gouezel2023exponential} in fact has:
    \begin{itemize}
       \item For all $g,h\in G$, $\#\{s\in S | (g,sh)\leq K\}\geq \#S-2$.
        \item For all $g,h\in G$, $\#\{s\in S | (g,s^{-1}h)\leq K\}\geq \#S-2$. 
    \end{itemize}
    Moreover, for $\eta=2^{-m}$, one can choose the Schottky set so that $\#S=2^{m+1}$.
\end{remark}

 Let $G$ be a group acting on a metric space $(Y,d)$ properly. For any $y\in Y$ we define:
    \begin{align*}
        B(y,L):=\{g\in G| d(y,g\cdot y)\leq L\}
    \end{align*}

\begin{theorem}[\cite{coornaert1993mesures}]%
    \label{thm: ballgrowth}
    Let $G$ be a group acting properly, co-compactly by isometries on a proper, geodesic $\delta$-hyperbolic metric space $Y$ whose boundary contains at least three points. Then there exist positive constants $C_1,C_2,h$ such that:
    \begin{equation*}
        \textnormal{for all } y\in Y,\textnormal{ for all } L\geq 0,\quad C_1e^{hL}\leq B(y,L)\leq C_2e^{hL} 
    \end{equation*}
\end{theorem}

\section{Quasigeodesic substitutions to the Pre-Wicks form}\label{sec: wicks form moves}
    
\noindent Fix a positive integer $n$. Let $G$ be a hyperbolic group and $S$ be a finite symmetric generating set. Equip $G$ with the word metric induced by $S$, so that it is $\delta$-hyperbolic. Define $l_0^{+}(n)\coloneq 2\delta\lceil\log_2(12n-1)\rceil$ and $l_0^{-}(n)\coloneq 2\delta\lceil\log_2(6n-1)\rceil$. 

\begin{definition}
    \label{def: substitution} Given a word $W$ in a set of unknowns $\mathcal{X}$, a \emph{$G$-valued substitution} $\mathfrak{g}$ of $W$ assigns to each unknown $X$ appearing as a letter in $W$ an element $g_{X} \in G \setminus \{1\}$. Now suppose $W=X_1\cdots X_m$, where each $X_i\in \mathcal{X}\cup\mathcal{X}^{-1}$. For each $i$, $g_{X_i}\coloneq g_{X}$ if $X_i=X$ and $g_{X_i}\coloneq g_{X}^{-1}$ if $X_i=X^{-1}$, where $X\in\mathcal{X}$.
    Any substitution $\mathfrak{g}$ naturally defines an element $\mathfrak{g}(W)\coloneq g_{X_1}\cdots g_{X_m}$. We will say that $\mathfrak{g}(W)$ is obtained from the word $W$ by the substitution $\mathfrak{g}$. In particular, when $G$ is a free group, we say that $\mathfrak{g}$ is a \emph{non-cancelling substitution} for $W$ if for every subword $X_iX_{i+1}$ of $W$, there is no cancellation between $g_{X_i}$ and $g_{X_{i+1}}$. $\mathfrak{g}$ is said to be a \emph{cyclically non-cancelling substitution} if for every cyclic permutation of $W$, $\mathfrak{g}$ is a non-cancelling substitution.  We also define the \textit{length of a substitution}, $l(\mathfrak{g})\coloneq\sum_{i=1}^ml(g_{X_i})$.
\end{definition}

\begin{definition}\label{def: path}
    Let $W=X_1\dots X_m$ be a word in a set of unknowns $\mathcal{X}$ and $\mathfrak{g}$ be a $G$-valued substitution. We associate a piece-wise geodesic path $\gamma(\mathfrak{g},W)$ to the tuple $(W,\mathfrak{g})$ as follows: for each unknown $X\in\mathcal{X}$, choose a geodesic word representing $g_X$. This defines a choice of a geodesic word representing $g_{X_i}$ for each $i$ in the obvious way. $\gamma(\mathfrak{g},W)$ is defined as the concatenation of the paths $\gamma_{X_1},\dots,\gamma_{X_m}$, where each $\gamma_{X_i}$ is a geodesic path joining $g_{X_1}\cdots g_{X_{i - 1}}$ to $g_{X_1}\cdots g_{X_{i}}$ corresponding to the geodesic word representing $g_{X_i}$. The argument $W$ in $\gamma(\mathfrak{g},W)$ will be skipped whenever $W$ is clear from the context. 
\end{definition}

Suppose $g$ is a group element having either orientable genus $n$ or non-orientable genus $n$, for some $n\geq 1$. Then, there exists a Wicks form $(W_0,\mathfrak{g}_0)$ such that $\mathfrak{g}_0(W_0)=g$, where $W_0$ is orientable if $g$ is orientable, and non-orientable otherwise. We illustrate this for the orientable case. Suppose $g$ is a group element of orientable genus $n$. We can write $g$ as a product of $n$ non-trivial commutators, $[a_1,b_1],\cdots,[a_n,b_n]$. Thus, we have an orientable Wicks form of genus $n$, namely $W_0\colonequals [A_1,B_1]\cdots[A_n,B_n]$ together with the substitution $\mathfrak{g}_0$ defined by: for all $1\leq i\leq n$, $g_{A_i}=a_i, g_{B_i}=b_i$, and $\mathfrak{g}_0(W_0)=g$. A piece-wise geodesic $\gamma(\mathfrak{g}_0)$ associated with $(W_0,\mathfrak{g}_0)$ can be wild, and in particular, there could be significantly long fellow-travelling between the geodesic pieces comprising it. We describe how to alter $(W_0, \mathfrak{g}_0)$ to a tuple $(W_t, \mathfrak{g}_t)$ such that $\mathfrak{g}_t(W_t) = \mathfrak{g}(W) = g$, and the geodesic pieces $\gamma_{t,x_1}, \dots, \gamma_{t,x_m}$ have no long fellow-travelling between each other. Consequently, $\gamma(\mathfrak{g}_t)$ is a quasigeodesic with uniform constants. We obtain $(W_t, \mathfrak{g}_t)$ in a sequence of steps starting at $(W_0, \mathfrak{g}_0)$, by applying certain moves to the tuple at hand, whenever possible. These moves are defined in the following subsection.

\subsection{Pre-Wicks Form Moves}\label{subsec: wicksformmoves}
We will now define three moves that take as input a tuple $(W,\mathfrak{g})$, where $W$ is a pre-Wicks form and $\mathfrak{g}$ is a $G$-valued substitution, and output another tuple $(W',\mathfrak{g}')$ such that $\mathfrak{g}'(W')=\mathfrak{g}(W)$. More details on what the moves do to the orientation of the pre-Wicks form and its genus are proved in \Cref{lemma: moveproperties}. We write the input pre-Wicks form as $W=X_1\dots X_{2m}$. Recall from \Cref{def: quadratic sign orientable non-orientable}, that for every $X_i$, there is a unique letter $X_{\bar{i}}$, the letter paired with $X_i$, such that $X_{\bar{i}}=X_i^{-1}$ or $X_i$, depending on the sign of $X_i$. There are two types of fellow travelling we need to address.:

\begin{itemize}
    \item Fellow travelling between $X_i$ and $X_k$ where $k\neq \bar{i}$: Move $1$ applies in this case.
    \item Fellow travelling between $X_i$ and $X_{\bar{i}}$: When $\text{sign}(X_i)=-1$, Move $2$ applies. Move $3$ applies otherwise.
\end{itemize}

\noindent\Cref{lemma: moveproperties} then ensures that starting from an initial tuple $(W_0,\mathfrak{g}_0)$ as before, we arrive at a terminal tuple $(W_t, \mathfrak{g}_t)$ to which none of the three moves apply.
\\In each Move, if $W$ is orientable, we assume it has orientable genus $n$, otherwise we assume it has non-orientable genus $n$. Moreover, the constant $l_0=l_0^+$ if $W$ is orientable, and $l_0=l_0^-$ otherwise.

\noindent \textbf{Move 1: Removing bottlenecks}

\begin{figure}[h]
    \centering
    \includegraphics[width=0.55\columnwidth]{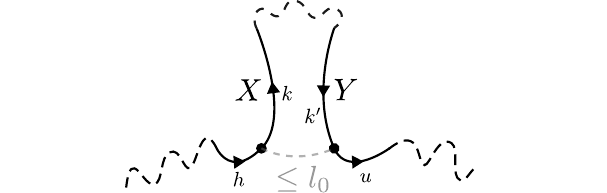}

    \caption{A bottleneck}
    \label{fig: bottleneck0}
\end{figure}

\begin{figure}[h]
    \centering
    \begin{subfigure}{\columnwidth}
        \centering
        \includegraphics[width=0.55\columnwidth]{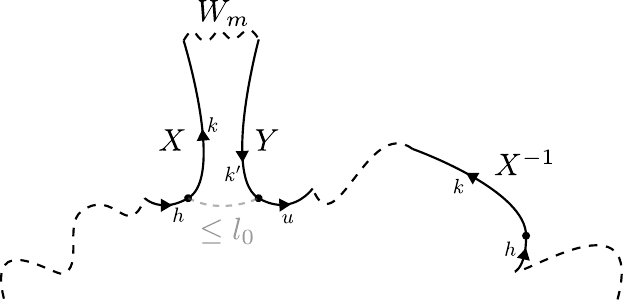}
        \caption{Input}
        \label{fig: bottleneck1 input}
    \end{subfigure}

    \begin{subfigure}{\columnwidth}
        \centering
        \includegraphics[width=0.55\columnwidth]{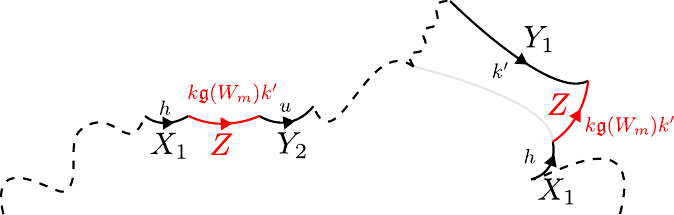}
        \caption{Output}
        \label{fig: bottleneck1 output}
    \end{subfigure}

    \caption{Removing bottlenecks: the $1^{\text{st}}$ case with $\text{sign}(X)=-1$}
    \label{fig: bottleneck1}
\end{figure}

\linkdest{Move 1}
\noindent\textbf{Input:} This move is defined for both orientable and non-orientable pre-Wicks forms, when the input $(W,\mathfrak{g})$ satisfies the following condition: There are letters $X_{j_1},X_{j_2}$ in $W$ with $j_2\neq j_1,\overline{j_1}$, and there are subsegments of $\gamma_{X_{j_1}}$ and $\gamma_{X_{j_2}}$ which $l_0$-fellow-travel. More precisely, for some choice of  $\gamma(\mathfrak{g})$, we have either of the following cases depending on whether $X_{j_1}$ appears before $X_{j_2}$ or not: (see \Cref{fig: bottleneck0}):

\begin{enumerate}[label=(\Alph*)]
    \item $W=W_iX_{j_1}W_mX_{j_2}W_f$ and $g_{X_{j_1}}=h \cdot k$, $g_{X_{j_2}}=k' \cdot u$, $\gamma_{X_{j_1}}=\gamma_h  \gamma_k$, $\gamma_{X_{j_2}}=\gamma_{k'} \gamma_u$. Writing $\tau(\gamma_h)=p$ and $\tau(\gamma_{k'})=q$, we have:
    \begin{enumerate}[label=(A\arabic*)]
        \item If $X_{\overline{j_1}}$ does not appear in $W_m$, then $l(k)> l_0\text{ and }d(p,q)\leq  l_0$.
        \item If $X_{\overline{j_1}}$ appears in $W_m$, then $l(h)> l_0\text{ and }d(p,q)\leq  l_0$.
    \end{enumerate}
    \item $W=W_iX_{j_2}W_mX_{j_1}W_f$ and $g_{X_{j_1}}=h \cdot k$, $g_{X_{j_2}}=k' \cdot u$, $\gamma_{X_{j_1}}=\gamma_h  \gamma_k$, $\gamma_{X_{j_2}}=\gamma_{k'} \gamma_u$. Writing $\tau(\gamma_h)=p$ and $\tau(\gamma_{k'})=q$, we have:
    \begin{enumerate}[label=(B\arabic*)]
        \item If $X_{\overline{j_1}}$ does not appear in $W_m$, then $l(h)> l_0\text{ and }d(p,q)\leq  l_0$.
        \item If $X_{\overline{j_1}}$ appears in $W_m$, then $l(k)> l_0\text{ and }d(p,q)\leq  l_0$.
    \end{enumerate}
\end{enumerate}
The output for case (B) is the symmetric counterpart of (A), so we only write down the output for (A). We assume that $X_{j_1}=X$ and $X_{j_2}=Y$, for some unknowns $X$ and $Y$. Again, the output is easily defined in the cases when $X_{j_1}=X^{-1}$ or $X_{j_2}=Y^{-1}$ analogously. Set $\varepsilon=\text{sign}(X_{j_1})$.
\\\textbf{Output:}  \linkdest{Move 1 Case A1}\textbf{Case A1: $X_{\overline{j_1}}$ does not appear as a letter in $W_m$} See \Cref{fig: bottleneck1}. 

We define a new word $W'$ as follows:
\begin{enumerate}
    \item Perform the substitution:
    \begin{equation}
        X\mapsto X'_1ZY_1^{-1}W_{m}^{-1},\quad Y\mapsto Y_1Y_2
    \end{equation}
    where $Z,X'_1,Y_1,Y_2$ are unknowns that do not appear in $W$.
    \item Perform reductions.
\end{enumerate}
Now we define a $G$-valued substitution $\mathfrak{g'}$ for $W'$:
\begin{enumerate}[resume]
    \item 
    \begin{equation}
        g'_{X'_1} \coloneq h,\quad g'_{Z}\coloneq k\mathfrak{g}(W_m)k',\quad g'_{Y_1}\coloneq k',\quad g'_{Y_2}\coloneq u
    \end{equation}
    and set $g'_V \coloneq g_V$ for every other letter $V$ appearing in $W'$.
\end{enumerate}
The new word $W'$ may have some unknowns that get assigned the identity element. There also may be some redundancies. Thus, as a final step:
\begin{enumerate}[resume]
    \item Remove redundancies from $(W',\mathfrak{g'})$ and remove any unknowns which get assigned identity by \(\mathfrak{g}'\).
\end{enumerate}

\begin{remark}[Removing backtracking in \(\gamma_{\mathfrak{g}}\)]\label{rmk: sharp corners}
    When \(X\) and \(Y\) above are consecutive letters in \(W\), that is, when \(W_m\) is the trivial word, then applying the bottleneck move removes backtracking between \(\gamma_X\) and \(\overline{\gamma}_Y\) beyond a certain length. See \Cref{fig: sharp corner} below: 
    % {\color{red}{$Y = Y_2 Y_1$ in the figure, inconsistent with proof.}}
    % {\color{gray}\(Y\) is sent, between different moves, to both \(Y_1Y_2\) and \(Y_2Y_1\). This can be made consistent.}

    \begin{figure}[h]
        \begin{subfigure}{\columnwidth}
        \centering
        \includegraphics[width=0.55\columnwidth]{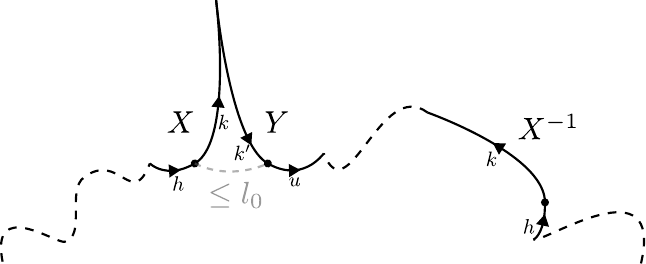}
        \caption{Input}
        \label{fig: sharp corner input}
        \end{subfigure}\\
        \begin{subfigure}{\columnwidth}
        \centering
        \includegraphics[width=0.55\columnwidth]{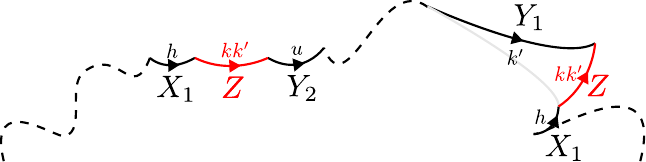}
        \caption{Output}
        \label{fig: sharp corner output}
        \end{subfigure}
        
        \caption{Removing backtracking between consecutive geodesic segments}
        \label{fig: sharp corner}
    \end{figure}

    When the underlying group $G$ is free, we can always choose to apply \linkref{Move 1}{Move 1} in only this fashion for removing backtracking. Furthermore, in this case, it will suffice to use only Move $1$ out of all three available Moves. 
\end{remark}

% {\color{gray}The assumption of the following case is when \(X^{-1}\) appears as a letter in \(W_m\). However, note that then there exists a letter \(Y'\) in \(W_m\) such that \(\gamma_{Y'}\) is both long and fellow-travels with \(\gamma_{X^{-1}}\). Therefore we can safely apply Case A1 to the pair \(X^{-1}\) and \(Y'\) instead, just as in the above note for free groups.}

\linkdest{Move 1 Case A2}\noindent\textbf{Case A2: $X_{\overline{j_1}}$ appears as a letter in $W_m$}\\In this case, we have $W_m=W_m^i X_{\overline{j_1}}W_m^f$ (see \Cref{fig: bottleneck2}), and 
we construct a new word $W'$ as follows:

\begin{figure}[h]
    \begin{subfigure}{.49\columnwidth}
        \centering
        \includegraphics[width=.8\columnwidth]{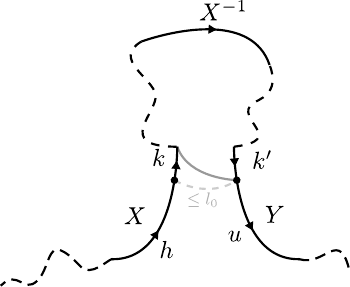}
        \caption{Input}
        \label{fig: bottlneck 2 input}
    \end{subfigure}
    \begin{subfigure}{.49\columnwidth}
        \centering
        \includegraphics[width=.8\columnwidth]{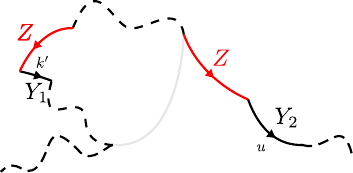}
        \caption{Output}
        \label{fig: bottleneck 2 output}
    \end{subfigure}

    \caption{Removing bottlenecks: the $2^{\text{nd}}$ case with $\text{sign}(X)=-1$}
    \label{fig: bottleneck2}
\end{figure}
\begin{enumerate}
    \item Perform the substitution:
    \begin{equation}
        X\mapsto \begin{cases}
            W_m^f Y_1Z^{-1}W_m^i &;\varepsilon=-1\\
            (W_m^i)^{-1}Z(W^f_mY_1)^{-1} &;\varepsilon=+1
        \end{cases} ,\quad Y\mapsto Y_1Y_2
    \end{equation}
    where $Z,Y_1,Y_2$ are unknowns that do not appear in $W$.
    \item Perform reductions.
\end{enumerate}
Now we define a $G$-valued substitution $\mathfrak{g'}$ for $W'$:
\begin{enumerate}[resume]
    \item 
    \begin{equation}
        \quad g'_{Z}\coloneq
            \mathfrak{g}(W_{m}^{i})(hk)^{\varepsilon}\mathfrak{g}(W_{m}^{f})k'
        ,\quad g'_{Y_1}\coloneq k',\quad g'_{Y_2}\coloneq u
    \end{equation}
    and set $g'_V \colonequals g_V$ for every other letter $V$ appearing in $W'$.
    \item Remove redundancies and unknowns which get assigned identity from $(W',\mathfrak{g'})$.
\end{enumerate}

\noindent \textbf{Move 2: Lowering the Rise
}
\linkdest{Move 2}

\begin{figure}[h]
    \begin{subfigure}{.49\columnwidth}
        \centering
        \includegraphics[width=.8\columnwidth]{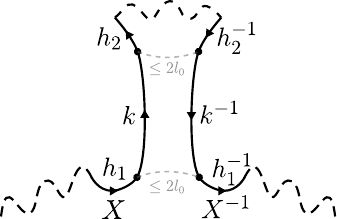}
        \caption{Input}
        \label{fig: rise 1 input}
    \end{subfigure}
    \begin{subfigure}{.49\columnwidth}
        \centering
        \includegraphics[width=.8\columnwidth]{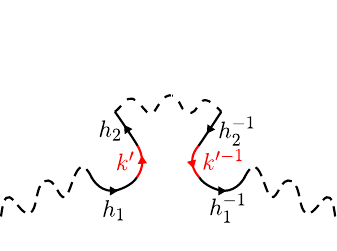}
        \caption{Output}
        \label{fig: rise 1 output}
    \end{subfigure}

    \caption{\protect\linkref{Move 2}{Move 2:} Lowering the rise, the $\text{sign}(X) = -1$ case.}
    \label{fig: rise 1}
\end{figure}

\noindent \textbf{Input:} This move is defined when \((W, \mathfrak{g})\) satisfies the following conditions (see \Cref{fig: rise 1}): There is a letter $X_j$ with $\text{sign}(X_j)=-1$  and there are matching subsegments of  $\gamma_{X_j}$ and $\gamma_{X_{\overline{j}}}$ that $l_0$-fellow travel. More precisely, for some choice of  $\gamma(\mathfrak{g})$ we have, $W=W_iX_jW_mX_{\bar{j}}W_f$, $g_{X_j}=h_1kh_2$, $\gamma_{X_j}=\gamma_{h_1}\gamma_k\gamma_{h_2}$, $\gamma_{X_{\bar{j}}}=\gamma_{h_2^{-1}}\gamma_{k^{-1}}\gamma_{h_{1}^{-1}}$. Writing $i(\gamma_k)=p_1$, $\tau(\gamma_{k^{-1}})=p_2$, $\tau(\gamma_k)=q_1$, $i(\gamma_{k^{-1}})=q_2$, we have:
    \begin{align*}
        l(k)> 2l_0+M  \text{ and }d(p_1,p_2)\leq 2l_0, d(q_1,q_2)\leq 2l_0
    \end{align*}
    where $M$ is the constant appearing in \Cref{lemma: conjugationbound}.

Given this input, we set $a_1\coloneq h_2\mathfrak{g}(W_m)h_{2}^{-1}$ and $a_2\coloneq (\mathfrak{g}(W_i)h_1)^{-1}\mathfrak{g}(W)(h_1^{-1}\mathfrak{g}(W_f))^{-1}$. Observe that we have:
    \begin{align*}
        ka_1k^{-1}=a_2
    \end{align*}
    Because $l(a_1)\leq 2l_0$ and $l(a_2)\leq 2l_0$, by \Cref{lemma: conjugationbound}, there exists a $k'\in G$ with $l(k')\leq 2l_0 +M$, such that: $k'a_1(k')^{-1}=a_2$. The output will be constructed using $k'$.
    \\\\\noindent\textbf{Output:} To begin with, we leave $W$ untouched, that is, $W'\coloneq W$, but we define a new $G$-valued substitution $\mathfrak{g'}$ as follows: 
    \begin{align*}
        g'_{X}\coloneq h_1k'h_2
    \end{align*}
    and set $g'_V\coloneq g_V$ for every other letter $V$ appearing in $W'$. Next, we remove redundancies and unknowns which get assigned identity from $(W',\mathfrak{g'})$.
\\\\\noindent\textbf{Move 3: Lowering the Rise (contd.)}
\linkdest{Move 3}
\begin{figure}[h]
    \begin{subfigure}{.49\columnwidth}
        \centering
        \includegraphics[width=.8\columnwidth]{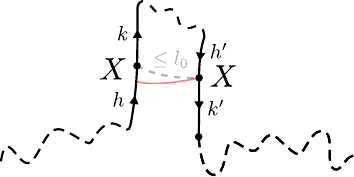}
        \caption{Input}
        \label{fig: move 3 input}
    \end{subfigure}
    \begin{subfigure}{.49\columnwidth}
        \centering
        \includegraphics[width=.8\columnwidth]{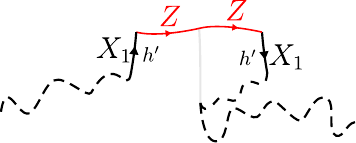}
        \caption{Output}
        \label{fig: move 3 output}
    \end{subfigure}

    \caption{\protect\linkref{Move 3}{Move 3}: Lowering the rise, the $\text{sign}(X)= +1$ case.}
    \label{fig: move 3}
\end{figure}
\\\noindent\textbf{Input:} This move is defined when $(W,\mathfrak{g})$ satisfies the following conditions: There is a letter $X_j$ with $\text{sign}(X_j)=1$, and there are subsegments of $\gamma_{X_{j}}$ and $\gamma_{X_{\bar{j}}}$ which come $l_0$-close. More precisely, $W=W_iX_jW_mX_{\bar{j}}W_f$ or $W=W_iX_{\bar{j}}W_mX_jW_f$, $X_{\bar{j}}=X_j$, and $g_{X_j}=hk=h'k'$, $\gamma_{X_j}=\gamma_h\gamma_k$, $\gamma_{X_{\bar{j}}}=\gamma_{h'}\gamma_{k'}$. Writing $\tau(\gamma_h)=p$ and $\tau(\gamma_{h'})=q$, we have:
    \begin{align*}
        l(h), l(k)>l_0\text{ and }d(p,q)\leq l_0
    \end{align*}
\noindent\textbf{Output:} We describe the output for the case when $W=W_iX_jW_mX_{\bar{j}}W_f$, and the other case is completely analogous. Note that $l(h)+l(k)=l(h')+l(k')=l(g_X)$. So either $l(k)\leq l(k')$ or $l(h)<l(h')$. Accordingly there are two cases for the output.\\\\
\linkdest{Move 3 Case 1}\noindent\textbf{Case 1: $\boldsymbol{{l(k)\leq l(k')}}$.} In this case we construct a new word $W'$ as follows (see \Cref{fig: move 3}): 
\begin{enumerate}
    \item Perform the substitution:
    \begin{equation}
        X\mapsto X'_1Z(X'_1)^{-1}W_m^{-1}
    \end{equation}
    where $Z$ and $X'_1$ are unknowns that do not appear in $W$.
    \item Perform reductions.
\end{enumerate}
Now we define a $G$-valued substitution for $W'$:
\begin{enumerate}[resume]
    \item 
    \begin{equation}
        g'_{X'_1}\coloneq h' ,\quad g'_Z\coloneq k'\mathfrak{g}(W_m)h'
    \end{equation}
    and set $g'_V\coloneq g_V$ for every other letter appearing in $W'$.
    \item Remove redundancies from $(W',\mathfrak{g}')$.
\end{enumerate}
\linkdest{Move 3 Case 2}\noindent\textbf{Case 2: $\boldsymbol{l(h)< l(h')}$.} In this case we construct a new word $W'$ as follows: 
\begin{enumerate}
    \item Perform the substitution:
    \begin{equation}
        X\mapsto W_m^{-1}(X'_2)^{-1}ZX'_2
    \end{equation}
    where $Z$ and $X'_2$ are unknowns that do not appear in $W$.
    \item Perform reductions.
\end{enumerate}
Now we define a $G$-valued substitution for $W'$:
\begin{enumerate}[resume]
    \item 
    \begin{equation}
        g'_{X'_2}\coloneq k',\quad g'_Z\coloneq k'\mathfrak{g}(W_m)h'
    \end{equation}
    and set $g'_V\coloneq g_V$ for every other letter appearing in $W'$.
    \item Remove redundancies from $(W',\mathfrak{g}')$ and unknowns which get assigned identity by \(\mathfrak{g}'\).
\end{enumerate}

Now that we have defined all the moves, we list down some basic properties that they share in the following lemma. 

\begin{lemma}\label{lemma: moveproperties}
    Let $g\in[G,G]$ be an element of orientable genus $n$. Let $(W,\mathfrak{g})$ be an orientable pre-Wicks form of orientable genus $n$ and a $G$-substitution with $\mathfrak{g}(W)=g$. For any $1\leq i\leq 2$, suppose we can apply Move \(i\) to $(W,\mathfrak{g})$ with the corresponding output being \((W', \mathfrak{g}')\). Then the following is true:
    \begin{enumerate}[font=\normalfont]
        \item $l(\mathfrak{g}')<l(\mathfrak{g})$
        \item the output $(W',\mathfrak{g'})$ is an orientable pre-Wicks form of orientable genus $n$ with $\mathfrak{g}'(W')=g$.
    \end{enumerate}
    Let $x\in G^2$ be an element of non-orientable genus $n$. Let $(U,\mathfrak{h})$ be a non-orientable pre-Wicks form of non-orientable genus $n$ and a $G$-substitution $\mathfrak{h}$ with $\mathfrak{h}(U)=x$. For any $1\leq i\leq 3$, suppose we can apply Move \(i\) to $(U,\mathfrak{h})$ with the corresponding output being \((U', \mathfrak{h}')\). Then the following is true:
    \begin{enumerate}[resume, font=\normalfont]
        \item $l(\mathfrak{h}')<l(\mathfrak{h})$
        \item $\mathfrak{h}'(U')=x$. $U'$ is either a non-orientable pre-Wicks form of non-orientable genus $n$ or an orientable pre-Wicks form of orientable genus $(n-1)/2$. The latter is possible only when $n$ is odd and $\textnormal{genus}^+(x)=(n-1)/2$.
    \end{enumerate}
\end{lemma}

\begin{proof}
 To begin with, we prove $(1)$ and $(3)$ for each Move, which say that whenever a Move is applied to a tuple $(W,\mathfrak{g})$, the length of the substitution (see \Cref{def: substitution}) must decrease. It suffices to prove this for the outputs obtained from each of the Moves, right before the last step is performed, that is, before removing redundancies. This is because removing redundancies does not increase the length of the substitution (directly follows from the triangle inequality).
        
        Moves $1$ and $2$ apply to both orientable and non-orientable pre-Wicks forms, and the arguments required to prove $(1)$ and $(3)$ for these moves are identical. So we only write down proofs for the orientable case, i.e. $(1)$.\\\noindent\linkref{Move 1}{Move 1}: We first consider \linkref{Move 1 Case A1}{Case A1}. It suffices to look at the lengths of elements corresponding to the unknowns $X, Y, X'_1, Y_1, Y_2$ and $Z$ (See \Cref{fig: bottleneck1}). So:
                \begin{align*}
                    l(\mathfrak{g}) - l(\mathfrak{g}')  &= 2 l(g_X) + 2l(g_Y) -2l(g'_{X'_1})- 2l(g'_{Z}) - 2l(g'_{Y_1}) - 2l(g'_{Y_2})\\
                    & = 2 (l(h) + l(k)) + 2(l(k') + l(u)) - 2 l(h) - 2l(k\mathfrak{g}(W_m)k') - 2l(k') - 2l(u)\\
                    &=  2l(k) - 2l(k\mathfrak{g}(W_m)k') > 0
                \end{align*}
                The inequality in the last step follows from the hypothesis on \((W, \mathfrak{g})\), that is, $l(k)>l_0\geq l(k\mathfrak{g}(W_m)k')$.\\Next, consider \linkref{Move 1 Case A2}{Case A2}. We only need to consider the unknowns $X, Y, Y_1, Y_2$ and $Z$. The remaining unknowns, including the ones appearing in $W_m^i$ and $W_m^f$, for example, have their assignments, and thus lengths, unchanged. See \Cref{fig: bottleneck2}. We also note that $l(g'_Z)\leq l(k)+d(p,q)=l(k)+l(k\mathfrak{g}(W_m)k')$.
                \begin{align*}
                l(\mathfrak{g}) - l(\mathfrak{g}') & = 2 l(g_X) + 2l(g_Y) - 2l(g'_{Y_1}) - 2l(g'_{Y_2}) - 2 l(g'_Z)\\
                & = 2 ( l(h) + l(k) ) + 2(l(k') + l(u)) - 2l(k') - 2l(u) - 2 l(g'_Z)\\
                &=2(l(h)+l(k))-2l(g'_Z)\\
                &\geq 2(l(h)+l(k))-2(l(k)+l(k\mathfrak{g}(W_m)k'))\\
                &=2(l(h)-l(k\mathfrak{g}(W_m)k'))>0
            \end{align*}
            where the inequality in the last step follows from the fact that $l(h) > l_0 \geq l(k\mathfrak{g}(W_m)k'))$.\\\noindent \linkref{Move 2}{Move 2}: Notice that in this move, only the substitution for $X$ has been changed. So we have:
            \begin{align*}
                l(\mathfrak{g}) - l(\mathfrak{g}') & = 2 l(g_X) - 2 l(g'_X)\\
                &=2(l(h_1)+l(k)+l(h_2))-2 l(g'_X)\\
                &\geq 2(l(h_1)+l(k)+l(h_2))-2(l(h_1)+l(k')+l(h_2))\\
                &=2(l(k)-l(k'))>0
            \end{align*}
            Once again, the final inequality holds from the hypothesis, because $l(k)>2l_0+M\geq l(k')$.
            \\\noindent \linkref{Move 3}{Move 3}: This move only applies to non-orientable pre-Wicks forms and it has two cases, which are symmetric, so we only argue for \linkref{Move 3 Case 1}{Case 1}, that is, when $U=U_iX_jU_mX_{\bar{j}}U_f$ and $l(k)\leq l(k')$. It suffices to look at the letters $X,X'_1$ and $Z$. See \Cref{fig: move 3}.
            \begin{align*}
                l(\mathfrak{h})-l(\mathfrak{h}')&=2l(g_X)-2l(g'_{X'_1})-2l(g'_Z)\\
                &=2(l(h')+l(k'))-2l(h')-2l(g'_Z)\\
                &=2(l(k')-l(k'\mathfrak{g}(W_m)h'))\\
                &\geq 2(l(k')-l(k'k^{-1})-l(k\mathfrak{g}(W_m)h'))\\
                &=2(l(k)-l(k\mathfrak{g}(W_m)h'))>0
            \end{align*}
            where we have used the fact that $hk=h'k'$ and that $l(k)>l_0\geq l(k\mathfrak{g}(W_m)h')$.

         Next, we prove $(2)$ and $(4)$. The fact that the output substitution gives the same element as in the input is a routine check. Given a tuple $(W,\mathfrak{g})$ with $\mathfrak{g}(W)=x$, on which say Move $i$ is applicable, let $(W_1,\mathfrak{g}_1)$ denote the pre-Wicks form and substitution obtained right before any unknowns getting assigned the identity element are removed. Then: 
         \begin{quote}
             If $W$ is orientable, so is $W_1$, and $\textnormal{genus}^+(W)=\textnormal{genus}^+(W_1)$. If $W$ is non-orientable, so is $W_1$, and $\textnormal{genus}^-(W)=\textnormal{genus}^-(W_1)$.
         \end{quote}
         To see this, note that performing reductions and removing redundancies do not change the quadraticity, orientation, and genus of the pre-Wicks form. Moreover, the substitutions performed in the first step of the Move can be realized by a cut-paste process done on the polygon $P_{W}$, which doesn't change the associated surface $\Sigma_W$, and therefore, both orientation and genus are preserved. The fact that quadraticity is preserved is again a consequence of the cut-paste process: cutting involves breaking up one (or more) unknown(s) as a product of new unknowns (of the same sign) and introducing a new unknown (corresponding to the cut) while pasting involves eliminating an unknown from the word.

         However, a change in orientation can occur when an unknown is removed. Suppose $W_2$ is the pre-Wicks form obtained by removing an unknown $X$ which is assigned the identity element by $\mathfrak{g}_1$. If $W_1$ is orientable, so is $W_2$, and $\textnormal{genus}^+(W_2)\leq \textnormal{genus}^+(W_1)=\textnormal{genus}^+(x)$. But the substitution on $W_2$ also evaluates to $x$, so $\textnormal{genus}^+(W_2)= \textnormal{genus}^+(W_1)$.
         Now suppose $W_1$ is non-orientable. The rest of the proof breaks into two cases:

         If $W_2$ is also non-orientable, then just as argued before, the non-orientable genus is preserved. If $W_2$ is orientable, then $\textnormal{genus}^+(W_2)\leq (\textnormal{genus}^-(W_1)-1)/2$. This follows from the following fact (cf. \cite[Proposition {\RomanCaps{1}}.6.10]{LyndonSchupp2001}):

         \begin{quote}
             \textit{If $W_0$ is a non-orientable quadratic word and $W'_0$ is an orientable quadratic word which is a homomorphic image of $W_0$, then $\textnormal{genus}^+(W'_0)\leq (\textnormal{genus}^-(W_0)-1)/2$.}
         \end{quote}

         Thus, it follows from \Cref{prop: square length upper bound}, that in this case, we must have $\textnormal{genus}^+(W_2)= (\textnormal{genus}^-(W_1)-1)/2$. This completes the proof.
\end{proof}

\begin{definition/proposition}\label{def: tracking pair}
    Let $g$ be an element of orientable (resp. non-orientable) genus $n$, and let $(W,\mathfrak{g})$ be an orientable (resp. non-orientable) pre-Wicks form of orientable (resp. non-orientable) genus $n$ and a $G$-valued substitution such that $\mathfrak{g}(W)=g$. Then, Moves $1$, $2$, and $3$ can only be applied finitely many times to $(W,\mathfrak{g})$. We denote a terminal pre-Wicks form and $G$-substitution tuple representing $g$, by $(W_t,\mathfrak{g}_t)$ and we call it a \emph{tracking pair} for $g$.
\end{definition/proposition}

\begin{corollary}\label{corollary: terminalwordexistence}
    For every group element of $g$ with orientable or non-orientable genus $n$, there is a tracking pair $(W_t,\mathfrak{g}_t)$ with $\mathfrak{g}_t(W_t)=g$.
    If $g$ has orientable genus $n$, then $W_t$ is an orientable pre-Wicks form with $\textnormal{genus}^+(W_t)=n$. If $g$ has non-orientable genus $n$, then $W_t$ is either a non-orientable pre-Wicks form of non-orientable genus $n$ or an orientable pre-Wicks form of orientable genus $(n-1)/2$. The latter is possible only when $n$ is odd and $\textnormal{genus}^+(g)=(n-1)/2$.
\end{corollary}
\begin{proof}
    Due to \Cref{def: tracking pair} and \Cref{lemma: moveproperties}, it suffices to show that there is an initial pair $(W_0,\mathfrak{g}_0)$ where $W_0$ is a genus $n$ pre-Wicks form and $\mathfrak{g}_0$ is a $G$-substitution such that $\mathfrak{g}_0(W_0)=g$. But this is clear. Suppose $g$ has orientable genus $n$. We can write $g$ as a product of $n$ non-trivial commutators, $[a_1,b_1],\cdots,[a_n,b_n]$. Thus, we have an orientable Wicks form of genus $n$, namely $W_0\colonequals [A_1,B_1]\cdots[A_n,B_n]$ together with the substitution $\mathfrak{g}_0$ defined by: $g_{A_i}=a_i, g_{B_i}=b_i$ for all $1\leq i\leq n$.
    If $g$ has non-orientable genus $n$, we can write $g$ as a product of $n$ non-trivial squares, $e_1^2,\dots, e_n^2$. Thus we have a non-orientable Wicks form of genus $n$, namely $U_0=E_1^2\dots E_n^2$, with the substitution $\mathfrak{h}_0$ defined by: $h_{E_i}\coloneq e_i$ for all $1\leq i\leq n$.
\end{proof}

\begin{remark}
    \leavevmode
    \begin{enumerate}[label=(\roman*)]
        \item From now on, whenever we say $g$ is an element of orientable genus $n$ and $(W,\mathfrak{g})$ is a tuple so that $\mathfrak{g}(W)=g$, it will be understood that $W$ is an orientable pre-Wicks form of orientable genus $n$. We mean similarly, when $g$ has non-orientable genus $n$, unless otherwise specified.
        \item A tracking pair $(W_t,\mathfrak{g}_t)$ depends on the order in which the moves are applied and on the initial tuple from which it is obtained. But for our purposes, this does not matter. The entire discussion in \Cref{subsec: properties of the tracking pair} applies to \emph{any} tracking pair representative. But, for the sake of convenience, the reader can fix the tracking pair to be the most efficient one, that is, $(W_t,\mathfrak{g}_t)$, such that both $l(W_t)$ and $l(\mathfrak{g}_t)$ are minimal.
    \end{enumerate}  
\end{remark}

\subsection{Properties of the tracking pair}\label{subsec: properties of the tracking pair}
Throughout this subsection, we set $D^{\pm}(n)=4l_0^{\pm}(n)+M+2$, where $M$ is the constant appearing in \Cref{lemma: conjugationbound}. All observations made in this subsection are about a tracking pair $(W_t,\mathfrak{g}_t)$ for an element $g$ of orientable/non-orientable genus $n$. In what follows, we will not explicitly mention the orientation, and $D=D^+(n)$, $l_0=l_0^+(n)$ in the orientable case, while $D=D^-(n)$, $l_0=l_0^-(n)$ in the non-orientable case.
\begin{definition}\label{def: longedges}
    Given an element $g$ of genus $n$ and a tracking pair $(W_t,\mathfrak{g}_t)$, a letter (or an unknown) $X$ that appears in $W_t$ will be said to be \emph{long} if $l(g_{X})\geq D$, otherwise we call it \emph{short}. Moreover, for any $c>0$, given a long letter $X$, and a choice of $\gamma(\mathfrak{g}_t)$, a vertex $p$ of $\gamma_X$ will be called a \emph{$c$-interior} vertex if it is at least a distance $c$ away from its initial and and terminal vertices. The subsegment comprising of the $c$-interior points of $\gamma_X$ will be denoted by $\gamma_X^{\text{c-int}}$, or just $\gamma_X^{\text{int}}$, when the value of $c$ is understood in context.
\end{definition}
Our moves have been designed so that $(W_t,\mathfrak{g}_t)$ $(3l_0+M+2,D,l_0)$-tracks $\gamma_g$. We will prove this shortly. 

\begin{lemma}\label{lemma: longedgepairing}
Let $g$ be an element of genus $n$ and let $(W_t,\mathfrak{g}_t)$ be a tracking pair for $g$, where $W_t=X_1\dots X_{2m}$. Let $\gamma_g$ be a geodesic joining the identity element to $g$ and consider any choice of $\gamma(\mathfrak{g_t})$. Let $X_j$ be a long letter in $W_t$ and consider any $l_0$-interior vertex $p$ of $X_j$. Consider any vertex point $q$ of $\gamma(\mathfrak{g}_t)\cup\gamma_g\backslash\gamma_{X_j}$ such that $p$ is $l_0$-close to $q$. Then: 
\begin{enumerate}[label=(\alph*), font=\normalfont]
    \item If $\text{sign}(X_j)=-1$, $q$ is either a point of $X_{\bar{j}}$ or it is a point of $\gamma_g$.
    \item If $\text{sign}(X_j)=1$, it is a point of $\gamma_g$.
\end{enumerate}
\end{lemma}
\begin{proof}
    To begin with, note that by \Cref{lem: length of a preWicks form}, $l(W_t)\leq 12n-2$ in the orientable case and $l(W_t)\leq 6n-2$ in the non-orientable case. Consequently, $\gamma_{X_j}$ is contained in the $l_0$-neighborhood of $\gamma(\mathfrak{g}_t)\cup\gamma_g\backslash\gamma_{X_j}$. Therefore, when $\text{sign}(X_j)=-1$, (a) follows from the fact that \hyperlink{Move 1}{Move 1} cannot be applied to the tracking pair. When $\text{sign}(X_j)=1$, (b) follows from the fact that \hyperlink{Move 1}{Move 1} and \hyperlink{Move 3}{Move 3} cannot be applied to the tracking pair.
\end{proof}

Informally, Lemma~\ref{lemma: longedgepairing} says that a long letter in a tracking pair can only fellow-travel with its paired letter (in case it has positive sign),  or a geodesic representing $g$. Now, we combine this observation with Lemma~\ref{lem: polygonpairing} to write down precisely the way in which such fellow-traveling happens. 

\begin{figure}[h]
    \centering
    \includegraphics[width=0.55\columnwidth]{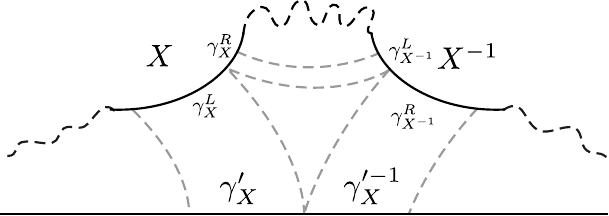}
    \caption{Parallel edge fellow travelling for a tracking pair}
    \label{fig: tracking for tracking pair}
\end{figure}

\begin{lemma}\label{lem: longedgeparallelfellow-travelling}
    Let $g$ be an element of orientable/non-orientable genus $n$ and let $(W_t,\mathfrak{g}_t)$ be a tracking pair for $g$. Consider the polygon $\gamma(\mathfrak{g}_t)\cup\gamma_g$. By Lemma~\ref{lem: polygonpairing}, there is a decomposition of each side of this polygon into subsegments, and moreover there is a pairing on the full collection of these subsegments such that paired subsegments $l_0$-fellow-travel (recall the length bound in \Cref{lem: length of a preWicks form}).
    Suppose $X_j$ is a long letter of $W_t$ and restrict the decomposition to $\gamma_{X_j}^{\text{int}}$ and $\gamma_{X_{\bar{j}}}^{\text{int}}$. Without any loss of generality, we assume that the letter $X_j$ appears before $X_{\bar{j]}}$ in $W_t$. 
    If $\text{sign}(X_j)=-1$, either of the following is true:
    \begin{enumerate}
        \item[\textbf{Case 1.}]\linkdest{longedgeparallelfellowtraveling case 1} \emph{There is fellow-traveling between subsegments of $\gamma^{\text{int}}_{X_j}$ and $\gamma^{\text{int}}_{X_{\bar{j}}}$}.  \\$\gamma_{X_j}^{\text{int}}=\gamma_{X_j}^L\gamma_{X_j}^R$, $\gamma_{X_{\bar{j}}}^{\text{int}}=\gamma_{X_{\bar{j}}}^L\gamma_{X_{\bar{j}}}^R$, where $\gamma_{X_j}^L$, $\gamma_{X_j}^R$, $\gamma_{X_{\bar{j}}}^L$, $\gamma_{X_{\bar{j}}}^R$ are subsegments coming from the aforementioned decomposition. There are subsegments $\gamma'_{X_j}$, $\gamma'_{X_{\bar{j}}}$ of $\gamma_g$, again coming from the decomposition, so that (See \Cref{fig: tracking for tracking pair}):
        \begin{enumerate}[font=\normalfont]
            \item $\gamma_{X_j}^R$ $l_0$-fellow-travels with a subsegment of $\gamma_{X_{\bar{j}}}$ and $\gamma_{X_{\bar{j}}}^L$ $l_0$-fellow-travels with a subsegment of $\gamma_{X_j}$, and:
            \begin{align*}
                l(\gamma^R_{X_j}), l(\gamma^L_{X_{\bar{j}}})\leq 2l_0 + M+2
            \end{align*}
            \item $\gamma_{X_j}^L$ and $\gamma'_{X_j}$ $l_0$-fellow-travel, $\gamma_{X_{\bar{j}}}^R$ and $\gamma'_{X_{\bar{j}}}$ $l_0$-fellow-travel.
        \end{enumerate}
        \item[\textbf{Case 2.}]\linkdest{longedgeparallelfellowtraveling case 2} \emph{There is no fellow-traveling between subsegments of $\gamma^{\text{int}}_{X_j}$ and $\gamma^{\text{int}}_{X_{\bar{j}}}$}. \\There are subsegments $\gamma'_{X_j}$, $\gamma'_{X_{\bar{j}}}$ of $\gamma_g$ so that:
        \begin{enumerate}[font=\normalfont]
            \item $\gamma_{X_j}^{\text{int}}$ and $\gamma'_{X_j}$ $l_0$-fellow-travel, $\gamma_{X_{\bar{j}}}^{\text{int}}$ and $\gamma'_{X_{\bar{j}}}$ $l_0$-fellow-travel.
        \end{enumerate}
    \end{enumerate}
    If $\text{sign}(X_j)=1$, then:
    \begin{enumerate}
        \item[\textbf{Case 3.}]\linkdest{longedgeparallelfellowtraveling case 3} $\gamma_{X_j}^{\text{int}}$ and $\gamma_{X_{\bar{j}}}^{\text{int}}$ are paired with subsegments $\gamma'_{X_j}$ and $\gamma'_{X_{\bar{j}}}$ of $\gamma_g$, and therefore $l_0$-fellow-travel with them.
    \end{enumerate}
    The collection of the interiors of the subsegments $\gamma'_{X_j}$ and $\gamma'_{X_{\bar{j}}}$ for all long letters $X_j$, is mutually disjoint. 
\end{lemma}

\begin{proof}
Suppose $\text{sign}(X_j)=-1$. Without any loss of generality, suppose $X_j=X$ and $X_{\bar{j}}=X^{-1}$ for some unknown $X$. To begin with, restrict the side decomposition to $\gamma_X^{\text{int}}$ and $\gamma_{X^{-1}}^{\text{int}}$.
Due to part $(a)$ of \Cref{lemma: longedgepairing}, the subsegments of $\gamma_X^{\text{int}}$ come in exactly two types: ones that $l_0$-fellow-travel with their paired counterparts in $\gamma_{X^{-1}}$ or ones that $l_0$-fellow-travel with their paired counterparts in $\gamma_{g}$. In what follows, whenever we say subsegment we always mean a subsegment coming from the aforementioned decomposition. 
\\Consider the rightmost subsegment of $\gamma_X^{\text{int}}$, that is the subsegment containing $\tau(\gamma_X^{\text{int}})$. The rest of the proof splits into two cases. Suppose this subsegment fellow travels with $\gamma_g$. 
    \begin{claim}\label{Claim: parallelfellow-travel2}
        All of $\gamma_{X}^{\text{int}}$ $l_0$ fellow-travels with a subsegment of $\gamma_g$.
    \end{claim}
    \begin{proof}
        Otherwise, we can write $\gamma_{X}^{\text{int}}$ as $\gamma^L_X\gamma_X^R$ where $\gamma_X^R$ is a (maximal) subsegment of $\gamma_X^{\text{int}}$ containing $\tau(\gamma_{X}^{\text{int}})$ which $l_0$-fellow-travels with a subsegment $\gamma'_X$ of $\gamma_g$. Removing the interiors of $\gamma_{X}^R\cup\gamma'_X$ from the polygon $\gamma(\mathfrak{g}_t)\cup \gamma_g$ gives us two components. $\gamma_X^L$ and $\gamma_{X^{-1}}$ belong to different components, so by \Cref{lem: polygonpairing} no subsegment of $\gamma_{X}^L$ can $l_0^+$-fellow-travel with a subsegment of $\gamma_X^{-1}$. We have arrived at a contradiction and so, indeed, all of $\gamma_{X}^{\text{int}}$ $l_0$ fellow-travels with a subsegment of $\gamma_g$, which we denote as $\gamma'_{X}$.
    \end{proof}
    Similarly, consider the leftmost subsegment of $\gamma_{X^{-1}}^{\text{int}}$, that is, the subsegment containing $i(\gamma_{X^{-1}}^{\text{int}})$. If it $l_0$ fellow-travels with $\gamma_g$, then just as we argued above, all of $\gamma_{X^{-1}}^{\text{int}}$ $l_0$ fellow-travels with a subsegment of $\gamma_g$, which we denote as $\gamma'_{X^{-1}}$. This proves \linkref{longedgeparallelfellowtraveling case 2}{Case 2} in the lemma. However, it is possible that the initial subsegment of $\gamma_{X^{-1}}^{\text{int}}$ is paired with a subsegment of $\gamma_X$, in which case, the latter must be a subsegment of $[\tau(\gamma_{X}^{\text{int}}),\tau(\gamma_{X})]$, and thus have length at most $l_0$. Therefore, $\gamma_{X^{-1}}^{\text{int}}=\gamma_{X^{-1}}^L\gamma_{X^{-1}}^R$, where $\gamma_{X^{-1}}^R$ $l_0$ fellow-travels with a subsegment of $\gamma_g$. 

Now suppose that the rightmost subsegment of $\gamma_X^{\text{int}}$, say $\gamma_X^R$ is paired with a subsegment of $\gamma_{X^{-1}}$, say $\beta$. 
\begin{itemize}
    \item If $\beta\subset[i(\gamma_{X^{-1}}), i(\gamma_{X^{-1}}^{\text{int}})]$, then we have $l(\gamma_X^R)\leq l_0$, and $\gamma_X^{\text{int}}=\gamma_X^L\gamma_X^R$, where $\gamma_X^L$ is paired with a subsegment of $\gamma_g$. In this case, all of $\gamma_{X^{-1}}^{\text{int}}$ is paired with a subsegment of $\gamma_g$ and therefore $l_0$ fellow-travels with it.
    \item Suppose $\beta\not\subset [i(\gamma_{X^{-1}}), i(\gamma_{X^{-1}}^{\text{int}})]$. In this case, we have a subsegment $\alpha$ of $\gamma_X$, which is paired with $\beta$. We can decompose $\gamma^{\text{int}}_X=\gamma_X^L\gamma_X^R$ and $\gamma^{\text{int}}_{X^{-1}}=\gamma_{X^{-1}}^L\gamma_{X^{-1}}^R$, where $\gamma_{X}^R=\alpha\cap\gamma_{X}^{\text{int}}$ and $\gamma_{X^{-1}}^L=\beta\cap\gamma_{X^{-1}}^{\text{int}}$. As argued before, \Cref{lem: polygonpairing} implies that $\gamma_{X}^L$ and $\gamma_{X^{-1}}^R$ are paired with subsegments of $\gamma_g$. We will now show that $l(\alpha), l(\beta)\leq 2l_0^++M+2$.
    
    Assume that the endpoints of $\alpha$ are vertices (there is an additive error of at most $2$ in doing so by possibly increasing the length of $\alpha$). Let $\tilde{\alpha}$ be the matching subsegment of $\alpha$ in $\gamma_{X^{-1}}$. More precisely, it is characterized by the following: $d(i(\tilde{\alpha}),i(\gamma_{X^{-1}}))=d(\tau(\alpha),\tau(\gamma_X))$ and $d(\tau(\tilde{\alpha}),i(\gamma_{X^{-1}}))=d(i(\alpha),\tau(\gamma_X))$. Because $\alpha$ and $\beta$ $l_0$ fellow-travel, and $d(\tau(\gamma_X^{\text{int}}),\tau(\gamma_X))=d(i(\gamma_{X^{-1}}),i(\gamma_{X^{-1}}^{\text{int}}))=l_0$, we have: 
    \begin{align*}
        d(i(\alpha),\tau(\tilde{\alpha}))\leq 2l_0,\quad d(\tau(\alpha),i(\tilde{\alpha}))\leq 2l_0 
    \end{align*}
    The fact that we cannot apply \hyperlink{Move 2}{Move 2} to $(W_t,\mathfrak{g}_t)$ implies that $l(\alpha)\leq 2l_0+M$. Similarly, we can argue that $l(\beta)\leq 2l_0+M$. This proves \linkref{longedgeparallelfellowtraveling case 1}{Case 1}.
    Suppose $\text{sign}(X_j)=1$. Restrict the side decomposition to $\gamma_{X_j}^{\text{int}}$ and $\gamma_{X_{\bar{j}}}^{\text{int}}$. Due to part $(b)$ of \Cref{lemma: longedgepairing}, all of $\gamma_{X_j}^{\text{int}}$ is paired with some subsegment $\gamma_{X_j}$ of $\gamma_g$, and all of $\gamma_{X_{\bar{j}}}^{\text{int}}$ is paired with some subsegment $\gamma'_{X_{\bar{j}}}$ of $\gamma_g$. This proves \linkref{longedgeparallelfellowtraveling case 3}{Case 3}. The fact that the collection of the interiors of the subsegments $\gamma'_{X_j}$ and $\gamma'_{X_{\bar{j}}}$ for all long letters $X_j$ is mutually disjoint follows directly from \Cref{lem: polygonpairing}.
\end{itemize}   
\end{proof}

\begin{corollary}\label{prop: trackingpairtracks}
    Suppose $g$ is an element of genus $n$,  $(W_t,\mathfrak{g}_t)$ is a tracking pair for $g$. Let $\gamma_g$ be a geodesic joining the identity element to $g$. Then, $\gamma(\mathfrak{g}_t)$ $(3l_0+M+2,D,l_0)$-tracks $\gamma_g$. 
\end{corollary}
\begin{proof}
    This readily follows from \Cref{lem: longedgeparallelfellow-travelling}.  
\end{proof}

\section{Element count: Upper bounds}\label{sec: Element count: Upper bounds}

\subsection{Counting in free groups}\label{subsec: Element count: Upper bounds: Free groups}

Fix a positive integer $r\geq 2$ and let $F_r$ denote the free group of rank $r$, equipped with the standard symmetric generating set of $2r$ elements. 

\begin{lemma}[\cite{culler1981using},\cite{ComerfordEdmunds1994}]\label{lem: canonical forms for group elements}
    Let $n$ be a positive integer. 
    \begin{enumerate}[font=\normalfont]
        \item If $g\in [G,G]$ has orientable genus $n$, then there exists an orientable pre-Wicks form $W$ of orientable genus $n$ and a non-canceling substitution $\mathfrak{g}$ such that $\mathfrak{g}(W)=g$.
        \item If $g\in G^2$ has non-orientable genus $n$, then either of the following holds true:
        \begin{enumerate}[font=\normalfont]
            \item There exists a non-orientable pre-Wicks form $W$ of non-orientable genus $n$ and a non-canceling substitution $\mathfrak{g}$ such that $\mathfrak{g}(W)=g$.
            \item there exists an orientable pre-Wicks form $W'$ of orientable genus $(n-1)/2$ and a non-canceling substitution $\mathfrak{g}'$ such that $\mathfrak{g}'(W')=g$.
        \end{enumerate}
        The latter is possible only when $n$ is odd and $\textnormal{genus}^+(g)=(n-1)/2$.
    \end{enumerate}
\end{lemma}
\begin{proof}
    In either case, we set $W=W_t(\coloneq X_1\cdots X_{2m})$ and $\mathfrak{g}_t$, where $(W_t,\mathfrak{g}_t)$ is the tracking pair associated with $g$ (\Cref{corollary: terminalwordexistence}). It remains to argue that the substitution $\mathfrak{g}_t$ is non-canceling. Consider two consecutive letters $X_i$ and $X_{i+1}$. The fact that \linkref{Move 1}{Move 1} and \linkref{Move 3}{Move 3} cannot be applied to $(W_t,\mathfrak{g}_t)$ implies that there cannot be any cancellation between $g_{X_{i}}$ and $g_{X_{i+1}}$.
\end{proof}

It follows from \Cref{lem: canonical forms for group elements} that the set $\partial B^{\pm}_n(L)$ is empty when $L$ is odd. Therefore in the rest of this subsection, the parameter $L$ is always assumed to be even.

\begin{definition}\label{def: setofnon-cancellingsubs}
    Let $W$ be a pre-Wicks form and let $\vec{l}$ be a length type for $W$. Let $G_W(\vec{l})$ be the set of group elements coming from non-cancelling substitutions $\mathfrak{g}$ for $W$ of length type $\vec{l}$, that is:
    \[
        G_W(\vec{l})\coloneq\{\mathfrak{g}(W)|\mathfrak{g}\text{ is a non-cancelling substitution for \(W\) of length type }\vec{l}\}
    \]
\end{definition}

Let us note that the number of reduced words in $F_r$ of length at most $m$ is $2r(2r-1)^{m-1}$ and the volume entropy $h$ equals $\log(2r-1)$. For the rest of this section, $h_r\coloneq \log(2r-1)$. 

\begin{definition}\label{def: the numbers tau}
    Given any pre-Wicks form $W=X_1\cdots X_{2m}$, define $\tau_W=\tau_W(r)$ to be the number of distinct non-canceling choices of starting and ending letters of $X_1,\dots, X_{2m}$. That is, $\tau_W$ is the number of distinct choices of $\{(s_j, t_j)\}_{j=1}^{2m}$, with each $s_j, t_j \in S$ such that:
         \begin{enumerate}
         \item For all $1\leq j\leq 2m$: 
         \begin{align*}
             (s_{\overline{j}},t_{\overline{j}})=\begin{cases}
             (t_j^{-1},s_j^{-1}) &\text{if sign}(X_j)=-1\\
             (s_j,t_j) &\text{if sign}(X_j)=+1
             \end{cases}
         \end{align*}
        \item For all $1\leq j\leq 2m-1$, $t_j \cdot s_{j+1}\neq 1$.
    \end{enumerate}
    Given any positive integers $n\geq 1$ and $r\geq 2$, define:
    \begin{align*}
        \tau^+_n&=\tau^+_n(r)\coloneq 2r(2r(2r-1)(2r-2))^{4n-1}\\ \tau^-_n&=\tau^-_n(r)\coloneq 2r(2r(2r-1)(2r-2))^{2n-1} 
    \end{align*}
\end{definition}

\begin{lemma}\label{lem: tau for maximal pre-Wicks forms}
    For all maximal orientable pre-Wicks forms $W$ of genus $n$, $\tau_W=\tau^+_n$. For all maximal non-orientable pre-Wicks forms $W$ of genus $n$, $\tau_W=\tau^-_n$. 
\end{lemma}
% \begin{proof}
    
% \end{proof}

Given $s,t\in S$, we define $B(s,t,L)$ to be the set of reduced words in $F_r$ of length $L$ that begin in $s$ and end in $t$. Our next claim says that for any choice of \(s\) and \(t\), the size of \(B(s, t, L)\) is approximately the same.

% {\color{gray}To check: is \cite{park2019conjugacy} the primary reference for the following lemma?}

\begin{lemma}\cite[Lemma 2.1]{park2019conjugacy}\label{lem: countingpinnedreducedwords}
    There exist constants $C_1,C_2$ depending only on $r$ such that, for all $L\geq 2$ and for all $s,t\in S$, we have:
    \begin{equation*}
        \frac{(2r - 1)^{L-1}}{2r} \left(1 + \frac{C_1}{r^{L-2}}\right) \leq |B(s,t,L)| \leq \frac{(2r - 1)^{L-1}}{2r} \left(1 + \frac{C_2}{r^{L-2}}\right)
    \end{equation*}
\end{lemma}

The next claim follows readily from \Cref{lem: countingpinnedreducedwords}, and from the fact that for free groups, any two distinct substitutions of \(W\) with the same underlying length type \(\vec{l}\) always give rise to distinct elements.

\begin{lemma}\label{lem: non-cancellingcount1}
    Let $W$ be a pre-Wicks form of length $2m$. Let $\vec{l}$ be a length type for $W$. There are constants $C'_1,C'_2>0$ depending only on $r$ and $m$, such that:
    \begin{align*}
        C'_1e^{\frac{h_r||\vec{l}||}{2}}\leq |G_W(\vec{l})|\leq C'_2e^{\frac{h_r||\vec{l}||}{2}}
    \end{align*}

\end{lemma}

We can strengthen the above claim as follows:

\begin{lemma}\label{lem: substitution count asymptote}
    Fix \(k > 0\). Let $W$ be a pre-Wicks form of genus $n$ and let $\vec{l}$ be a length type for $W$ such that each entry of \(\vec{l}\) is at least \(k \log L \). Let $C_1$, $C_2$ be the constants in \Cref{lem: countingpinnedreducedwords}. Then:
    \begin{equation*}
        \left (1+ \frac{C_1}{r^{k\log L-2}} \right)^m\leq \frac{|G_W(\vec{l})|}{\displaystyle{\tau_W \cdot\frac{(2r - 1)^{\frac{||\vec{l}||}{2} - m}}{\displaystyle (2r)^m}}}\leq \left(1+\frac{C_2}{r^{k\log L-2}}\right)^{m}
    \end{equation*}
\end{lemma}
\begin{proof}
     Let $W=X_1\cdots X_{2m}$. Consider a collection \(\{(s_j, t_j)\}_{j=1}^{2m}\), with each \(s_j, t_j \in S\) such that:
      \begin{enumerate}
         \item For all $1\leq j\leq 2m$: 
         \begin{align*}
             (s_{\overline{j}},t_{\overline{j}})=\begin{cases}
             (t_j^{-1},s_j^{-1}) &\text{if sign}(X_j)=-1\\
             (s_j,t_j) &\text{if sign}(X_j)=+1
             \end{cases}
         \end{align*}
        \item For all $1\leq j\leq 2m-1$, $t_j \cdot s_{j+1}\neq 1$.
    \end{enumerate}
    Since each entry of \(\vec{l}\) is at least 2, every non-cancelling substitution of length type \(\vec{l}\), that is, every element of \(G_W(\vec{l})\), corresponds to such a non-cancelling collection of starting and ending letters.
    
    Given a choice \(\{(s_j, t_j)\}_{j=1}^{2m}\) as above, there are \(\prod_j B(s_j, t_j, l_j)\) many elements of \(G_W(\vec{l})\) restricting to this choice. By applying \Cref{lem: countingpinnedreducedwords} to each \(B(s_j, t_j, l_j)\) in the product and then summing over the \(\tau\) many possibilities of \(\{(s_j, t_j)\}_{j=1}^{2m}\), we have:
    \begin{equation}\label{eq: bound for substitution count}
        \tau \cdot\frac{(2r - 1)^{L - m}}{\displaystyle (2r)^m} \prod_j \left(1 + \frac{C_1}{r^{l_j - 2}}\right) 
        \leq |G'(\vec{l})| \leq 
        \tau \cdot\frac{(2r - 1)^{L - m}}{\displaystyle (2r)^m} \prod_j \left(1 + \frac{C_2}{r^{l_j - 2}}\right)
    \end{equation}
    Now the lema follows from the fact that $k\log L <l_j<\frac{L}{2}$ for each $1\leq j\leq 2m$.
\end{proof}

For every positive integer $L$, let $\partial B_n^{+}(L)$ denote the set of group elements $g\in [G,G]$ having orientable genus $n$ and length equal to $L$. Let $\partial B_n^{-}(L)$ denote the set of group elements $g\in G^2$ having non-orientable genus $n$ and length equal to $L$. 

\begin{proposition}\label{thm: free grp shell upper bound asymptote}
    Let $n$ and $r$ be positive integers with $r\geq 2$.
    \begin{equation*}
        \limsup_{\substack{L\to\infty\\ L\textnormal{ even}}}\frac{\partial B^{+}_n(L)}{L^{6n-2}e^{\frac{h_rL}{2}}}\leq \frac{\tau^+_n(r) \mathcal{N}_{\text{pre}}^{+}(n)}{2^{6n-2}(6n-2)!(2r)^{6n-1}(2r-1)^{6n-1}}
    \end{equation*}
\end{proposition}
\begin{proof}
    Given any pre-Wicks form $W$ of length $2m$ and $L>0$ (even) we define two families of length types for $W$ as follows:
    \begin{align*}
        &\mathcal{L}_1(W,m,r,L)\coloneq \left\{\vec{l}\in\mathbb{Z}_{\geq 0}^{m} \middle |\; \sum_{i=1}^{m}l_i= \frac{L}{2}, l_i>2\;(1\leq i\leq m)\right\}\\
        &\mathcal{L}_2(W,m,r,L)\coloneq \left\{\vec{l}\in\mathbb{Z}_{\geq 0}^{m} \middle |\; \sum_{i=1}^{m}l_i= \frac{L}{2},\, l_i\geq \log L \; (1\leq i\leq m) \right\}
    \end{align*}
    We will suppress the arguments of $\mathcal{L}_1,\mathcal{L}_2$ below whenever they are clearly understood.
    \begin{claim}\label{claim: claim in free grp upper bound asymptote proof}
    \leavevmode
        \begin{enumerate}
            \item For any non-maximal orientable pre-Wicks form $W$, we have:
            \begin{equation*}
                \lim_{\substack{L\to\infty\\ L\textnormal{ even}}} \frac{\sum_{\vec{l}\in \mathcal{L}_1} |G_W(\vec{l}\,)|}{L^{6n-2}e^{\frac{h_rL}{2}}} =0 
            \end{equation*}
            \item For any maximal orientable pre-Wicks form $W$:
            \begin{equation*}
                \lim_{\substack{L\to\infty\\ L\textnormal{ even}}}\frac{\sum_{\vec{l}\in\mathcal{L}_1-\mathcal{L}_2}|G_W(\vec{l})|}{L^{6n-2}e^{\frac{h_rL}{2}}}=0
            \end{equation*}
            \item For any maximal orientable pre-Wicks form $W$:
            \begin{equation*}
                \lim_{\substack{L\to\infty\\ L\textnormal{ even}}}\frac{\sum_{\vec{l}\in \mathcal{L}_2}|G_{W}(\vec{l})|}{L^{6n-2}e^{\frac{h_rL}{2}}}=\frac{\tau_n^+(r)}{2^{6n-2}(6n-2)!(2r)^{6n-1}(2r-1)^{6n-1}}
            \end{equation*}
            where the number $\tau^+_n(r)$ is the constant defined in \Cref{def: the numbers tau}.
        \end{enumerate}
    \end{claim}
    Before we prove the claim, note that the proposition follows from the claim, the fact that (\Cref{lem: canonical forms for group elements}):
    \[\partial B^+_n(L)\subset\bigcup_{W}\bigcup_{\vec{l}\in \mathcal{L}_1(W,L)}G_W(\vec{l}\,)\]
    and the fact that for any two maximal orientable pre-Wicks forms of the same genus, $|\mathcal{L}_2(W,n,r,L)|=|\mathcal{L}_2(W',n,r,L)|$.
    \begin{proof}[Proof of \Cref{claim: claim in free grp upper bound asymptote proof}]
         
    Let $C_2$ be the constant in \Cref{lem: countingpinnedreducedwords}. It follows from \Cref{lem: countingpinnedreducedwords} and \Cref{def: the numbers tau} that for all $\vec{l}\in\mathcal{L}_1(W,m,r,L)$:
    \begin{equation}\label{eq: free group substitution upper bound}
         |G_W(\vec{l})|\leq \frac{\tau_W}{(2r)^m(2r-1)^m}\left(1+\frac{C_2}{r}\right)^m \times e^{\frac{h_rL}{2}}
    \end{equation}
    
    \begin{enumerate}
        \item Let $W$ be any non-maximal orientable pre-Wicks form of length $2m$ and genus $n$. It follows from \Cref{eq: free group substitution upper bound} and the fact that $\tau_W< (2r)^{2m}$ that:
        \begin{align*}
            \sum_{\vec{l}\in\mathcal{L}_1(W,m,L)}|G_W(\vec{l})|&\leq |\mathcal{L}_1(W,m,L)|\times \left(\frac{2r}{2r-1}\right)^m\left(1+\frac{C_2}{r}\right)^m \times e^{\frac{h_rL}{2}}\\
            &\leq \left(\frac{2r}{2r-1}\right)^m\left(1+\frac{C_2}{r}\right)^m \times \frac{L^{m-1}}{(m-1)!}\times e^{\frac{h_rL}{2}}
        \end{align*}
    Because $W$ is non-maximal, $m<6n-1$, and $(1)$ follows from the above inequality.
    \item Let $W$ be any orientable pre-Wicks form of length $2m$ and genus $n$. Note that:
    \begin{align*}
        &\mathcal{L}_1(W,m,L)-\mathcal{L}_2(W,m,L)=\bigcup_{j=1}^m \left \{ \vec{l}\in  \mathbb{Z}_{\geq 0}^m \middle | \; \sum_{i=1}^m l_i = \frac{L}{2}, l_j\leq \log L, l_i>2 \;(1\leq i\leq m) \right \}\\
        &\left| \left \{ \vec{l}\in  \mathbb{Z}_{\geq 0}^m \middle | \; \sum_{i=1}^m l_i = \frac{L}{2}, l_j\leq \log L, l_i>2 \;(1\leq i\leq m) \right \} \right | \leq \frac{(\log L) L^{m-2}}{(m-2)!}
    \end{align*}
        Therefore, using the fact that $\tau_W<(2r)^{2m}$ and \Cref{eq: free group substitution upper bound}, we get:
        \begin{align*}
            \frac{1}{L^{6n-2}e^{\frac{h_rL}{2}}}\sum_{\vec{l}\in \mathcal{L}_1-\mathcal{L}_2}|G_W(\vec{l}\,)|\leq \frac{m}{(m-2)!}\left(\frac{2r}{2r-1}\right)^m\left(1+\frac{C_2}{r}\right)^m\times \frac{(\log L)}{L^{6n-m}}
        \end{align*}
        When $W$ is maximal, $m=6n-1$, and now $(2)$ follows. 
        \item Let $W$ be a maximal orientable pre-Wicks form of genus $n$. $(3)$ follows from: the fact that $\lim_{L\to\infty}\frac{|\mathcal{L}_2(W,n,r,L)|}{L^{6n-2}}=\frac{1}{2^{6n-2}(6n-2)!}$, the (uniform) bounds for $|G_W(\vec{l\,})|$ for each $\vec{l}\in\mathcal{L}_2(W,L)$ in \Cref{lem: substitution count asymptote}, and the fact that, $\tau_W=\tau^+_n(r)$ (\Cref{lem: tau for maximal pre-Wicks forms}). \qedhere
    \end{enumerate}
    \end{proof}
    This ends the proof of Proposition \ref{thm: free grp shell upper bound asymptote}.
\end{proof}
Similarly, we have:
\begin{proposition}\label{thm: free grp shell upper bound asymptote non-orientable}
    Let $n>1$ and $r\geq 2$ be positive integers.
    \begin{equation*}
        \limsup_{\substack{L\to\infty\\ L\textnormal{ even}}}\frac{\partial B^{-}_n(L)}{L^{3n-2}e^{\frac{h_rL}{2}}}\leq \frac{\tau^-_n(r) \mathcal{N}_{\text{pre}}^{-}(n)}{2^{3n-2}(3n-2)!(2r)^{3n-1}(2r-1)^{3n-1}}
    \end{equation*}
\end{proposition}
\begin{proof}
    The proof follows using the same arguments used in the proof of \Cref{thm: free grp shell upper bound asymptote}, with the estimates in \Cref{claim: claim in free grp upper bound asymptote proof} replaced by their non-orientable analogues. However there is one subtlety that needs to be accounted for. By \Cref{lem: canonical forms for group elements}, there are two types of group elements in $\partial B^-_n(L)$: those that arise from non-canceling substitutions of non-orientable pre-Wicks forms $W$ and $\textnormal{genus}^-(W)=n$ and those that arise from non-canceling substitutions of orientable pre-Wicks forms $W'$ with $\textnormal{genus}^+(W')=(n-1)/2$. However, asymptotically as $L\to\infty$ the second type of elements are rare. More precisely, the number of group elements of the second type is $O(L^{3n-5}e^{\frac{h_rL}{2}})$ and therefore do not contribute to the estimate in the statement of the proposition.
\end{proof}

\subsection{Counting in hyperbolic groups}\label{subsec: Element count: Upper bounds: Hyperbolic groups}
We use the constants $D^{\pm}=D^{\pm}(n)$ and $l_0^{\pm}=l_0^{\pm}(n)$ defined in \Cref{sec: wicks form moves} and \Cref{subsec: properties of the tracking pair}.
\begin{theorem}\label{thm: element count upper bound orientable hyperbolic}
    Let $G$ be a hyperbolic group equipped with the word metric induced by some finite generating set so that it is $\delta$-hyperbolic. Let $h$ denote the volume entropy of $G$. There exist positive constants $C=C(\delta,n)$ and $L_0=L_0(h,\delta,n)$ so that:
    \begin{equation*}
        \forall L\geq L_0, |B_n^+(L)|\leq CL^{6n-2}e^{hL/2}
    \end{equation*}
\end{theorem}
\begin{proof}
Consider any group element $g$ of genus $n$ with length at most $L$. By \Cref{corollary: terminalwordexistence}, there is a tracking pair (cf. \Cref{def: tracking pair}) $(W_t,\mathfrak{g}_t)$ for $g$. $W_t$ can be written as: $W_t=S_1L_1S_2\dots S_{m}L_{m}S_{m+1}$, where $L_1,\dots,L_m$ are the long letters of $W_t$ (so $m$ is necessarily even) and $S_1,\dots,S_{m+1}$ are products of short letters (allowed to be the empty word) of $W_t$. 
\begin{claim}\label{claim: upper bound on length of tracking pair substitution}
    \begin{equation}\label{eqn: sum of long edges}
        \sum_{i=1}^m l(g_{L_i})\leq L+mD^+\leq L+(12n-2)D^+
    \end{equation}
\end{claim}
\begin{proof}
    Let $\gamma_g$ be any geodesic joining the identity element to $g$ and consider any choice of $\gamma(\mathfrak{g}_t)$. 
    Recall \Cref{lem: longedgeparallelfellow-travelling} and note that for every long edge $L_i$, $\gamma_{L_i}^{\text{int}}$ can be written as $\gamma_{L_i}^L\gamma_{L_i}^R$ and there is a subsegment $\gamma'_{L_i}$ of $\gamma_g$ so that $\gamma_{L_i}$ $l_0^+$-fellow-travels with $\gamma'_{L_i}$ and $l(\gamma_{L_i}^R)\leq 2l_0^++M+2$. Thus, the claim follows from the following estimate and the fact that $m\leq 12n-2$ (see \Cref{lem: length of a preWicks form}):
    \begin{align*}
        l(\gamma_{L_i})=l(\gamma_{L_i}^{\text{int}})+2l_0^+=l(\gamma_{L_i}^L)+l(\gamma_{L_i}^R)+2l_0^+
        \leq l(\gamma'_{L_i})+(2l_0^++M+2)+2l_0^+
        = l(\gamma'_{L_i})+ D^{+}. \,\,\,\,\,\,\qedhere
    \end{align*}
\end{proof}
Before we begin our counting arguments, let us note that it follows from \Cref{thm: ballgrowth} applied for $G$ acting on its Cayley graph, that there exist constants $C_2, h>0$ so that: $|B(L)|\leq C_2 e^{hL}$, where $B(L)$ denotes the number of group elements of length at most $L$. Consider the following condition:
\begin{itemize}
    \item[\textbf{(C1)}] $W$ is a pre-Wicks form of orientable genus $n$, $W$ admits a decomposition of the form $W=S_1L_1S_2\dots S_mL_mS_{m+1}$ where $L_i (1\leq i\leq m)$ are long letters, $S_j=S_j^{(1)}\dots S_j^{(k_j)}$ where $S_j^{(k)} (1\leq j\leq m+1, 1\leq k\leq k_j)$ are short letters. Moreover, \Cref{eqn: sum of long edges} holds and:
    \begin{align*}
        m+\sum_{j=1}^{m+1}k_j\leq 12n-2
    \end{align*}
\end{itemize}
Given $4n \leq m'\leq 12n-2$, $0\leq m\leq m'$, and a pre-Wicks form $W$ of orientable genus $n$ with $l(W)=m'$, let $B'_n(W,L, m, m')$ be the set of substitutions $\mathfrak{g}$ for $W$ satisfying \textbf{(C1)} such that $W$ has exactly $m$ many long letters. It follows from \Cref{claim: upper bound on length of tracking pair substitution} that:
\begin{align*}
    |B_n^+(L)|&\leq \sum_{\substack{W\\\text{genus}(W)=n}}|B'_n(W,L,m,m')|
\end{align*}
Note that on the right hand side we are summing over all non-equivalent pre-Wicks forms $W$ with orientable genus $n$. Each term on the right hand side has the following upper bound:
\begin{align*}    
  |B''_n(W,L,m,m')|&\leq \binom{\frac{m'}{2}}{\frac{m}{2}} \times C_2^{\frac{m'-m}{2}}e^{\frac{hD^+(m'-m)}{2}}\times \sum_{\substack{(l_1,\dots,l_{m/2})}} C_2^{m/2}e^{h\sum_{i=1}^{m/2}l_i}
\end{align*}
The first factor comes from choosing long letters: $m/2$ many long unknowns from $m'/2$ many unknowns appearing in $W$. The second factor comes from choosing group elements which substitute $\frac{m'-m}{2}$ many short unknowns while the third factor comes from choosing group elements which substitute for the long unknowns. Note that in this factor, we are summing over solutions to 
\begin{align}\label{eqn: countingargumentlengthsumeqn}
    \sum_{i=1}^{m/2}l_i\leq \frac{L}{2}+\frac{mD^+}{2}\text{ and }l_i\geq D^+\text{ for all }1\leq i\leq \frac{m}{2}
\end{align}
But most of the contribution to the third factor comes from a  thin ``shell" for large values of $L$, that is, there is a constant $c'$ depending only on $h$ such that:
\begin{align}\label{eqn: boundparitionlength}
    \sum_{\substack{(l_1,\dots,l_{m/2})}} e^{h\sum_{i=1}^{m/2}l_i}\leq c'{\lceil \frac{L}{2}+\frac{mD^+}{2}\rceil\choose \frac{m}{2}-1}e^{\frac{hL}{2}+\frac{mhD^+}{2}}
\end{align}
Using \Cref{eqn: countingargumentlengthsumeqn} and \Cref{eqn: boundparitionlength}, we have:
\begin{align*}
    |B'_n(W,L,m,m')|\leq c'C_2^{\frac{m'}{2}}e^{\frac{hm'D^+}{2}} {\frac{m'}{2} \choose \frac{m}{2}} {\lceil \frac{L}{2}+\frac{m'D^+}{2}\rceil\choose \frac{m'}{2}-1}e^{\frac{hL}{2}}
\end{align*}
\noindent There exists $L_0=L_0(h,\delta,n)$ so that for $L\geq L_0$, the terms corresponding to maximal pre-Wicks forms $B''_n(W,L,12n-2,12n-2)$ dominate the rest. So for large $L$, we have:
\begin{align*}
    |B_n^+(L)|\leq \mathcal{N}_{\text{pre}}^+(n)c'C_2^{6n-1}e^{(6n-1)hD^+} {\lceil \frac{L}{2}+(6n-1)D^+\rceil\choose 6n-2}e^{\frac{hL}{2}}
\end{align*}
where $\mathcal{N}_{\text{pre}}^+(n)$ is the number of non-equivalent pre-Wicks forms of orientable genus $n$. If necessary, choose $L_0$ even larger so that, for $L\geq L_0$, $(12n-2)D^+<L$. Thus, for $L\geq L_0$, we have:
\begin{align*}
    |B_n^+(L)|\leq \frac{c'C_2^{6n-1}e^{(6n-1)hD^+}}{(6n-2)!} \mathcal{N}_{\text{pre}}^+(n)L^{6n-2}e^{hL/2 }. \quad \quad\qedhere
\end{align*}
\end{proof}

Let $B_n^-(L)$ denote the set of elements in $G$ of non-orientable genus $n$ having length at most $L$.
\begin{theorem}\label{thm: element count upper bound non-orientable hyperbolic}
    Let $G$ be a hyperbolic group equipped with the word metric induced by some finite generating set $S$ so that it is $\delta$-hyperbolic. Let $h$ denote the volume entropy of $G$. There exist positive constants $C=C(\delta,n)$ and $L_0=L_0(h,\delta,n)$ so that:
    \begin{equation}
        \forall L\geq L_0, |B_n^-(L)|\leq CL^{3n-2}e^{hL/2}
    \end{equation}
\end{theorem}
\begin{proof}
    The proof is essentially the same as that for the orientable genus case, mutatis mutandis. All occurrences of $l_0^+$ and $D^+$ are replaced with $l_0^-$ and $D^-$ respectively, and those of $n$ are replaced with $n/2$. There is only one nuance: there are two types of group elements in $ B^-_n(L)$ (cf. \Cref{corollary: terminalwordexistence}) - those that arise from a tracking pair where $W_t$ is a non-orientable pre-Wicks form with $\textnormal{genus}^-(W_t)=n$ and those that arise from a tracking pair where $W_t$ is an orientable pre-Wicks form with $\textnormal{genus}^+(W_t)=(n-1)/2$. However, the number of group elements of the second type is $O(L^{3n-5}e^{\frac{hL}{2}})$ and therefore for large $L$ their contribution is negligible in comparison to the first type.
\end{proof}

\section{Element count: Lower bounds}\label{sec: Element count: Lower bounds}

\begin{definition}\label{def: lengthtype}
    Let $W$ be a pre-Wicks form. A \textit{length type} for $W$ is an assignment of non-negative integer values to each unknown appearing in $W$. Equivalently, writing $W=X_1\dots X_{2m}$, a length type for $W$ is an assignment of non-negative integer values to the $X_i$'s such that paired letters get assigned the same value. That is, it is a $2m$-tuple $\vec{l}\in \mathbb{Z}_{\geq 0}^{2m}$, satisfying: for all $1\leq i,j\leq 2m$,  if $ j=\bar{i}$ then $l_j=l_{\bar{i}}$. Define the length of $\vec{l}$ to be its $l_1$-norm, that is, $\|\vec{l}\|\coloneq\sum_{i=1}^{2m}l_i$. A substitution $\mathfrak{g}$ is said to be of length type $\vec{l}$ if for every letter $X_j$, $l(g_{X_j})=l_j$.
\end{definition}
\begin{remark}\label{rmk: length type as an m-tuple}
    There will be times when we will think of a length type as an $m$-tuple. To do this, we will choose a total order on the set of unknowns appearing in $W$.
\end{remark}

\Cref{subsec: combinatorialreduction} contains some of the main combinatorial results that go into the lower bound proofs. The reader can choose to read through the basic definitions and the statement of \Cref{thm: the main reduction theorem}, skip \Cref{subsubsec: tools for the proof of the reduction theorem}  and \ref{subsubsec: Proof of the reduction theorem}, and proceed to element count lower bound proofs.

\subsection{Combinatorial Reduction}\label{subsec: combinatorialreduction}

Let $A_0$ and $H$ be non-negative integers.

\begin{definition}[Matching Positions]\label{def: matching position}
    Let $W=X_1\dots X_{2m}$ be a pre-Wicks form. Let $\vec{l}$ be a length type for $W$ and let $\vec{b}=(b_0,\dots,b_{2m})$ be a sequence of integers in $[A_0,A_0+H]$. Given this data, define $L(W,\vec{l},\vec{b})\coloneq ||\vec{l}||+||\vec{b}||$. By a \textit{position}, we mean an integer $1\leq i\leq L(W,\vec{l},\vec{b})$. The data $(W,\vec{l},\vec{b})$ determines a \emph{canonical partition} of the set of positions $\{1,\dots,L(W,\vec{l},\vec{b})\}$ corresponding to letters and junctions, defined as follows:

    \begin{figure}[h]
        \centering
        \includegraphics[width=0.9\columnwidth]{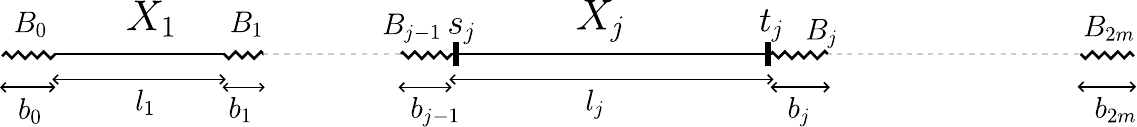}
        \caption{Length types, junctions, starting and ending positions}
        \label{fig: positions}
    \end{figure}
    
    For each $1\leq j\leq 2m$, the positions $\{\sum_{k=0}^{j-1}b_k+\sum_{k=1}^{j-1}l_k+1,\dots, \sum_{k=0}^{j-1}b_k+\sum_{k=1}^{j}l_k \}$ will be said to correspond to the letter $X_j$. For each $1\leq j\leq 2m-1$, the positions $\{\sum_{i=0}^{j-1}b_i+\sum_{k=1}^{j}l_k+1,\dots, \sum_{i=0}^{j}b_i+\sum_{k=1}^{j}l_k \}$ will be said to comprise the $j^{\text{th}}$ \emph{junction}, denoted by $B_j$. In addition, the initial \(b_0\) many positions and the terminal \(b_{2m}\) many positions of \(L(W, \vec{l}, \vec{b})\) are respectively the \(0\)-th and \(2m\)-th junctions. We will say that a letter $X_j$ or a junction $B_j$ \emph{contains} the position $i$ with respect to $(W, \vec{l}, \vec{b})$, if $i$ is a position belonging to the corresponding set of positions with respect to $\vec{l}$. 
    
    Accordingly, given any letter $X_j$, the \emph{starting} and \emph{ending positions} of $X_j$ refer to the smallest and largest position contained in $X_j$ respectively, and they will be denoted by $s_{X_j}$ and $t_{X_j}$ respectively. For any position $i$ contained in a letter \(X_j\) (equivalently, \(s_{X_j} \leq i \leq t_{X_j}\)), its \textit{matching position} $\underline{i}$, with respect to the length type $\vec{l}$, is defined as:
    \begin{equation}\label{eq: matching position 1}
        \underline{i}\coloneq\begin{cases}
            t_{X_{\bar{j}}} + s_{X_j} - i & \text{; sign}(X_j)=-1\\
            s_{X_{\overline{j}}} + i - s_{X_{j}} & \text{; sign}(X_j)=+1
        \end{cases}
    \end{equation} 
    (where $X_{\bar{j}}$ is the letter paired with $X_j$ in $W$). Similarly, given a sequence of consecutive positions, $\mathbf{i}=(a, \dots,a+(k-1))$, we will say that the letter $X_j$ (or the junction $B_j$) \emph{contains} $\mathbf{i}$ if all the positions comprising the sequence belong to $X_j$ (or $B_j$ respectively). Suppose the letter $X_j$ contains $\mathbf{i}$; we then define its sequence of matching positions $\underline{\mathbf{i}}$:
    \begin{equation*}
        \underline{\mathbf{i}}\coloneq\begin{cases}
            (\underline{a},\dots,\underline{a+(k-1)}) &\text{; sign}(X_j)=+1\\
            (\underline{a+(k-1)},\dots,\underline{a}) &\text{; sign}(X_j)=-1
        \end{cases} 
    \end{equation*}
    Note that when \(\mathbf{i}\) is contained in \(X_j\), \(\underline{\mathbf{i}}\) is contained in \(X_{\overline{j}}\) by definition.
\end{definition}

\begin{definition}[Reduction in degrees of freedom]\label{def: reduction in degrees of freedom}
     Let $W$ and $W'$ be pre-Wicks forms of lengths $2m$ and $2m'$ respectively. Let $\vec{l},\vec{l}'$ be two length types for $W$ and $W'$ respectively. Let $\vec{b}\in[A_0,A_0+H]^{2m+1}$ and $\vec{b}'\in[A_0,A_0+H]^{2m'+1}$. Suppose we have $L(W,\vec{l},\vec{b})=L(W',\vec{l}',\vec{b}')= L$. Then, any position \(i\) gives rise to two matching positions \(\underline{i}\) and \(\underline{i}'\) with respect to \(\vec{l}\) and \(\vec{l}'\) respectively. We will say that a sequence of positions $\mathbf{i}$ is \emph{admissible} if there are letters $X$ and $X'$ such that $\mathbf{i}$ is contained in $X$ with respect to the length type $\vec{l}$ and $\mathbf{i}$ is contained in $X'$ with respect to the length type $\vec{l}'$.
\begin{figure}[h]
    \centering
    \includegraphics[width=0.7\columnwidth]{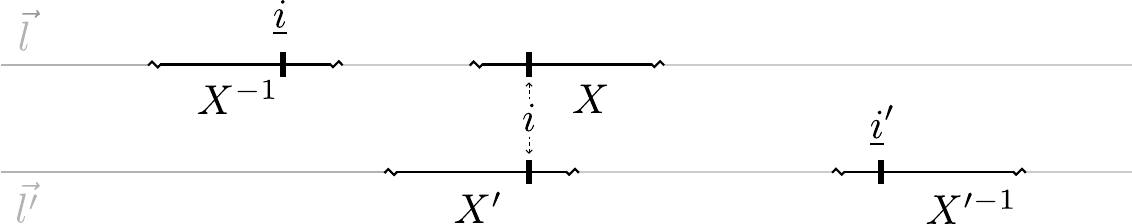}
    \caption{Matching positions \(\underline{i}\) and \(\underline{i}'\) when both $X$ and $X'$ have sign $-1$.}
    \label{fig: matching position}
\end{figure}
    Let $N$ and $D_0$ be positive integers. We will say that there is a $D_0$-\emph{reduction in the degrees of freedom by} $N$ if there is a sequence of $N$ consecutive admissible positions $\mathbf{i}=(i_1,\dots,i_N)\subset\{1,\dots,L\}$ such that both $\underline{\mathbf{i}}$ and $\underline{\mathbf{i}}'$ are admissible, and $|\underline{i_1}- \underline{i_1}'|\geq D_0$. In the simplest case of $D_0 = 1$, this condition is equivalent to $\underline{i_1} \neq \underline{i_1}'$. 
\end{definition}

The motivation for this notion comes from \Cref{lem: element intersection count for reduction in degrees of freedom} and \Cref{lem: hyperbolic element intersection count for reduction in degrees of freedom}. Moreover, the condition $|\underline{i_1}- \underline{i_1}'|\geq D_0$ is included in the previous definition without loss of generality: we soon prove in \Cref{lem: matching position extended} that given $\text{sign}(X)=\text{sign}(X')$, $|\underline{i_1}- \underline{i_1}'|\geq D_0$ is equivalent to $|\underline{i_k}- \underline{i_k}'|\geq D_0$ for any $1 \leq k \leq N$. We are now ready to state the main result of this section. 

\begin{proposition}\label{thm: the main reduction theorem}
    Let $c, A_0$,$D_0$, and $H$ be positive integers. Let $W$ and $W'$ be maximal pre-Wicks forms of lengths $2m$ and $2m'$ respectively. Let $\vec{l}$ and $ \vec{l}'$ be the length types for $W$ and $W'$ respectively. Furthermore, let $\vec{b}\in[A_0,A_0+H]^{2m+1}$ and $\vec{b}'\in[A_0,A_0+H]^{2m'+1}$. Now, suppose the following holds:
    \begin{enumerate}[font=\normalfont, label=(L\arabic*)]
        \item\label{item: L1} $L(W,\vec{l},\vec{b})=L(W',\vec{l}',\vec{b}')$
        \item\label{item: L2} 
         For all $1\leq i\leq 2m$, $1\leq j\leq 2m'$, $l_i,l'_j>6c+5(A_0+H)+4D_0$.
    \end{enumerate}
    Then, either:
    \begin{enumerate}[font=\normalfont]
        \item There is $D_0$-reduction in degrees of freedom by $c$, or
        \item $W=W'$, therefore \(m = m'\), and moreover for all $1\leq i\leq 2m$, $|l_i-l'_i|\leq 3(D_0-1)+4H$.
    \end{enumerate}
\end{proposition}

We note down the special case when $A_0=H=0$ and $D_0=1$ for future reference.

\begin{corollary}\label{cor: the main reduction theorem for free groups}
    Let $c\geq 1$ be a positive integer. Let $W$ and $W'$ be maximal pre-Wicks forms of lengths $2m$ and $2m'$ respectively. Let $\vec{l}$ and $ \vec{l}'$ be length types for $W$ and $W'$ respectively. Now, suppose the following holds:
    \begin{enumerate}[font=\normalfont, label=(\alph*)]
        \item\label{item: L1 free} $||\vec{l}||=||\vec{l}'||$
        \item\label{item: L2 free} 
         For all $1\leq i\leq 2m$, $1\leq j\leq 2m'$, $l_i,l'_j>6c+4$.
    \end{enumerate}
    Then, either:
    \begin{enumerate}[font=\normalfont]
        \item There is 1-reduction in the degrees of freedom by $c$, or
        \item $W=W'$, therefore \(m = m'\), and moreover $\vec{l} = \vec{l'}$.
    \end{enumerate}
\end{corollary}

We will prove \Cref{thm: the main reduction theorem} through propositions \ref{prop: reduction or shifted matching for different wicks forms}, \ref{prop: reduction or shifted matching}, and \ref{prop: no reduction implies length types are equal}, forming the core of the subsection that follows immediately below. However, the reader can skip to the later subsections to see how \Cref{thm: the main reduction theorem} is used to obtain the element count lower bounds with \Cref{subsec: Element count: Lower bounds: Free groups} containing the discussion for free groups, and which is then generalized in the more involved \Cref{subsec: Element count: Lower bounds: Hyperbolic groups} for arbitrary hyperbolic groups. 

\subsubsection{Tools for the proof of \Cref{thm: the main reduction theorem}}\label{subsubsec: tools for the proof of the reduction theorem}
Let $c,D_0,A_0,H$, $W,m,\vec{b},\vec{l}$, $W',m',\vec{b}',\vec{l}'$ satisfy \ref{item: L1} as in \Cref{thm: the main reduction theorem}, and the following condition which is weaker than \ref{item: L2}:

\begin{enumerate}[font=\normalfont, label=(L\arabic*), start=3]
    \item\label{item: L3} For all $1\leq i\leq 2m$, $1\leq j\leq 2m'$, $l_i,l'_j>3c+2(A_0+H)+2D_0$ 
\end{enumerate}

These will be the standing assumptions for the rest of this subsection. 

The absence of any $D_0$-reduction in degrees of freedom by $c$ between $(W,\vec{l},\vec{b})$ and $(W',\vec{l}',\vec{b}')$ has the following consequences. Informally:
\begin{enumerate}
    \item Suppose a letter $X$ of $W$ has a \emph{large} overlap with a letter $X'$ of $W'$, that is, there is a long sequence of consecutive positions that belong both to $X$ with respect to $\vec{l}$ and to $X'$ with respect to $\vec{l}'$. Then for each position in the overlap, the corresponding matching positions with respect to the length types $\vec{l}$ and $\vec{l}'$ have to be at most $(D_0-1)$-separated, and moreover, $\text{sign}(X)=\text{sign}(X')$.
    \item A letter $X'$ of $W'$ cannot have \emph{large} overlap with two consecutive letters $X_j$, $X_{j+1}$ of $W$.
\end{enumerate}

\noindent We prove these observations below and then utilize them frequently in the course of proving \Cref{thm: the main reduction theorem}. Observation $(1)$ is proved in \Cref{lem: common positions} while $(2)$ is proved in \Cref{lem: reduction or redundancy}.

We introduce some notation that will be used throughout \Cref{subsec: combinatorialreduction} and \Cref{subsec: Combreduction revisited}. For any three integers $a,b,c$ satisfying $|a-b|\leq c$, we write $a\simeq_{c} b$. Note that $a\simeq_c b$ if and only if $b\simeq_c a$. Also, if $a\simeq_c b$ and $b\simeq_d e$, then $a\simeq_{c+d}e$.

\begin{lemma}\label{lem: matching position extended}
    Let $i$ be a position contained both in a letter \(X_j\) with respect to $(W, \vec{l}, \vec{b})$ and in \(X'_k\) with respect to $(W', \vec{l'}, \vec{b'})$. Moreover, suppose that  $\text{sign}(X_j)=\text{sign}(X'_k)$. Then:
    \begin{enumerate}[font=\normalfont]
        \item If $\underline{i}\simeq_{D_0-1}\underline{i}'$, then for every position $p$ common to $X_j$ and $X'_k$, we have $\underline{p}\simeq_{D_0-1}\underline{p}'$.
        \item The condition $\underline{i}\simeq_{D_0-1}\underline{i}'$ is equivalent to:
        \begin{itemize}
            \item $t_{X_{\overline{j}}} - t'_{X'_{\overline{k}'}} \simeq_{D_0-1} s'_{X'_{k}} - s_{X_j}$, when $\text{sign}(X_j)=\text{sign}(X'_k)=-1$.
            \item $s_{X_{\overline{j}}} - s'_{X'_{\overline{k}'}} \simeq_{D_0-1} s_{X_j} - s'_{X'_{k}}$, when $\text{sign}(X_j)=\text{sign}(X'_k)=+1$.
        \end{itemize}
    \end{enumerate}
\end{lemma}

\begin{proof}
    \leavevmode
    \begin{enumerate}[font=\normalfont]
        \item Recalling the expressions for matching positions from \Cref{eq: matching position 1}, observe that for all $l\geq 0$:
        \begin{align}\label{eq: formula for matching position of shift}
            i\pm l \text{ is contained in }X_j&\implies \underline{i\pm l}=\underline{i}\pm \varepsilon l\\
            i\pm l \text{ is contained in }X'_k&\implies \underline{(i\pm l)}'=\underline{i}'\pm \varepsilon l
        \end{align}
        where $\varepsilon\coloneq\text{sign}(X_j)=\text{sign}(X'_k)$. The claim follows directly from this observation.

        \item  Follows directly from \Cref{eq: matching position 1}.
    \end{enumerate}
\end{proof}

\begin{lemma}\label{lem: common positions}
    Suppose there is no $D_0$-reduction in degrees of freedom by $c$ between $(W,\vec{l},\vec{b})$ and $(W',\vec{l}',\vec{b}')$. Moreover, suppose there are $c_0\geq 3c+2(A_0+H)+D_0$ many consecutive positions $p+1,\dots,p+c_0$ common between $X_j$, a letter of $W$, and $X'_k$, a letter of $W'$. Then:
    \begin{enumerate}[font=\normalfont]
        \item $\text{sign}(X_j)=\text{sign}(X'_k)$.
        \item For all $1\leq i
        \leq c_0$, $\underline{p+i}\simeq_{D_0-1}\underline{(p+i)}'$.
    \end{enumerate}
\end{lemma}

    \begin{figure}[h]
        \centering
        \includegraphics[width=0.85\columnwidth]{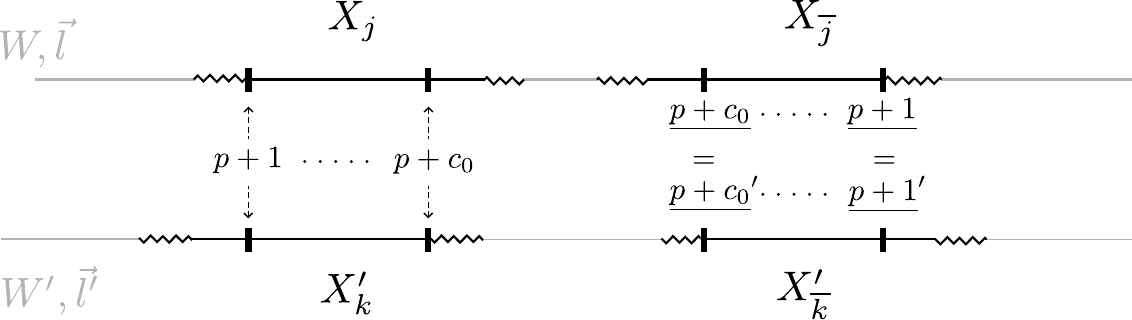}
        \caption{For $D_0=1$, when there are at least \(c_0\) many positions in common between two letters with respect to two different length types, the two corresponding matching positions of each of these common positions are equal.}
        \label{fig: common positions}
    \end{figure}

\begin{proof}
    \leavevmode
    \begin{enumerate}
        \item[(1)] Suppose, on the contrary, that $\text{sign}(X_j)\neq\text{sign}(X'_k)$. In terms of \Cref{fig: common positions}, this means that the segments corresponding to the letters $X_{\bar{j}}$ and $X'_{\bar{k}}$ are oriented opposite to each other. 
        If $|\underline{p+1}-\underline{p+1}'|\geq D_0$, set $p_0\coloneq p+1$. Otherwise, if  $\underline{p+1}\simeq_{D_0-1}\underline{p+1}'$, then set $p_0\coloneq p+1+D_0$. In either case, $|\underline{p_0}-\underline{p_0}'|\geq D_0$.
        
        Without any loss of generality, suppose $\text{sign}(X_j)=+1$, and $\text{sign}(X'_k)=-1$.
        \begin{claim}\label{claim: consecutiveadmissiblepositions}
            There is a sequence of $c$ many consecutive positions $\mathbf{i}\subset \{p_0,\dots, p_0+3c+2(A_0+H)-1\}$ such that both $\underline{\mathbf{i}}$ and $\underline{\mathbf{i}}'$ are admissible.
        \end{claim} 

        Before we prove the claim, note that it implies a $D_0$-reduction in degrees of freedom by $c$, which gives a contradiction.
        \begin{proof}
            There are three cases, depending on whether $\underline{p_0}$ and $\underline{p_0}'$ are admissible or not. In each case, using the bounds on the junction sizes and \ref{item: L3}, we argue that if at all there is any obstruction to admissibility, by skipping a maximum of $2c+2(A_0+H)$ many positions from the left in $\{p_0,\dots, p_0+3c+2(A_0+H)-1\}$, we can find a sequence $\mathbf{i}$ of $c$ many consecutive positions such that both $\underline{\mathbf{i}}$ and $\underline{\mathbf{i}}'$ are admissible:

                \emph{Case 1:}
                Suppose that neither $\underline{p_0}$ nor $\underline{p_0}'$ is admissible: in other words, $\underline{p_0}$ and $\underline{p_0}'$ are contained in a junction \(B_r'\) of \((W', \vec{l}', \vec{b}')\) and \(B_s\) of \((W, \vec{l}, \vec{b})\) respectively. Because the size of these junctions is at most $A_0+H$, the positions $\underline{p_0+A_0+H}$ and $\underline{(p_0+A_0+H)}'$ are both admissible. Consider the sequence $\mathbf{i}=(p_0+A_0+H,\dots, p_0-1+A_0+H+c )$. Because $l'_{X'_{r+1}}, l_{X_{s}}> c+(A_0+H)$, $\underline{\mathbf{i}}$ belongs to $X'_{r+1}$ while $\underline{\mathbf{i}}'$ belongs to $X_{s}$.
                
                \emph{Case 2:} Next, suppose that exactly one among $\underline{p_0}$ and $\underline{p_0}'$ is admissible. We argue for the case when $\underline{p_0}$ is admissible. The other case can be argued similarly. Since $\underline{p_0}'$ is not admissible, it belongs to some junction $B_s$. The size of $B_s$ is at most $A_0+H$, and $l_{X_{s}}>2c +2(A_0+H)$. Therefore, the $(2c+A_0+H)$ many consecutive positions $\underline{(p_0+A_0+H)}',\dots,\underline{(p_0-1+2c+2(A_0+H))}'$ are contained in the letter $X_{s}$. 
                
                Now, if $(\underline{p_0+A_0+H},\dots, \underline{p_0-1+A_0+H+c})$ is admissible, setting \(\mathbf{i}\) to be the same accomplishes our claim. Otherwise, we argue in two further sub-cases. First, suppose that $\underline{p_0+A_0+H}$ is admissible, and say it belongs to the letter $X'_r$.
                Since the size of $B'_{r}$ is at most $A_0+H$ and $l'_{X'_{r+1}}>2c+A_0+H$, the sequence of $c$ many consecutive positions $(\underline{p_0+2(A_0+H)+c},\dots,\underline{p_0-1+2(A_0+H)+2c})$ belongs to $X'_{r+1}$. Thus, setting $\mathbf{i}=(p_0+2(A_0+H)+c,\dots,p_0-1+2(A_0+H)+2c)$ completes this case. In the second case, $\underline{p_0+A_0+H}$ is not admissible and it belongs to some junction $B'_r$. As before, the size of $B'_r$ is at most $A_0+H$ and $l'_{X'_{r+1}}>c+A_0+H$. Therefore the sequence of $c$ many consecutive positions $(\underline{p_0+2(A_0+H)},\dots,\underline{p_0-1+2(A_0+H)+c} )$ is contained in $X'_{r+1}$. Setting $\mathbf{i}$ to be the tuple $(p_0+2(A_0+H),\dots,p_0-1+2(A_0+H)+c )$  completes this case too.
                
                \emph{Case 3:} Finally, suppose that both $\underline{p_0}$ and $\underline{p_0}'$ are admissible: say $\underline{p_0}$ belongs to $X'_r$ and $\underline{p_0}'$ belongs to $X_s$. Set $\mathbf{i}_0=(p_0,\dots, p_0-1+c)$. If both $\underline{\mathbf{i}_0}$ and $\underline{\mathbf{i}_0}'$ are admissible, $\mathbf{i}=\mathbf{i}_0$ works and we're done. Otherwise at least one of $\underline{\mathbf{i}_0}$ and $\underline{\mathbf{i}_0}'$ is not admissible. We argue the case when $\underline{\mathbf{i}_0}'$ is not admissible, the other case being entirely similar. Because the size of the junction $B_{s-1}$ is at most $A_0+H$, it follows that $\tilde{c}=\underline{p_0}'-t_{X_{s-1}}+1<c+(A_0+H)$. Therefore, since $l_{X_{s-1}}>2c+(A_0+H)$, the sequence of $(2c+A_0+H)$ many consecutive positions $\underline{(p_0+\tilde{c})}',\dots,\underline{(p_0-1+\tilde{c}+2c+A_0+H)}'$ belongs to $X_{s-1}$ and is therefore admissible. Set $\mathbf{i}_1=(p_0+\tilde{c},\dots,p_0-1+\tilde{c}+2c+A_0+H)$. If $\underline{\mathbf{i}_1}$ is admissible, then $\mathbf{i}=\mathbf{i}_1$ works, and we're done. Otherwise, if $\underline{p_0+\tilde{c}}$ is admissible but $\underline{\mathbf{i}_1}$ is not, then like we've argued before in Case 2, the sequence of $c$ many consecutive positions  $\underline{p_0+\tilde{c}+c+(A_0+H)},\dots,\underline{p_0-1+\tilde{c}+2c+(A_0+H)}$ is admissible, and setting $\mathbf{i}=(p_0+\tilde{c}+c+(A_0+H),\dots,p_0-1+\tilde{c}+2c+(A_0+H))$ $\mathbf{i}=(p_0+2c+2(A_0+H),\dots, p_0-1+3c+2(A_0+H))$, we're done. Otherwise, $\underline{p_0+\tilde{c}}$ is not admissible, and setting $\mathbf{i}=(p_0+\tilde{c}+(A_0+H),\dots, p_0-1+\tilde{c}+c+(A_0+H))$ works.
            \end{proof}
        \item[(2)] We have $\text{sign}(X_j)=\text{sign}(X'_k)$. If $|\underline{(p+1)}-\underline{(p+1)}'|\geq D_0$, it follows from \Cref{lem: matching position extended} that $|\underline{(p+i)}-\underline{(p+i)}'|\geq D_0$ for all $i\in\{1,\dots,c_0\}$. Now, the same arguments made in the proof \Cref{claim: consecutiveadmissiblepositions} show that there is a sequence of $c$ many consecutive positions $\mathbf{i}\subset\{p+1,\dots, p+3c+2(A_0+H)\}$ such that both $\underline{\mathbf{i}}$ and $\underline{\mathbf{i}}'$ are admissible. This means that there is a $D_0$-reduction in degrees of freedom by $c$, which contradicts our hypothesis.
    \end{enumerate}
\end{proof}

%We need l_i>2(D_0-1)+A_0+H for this lemma
\begin{lemma}\label{lem: reduction or redundancy}
    Suppose there is no $D_0$-reduction in degrees of freedom by $c$ between $(W,\vec{l},\vec{b})$ and $(W',\vec{l}',\vec{b}')$. Then it cannot be the case that a letter $X'_k$ of $W'$ has $3c+2(A_0+H)+D_0$ or more positions in common with each of two consecutive letters $X_j$, and $X_{j+1}$ of $W$.
\end{lemma}

\begin{figure}[h]
    \centering
    \includegraphics[width=0.85\columnwidth]{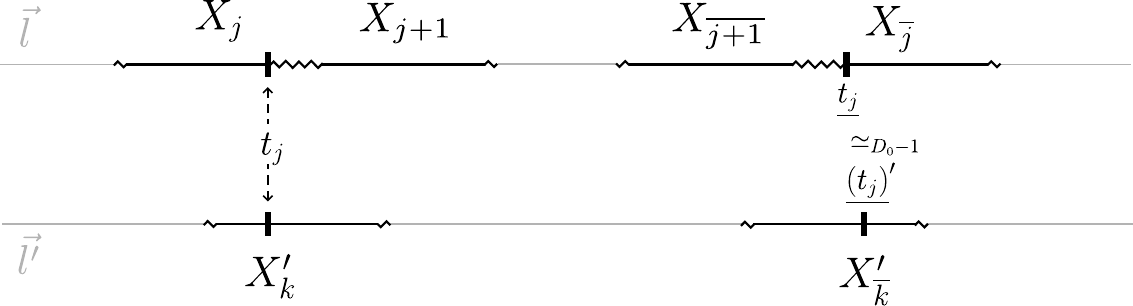}
    \caption{When \(\underline{(t_{j})} \simeq_{D_0 - 1} \underline{(t_j)}'\), \(W\) contains the redundancy \(X_{j}X_{j+1}\), a contradiction.}
    \label{fig: no redundancy}
\end{figure}

\begin{proof}
    Suppose the contrary. Our hypotheses imply that both $t_{X_j}$, the ending position of $X_j$, and $s_{X_{j+1}} = t_{X_j} + b_j+1$, the starting position of $X_{j+1}$, must be contained in $X'_k$. Moreover, there are at least $3c+2(A_0+H)+D_0$ consecutive positions to the left of $t_{X_j}$ which are contained in $X'_k$, and there are at least $3c+2(A_0+H)+D_0$ consecutive positions to the right of $s_{X_{j+1}}$, contained in $X'_k$. Thus, it follows from \Cref{lem: common positions} that:
    \begin{enumerate}[label=(\roman*)]
        \item $\text{sign}(X_j)=\text{sign}(X'_k)=\text{sign}(X_{j+1})(\eqcolon\varepsilon\text{, say})$.
        \item $\underline{t_{X_j}}\simeq_{D_0-1}\underline{(t_{X_j})}'$ and $\underline{s_{X_{j+1}}}\simeq_{D_0-1}\underline{(s_{X_{j+1}})}'$. 
    \end{enumerate}
    Irrespective of whether $\varepsilon=-1$ or $+1$, $|\underline{(t_{X_j})}'-\underline{(s_{X_{j+1}})}'|=|s_{X_{j+1}}-t_{X_j}|=b_j+1$. Therefore, using ({\romannumeral 2}), we have:
    \[ |\underline{t_{X_j}}- \underline{s_{X_{j+1}}}|\leq 2(D_0-1)+A_0+H+1\]
    
    \noindent If $\varepsilon=-1$, then: $\underline{t_{X_j}}=s_{X_{\bar{j}}}$, $\underline{s_{X_{j+1}}}=t_{X_{\overline{j+1}}}$. Therefore, $|s_{X_{\bar{j}}}-t_{X_{\overline{j+1}}}|\leq 2(D_0-1)+A_0+H+1$. Note that $l_i>2(D_0-1)+A_0+H$ for all $1\leq i\leq 2m$. So it follows that $\bar{j}=\overline{j+1}+1$.

    \noindent If $\varepsilon=1$, then: $\underline{t_{X_j}}=t_{X_{\bar{j}}}$, and $\underline{s_{X_{j+1}}}=s_{X_{\overline{j+1}}}$. Therefore, $|s_{X_{\overline{j+1}}}-t_{X_{\overline{j}}}|\leq 2(D_0-1)+A_0+H+1$. So, as before, it follows that $\overline{j+1}=\overline{j}+1$.
    
    In any case, we infer that $(X_jX_{j+1})^{\varepsilon}$ is also a subword of $W$, that is, $W$ has a redundancy, which is not possible.

\end{proof}

% {\color{gray} If it's appropriate, we can include the gist of Inhyeok's proof sketch for the lower bound applied to free groups here. The idea is quite straightforward, and if we forget about proving \Cref{intro thm: free groups asymptotics} for now, then the argument allows the reader to skip several of the lengthy proofs that follow.}

\subsubsection{Proof of \Cref{thm: the main reduction theorem}}\label{subsubsec: Proof of the reduction theorem}

Let us first recall some notation that will be relevant to the proof of this theorem. For the rest of this subsection, $W=X_1\cdots X_{2m}$ and $W'=X'_1\cdots X'_{2m'}$ are two pre-Wicks forms having length $2m$ and $2m'$ respectively. $\vec{l}$ and $\vec{l}'$ are two length types for $W$ and $W'$ respectively. We also have two sequences $\vec{b}\in[A_0,A_0+H]^{2m+1}$ and $\vec{b}'\in[A_0,A_0+H]^{2m'+1}$. $(W,\vec{l},\vec{b})$ and $(W',\vec{l}',\vec{b}')$ satisfy conditions \ref{item: L1} and \ref{item: L2} from \Cref{thm: the main reduction theorem}.
For a letter $X_j$ of $W$, recall that $s_{X_j}$ and $t_{X_j}$ denote the starting and ending positions of $X_j$ with respect to $\vec{l}$. Similarly, for a letter $X'_k$ of $W'$, let $s'_{X'_k}$ and $t'_{X'_k}$ be defined analogously with respect to $\vec{l}'$.

\begin{proposition}\label{prop: reduction or shifted matching for different wicks forms}
    Suppose that there is no $D_0$-reduction in degrees of freedom by $c$ for $(W,\vec{l},\vec{b})$ and $(W',\vec{l}',\vec{b}')$. Then, $W$ is equivalent to $W'$ as a pre-Wicks form. In particular, $m=m'$ and for all $1\leq j\leq 2m$, $X_j$ and $X'_j$ have at least $3c+2(A_0+H)+D_0$ many positions in common.
\end{proposition}

\begin{proof}
We begin by making the following observation.

\begin{claim}\label{claim: reduction or controlled covering of letters}
    Given any letter $X_i$ of $W$, consider the sequence of letters in $W'$ \textit{covering} it, that is, the largest collection of consecutive letters $X'_j,\dots,X'_{j+s}$ of $W'$ with each of them having positions common with $X_i$. Then by \Cref{lem: reduction or redundancy}, either: 
    \begin{enumerate}[font=\normalfont]
        \item $s=0$, $X_i$ is exactly covered by one letter $X'_j$ of $W'$.
        \item $s=1$, $X_i$ is covered by two letters, $X'_j$ and $X'_{j+1}$, with exactly one of them having more than $3c+2(A_0+H)+2D_0$ positions in common with $X_i$.
        \item $s=2$, $X_i$ is covered by three letters, $X'_j$, $X'_{j+1}$, and $X'_{j+2}$, with only $X'_{j+1}$ having more than $3c+2(A_0+H)+2D_0$ positions in common with $X_i$.
    \end{enumerate}
\end{claim}

\begin{proof}
    To begin with, because $l_{X_1},l_{X_{2m}},l'_{X'_1},l'_{X'_{2m'}}> 3c+3(A_0+H)+2D_0$, $b_1,b'_1,b_{2m+1},b'_{2m'+1}\leq A_0+H$, and $L(W,\vec{l},\vec{b})=L(W',\vec{l'},\vec{b'})$, $X_1$ and $X'_1$, and $X_{2m}$ and $X'_{2m'}$ respectively, have at least $3c+2(A_0+H)+2D_0$ many positions in common. It follows from \Cref{lem: reduction or redundancy} that $X_1$ and $X_{2m}$ are covered by at most $2$ consecutive letters of $W'$. Moreover, by \Cref{lem: common positions} it follows that $\text{sign}(X_1)=\text{sign}(X'_1)$, $\text{sign}(X_{2m})=\text{sign}(X'_{2m'})$, and $X_{\overline{1}}$ and $X'_{\overline{1}'}$, and $X_{\overline{2m}}$ and $X'_{\overline{2m'}}$ respectively, have at least $3c+2(A_0+H)+D_0$ many positions in common.

    Next, consider any letter $X_i$ where $i\neq 1,2m, \overline{1},\overline{2m}$. If $X_i$ is covered by exactly one letter $X'_j$, then note that they have at least $3c+2(A_0+H)+D_0$ many positions in common because $l_{X_i},l'_{X'_j}>3c+2(A_0+H)+D_0$. If $X_i$ is covered by exactly two letters $X'_j$, $X'_{j+1}$, because $l_{X_i}>6c+5(A_0+H)+4D_0$ and \Cref{lem: reduction or redundancy}, $X_i$ has at least $3c+2(A_0+H)+2D_0$ many positions in common with exactly one of $X'_j$ and $X'_{j+1}$. Finally, because all letters of $W'$ have length greater than $3c+2(A_0+H)+D_0$, $X_i$ cannot be covered by more than three consecutive letters of $W'$, due to \Cref{lem: reduction or redundancy}. In the case when $X_i$ is covered by exactly three consecutive letters $X'_j,X'_{j+1},X'_{j+2}$, as before, it follows from \Cref{lem: reduction or redundancy} that only $X'_{j+1}$ has at least $3c+2(A_0+H)+2D_0$ many positions in common with $X_i$.
\end{proof}

Thus, for every letter $X_i$, there is a unique letter $X'_{\varphi(i)}$ such that $X_i$ and $X'_{\varphi(i)}$ have more than $3c+2(A_0+H)+D_0$ positions in common. Due to \Cref{lem: common positions}:

\begin{enumerate}[label=(\roman*)]
    \item $\text{sign}(X_i)=\text{sign}(X'_{\varphi(i)})$.
    \item $\varphi(\bar{i})=\overline{\varphi(i)}$, i.e., $X_{\bar{i}}$, the letter in $W$ paired with $X_i$, has at least $3c+2(A_0+H)+D_0$ positions in common with $X'_{\overline{\varphi(i)}}$, the letter in $W'$ paired with $X'_{\varphi(i)}$.
\end{enumerate}
 In other words, we have an injective, pairing-preserving map from the letters of $W$ to the letters of $W'$. But because $L(W, \vec{l}, \vec{b})=L(W', \vec{l}', \vec{b}')$ (condition $(\mathbf{L1})$), by switching the roles of $W$ and $W'$ in the above argument, we conclude that this map also has to be surjective. Therefore, $W'$ is equivalent to $W$. It follows that $m=m'$, and that $X_j$ and $X'_j$ have at least $3c+2(A_0+H)+D_0$ positions in common for all $1\leq j\leq 2m$.
\end{proof}

\begin{remark}\label{rmk: reduction or equality strengthened}
    In the proof of \Cref{prop: reduction or shifted matching for different wicks forms}, we do not require the full strength of the assumption \ref{item: L2}. In fact, it suffices if $(W,\vec{l},\vec{b})$ and $(W',\vec{l}',\vec{b}')$ satisfy $(\mathbf{L1})$ and the following condition: $l_{X_1},l_{X_{2m}},l'_{X'_1},l'_{X'_{2m'}}> 3c+3(A_0+H)+2D_0$, and for all $1\leq j\leq 2m$, $1\leq j'\leq 2m'$ with $j\neq 1,2m,\overline{1},\overline{2m}$, $j'\neq 1,2m',\overline{1}',\overline{2m'}'$, $l_{X_j},l_{X'_{j'}}>6c+5(A_0+H)+4D_0$.
\end{remark}

\Cref{prop: reduction or shifted matching for different wicks forms} reduces the proof of \Cref{thm: the main reduction theorem} to the case when $W=W'$. We will now introduce a definition:

\begin{definition}[Shifts]\label{def: shifts}
    Given two length types $\vec{l}, \vec{l}'$ for a pre-Wicks form $W=X_1\dots X_{2m}$, define the \textit{shift} of $X_j$, $\text{shift}(j)\coloneq t_{X_j}-t'_{X_j}$. Given two sequences $\vec{b},\vec{b}'\in [A_0,A_0+H]^{2m+1}$, we define $\Delta(j)=b_j-b'_j$ for all $0 \leq j \leq 2m$.
\end{definition}

We will also make a simplification in notation in view of \Cref{prop: reduction or shifted matching for different wicks forms}. For $W=X_1\dots X_{2m}$ and two length types $\vec{l},\vec{l}'$: $s_j \coloneq s_{X_j}, s'_j \coloneq s'_{X_j}, t_j\coloneq t_{X_j}$ and $t'_j\coloneq t'_{X_j}$.

\begin{corollary}\label{cor: shift equations}
    Consider $(W,\vec{l},\vec{b})$ and $(W,\vec{l}',\vec{b}')$. Suppose there is a position $i$, contained in $X_j$, with respect to both $\vec{l}$ and $\vec{l}'$. Then the condition $\underline{i}\simeq_{D_0-1}\underline{i}'$ implies the following:
    \begin{enumerate}[font=\normalfont]
        \item If $\text{sign}(X_j)=-1$, then $\text{shift}(\overline{j})\simeq_{D_0+H-1}-\text{shift}(j-1)$.
        \item If $\text{sign}(X_j)=1$, then $\text{shift}(j-1)\simeq_{D_0+2H-1}\text{shift}(\overline{j}-1)$ and $\text{shift}(j)\simeq_{D_0-1}\text{shift}(\overline{j})$.
    \end{enumerate}
\end{corollary}

\begin{proof}
    This directly follows from \Cref{lem: matching position extended}. Here's an explicit check, for the sake of clarity. Note that for all $1\leq j\leq 2m$, $s_j=t_{j-1}+b_{j-1}+1$ and $s'_j=t'_{j-1}+b'_{j-1}+1$. Therefore:
    \begin{equation*}
        s_j-s'_j=\text{shift}(j-1)+\Delta(j-1)
    \end{equation*}
    We use this, together with the fact that for all $1\leq j\leq 2m$, $|\Delta(j)|\leq H$.
    \begin{enumerate}
        \item If $\text{sign}(X_j)=-1$, by \Cref{lem: matching position extended}, we get:
        \begin{align*}
            \underline{i}\simeq_{D_0-1}\underline{i}' &\implies t_{\overline{j}}-t'_{\overline{j}}\simeq_{D_0-1}-(s_j-s'_j)=-(\text{shift}(j-1)+\Delta(j-1))\\
            &\implies \text{shift}(\overline{j})\simeq_{D_0+H-1}-\text{shift}(j-1)
        \end{align*}
        \item If $\text{sign}(X_j)=1$, by \Cref{lem: matching position extended}, we get:
        \begin{align*}
            \underline{i}\simeq_{D_0-1}\underline{i}' &\implies s_j-s'_j\simeq_{D_0-1}s_{\overline{j}}-s'_{\overline{j}} \implies \text{shift}(j-1)+\Delta(j-1)\simeq_{D_0-1}\text{shift}(\overline{j}-1)+\Delta(\overline{j}-1)\\
             &\implies \text{shift}(j-1)\simeq_{D_0+2H-1}\text{shift}(\overline{j}-1)
        \end{align*}
        Now note that $l_{X_j}=t_j-s_j+1$, thus we also get:
        \begin{align*}
            \underline{i}\simeq_{D_0-1}\underline{i}' &\implies t_j-t'_j\simeq_{D_0-1}t_{\overline{j}}-t'_{\overline{j}}
            \implies \text{shift}(j)\simeq_{D_0-1}\text{shift}(\overline{j})
        \end{align*}
    \end{enumerate}
\end{proof}

\begin{proposition}[Shifted matching]\label{prop: reduction or shifted matching}
    Let $W$ be a pre-Wicks form and let $\vec{l},\vec{l}'$ be two length types for $W$ satisfying conditions \ref{item: L1} and \ref{item: L2}. Suppose there is no $D_0$-reduction in degrees of freedom by $c$. Then for all $1\leq j\leq 2m$:
    \begin{enumerate}[font=\normalfont]
            \item If \(\text{sign}(X_j) = -1\), then
            \begin{enumerate}[font=\normalfont,label=(1\alph*)]
                \item\label{item: shift equation 1a} $\text{shift}(\overline{j})\simeq_{D_0+H-1}-\text{shift}(j-1)$, and
                \item\label{item: shift equation 1b} $\text{shift}(j) \simeq_{D_0+H-1} -\text{shift}(\overline{j} - 1)$.
            \end{enumerate}
            
            \item If \(\text{sign}(X_j) = 1\), then
            \begin{enumerate}[font=\normalfont,label=(2\alph*)]
                \item\label{item: shift equation 2a} $\text{shift}(j-1)\simeq_{D_0+2H-1}\text{shift}(\overline{j}-1)$, and
                \item\label{item: shift equation 2b} $\text{shift}(j)\simeq_{D_0-1}\text{shift}(\overline{j})$.
            \end{enumerate}
    \end{enumerate}
\end{proposition}

\begin{proof}
    By \Cref{prop: reduction or shifted matching for different wicks forms}, for all $1\leq j\leq 2m$, $X_j$ has more than $3c+2(A_0+H)+D_0$ positions in common with respect to $\vec{l}$ and $\vec{l}'$ and it follows from \Cref{lem: common positions} that for any such common position $i$, we have $\underline{i}\simeq_{D_0-1}\underline{i}'$. The conclusions now follow from \Cref{cor: shift equations}: $(1)$ of \Cref{cor: shift equations} proves $1(a)$, $1(b)$ follows by swapping $j$ and $\overline{j}$ in $1(a)$, and $(2)$ of \Cref{cor: shift equations} proves both $2(a)$ and $2(b)$.

\end{proof}

\begin{proposition}\label{prop: no reduction implies length types are equal}
    Let $c$ and $D_0$ be positive integers. Let $W$ be a maximal pre-Wicks form of genus $n$ and let $\vec{l}, \vec{l}'$ be two length types for $W$, and let $\vec{b},\vec{b}'\in [A_0,A_0+H]^{l(W)+1}$ satisfying \ref{item: L1} and \ref{item: L2}.
    If there is no $D_0$-reduction in degrees of freedom by $c$, then for all $1\leq j\leq l(W)$, $|l_j-l'_j|\leq 3(D_0-1)+ 4H$.
\end{proposition}

\begin{proof}
    $W$ is a maximal pre-Wicks form, so by \Cref{prop: what does a preWicks form look like}, either:
    \begin{itemize}
        \item $W=DAU_1A^{-1}B^{-1}U_2BD^{-1}$, and the word $U_1X^{-1}U_2X$ is a Wicks form, or
        \item $W=DAU_1BAU_2BD^{-1}$, and the word $U_1XU_2X$ is a Wicks form.
    \end{itemize}
    where \(U_1, U_2\) are subwords of \(W\), and \(A, B, D\) and $X$ are unknowns that do not appear in $U_1$ or $U_2$. $l(W)=12n-2$ or $l(W)=6n-2$ according to whether $W$ is orientable or non-orientable respectively.

    Note that for all \(1\leq j \leq l(W)\), we have:
    \begin{align}\label{eq: shifted length match}
        l_j = l_j' + \text{shift}(j) - \text{shift}({j-1}) - \Delta(j-1)
    \end{align}

    Since $|\Delta(j)| \leq H$ for all $0 \leq j \leq l(W)$, in order to prove the proposition, it suffices to show that: 
    
    \begin{claim}\label{eq: shift bound in element count}
        For all $0\leq j \leq l(W)$, $|\text{shift}(j)| \leq 3(D_0+H-1)/2$.
    \end{claim}
    \begin{proof}
        There is nothing to check for $j=0$, as $\text{shift}(0) = 0$. For $j=l(W)$, $t_{l(W)} + b_{l(W)} = t'_{l(W)} + b'_{l(W)}$ by \textbf{(L1)}, and therefore $|t_{l(W)}-t'_{l(W)}|=|b_{l(W)}-b'_{l(W)}|\leq H$. So, suppose $0 < j< l(W)$. 
        Recall the polygon $P_W$ and the graph $\Gamma_W$ associated with $W$ (cf. \Cref{def: wicks form genus}).
        The cycle $\partial P_W$ based at $v$ (corresponding to reading the letters of $W$ from left to right) projects to a cycle $\alpha$ in $\Gamma_W$. We think of $\alpha$ as a vertex path $(v_0,\dots,v_{l(W)})$ and observe that it satisfies the following:
        
    \begin{enumerate}
        \item Every edge of $\Gamma_W$ is traversed exactly twice: both in the same direction if the sign of the letter labelling the edge is $+1$, and in opposite directions otherwise.
        \item There is no backtracking.
    \end{enumerate}
    The following observation plays a key role in the proof:
    \begin{quote}
        There is exactly one vertex of degree $1$ in $\Gamma_W$, and the rest of the vertices have degree $3$. 
    \end{quote}
    In particular, observe that $\text{degree}(v_0)=1$, and $v_i\neq v_0$ for all $i$ except $i=0 \text{ and }l(W)$. Therefore, $v_j$ is a trivalent vertex. It follows from $(1)$ and $(2)$ above that the turns made by $\alpha$ in a $1$-neighbourhood of a trivalent vertex can be exactly one among five types, as depicted in \Cref{fig: circuit in 1-neighbourhood}.

    \begin{figure}[h]
            \centering
            \includegraphics[width=0.85\columnwidth]{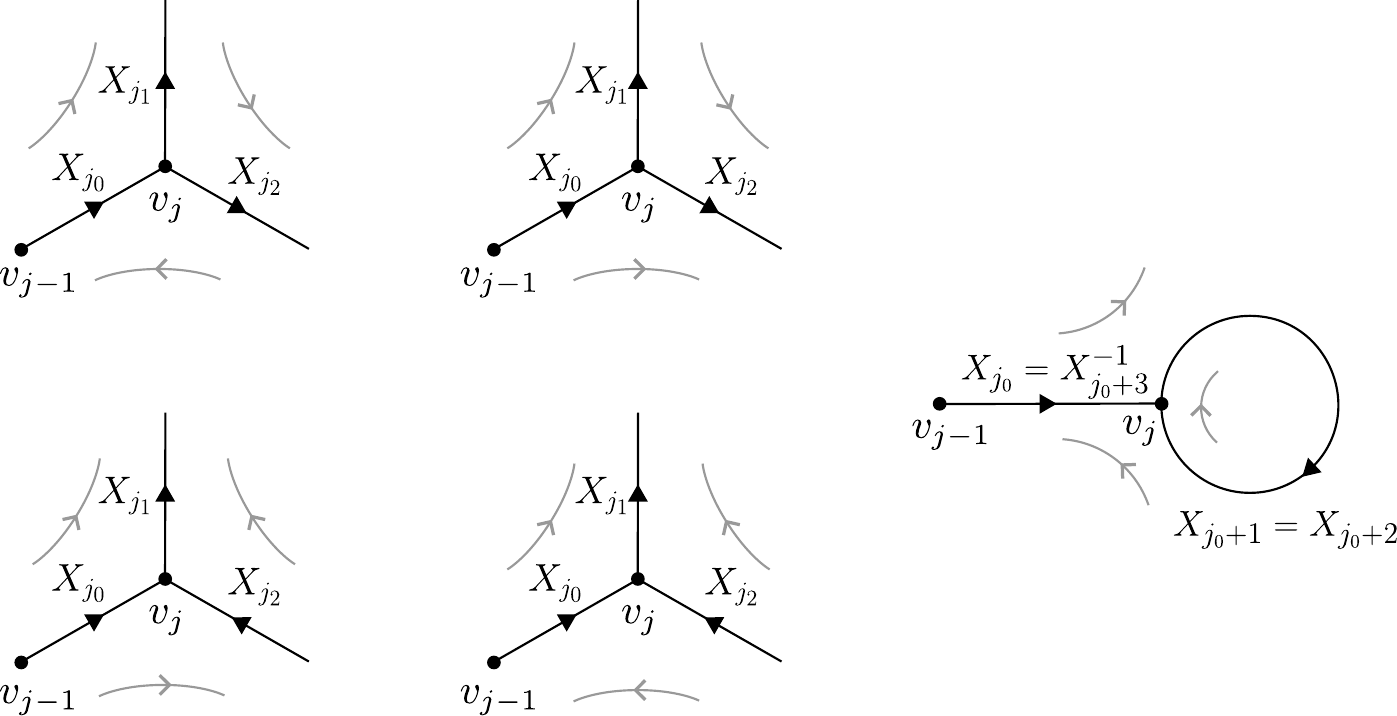}
            \caption{The different cases of the circuit corresponding to $W$ intersecting the 1-neighbourhood of $v_j$ in $\Gamma_W$. The grey arrows denote the passing of the circuit through $v_j$. More precisely, for example, in the top left figure the circuit contains the subpaths $e_{j_0} e_{j_1}$, $\overline{e}_{j_1}e_{j_2}$, and $\overline{e}_{j_2}\overline{e}_{j_0}$. Equivalently, $W$, when considered as a cyclic word, contains the subwords $X_{j_0} X_{j_1}$, $X_{\overline{j_1}} X_{j_2}$, and $X_{\overline{j_2}}X_{\overline{j_0}}$}. 
            \label{fig: circuit in 1-neighbourhood}
    \end{figure}
    In order to read out the turns made by $\alpha$ as it traverses the vertex $v_j$ in $\Gamma_W$ during its full run, we define indices $j_0, j_1$, and $j_2$ as follows: $j_0\coloneq j$, $j_1\coloneq j+1$, and $X_{j_2}$ is the unique letter adjacent to $X_{\overline{j_1}}$ with a vertex (projecting down to) $v_j$.
    \\\\More explicitly, the following is an algebraic description of the relative positions of $X_{j_0}, X_{j_1},$ and $X_{j_2}$ in $W$ depending on the possible signs of each of these letters. Here we have set $X_{j_0}=X, X_{j_1}=Y$, and $X_{j_2}=Z$. These cases are also shown in \Cref{fig: circuit in 1-neighbourhood}:
    \begin{enumerate}[label=(\alph*)]
        \item \underline{Top left, without loop:} $\text{sign}(X_{j_0})=\text{sign}(X_{j_1})=\text{sign}(X_{j_2})=-1$ and $W$ contains the subwords $XY, Y^{-1}Z, \text{and } Z^{-1}X^{-1}$.
        \item \underline{Top right, without loop:} $\text{sign}(X_{j_0})=\text{sign}(X_{j_2})=+1$, $\text{sign}(X_{j_1})=-1$, and $W$ contains the subwords $ XY, Y^{-1}Z, \text{and } XZ$.
        \item \underline{Bottom left, without loop:} $\text{sign}(X_{j_0})=\text{sign}(X_{j_1})=+1$, $\text{sign}(X_{j_2})=-1$, and $W$ contains the subwords $ XY, ZY, \text{and } XZ^{-1}$.
        \item \underline{Bottom right, without loop:} $\text{sign}(X_{j_0})=-1$, $\text{sign}(X_{j_1})=\text{sign}(X_{j_2})=+1$, and $W$ contains the subwords $ XY, ZY, \text{and } ZX^{-1}$.
        \item \underline{with loop:} This splits into two sub-cases, although they essentially have the same algebraic form:
        \begin{enumerate}[label=(\roman*)]
            \item If $v_{j-1}\neq v_j$, then $\text{sign}(X_{j_0})=\text{sign}(X_{j_2})=-1$, and $\text{sign}(X_{j_1})=+1$. In this case:  $j_0\coloneq j$, $j_1\coloneq j+1$, and $j_2\coloneq\overline{j_0}$. $W$ contains $XY^2X^{-1}$ as a subword.
            \item If $v_{j-1}=v_j$, then $\text{sign}(X_{j_0})=\text{sign}(X_{j_1})=+1$, and $\text{sign}(X_{j_2})=-1$. In this case: $j_0\coloneq j$, $j_1\coloneq j+1$, and $j_2\coloneq j_0-1$. $W$ contains $ ZX^2Z^{-1}$ as a subword.
        \end{enumerate}
        
    \end{enumerate}
    Note that $(a)$ is the only relevant neighbourhood type when $W$ is orientable. However, when $W$ is non-orientable, all five types need to be taken into account. 
    
    We can now use the shift equations from \Cref{prop: reduction or shifted matching} to prove the claim. We argue for each of the cases in the algebraic description above:

    \begin{enumerate}[label=(\alph*)]
        \item \underline{Top left, without loop:} Using the algebraic description, we get: $j_2=\overline{j_1}+1$ and $\overline{j_2}+1=\overline{j_0}$. Now, applying \ref{item: shift equation 1a} from \Cref{prop: reduction or shifted matching} thrice:
        
        \begin{multline*}
            \text{shift}(j_0)\overset{1(a)}{\simeq}_{D_0+H-1}-\text{shift}(\overline{j_0+1})=-\text{shift}(\overline{j_1})\overset{1(a)}{\simeq}_{D_0+H-1}\text{shift}(\overline{\overline{j_1}+1})\\
            =\text{shift}(\overline{j_2})\overset{1(a)}{\simeq}_{D_0+H-1}-\text{shift}(\overline{\overline{j_2}+1})=-\text{shift}(\overline{\overline{j_0}})=-\text{shift}(j_0)
        \end{multline*}

        Therefore, $\text{shift}(j_0) \simeq_{3(D_0+H-1)} -\text{shift}(j_0)$, and hence $|\text{shift}(j)|=|\text{shift}(j_0)|\leq 3(D_0+H-1)/2$.
        \item \underline{Top right, without loop:} Using the algebraic description, we get: $j_2=\overline{j_1}+1$ and $\overline{j_2}-1=\overline{j_0}$. Applying \ref{item: shift equation 1a} followed by \ref{item: shift equation 2a} from \Cref{prop: reduction or shifted matching}:
        \begin{multline*}
            \text{shift}(j_0)\overset{1(a)}{\simeq}_{D_0+H-1}-\text{shift}(\overline{j_0+1})=-\text{shift}(\overline{j_1})\\
            =-\text{shift}(j_2-1)\overset{2(a)}{\simeq}_{D_0+2H-1}-\text{shift}(\overline{j_2}-1)=-\text{shift}(\overline{j_0})
        \end{multline*}
        Further, applying \ref{item: shift equation 2b} gives: $\text{shift}(j_0)\simeq_{D_0-1}\text{shift}(\overline{j_0})$. Therefore, $|\text{shift}(j)|=|\text{shift}(j_0)|\leq3(D_0+H-1)/2$.
        \item \underline{Bottom left, without loop:} Using the algebraic description, we get: $j_2=\overline{j_1}-1$ and $\overline{j_0}=\overline{j_2}-1$. Applying \ref{item: shift equation 2a} followed by \ref{item: shift equation 1b} and \ref{item: shift equation 2b}, we get:
        \begin{multline*}
            \text{shift}(j_0)\overset{2(a)}{\simeq}_{D_0+2H-1}\text{shift}(\overline{j_1}-1)=\text{shift}(j_2)\\
            \overset{1(b)}{\simeq}_{D_0+H-1}-\text{shift}(\overline{j_2}-1)=-\text{shift}(\overline{j_0})\overset{2(b)}{\simeq}_{D_0-1}-\text{shift}(j_0)
        \end{multline*}
        Therefore, $|\text{shift}(j)|=|\text{shift}(j_0)|\leq 3(D_0+H-1)/2$.
        \item \underline{Bottom right, without loop:} Using the algebraic description, we get: $j_2=\overline{j_1}-1$ and $\overline{j_0}=\overline{j_2}+1$. Applying \ref{item: shift equation 2a}, \ref{item: shift equation 2b}, and \ref{item: shift equation 1b}, we get:
        \begin{align*}
            \text{shift}(j_0)\overset{2(a)}{\simeq}_{D_0+2H-1}\text{shift}(\overline{j_1}-1)=\text{shift}(j_2)\overset{2(b)}{\simeq}_{D_0-1}\text{shift}(\overline{j_2})=\text{shift}(\overline{j_0}-1)\overset{1(b)}{\simeq}_{D_0+H-1}-\text{shift}(j_0)
        \end{align*}
        Therefore, $|\text{shift}(j)|=|\text{shift}(j_0)|\leq 3(D_0+H-1)/2$.
        \item \underline{with loop:}
        \begin{enumerate}[label=(\roman*)]
            \item Using the algebraic description, we get: $\overline{j_1}=j_1+1=j_0+2$ and $j_2=\overline{j_1}+1=j_0+3=\overline{j_0}$. Applying \ref{item: shift equation 2a}, \ref{item: shift equation 2b}, and \ref{item: shift equation 1b},  we get:
            \begin{multline*}
                \text{shift}(j_0)\overset{2(a)}{\simeq}_{D_0+2H-1}\text{shift}(\overline{j_1}-1)=\text{shift}(j_1)\overset{2(b)}{\simeq}_{D_0-1}\text{shift}(\overline{j_1})\\
                =\text{shift}(j_2-1)=\text{shift}(\overline{j_0}-1)\overset{1(b)}{\simeq}_{D_0+H-1}-\text{shift}(j_0)
            \end{multline*}
            Therefore, $|\text{shift}(j)|=|\text{shift}(j_0)|\leq 3(D_0+H-1)/2$.
            \item The former case provides us with
            \[
                \text{shift}(j_0) {\simeq}_{D_0+2H-1} \text{shift}(j_1) \simeq_{2D_0 + H - 2} - \text{shift}(j_0) {\simeq}_{D_0+2H-1} -\text{shift}(j_1),
            \]            
            implying $|\text{shift}(j_1)|=|\text{shift}(\overline{j_1})|\leq 3(D_0+H-1)/2$. This completes the proof, since the current case (corresponding to the subword $ZX^2Z^{-1}$) includes as a special version of the former (corresponding to $XY^2X^{-1}$) with $j_0$ here being included as $j_1$.
        \end{enumerate}
    \end{enumerate}
    \end{proof}
    This completes the proof of \Cref{prop: no reduction implies length types are equal} and therefore that of \Cref{thm: the main reduction theorem}.
    \end{proof}
    
    \begin{remark}\label{rmk: injectivity fails when there is degree four vertex}
        %{\color{gray}(This remark was originally written for free groups: how should it be included now? Recall the fact that  no reduction by $k \log L$ implies the maximum length of $\gamma$ is bounded by $k \log L$), which affects the lower bound $L$-exponent for all hyperbolic groups.}
        
        The assumption of \(W\) being maximal is important in the arguments of \Cref{subsec: combinatorialreduction}. For non-maximal Wicks forms and length types (satisfying conditions \Cref{item: L1 free} and \ref{item: L2 free}), there might still by some reduction in the degrees of freedom, but not necessarily by $N$.
        
        For instance, suppose \(W\) has subwords of the form \( A  B^{-1} \),\( B C^{-1}  \),\( C D^{-1}\), and \( D A^{-1}\). Note that \(W\) is not maximal as \(\Gamma_W\) contains a vertex of degree at least four. 
        For a concrete example, consider $W=AB^{-1}XY^{-1}  BC^{-1}Q^{-1}R CD^{-1}X^{-1}Y DA^{-1}QR^{-1}$, a Wicks form of orientable genus $2$. 
        Fix $N\geq 1$. Consider two non-cancelling substitutions \(\mathfrak{g}\) and \(\mathfrak{g}'\) of \(W\) which agree on each letter of \(W\) except \(A, B, C, D\), and

        \begin{alignat*}{4}
            \mathfrak{g}_A &= a \quad& \mathfrak{g}_B &= b\gamma \quad& \mathfrak{g}_C &= c \quad& \mathfrak{g}_D &= d \gamma \\
            \mathfrak{g}'_A &= a\gamma^{-1} \quad& \mathfrak{g}'_B &= b \quad& \mathfrak{g}'_C &= c \gamma^{-1} \quad& \mathfrak{g}'_D &= d 
        \end{alignat*}
        with \(a, b, c, d, \gamma\) being elements in the free group, with $l(\gamma)<N$ and $l(a),l(b),l(c),l(d)\geq 10N$. With these assumptions, note that \(\mathfrak{g}(W) = \mathfrak{g}'(W)\) and for $\gamma$ non-trivial, the corresponding length types are also different. However, there is no ($1$-)reduction in the degrees of freedom by $N$. 
        
    \end{remark}

We now move on to applying the main reduction theorem to prove the desired element count lower bounds. In the following subsection, we show this for free groups, where the argument is relatively straightforward. The proof for arbitrary hyperbolic groups is identical in spirit, but the explicit details become more involved, especially surrounding the inclusions of various fellow-travelling constants. We therefore deal with this case in a separate subsection.

\subsection{Counting in Free groups}\label{subsec: Element count: Lower bounds: Free groups}

Fix a positive integer $r\geq 2$ and $F_r$ denote the free group of rank $r$, with a symmetric generating set $S$ of \(2r\) elements. Fix $n\geq 1$. In this subsection, we always set $A_0 = H = 0$, and $D_0 = 1$ while applying the results of \Cref{subsec: combinatorialreduction}. This also allows us to ignore $\vec{b}$ here, since now it must always be the zero vector. Therefore, we only refer to the tuple $(W,\vec{l})$, with $W$ being a pre-Wicks form and $\vec{l}$ a length type for $W$, as usual. We will use definitions and notation introduced in \Cref{subsec: Element count: Upper bounds: Free groups}. In particular, the reader should recall \Cref{def: setofnon-cancellingsubs} and \Cref{def: the numbers tau}. Due to \Cref{lem: canonical forms for group elements}, the set $\partial B^{\pm}_n(L)$ is empty when $L$ is odd. Therefore in the rest of this subsection, the parameter $L$ is always assumed to be even.

When two different non-cancelling substitutions to a pre-Wicks form $W$ give rise to the same group element, this forces both the substitutions to have degeneracies, in the sense that parts of the substitution are completely determined by some other parts. This, in turn, implies that such an event is not generic. More precisely:

\begin{lemma}\label{lem: element intersection count for reduction in degrees of freedom}
    Let $N\geq 1$ be a positive integer. Let $W$ and $W'$ be pre-Wicks forms of genus $n$, let $\vec{l}$ and $\vec{l}'$ be respective distinct length types with $||\vec{l}||=||\vec{l}'||=L$. If there's a reduction in degrees of freedom by \(N\) between two length types \(\vec{l}\) and \(\vec{l}'\), then there is a constant $C_3$ such that:
    \[
        \frac{|G_W(\vec{l})\cap G_{W'}(\vec{l}')|}{|G_{W}(\vec{l})|} \leq C_3 e^{-h_r N}
    \]
    where $h_r = \log{(2r-1)}$.
\end{lemma}

In order to see this, we will need to relate the abstract combinatorial reduction developed in \Cref{subsec: combinatorialreduction} to non-cancelling substitutions of a given length type, which we do in the following remark.

\begin{remark}[Relating combinatorial reduction to degeneracy in substitutions]\label{rmk: Relating combinatorial reduction to degeneracy in substitutions}
Recall \Cref{def: matching position} and \Cref{def: reduction in degrees of freedom}. Given $(W,\vec{l})$, and a non-cancelling substitution $\mathfrak{g}$ for it, the path $\gamma(\mathfrak{g}(W))$ is the geodesic segment of length $L$($=\|\vec{l}\|$) joining the identity element to $\mathfrak{g}(W)$. We also have the set of positions $\{1,\dots,L\}$ associated to $(W,\vec{l})$, which comes with a natural partition induced by $W$. Therefore, there is a natural one-to-one correspondence between the positions $\{1,\dots, L\}$ and $\gamma(\mathfrak{g}(W))$: for any position $i$, we denote by $x_i\in S$ the $i^{\text{th}}$ letter of $\mathfrak{g}(W)$.\\Now, suppose the letter $X_j$ contains a position $i$. Then, its matching position $\underline{i}$ is contained in the paired letter $X_{\bar{j}}$ and satisfies: 
    \begin{align}\label{eq: matching position equations}
        x_{\underline{i}} = x_i^{\epsilon_{X_j}}
    \end{align}
    where $\epsilon_{X_j} = \text{sign}(X_j)$. Recall that if $X_j$ has sign \(-1\), then the matching position $\underline{i}$ is defined as:
    \begin{equation}
        \underline{i} = t_{X_{\bar{j}}} + s_{X_j} - i
    \end{equation}
    In this case, $x_{\underline{i}} = x_i^{-1}$. And, if $X_j$ has sign $+1$, then the matching position $\underline{i}$ is defined as:
    \begin{equation}
        \underline{i} = s_{X_{\overline{j}}} + i - s_{X_{j}}
    \end{equation}
    In this case, we have $x_{\underline{i}} = x_i$. \\\\Next, suppose $\mathfrak{g},\mathfrak{g}'$ are two non-cancelling substitutions for $(W,\vec{l})$ giving rise to the same group element. So we have $||\vec{l}||=||\vec{l}'||(=L)$. For any position $i\in\{1,\dots,L\}$, let $x_i\in S$ and $x'_i\in S$ be the $i^{\text{th}}$ letters of $\mathfrak{g}(W)$ and $\mathfrak{g}'(W)$ respectively, just as defined above. \\We are now ready to see how the notion of reduction in degrees of freedom (cf. \Cref{def: reduction in degrees of freedom}) relates to the substitutions. Suppose there is a position $i$ contained in the letter $X$ with respect to $(W, \vec{l})$ and in $X'$ with respect to $(W', \vec{l}')$. Let $\varepsilon_X$, $\varepsilon_{X'}$ be the signs of $X$ and $X'$ respectively. 
Observe that if $\underline{i} \neq \underline{i}'$, we get a \textit{non-trivial equation} for $\mathfrak{g}$ and $\mathfrak{g}'$:
    \begin{equation}\label{eq: reduction equations - free}
        (x_{\underline{i}})^{\varepsilon_X}=(x_{\underline{i}'})^{\varepsilon_{X'}}, \quad (x'_{\underline{i}})^{\varepsilon_X}=(x'_{\underline{i}'})^{\varepsilon_{X'}}
    \end{equation}
thereby restricting the degree of freedom we have in choosing $\mathfrak{g}$ and \(\mathfrak{g}'\) by one letter. While these equations can be read off \Cref{fig: matching position}, let us also deduce them, formally. Firstly, we have: 
\begin{equation}\label{eq: matching position equations - free}
    x_{\underline{i}} = x_i^{\epsilon_X}\quad \text{ and }\quad  x'_{\underline{i}'} = (x'_i)^{\epsilon_{X'}}
\end{equation}
These equations are obtained by using \Cref{eq: matching position equations} for the position $i$ with respect to to each of the length types $\vec{l}$ and $\vec{l}'$. Here $\underline{i}$ and $\underline{i}'$ are the matching positions of $i$ with respect to $\vec{l}$ and $\vec{l}'$ respectively. Secondly, since we are in a free group, $\mathfrak{g}(W) = \mathfrak{g}'(W)$ implies that:

\begin{equation}\label{eq: admissibilty equations - free}
    x_i = x'_i,\quad x_{\underline{i}}=x'_{\underline{i}}, \quad x_{\underline{i}'} =  x'_{\underline{i}'}
\end{equation}
Now the desired equations follow.
\end{remark}

\begin{proof}[Proof of \Cref{lem: element intersection count for reduction in degrees of freedom}]
    Let $\mathfrak{g}$ and $\mathfrak{g}'$ be non-cancelling substitutions for $(W,\vec{l})$ and $(W,\vec{l}')$ respectively such that $\mathfrak{g}(W)=\mathfrak{g'}(W')$. As there is reduction in degrees of freedom by $N$, there is an admissible sequence of $N$ consecutive positions $\mathbf{i}=(i_1,\dots,i_N)$, such that:
    \begin{itemize}
        \item Both $\underline{\mathbf{i}}$ and $\underline{\mathbf{i}'}$ are admissible. 
        \item $\underline{i_1}\neq \underline{i_1}',\dots, \underline{i_N}\neq \underline{i_N}'$.
    \end{itemize}
    $\underline{\mathbf{i}}$ is admissible, so there are letters $X$ and $X'$ in $W$ and $W'$ respectively, such that $\underline{\mathbf{i}}$ is contained in $X$ with respect to $\vec{l}$, and in $X'$ with respect to  $\vec{l}'$. Following the discussion and notation in \Cref{rmk: Relating combinatorial reduction to degeneracy in substitutions}, we get $N$ \emph{non-trivial} equations for $\mathfrak{g}$ (cf. \Cref{eq: reduction equations - free}):
    \begin{equation*}
        (x_{\underline{i_1}})^{\varepsilon_X}=(x_{\underline{i_1}'})^{\varepsilon_{X'}},\dots, (x_{\underline{i_N}})^{\varepsilon_X}=(x_{\underline{i_N}'})^{\varepsilon_{X'}}
    \end{equation*}
    where $\varepsilon_X$ and $\varepsilon_{X'}$ denote the signs of the letters $X$ and $X'$ respectively. Note that the equations are completely determined by the sequence $\mathbf{i}$, the signs $\varepsilon_X$, $\varepsilon_{X'}$, and the length types $\vec{l}$, $\vec{l}'$. In particular, we have the same set of equations for the substitution $\mathfrak{g}'$.
    Observe that due to the equations, there exist $N$ distinct positions in $\underline{\mathbf{i}}\cup\underline{\mathbf{i}}'$ at which the substitution is completely determined by the values it takes at the rest of the positions. More precisely, the substitution $\mathfrak{g}$ must have one of the following degenerate algebraic forms:
    \begin{enumerate}[font=\normalfont]
        \item  There are two distinct unknowns $X$ and $Y$ of $W$ such that:
        \begin{equation*}
            g_X=uzv,\quad g_{Y}=u'z^{\varepsilon}v'
        \end{equation*}
        where $\varepsilon\in\{+1,-1\}$, and:
        \begin{equation*}
            l(z)=N,\quad l(u)+l(v)+N=l_X,\quad l(u')+l(v')+N=l_{Y}
        \end{equation*}
        \item  There is an unknown $X$ of $W$ such that:
        \begin{equation*}
            g_X=uzwz^{\varepsilon}v
        \end{equation*}
        where $\varepsilon\in\{+1,-1\}$, and:
        \begin{equation*}
            l(z)=N,\quad l(u)+l(w)+l(v)+2N=l_X
        \end{equation*}
        \item There is an unknown $X$ of $W$, and a positive integer $T\geq 1$ such that:
        \begin{equation*}
            g_X=uz_0^kz_1v
        \end{equation*}
        where $k=\left\lfloor\frac{N}{T}\right\rfloor+1$, $l(z_0) = T$, $l(z_1)<T$, and $z_{1}$ is a group element of length $N-(k-1)T$ which appears as the prefix in the reduced word representing $z_0$. Furthermore:
        \begin{equation*}
            l(u)+l(v)+l(z_1)+kT=l_X,\quad (k-1)T+l(z_1)=N
        \end{equation*}
    \end{enumerate}
    We derive these algebraic forms in the more general set-up of hyperbolic groups later (see proof of \Cref{prop: hyperbolic degenerate substitutions}).
    Now it follows from \Cref{lem: non-cancellingcount1}, that:
    \begin{align*}
        \frac{|G_W(\vec{l})\cap G_{W'}(\vec{l}')|}{|G_{W}(\vec{l})|}\leq \left(\frac{C'_2}{C'_1}\right)\frac{e^{h_r(L-N)}}{e^{h_rL}}=C'_3 e^{-h_r N}
    \end{align*}
    
\end{proof}

We are now in a position to complete the lower bound proofs for $\partial B^\pm _n(L)$, done in propositions \ref{thm: free grp shell lower bound asymptote} and \ref{thm: free grp shell lower bound asymptote non-orientable}. The proofs involve an inclusion-exclusion argument, the bounds for whose terms are in the two lemmas that follow. We use \Cref{thm: the main reduction theorem} along with the bounds derived in \Cref{lem: non-cancellingcount1} and \Cref{lem: element intersection count for reduction in degrees of freedom} to give proofs for these two lemmas.

\begin{lemma}\label{lem: non-cancellingcount2}
    Let $W$ be a maximal orientable pre-Wicks form of genus $n$. Then, for all $L\geq 1$, there exists a collection of length types $\mathcal{L}=\mathcal{L}(W,n,r,L)$ and there are constants $L_0,C'_3\geq 1$ depending only on $r$ and $n$ satisfying, for all $L\geq L_0$:
    \begin{enumerate}[font=\normalfont]
        \item For all $\vec{l}\in \mathcal{L}(W,n,r,L)$, $||\vec{l}||=L$ and \[\lim_{L\to\infty}\frac{|\mathcal{L}(W,n,r,L)|}{L^{6n-2}}=\frac{1}{2^{6n-2}(6n-2)!}\]
    \end{enumerate}
    Let $W'$ be another orientable pre-Wicks form of genus $n$, and let $\vec{l}\in \mathcal{L}(W,n,r,L)$, $\vec{l}'\in \mathcal{L}(W',n,r,L)$.
    \begin{enumerate}[font=\normalfont, resume]
        \item If $(W,\vec{l})\neq(W',\vec{l}')$, then:
        
        \[|G_W(\vec{l})\cap G_{W'}(\vec{l}')|\leq C'_3\frac{e^{\frac{h_rL}{2}}}{L^{6n}}\]
        
    \end{enumerate}
\end{lemma}

This lemma is proved together with its the non-orientable variant:

\begin{lemma}\label{lem: non-cancellingcount2 non-orient}
     Let $W$ be a maximal non-orientable pre-Wicks form of genus $n$. Then, for all $L\geq 1$, there exists a collection of length types $\mathcal{L}=\mathcal{L}(W,n,r,L)$ and there are constants $L_0,C'_3\geq 1$ depending only on $r$ and $n$ satisfying, for all $L\geq L_0$:
    \begin{enumerate}[font=\normalfont]
        \item For all $\vec{l}\in \mathcal{L}(W,n,r,L)$, $||\vec{l}||=L$ and \[\lim_{L\to\infty}\frac{|\mathcal{L}(W,n,r,L)|}{L^{3n-2}}=\frac{1}{2^{3n-2}(3n-2)!}\]
    \end{enumerate}
    Let $W'$ be another non-orientable pre-Wicks form of genus $n$, and let $\vec{l}\in \mathcal{L}(W,n,r,L)$, $\vec{l}'\in \mathcal{L}(W',n,r,L)$.
    \begin{enumerate}[font=\normalfont, resume]
        \item If $(W,\vec{l})\neq(W',\vec{l}')$, then:
        \begin{align*}
            |G_W(\vec{l})\cap G_{W'}(\vec{l}')|\leq C'_3\frac{e^{\frac{h_rL}{2}}}{L^{3n}}
        \end{align*}
        \end{enumerate}
\end{lemma}

\begin{proof}[Proof of \Cref{lem: non-cancellingcount2} and \Cref{lem: non-cancellingcount2 non-orient}]
    Suppose $W$ is a maximal orientable pre-Wicks form. Then $W$ is of length $12n-2$. Below, we think of length types for $W$ as $(6n-1)$-tuples (cf. \Cref{rmk: length type as an m-tuple}). Consider the following collection of length types for $W$:
    \begin{align*}
         \mathcal{L}=\mathcal{L}(W,n,r,L)\coloneq\left\{\vec{l}\in\mathbb{Z}^{6n-1}_+\middle|\quad  \sum_{i=1}^{6n-1}l_i = \frac{L}{2};\text{ for all }1\leq i\leq 6n-1, l_i\geq 6\left \lfloor \frac{6n\log L}{h_r}\right \rfloor+4\right\}
     \end{align*}
    Note that every pair of length types in $\mathcal{L}$ satisfies conditions \ref{item: L1 free} and \ref{item: L2 free}, with $c=\left \lfloor \frac{6n\log L}{h_r}\right \rfloor$ as stated in \Cref{cor: the main reduction theorem for free groups}. Also,
    \begin{align*}
         \lim_{L\to\infty}\frac{|\mathcal{L}|}{L^{6n-2}}&=\lim_{L\to\infty}\frac{1}{L^{6n-2}}\binom{\frac{L}{2} - (6n-1)\left (6\left \lfloor \frac{6n}{h_r}\log L\right \rfloor+4\right )}{6n-2}= \frac{1}{2^{6n-2}(6n-2)!}
     \end{align*}
    which proves $(1)$. Next, suppose $W'\neq W$ is another pre-Wicks form. Suppose $\vec{l}\in \mathcal{L}(W,n,r,L)$, $\vec{l}'\in \mathcal{L}(W',n,r,L)$. Suppose $\mathfrak{g}$ and $\mathfrak{g}'$ are two non-cancelling substitutions of length type $\vec{l}$ and $\vec{l}'$ respectively, which give rise to the same element $g$, i.e. $\mathfrak{g}(W)=\mathfrak{g}'(W')=g$. Now if either $\vec{l}\neq\vec{l}'$ or $W\neq W'$, using \Cref{cor: the main reduction theorem for free groups} for $c=\left \lfloor \frac{6n\log L}{h_r}\right \rfloor$, there has to be a reduction in degrees of freedom by $\left \lfloor \frac{6n\log L}{h_r}\right \rfloor$. Applying the bounds in \Cref{lem: non-cancellingcount1} and \Cref{lem: element intersection count for reduction in degrees of freedom} gives:
    \begin{align*}
        |G_W(\vec{l})\cap G_{W'}(\vec{l}')|\leq  C'_2e^{-h_r\left \lfloor \frac{6n\log L}{h_r}\right \rfloor }e^{\frac{h_rL}{2}}\leq C'_2\frac{e^{h_r\left (\frac{L}{2}+1 \right )}}{L^{6n}}= C'_2(2r-1)\frac{e^{\frac{h_rL}{2}}}{L^{6n}}
    \end{align*}
    which proves both $(2)$ and $(3)$.\\If instead, $W$ is a maximal non-orientable pre-Wicks form, then it has length $6n-2$ and the same arguments as above work for the collection:
    \begin{align*}
         \mathcal{L}'\coloneq\left\{\vec{l}\in\mathbb{Z}^{3n-1}_+\middle|\quad  \sum_{i=1}^{3n-1}l_i = \frac{L}{2};\text{ for all }1\leq i\leq 3n-1, l_i\geq 6\left \lfloor \frac{3n\log L}{h_r}\right \rfloor+4\right\}
     \end{align*}
\end{proof}

Here are the two final propositions. We only provide a proof of the former, since the latter's proof is nearly identical except for using \Cref{lem: non-cancellingcount2} with \Cref{lem: non-cancellingcount2 non-orient}.

\begin{proposition}\label{thm: free grp shell lower bound asymptote}
    Let $n$ and $r$ be positive integers with $r\geq 2$. Then:
    \begin{equation*}
        \liminf_{\substack{L\to\infty\\ L\textnormal{ even}}}\frac{|\partial B^{+}_n(L)|}{L^{6n-2}e^{\frac{h_rL}{2}}}\geq \frac{\tau^+_n(r) \mathcal{N}_{\text{pre}}^{+}(n)}{2^{6n-2}(6n-2)!(2r)^{6n-1}(2r-1)^{6n-1}}
    \end{equation*}
\end{proposition}

\begin{proposition}\label{thm: free grp shell lower bound asymptote non-orientable}
    Let $n$ and $r$ be positive integers with $r\geq 2$. Then:
    \begin{equation*}
        \liminf_{\substack{L\to\infty\\ L\textnormal{ even}}}\frac{|\partial B^{-}_n(L)|}{L^{3n-2}e^{\frac{h_rL}{2}}}\geq \frac{\tau^-_n(r) \mathcal{N}_{\text{pre}}^{-}(n)}{2^{3n-2}(3n-2)!(2r)^{3n-1}(2r-1)^{3n-1}}
    \end{equation*}
\end{proposition}

\begin{proof}[Proof of \Cref{thm: free grp shell lower bound asymptote}]
    Let $W$ be a maximal orientable pre-Wicks form of genus $n$, thus having length \(12n - 2\) and \(6n-1\) unknowns (by \Cref{lem: length of a preWicks form}). Let $\mathcal{L}(W,n,r,L), L_0$ and $C'_3$ be as in \Cref{lem: non-cancellingcount2}. For any $\vec{l}\in\mathcal{L}$, a group element in $\partial G_W(\vec{l}\,)$ has genus at most $n$. That is we have:
    \[\bigcup_{\substack{W\\\text{maximal}\\\text{orientable}\\\text{genus }n}}\bigcup_{\vec{l} \in \mathcal{L}}G_W(\vec{l}\,)\subset \partial B_{\leq n}^+(L)\]
    
    where \(\partial B_{\leq n}^+(L)\coloneq \bigcup_{i=1}^n\partial B^{+}_i(L)\) is the set of all elements in \(F_r\) having orientable genus lesser than or equal to \(n\) and length equal to $L$. In what follows, we will often hold back from writing all arguments in $\mathcal{L}(W,n,r,L)$ when they are clearly understood.
    So we have:
    \begin{align*}
        |\partial B_{\leq n}^+(L)|\geq \sum_W\sum_{\vec{l}\in \mathcal{L}(W)}|G_W(\vec{l})| - \frac{1}{2}\sum_{\substack{W,W'\\\vec{l}\in\mathcal{L}(W),\vec{l}'\in \mathcal{L}(W')\\(W,\vec{l})\neq(W',\vec{l}')}}|G_W(\vec{l})\cap G_W(\vec{l}')|
    \end{align*}
    
    Using \Cref{lem: non-cancellingcount2} and the fact that for any two maximal orientable pre-Wicks forms of the same genus, $|\mathcal{L}(W,n,r,L)|=|\mathcal{L}(W',n,r,L)| \coloneq |\mathcal{L}|$ :
    \begin{align*}
        \frac{1}{L^{6n-2}e^{\frac{h_rL}{2}}}\sum_{\substack{W,W'\\\vec{l}\in\mathcal{L}(W),\vec{l}'\in \mathcal{L}(W')\\(W,\vec{l})\neq(W',\vec{l}')}}|G_W(\vec{l})\cap G_{W'}(\vec{l}')|&\leq \frac{1}{L^{6n-2}e^{\frac{h_rL}{2}}}\cdot (\mathcal{N}_{\text{pre}}^+(n))^2|\mathcal{L}|^2C'_3 \frac{e^{\frac{h_rL}{2}}}{L^{6n}}\\
        &=C'_3(\mathcal{N}_{\text{pre}}^+(n))^2 \frac{|\mathcal{L}|^2}{L^{12n-2}}=\frac{C'_3(\mathcal{N}_{\text{pre}}^+(n))^2e^{-O(n\log n)}}{L^2}
    \end{align*}
    \begin{equation*}
        \text{Therefore, }\lim_{L\to\infty}\frac{1}{L^{6n-2}e^{\frac{h_rL}{2}}}\sum_{\substack{W,W'\\\vec{l}\in\mathcal{L}(W),\vec{l}'\in \mathcal{L}(W')\\(W,\vec{l})\neq(W',\vec{l}')}}|G_W(\vec{l})\cap G_{W'}(\vec{l}')|=0
    \end{equation*}
    Also, using \Cref{lem: substitution count asymptote}, \Cref{lem: tau for maximal pre-Wicks forms}, and the fact that $\lim_{L\to\infty}\frac{|\mathcal{L}(W,n,r,L)|}{L^{6n-2}}=\frac{1}{2^{6n-2}(6n-2)!}$, we have:

    \begin{equation*}
        \lim_{L\to\infty}\frac{1}{L^{6n-2}e^{\frac{h_rL}{2}}} \sum_W\sum_{\vec{l}\in \mathcal{L}(W)}|G_W(\vec{l})|= \frac{\tau^+_n(r) \mathcal{N}_{\text{pre}}^{+}(n)}{2^{6n-2}(6n-2)!(2r)^{6n-1}(2r-1)^{6n-1}}
    \end{equation*}
    \begin{align*}
        \text{Therefore, } \liminf_{L\to\infty}\frac{|\partial B^{+}_{\leq n}(L)|}{L^{6n-2}e^{\frac{h_rL}{2}}}\geq \frac{\tau^+_n(r) \mathcal{N}_{\text{pre}}^{+}(n)}{2^{6n-2}(6n-2)!(2r)^{6n-1}(2r-1)^{6n-1}}
    \end{align*}
    Note that by \Cref{thm: free grp shell upper bound asymptote},
    \begin{align*}
        \text{For }i<n, \quad \lim_{L\to\infty}\frac{|\partial B_i^+(L)|}{L^{6n-2}e^{\frac{h_rL}{2}}}=0
    \end{align*}
    As  \(\partial B_{\leq n}^+(L)\coloneq \bigcup_{i=1}^n\partial B^{+}_i(L)\), we have:
    \begin{align*}
        \liminf_{L\to\infty}\frac{|\partial B^{+}_n(L)|}{L^{6n-2}e^{\frac{h_rL}{2}}}\geq \frac{\tau^+_n(r) \mathcal{N}_{\text{pre}}^{+}(n)}{2^{6n-2}(6n-2)!(2r)^{6n-1}(2r-1)^{6n-1}}
    \end{align*}
\end{proof}

\subsection{Counting in Hyperbolic groups}\label{subsec: Element count: Lower bounds: Hyperbolic groups}
Throughout this section, $G$ will denote a hyperbolic group. We will equip $G$ with the word metric induced by some symmetric generating set $S$.

\begin{definition}\label{def: chainsubstitution}
    Let $C,D\geq 0$, $W=X_1\dots X_m$ be a word, and $\mathfrak{g}$ be a $G$-valued substitution. Let $x_0=1_{G}$, $x_i=g_{X_1}\dots g_{X_i}$ for $1\leq i\leq m$. $\mathfrak{g}$ will be called a \emph{$(C,D)$-chain substitution} for $W$ if the sequence $x_0,\dots,x_m$ is a $(C,D)$-chain (cf. \Cref{def: alignedchain}).
\end{definition}

\begin{lemma}\label{lem: goodsubstitutionsexist}
     Let $W$ be a pre-Wicks form of length $2m$. Denote the unknowns in $W$ by $X_1,\dots, X_m$. Then there exists a constant $K>0$, such that for every $R_0\geq3(K+\delta)$, there is a $(K,R_0)$-Schottky set $\mathcal{S}(R_0)$, such that for every $R_1\geq R_0$ and any $G$-valued substitution $\mathfrak{g}$ with $l(g_{X_i})\geq R_1$, there exist $l_i,r_i\in \mathcal{S}(R_0)$ for all $1\leq i\leq m$, so that the substitution $\mathfrak{g}'$ defined by:
    \begin{align*}
        \text{for all }1\leq i\leq m, \quad g'_{X_i}=l_ig_{X_i}r_i 
    \end{align*}
    is a $(K+2\delta,R_1)$-chain substitution. Moreover:
    \begin{align*}
        \text{for all }1\leq i\leq m, \quad l(g_{X_i})\leq l(g'_{X_i})\leq l(g_{X_i})+ 2K'
    \end{align*}
    where $K'=K'(R_0)>0$ is the length of the longest element in $\mathcal{S}(R_0)$.
\end{lemma}

\begin{proof}
    The lemma follows from a standard application of \Cref{prop: existenceofSchottkyset}, \Cref{rmk: rmktoexistenceofSchottkyset}, and \Cref{lem: alignmentandprogress}. The idea is as follows: given an arbitrary substitution $\mathfrak{g}$, there could be uncontrolled ``cancellation" between consecutive letters of the Wicks form. But this can be remedied by multiplying the group element corresponding to each letter on the left and right with Schottky elements. These Schottky elements need to satisfy a number of alignment conditions and we can arrange a large enough Schottky set so that this happens.\\ Recalling that the unknowns of $W$ are $X_1,\dots, X_m$, we write $W=Y_1\dots Y_{2m}$ where each $Y_i\in\{X_k,X_k^{-1}\mid 1\leq k\leq m \}$. For any choice of $\{l_i\}_{1\leq i\leq m},\{r_i\}_{1\leq i\leq m}$, define the sequence of points: $\{y_{i}\}_{0\leq i\leq 6m}$ as follows: Set $y_0=1_G$. For each $i\in\{1 ,\dots 2m\}$, $Y_i=X_k$ or $X_k^{-1}$. Say $Y_i=X_k$, then:
    \begin{equation*}
        y_{3i-2}\coloneq l_k\cdot y_{3i},\quad y_{3i-1}\coloneq l_kg_{X_k}\cdot y_{3i},\quad y_{3i+3}\coloneq l_kg_{X_k}r_k\cdot y_{3i} 
    \end{equation*}
    and if $Y_i=X_k^{-1}$, then:
    \begin{equation*}
        y_{3i-2}\coloneq r_k^{-1}\cdot y_{3i},\quad y_{3i-1}\coloneq r_k^{-1}g_{X_k}^{-1}\cdot y_{3i},\quad y_{3i}\coloneq r_k^{-1}g_{X_k}^{-1}l_k^{-1}\cdot y_{3i} 
    \end{equation*}
    
    \begin{claim}\label{claim: numberofgoodsubs}
        There exists $K>0$ such that for every $R_0>0$ there is a Schottky set $\mathcal{S}(R_0)$ such that for every substitution $\mathfrak{g}$ for $W$ satisfying $l(g_{X_i})\geq R_0$ for all $1\leq i\leq m$, there exist elements $\{l_i\}_{1\leq i\leq m},\{r_i\}_{1\leq i\leq m}\subset \mathcal{S}(R_0)$ such that the sequence $\{y_i\}_{0\leq i\leq 6m}$ is a $(K,R_0)$-chain.
    \end{claim}
    \begin{proof}
        By \Cref{prop: existenceofSchottkyset}, and \Cref{rmk: rmktoexistenceofSchottkyset}, there exists a constant $K>0$ (depending only on the group $G$ and the choice of generating set), such that for every $R_0>0$ we have a $(\frac{1}{8},K,R)$-Schottky set $\mathcal{S}(R_0)$ of size $16$. In choosing $l_i$ (and $r_i$ resp.), we need to satisfy at most $6$ alignment conditions. We demonstrate this as follows: Consider two consecutive letters of $W$, $Y_k$ and $Y_{k+1}$ and suppose $Y_k=X_i$ and $Y_{k+1}=X_j$. For a choice of $l_i,l_j,r_i,r_j\in\mathcal{S}$, the sequence $y_k,y_{k+1},\dots, y_{k+6}$ will be a $(K,R_0)$-chain if and only if the following families of alignment conditions are satisfied: 
        % (\Cref{fig: aligned points}):
        \linkdest{alignmentconditions}
    \begin{enumerate}
        \item[(F1)] $(l_i^{-1},g_{X_i})_{1_G}\leq K$, $(l_j^{-1},g_{X_j})_{1_G}\leq K$
        \item[(F2)] $(g_{X_i}^{-1},r_i)_{1_G}\leq K$, $(g_{X_j}^{-1},r_j)_{1_G}\leq K$
        \item[(F3)] $(r_i^{-1},l_j)_{1_G}\leq K$
    \end{enumerate}
        There are analogous families of conditions for the cases when: $X_i, X_j^{-1}$ are adjacent, $X_i^{-1},X_j$ are adjacent, and $X_i^{-1},X_j^{-1}$ are adjacent. Depending on the specific algebraic form of the Wicks form $W$, the appropriate conditions are assumed for each family. As $W$ is quadratic, evidently there are at most six alignment conditions that each $l_i$ and $r_i$ needs to satisfy. By \Cref{rmk: rmktoexistenceofSchottkyset}, it follows that we have at least $4$ choices each for $l_i$ and $r_i$ respectively.   
    \end{proof}
    Now let $R_0, R_1$ be positive integers such that $R_1\geq R_0\geq 3(K+\delta)$, and let $\mathfrak{g}$ be a substitution for $W$ such that $l(g_{X_i})\geq R_1$ for all $1\leq i\leq m$. We show that the choices of $\{l_i\}_{1\leq i\leq m},\{r_i\}_{1\leq i\leq m}\subset \mathcal{S}(R_0)$ made in \Cref{claim: numberofgoodsubs} ensure that the substitution $\mathfrak{g}'$ as defined in the statement of the lemma, is a $(K+2\delta,R_1)$-substitution.\\Suppose that the unknowns $X_i$ and $X_j$ are adjacent in the Wicks form (completely analogous arguments can be made in the other three cases), that is $Y_k=X_i$ and $Y_{k+1}=X_j$ for some $k$.

    % \begin{figure}[h]
    %     \centering
    %     \includegraphics[width=0.85\columnwidth]{aligned_points}
    %     \caption{}
    %     \label{fig: aligned points}
    % \end{figure}
    % {\color{red} Need to relabel this figure}
    
    The sequence $y_{3(k-1)},y_{3k-2},\dots, y_{3(k+1)}$ is a $(K,R_0)$-chain. So it follows from \Cref{lem: alignmentandprogress}, that:
        \begin{align*}
            l(g'_{X_i})=l(l_ig_ir_i)=d(y_{3(k-1)},y_{3k})&\geq \sum_{s=3k-3}^{3k-1}(d(y_{s},y_{s+1})-(2K+2\delta))\\
            &\geq (2R_0-3(2K+2\delta)) + l(g_{X_i})\\
            &\geq l(g_{X_i})\geq R_1 \\
            (y_{3(k-1)},y_{3(k+1)})_{y_{3k}}&\leq K+2\delta
        \end{align*}
    Also, by triangle inequality, we have $l(g'_{X_i})\leq  l(g_{X_i})+2K'$ where $K'=K'(R_0)$ is the length of the longest element in $\mathcal{S}(R_0)$.
\end{proof}

\begin{definition}[Good Substitutions]\label{def: goodsubs}
    Let $R_0$ and $R_1$ be positive integers satisfying $3(K+\delta)\leq R_0\leq R_1$. Let $K, K'(R_0),$ and $\mathcal{S}(R_0)$ be as in in \Cref{lem: goodsubstitutionsexist}. Let $W$ be a pre-Wicks form of length $2m$, and let $X_1,\dots,X_m$ be its unknowns. Let $\vec{l}=(l_1,\dots,l_m)$ be a length type for $W$ with $l_i\geq R_1$ for all $1\leq i\leq m$. Set $L=2\sum_{i=1}^ml_i$.\\We define a map $\Phi=\Phi(\vec{l}):\prod_{i=1}^m\partial B(l_i)\rightarrow B(L+4mK')$, where $\partial B(l_i)$ denotes the set of group elements having length exactly $l_i$ as follows: Any $(g_1,\dots,g_m)\in \prod_{i=1}^m\partial B(l_i)$, defines a substitution $\mathfrak{g}$ for $W$ of length type $\vec{l}$. By \Cref{lem: goodsubstitutionsexist}, there exist $\{l_i\}_{1\leq i\leq m},\{r_i\}_{1\leq i\leq m}\subset\mathcal{S}(R_0)$ which defines the substitution $\mathfrak{g}'$: for all $1\leq i\leq m$,  $g'_{X_i}=l_ig_{X_i}r_i$. Define $\Phi(\mathfrak{g})\coloneq \mathfrak{g'}(W)$ and $G_W(\vec{l})\coloneq \Phi(\prod_{i=1}^m\partial B(l_i))$. We will say that $\mathfrak{g'}$ is the \emph{good substitution} associated with $\mathfrak{g}$. $G_W(\vec{l})$ is thus the set of group elements obtained from good substitutions for the length type $\vec{l}$. Note that the good substitution $\mathfrak{g}'$ satisfies the following properties:
    \begin{enumerate}
        \item The alignment conditions \linkref{alignmentconditions}{(F1)-(F3)} and other analogous ones.
        \item $\mathfrak{g}'$ is a $(K+2\delta,R_1)$-chain substitution.
        \item $\text{for all }1\leq i\leq m$, $l(g_{X_i})\leq l(g'_{X_i})\leq l(g_{X_i})+ 2K'$.
    \end{enumerate}
\end{definition}

\begin{remark}\label{rmk: good subs for Wicks forms are cyclically good}
    Let us observe that the good substitutions for Wicks forms are in fact \emph{cyclically good}: Let $W$ be a Wicks form. Let $\mathfrak{g}$ be a substitution for $W$ of length type $\vec{l}$, and let $\mathfrak{g}'$ be the associated good substitution. Then, a quick inspection of the proof of \Cref{lem: goodsubstitutionsexist} helps observe that for any cyclic permutation $\sigma$ of $W$, $\mathfrak{g}'$ is a good substitution for $\sigma(W)$ as well, that is, the sequence: $x_0=1_{G}$, $x_{i}=g'_{X_{\sigma(1)}}\cdots g'_{X_{\sigma(i)}}$ is a $(K+2\delta,R_1)$ chain.
\end{remark}

\Cref{def: goodsubs} depends on the choice of a $R_0\geq 3(K+\delta)$. We will choose $R_0$ in \Cref{construction: geodesicencoding}.

\begin{lemma}\label{lem: lowerboundgoodsubs}
    Let $R_1\geq R_0$ be a positive integer. Let $W=X_1\dots X_{2m}$ be a pre-Wicks form of length $2m$ and let $\vec{l}$ be a length type for $W$ with $l_i\geq R_1$ for all $1\leq i\leq 2m$. Then, there exists $C'_1>0$ depending only on $m$, $\delta$, $G$ and the choice of the generating set, so that:
    \begin{align*}
        |G_W(\vec{l})|\geq C'_1e^{\frac{h||\vec{l}||}{2}}
    \end{align*}
    where $h$ is the volume entropy of $G$.
\end{lemma}

\begin{proof}
    In view of \Cref{rmk: length type as an m-tuple} and \Cref{def: goodsubs}, we have:
    \begin{align*}
        \prod_{i=1}^m\partial B(l_i)=\bigsqcup_{g\in G_W(\vec{l})} \left (\Phi|_{\vec{l}}\right )^{-1}(g)
    \end{align*}
    Therefore, in order to obtain the desired lower bound for $\left|G_W(\vec{l})\right|$, it suffices to obtain a uniform upper bound on the fibers, i.e.  $\left|\left (\Phi|_{\vec{l}}\right )^{-1}(g)\right|\leq C$ for some $C>0$. On that note, let $\mathfrak{g}_1$ and $\mathfrak{g}_2$ be two substitutions of length type $\vec{l}$ such that their associated good substitutions give rise to the same group element, i.e. $\mathfrak{g}'_1(W)=\mathfrak{g}'_2(W)=g$.  Fix a geodesic $\gamma_g$ joining the identity to $g$. $\mathfrak{g}_i$ is a $(K+2\delta, R)$-chain substitution, so by \Cref{prop: chainfellowtravelling}, $\gamma(\mathfrak{g}_i)$ is contained in a $(K+32\delta)$-neighbourhood of $\gamma_g$. This allows us to construct the following sequence of points: Let $x_0^{(i)}\coloneq 1_G, x_{2m}^{(i)}\coloneq \tau(\gamma_g)$. For all $0<j<2m$, let $x_j^{(i)}$ be a vertex on $\gamma_g$ such that $d(x_j^{(i)},\tau(\gamma^{(i)}_{X_j}))\leq K+32\delta + 1$.\\There is an upper bound for the maximum distance between $x^{(1)}_j$ and $x^{(2)}_j$.
    \begin{claim}
        $\max\{d(x_j^{(1)},x_j^{(2)})|1\leq j\leq 2m\}\leq (4m-2)(K+K'+32\delta+1)$
    \end{claim}
    \begin{proof}
        This follows directly from property $(3)$ of a good substitution (cf. \Cref{def: goodsubs}) and triangle inequalities.
    \end{proof}
    It follows that for any $X_j$, we have:
    \begin{align*}
        d(\tau(\gamma^{(1)}_{X_i}),\tau(\gamma^{(2)}_{X_i})), d(i(\gamma^{(1)}_{X_i}),i(\gamma^{(2)}_{X_i}))&\leq (4m-2)(K+K'+32\delta+1)+2(K+32\delta+1)\\
        &\leq 4m(K+K'+32\delta+1)(\coloneq K_0)
    \end{align*}
    So for every letter $X_j$, there exist $w_{j},w_{j}'\in B(K_0)$, such that $g^{(1)}_{X_j}=w_jg^{(2)}_{X_j}w'_j$. As the substitution is determined by $m$ unknowns, we have:
    \begin{align*}
        \forall g\in\Phi(G_W(\vec{l})),\quad \left|\left (\Phi|_{G_W(\vec{l})}\right )^{-1}(g)\right |\leq (|B(K_0)|)^{m}
    \end{align*}
     Using \Cref{thm: ballgrowth}, there are constants $C_1,C_2, h>0$ so that: $C_1 e^{hL}\leq |B(L)|\leq C_2 e^{hL}$, where $h$ is the volume entropy of $G$. Thus, we get:
     \begin{align*}
         \forall g\in\Phi(G_W(\vec{l})),\quad \left|\left (\Phi|_{G_W(\vec{l})}\right )^{-1}(g)\right |\leq C_2^{m}e^{mhK_0}
     \end{align*}
     Due to the non-amenability of $G$, there exists a constant $c_1$ depending only on $G$ and the choice of generating set, so that $|\partial B(L)|\geq c_1|B(L)|$, for all $L\geq 1$. So we have, $|\prod_{i=1}^m\partial B(l_i)|\geq (c_1C_1)^me^{h\sum_{i=1}^ml_i}$, and therefore:
     \begin{align*}
         |G_W(\vec{l})|\geq\left\{ \left(\frac{c_1C_1}{C_2}\right)^me^{-mhK_0}\right\}e^{\frac{h||\vec{l}||}{2}}
     \end{align*}
\end{proof}

\begin{remark}
    For the rest of this section, whenever we write the triple(s) $(W,\mathfrak{g}_1,\vec{l})$ and/or $(W',\mathfrak{g}_2,\vec{l}')$, we will mean that $W=X_1\cdots X_{2m}$ and $W'=X'_1\cdots X'_{2m'}$ are pre-Wicks forms of lengths $2m$ and $2m'$ respectively. $\vec{l}$ and $\vec{l}'$ are length types for $W$ and $W'$ respectively. $\mathfrak{g}_1$ and $\mathfrak{g}_2$ are substitutions for $W$ and $W'$ of length types $\vec{l}$ and $\vec{l}'$ respectively. Moreover, for the rest of this subsection we fix $R_0$ and $R_1$ so that: $R_0\leq R_1$, $R_0\geq 5(K+30\delta)$ and $R_1\geq 12(K+30\delta)$. The lower bounds on $R_0$ and $R_1$ are sufficient to ensure that \Cref{construction: geodesicencoding} can be done.
\end{remark}

\begin{lemma}\label{lem: hyperbolic element intersection count for reduction in degrees of freedom}
    Let $N$ be a positive integer. Let $W$ and $W'$ be maximal pre-Wicks forms of length $2m$ and $2m'$ respectively. There exists a positive integer $D_0=D_0(K,R_0,h,\delta,m,m')$ such that the following is true: let $\vec{l}$ and $\vec{l}'$ be any length types for $W$ and $W'$ respectively, such that:
    \begin{enumerate}[font=\normalfont]
        \item $||\vec{l}||=||\vec{l}'||(\coloneq L)$
        \item For all $1\leq i\leq m$, $1\leq j\leq m'$, $l_i,l'_j>\max\left\{6N+5D_2+2,\frac{R_1}{D_0}\right\}$
        \item $(W,\vec{l})\neq (W',\vec{l}')$
    \end{enumerate}
    Then, there is a constant $C'_2>0$ depending only on $m,\delta, K,h,D_0$ such that:
    \begin{equation*}
        \frac{|G_W(10D_0\cdot\vec{l})\cap G_{W'}(10D_0\cdot\vec{l}')|}{|G_W(10D_0\cdot\vec{l})|}\leq C'_2L^2e^{-\frac{hND_0}{2}}
    \end{equation*}
    where $h$ is the volume entropy of $G$.
\end{lemma}

\begin{lemma}\label{lem: hyperbolic element intersection count}
    Let $n$ and $L$ be positive integers. Let $D_0$ be the constant in \Cref{lem: hyperbolic element intersection count for reduction in degrees of freedom}. There exist positive constants $L_0$, $C'_3$ depending on $K,\delta,n,R_0,R_1,h$ such that for all $L\geq L_0$, for any maximal orientable pre-Wicks form $W$ of genus $n$, there is a collection of length types $\mathcal{L}=\mathcal{L}(W,n,L)$, satisfying:
    \begin{enumerate}[font=\normalfont]
        \item $|\mathcal{L}| = e^{-O(n\log n)}L^{6n-2}$ and $||\vec{l}||=2\left\lceil \frac{L}{2}\right\rceil$.
    \end{enumerate}
    Let $W'$ be another orientable pre-Wicks form of genus $n$, and let $\vec{l}\in \mathcal{L}(W,n,L)$, $\vec{l}'\in \mathcal{L}(W',n,L)$.
    \begin{enumerate}[font=\normalfont, resume]
        \item If $(W,\vec{l})\neq(W',\vec{l}')$, then for all $L\geq L_0$:
        \begin{align*}
            |G_W(10D_0\cdot\vec{l})\cap G_{W'}(10D_0\cdot\vec{l}')|\leq C'_3\frac{e^{\frac{h\cdot10D_0L}{2}}}{L^{6n}}
        \end{align*}
        where $h$ is the volume entropy of $G$.
    \end{enumerate}
\end{lemma}
\begin{proof}
     Suppose $W$ is a maximal orientable pre-Wicks form. Then $W$ is of length $12n-2$. Below, we think of length types for $W$ as $(6n-1)$-tuples (cf. \Cref{rmk: length type as an m-tuple}). Consider the following collection of length types for $W$:
    \begin{align*}
         &\mathcal{L}(W,n,L)\\
         &\coloneq\left\{\vec{l}\in\mathbb{Z}^{6n-1}_+\middle|\; \sum_{i=1}^{6n-1}l_i = \left\lceil\frac{L}{2}\right\rceil; l_i\geq \max\left\{ 6\left \lceil \frac{2(6n+2)\log L}{hD_0}\right \rceil+5D_2+2,\frac{R_1}{D_0}\right\},\;(1\leq i\leq 6n-1)\right\}
     \end{align*}
    We choose $L_0$ so that $\mathcal{L}(W,n,L)$ is non-empty, for instance the following condition suffices:
    \[\text{for all }L\geq L_0,\quad L>100(12n-2)\max\left\{ 6\left \lceil \frac{2(6n+2)\log L}{hD_0}\right \rceil+5D_2+2,\frac{R_1}{D_0}\right\}\]
    Then, we have:
    \begin{align*}
         |\mathcal{L}|&= \binom{\left\lceil\frac{L}{2}\right\rceil - (6n-1)\max\left\{ 6\left \lceil \frac{2(6n+2)\log L}{hD_0}\right \rceil+5D_2+2,\frac{R_1}{D_0}\right\}}{6n-2}\leq \frac{L^{6n-2}}{2^{6n-2}(6n-2)!}
     \end{align*}
    which proves $(1)$. Let $W'$ be another maximal orientable pre-Wicks form of genus $n$. Let $\vec{l}\in \mathcal{L}(W,n,L)$ and $\vec{l}'\in \mathcal{L}(W',n,L)$ such that $(W,\vec{l})\neq(W',\vec{l}')$. $(W,\vec{l})$ and $(W',\vec{l}')$ satisfy the hypotheses of \Cref{lem: hyperbolic element intersection count for reduction in degrees of freedom} for $N=\left \lceil \frac{2(6n+2)\log L}{hD_0}\right \rceil$. Therefore, using \Cref{lem: lowerboundgoodsubs} and \Cref{lem: hyperbolic element intersection count for reduction in degrees of freedom}:
    \[G_W(10D_0\cdot\vec{l})\cap G_{W'}(10D_0\cdot \vec{l})\leq (C'_1 e^{h\cdot 10D_0\cdot\lceil\frac{L}{2}\rceil})\left(C'_2L^2e^{-\frac{hD_0}{2}\times \left\lceil \frac{2(6n+2)\log L}{hD_0}\right \rceil }\right)\leq (C'_1C'_2e^{10hD_0})\frac{e^{\frac{h\cdot10D_0L}{2}}}{L^{6n}}\]
\end{proof}

\begin{theorem}\label{thm: lowerboundhyp}
    Let $G$ be a hyperbolic group equipped with the word metric induced by some finite symmetric generating set $S$. Let $h$ denote the volume entropy of $G$. There exist positive constants $C'_4,L_0$ depending only on $n$, $G$, and the choice of the generating set, such that:
    \begin{equation*}
        \forall L\geq L_0, |B^{+}_n(L)|\geq C'_4 L^{6n-2}e^{hL/2}
    \end{equation*}
\end{theorem}

\begin{proof}
    Let $D_0$ be the constant in \Cref{lem: hyperbolic element intersection count for reduction in degrees of freedom}. Let $L_0,C'_1,C'_3$ be as in \Cref{lem: lowerboundgoodsubs} and \Cref{lem: hyperbolic element intersection count}. For every $L\geq L_0$, and every maximal orientable pre-Wicks form $W$ of genus $n$, let $\mathcal{L}_W(L)$ denote the collection of length types $\mathcal{L}(W,n,\lceil L/10D_0\rceil)$ (as defined in \Cref{lem: hyperbolic element intersection count for reduction in degrees of freedom}). So for every $\vec{l}\in\mathcal{L}(W,L)$, $10D_0\cdot ||\vec{l}||=20D_0\lceil\frac{1}{2}\lceil\frac{L}{10D_0}\rceil\rceil$, and $L\leq 10D_0\cdot ||\vec{l}||\leq L+20D_0$. Thus, it follows from \Cref{lem: goodsubstitutionsexist} that:
    \[\bigcup_{\substack{W\\\text{maximal}\\\text{orientable}\\\text{genus }n}}\bigcup_{\vec{l}\in \mathcal{L}_W(L)} G_W(10D_0\cdot\vec{l}\,)\subseteq B^{+}_{\leq n}(L+2(12n-2)K'+20D_0)\]
    \begin{align*}
        &|B_{\leq n}^{+}(L+2(12n-2)K'+20D_0)|\\
        &\geq \sum_{W}\sum_{\vec{l}\in \mathcal{L}_W(L)}|G_W(10D_0\cdot\vec{l}\,)| - \frac{1}{2}\sum_{W}\sum_{\vec{l}\in \mathcal{L}_W(L)}\sum_{W'}\sum_{\vec{l}'\in \mathcal{L}_{W'}(L)}|G_W(10D_0\cdot\vec{l}\,)\cap G_{W'}(10D_0\cdot\vec{l}'\,)|
    \end{align*}
    Note that for all maximal orientable Wicks forms $W,W'$, we have $|\mathcal{L}_{W}(L)|=|\mathcal{L}_{W'}(L)|(\coloneq \mathcal{L})$. Using \Cref{lem: lowerboundgoodsubs} and \Cref{lem: hyperbolic element intersection count}, we have:
    \begin{align*}
        &\sum_{W}\sum_{\vec{l}\in \mathcal{L}_W(L)}|G_W(10D_0\cdot\vec{l}\,)|\geq C'_1\mathcal{N}_{\text{pre}}^{+}(n)\mathcal{L}e^{\frac{hL}{2}}\\
        &\sum_{W}\sum_{\vec{l}\in \mathcal{L}_W(L)}\sum_{W'}\sum_{\vec{l}'\in \mathcal{L}_{W'}(L)}|G_W(10D_0\cdot\vec{l}\,)\cap G_{W'}(10D_0\cdot\vec{l}'\,)|\leq (10D_0)^{6n}C'_3e^{10hD_0}(\mathcal{N}_{\text{pre}}^{+}(n))^2\mathcal{L}^2\frac{e^{\frac{hL}{2}}}{L^{6n}}\\
        \text{So, }\quad &\frac{\sum_{W,\vec{l},W',\vec{l}'}|G_W(10D_0\cdot\vec{l}\,)\cap G_{W'}(10D_0\cdot\vec{l}'\,)|}{\sum_{\vec{l}\in \mathcal{L}_W(L)}|G_W(10D_0\cdot\vec{l}\,)|}\leq \left(\frac{100C'_3D_0^2e^{10hD_0}}{2^{6n-2}C'_1}\frac{\mathcal{N}_{\text{pre}}^{+}(n)}{(6n-2)!}\right)\times\frac{1}{L^2}
    \end{align*}
    where $\mathcal{N}_{\text{pre}}^{+}(n)$ denotes the number of maximal orientable pre-Wicks forms of genus $n$. Thus, for large $L$, we have:
    \begin{align*}
        &|B_{\leq n}^+(L+2(12n-2)K'+20D_0)|\\
        &\geq \frac{1}{2}\sum_{W}\sum_{\vec{l}\in \mathcal{L}_W(L)}|G_W(10D_0\cdot\vec{l}\,)|\geq \frac{1}{2}C'_1\mathcal{N}_{\text{pre}}^{+}(n)\mathcal{L}e^{\frac{hL}{2}}\geq\left(\frac{C'_1\mathcal{N}_{\text{pre}}^{+}(n)e^{-O(n\log n)}}{(10D_0)^{6n-2}}\right)\times L^{6n-2}e^{\frac{hL}{2}}
    \end{align*}
    Note that by \Cref{thm: element count upper bound orientable hyperbolic}, \(B_i^+(L)\) for each \(i \in\{0, \cdots, n-1\}\) are smaller than \(B_{\leq n}^+(L)\) by a multiplicative factor which is polynomial in \(1/L\). Therefore the above estimate continues to hold for $B^+_n(L)$, up to possibly modifying \(L_0\).
\end{proof}

\begin{lemma}\label{lem: hyperbolic element intersection count non-orientable}
    Let $n$ and $L$ be positive integers. Let $D_0$ be the constant in \Cref{lem: hyperbolic element intersection count for reduction in degrees of freedom}. There exist positive constants $L_0=L_0(K,\delta,n,R_0,R_1)$, $C'_3=C'_3(K,\delta,n,R_0,h)$ such that for all $L\geq L_0$, for any maximal non-orientable pre-Wicks form $W$ of genus $n$, there is a collection of length types $\mathcal{L}=\mathcal{L}(W,n,L)$, satisfying:
    \begin{enumerate}[font=\normalfont]
        \item $|\mathcal{L}| = e^{-O(n\log n)}L^{3n-2}$ and $||\vec{l}||=2\left\lceil\frac{L}{2}\right\rceil$.
    \end{enumerate}
    Let $W'$ be another non-orientable pre-Wicks form of genus $n$, and let $\vec{l}\in \mathcal{L}(W,n,L)$, $\vec{l}'\in \mathcal{L}(W',n,L)$.
    \begin{enumerate}[font=\normalfont, resume]
        \item If $(W,\vec{l})\neq(W',\vec{l}')$, then for all $L\geq L_0$:
        \begin{align*}
            |G_W(10D_0\cdot\vec{l})\cap G_{W'}(10D_0\cdot\vec{l}')|\leq C'_3\frac{e^{\frac{h\cdot 10D_0L}{2}}}{L^{3n}}
        \end{align*}
        where $h$ is the volume entropy of $G$.
    \end{enumerate}
\end{lemma}

\begin{proof}
    A maximal non-orientable pre-Wicks form has length $6n-2$. Define
     \begin{align*}
         &\mathcal{L}(W,n,L)\\
         &\coloneq\left\{\vec{l}\in\mathbb{Z}^{3n-1}_+\middle|\; \sum_{i=1}^{3n-1}l_i = \left\lceil\frac{L}{2}\right\rceil; l_i\geq \max\left\{ 6\left \lceil \frac{2(3n+2)\log L}{hD_0}\right \rceil+5D_2+2,\frac{R_1}{D_0}\right\},\;(1\leq i\leq 3n-1)\right\}
     \end{align*}
      The lemma follows from \Cref{lem: hyperbolic element intersection count for reduction in degrees of freedom} in the same way in which we argued in the proof of \Cref{lem: hyperbolic element intersection count}. 
\end{proof}

\begin{theorem}\label{thm: lowerboundhyp nonorientable}
    Let $G$ be a hyperbolic group equipped with the word metric induced by some finite symmetric generating set $S$. Let $h$ denote the volume entropy of $G$. There exist positive constants $C'_4,L_0$ depending only on $n$, $G$, and the choice of the generating set, such that:
    \begin{equation*}
        \forall L\geq L_0, |B^{-}_n(L)|\geq C'_4 L^{3n-2}e^{hL/2}
    \end{equation*}
\end{theorem}

\subsection{Relating combinatorial reduction to degeneracy in substitutions}\label{subsec: relating combinatorial reduction to degeneracy in substitutions - hyperbolic}

\begin{construction}[Geodesic encoding]\label{construction: geodesicencoding}
    Let $R_0$ and $R_1$ be positive integers such that: $R_0\leq R_1$, $R_0\geq 5(K+30\delta)$ and $R_1\geq 12(K+30\delta)$. Let $W=X_1\dots X_{2m}$ be a pre-Wicks form and let $\vec{l}$ be a length type for $W$ where $l_i\geq R_1$ for all $1\leq i\leq m$. Let $\mathfrak{g}$ be a substitution for $W$ of length type $\vec{l}$, and let $\mathfrak{g}'$ be the associated good substitution. Set $g=\mathfrak{g}'(W)$, and let $\gamma_g$ be a geodesic segment joining identity to $g$. We will now proceed to construct a partition of $\gamma_g$, which encodes the data of the substitution. First, we set up some notation. 
    \\\\Recall that for each $X_j$, $X_{\overline{j}}$ is the letter in $W$ it is paired with and $X_{\overline{j}}=X_j^{\text{sign}(X_j)}$. In order to define $\mathfrak{g}'$, we chose appropriate Schottky elements $\{l_i\}_{1\leq i\leq m}$, $\{r_i\}_{1\leq i\leq m}$. For all $1\leq j\leq 2m$, $g'_{X_j}=l_{X_j}g_{X_j}r_{X_j}$, where: 
    \begin{align*}
        l_{X_{\overline{j}}}\coloneq\begin{cases}
            (r_{X_j})^{-1}& \text{if sign}(X_j)=-1\\
            l_{X_j}& \text{if sign}(X_j)=1
        \end{cases}
        \hspace{0.3cm} r_{X_{\overline{j}}}\coloneq\begin{cases}
            (l_{X_j})^{-1}& \text{if sign}(X_j)=-1\\
            r_{X_j}& \text{if sign}(X_j)=1
        \end{cases}
    \end{align*}
    Let $\alpha_{X_j}$ denote the piece-wise geodesic path $\gamma_{l_{X_j}}\cup\gamma_{g_{X_j}}\cup\gamma_{r_{X_j}}$, where $\gamma_{l_{X_j}}$ is a geodesic segment joining $\mathfrak{g}'(X_1\cdots X_{j-1})$ to $\mathfrak{g}'(X_1\cdots X_{j-1})l_{X_j}$, $\gamma_{g_{X_j}}$ is a geodesic segment joining $\mathfrak{g}'(X_1\cdots X_{j-1})l_{X_j}$ to $\mathfrak{g}'(X_1\cdots X_{j-1})l_{X_j}g_{X_j}$, and $\gamma_{r_{X_j}}$ is a geodesic segment joining $\mathfrak{g}'(X_1\cdots X_{j-1})l_{X_j}g_{X_j}$ to $\mathfrak{g}'(X_1\cdots X_{j-1})l_{X_j}g_{X_j}r_{X_j}$. Moreover, the geodesic segments $\{\gamma_{g_{X_j}}\}_{1\leq j\leq 2m}$ are chosen so that:
    \begin{enumerate}[font=\normalfont]
        \item[(\textbf{P1})]\linkdest{encoding: matching segments 1} If $s_1\cdots s_k$ is the word representing $\gamma_{g_{X_j}}$, then the word representing $\gamma_{g_{X_{\overline{j}}}}$ is $s_1\cdots s_k$ if $\text{sign}(X_j)=+1$, and it is $s_k^{-1}\cdots s_1^{-1}$ if $\text{sign}(X_j)=-1$.
    \end{enumerate}
    For all $1\leq j\leq 2m$, $\alpha_{X_j}$ is a $(K,R_0)$-chain  (\Cref{claim: numberofgoodsubs}), and $\gamma(\mathfrak{g})$ is $(K+2\delta,R_1)$-chain (\Cref{lem: goodsubstitutionsexist}). Using \Cref{prop: chainfellowtravelling}, there exist decompositions for $\{\gamma_{g_{X_j}}\}_{1\leq j\leq 2m}$ and $\gamma_g$ (see \Cref{fig: double tracking notations}):
    \begin{equation}\label{eq: geodesic encoding}
        \gamma_{g_{X_j}}=\beta_j^{(0)}\gamma_{X_j}\beta_j^{(1)} \;(1\leq j\leq 2m),\quad \gamma_g=\beta'_0\pi_{X_1}\beta'_1\dots\pi_{X_{2m}}\beta'_{2m}
    \end{equation}
    and there exist positive constants $D_1=D_1(K,R_0,\delta)$, $D_2=D_2(K,R_0,\delta)$, $D_3=D_3(K,\delta)$, and $\tilde{D}_3=\tilde{D}_3(K,R_0,\delta)$ such that:
    \begin{enumerate}[font=\normalfont, resume]
        \item[(\textbf{P2})]\linkdest{encoding: matching segments 2} For all $1\leq j\leq 2m$, the subsegments $\gamma_{X_j}$ and $\pi_{X_j}$ begin and end at vertices, and $l(\gamma_{X_j})=l(\pi_{X_j})$.
        \item[(\textbf{P3})]\linkdest{encoding: pairing} For all $1\leq j\leq 2m$, $l(\pi_{X_j})=l(\pi_{X_{\overline{j}}})$.
        \item[(\textbf{P4})]\linkdest{encoding: junctions} For all $1\leq j\leq  2m-1$, $D_1\leq l(\beta'_j)\leq D_2$. $\frac{D_1}{2}\leq l(\beta'_0)\leq \frac{D_2}{2}$ and $\frac{D_1}{2}\leq l(\beta'_{2m})\leq \frac{D_2}{2}$.
        \item[(\textbf{P5})]\linkdest{encoding: mismatch bound} For all $1\leq j\leq 2m$, $l(\beta^{(0)}_j)\leq D_3$ and $l(\beta^{(1)}_j)\leq D_3$.
        \item[(\textbf{P6})]\linkdest{encoding: fellow travel bound}  For all $1\leq j\leq 2m$, the subsegments $\gamma_{X_j}$ and $\pi_{X_j}$ $8(K+40\delta)$-fellow travel.
    \end{enumerate}
    Furthermore, for all $1\leq j\leq 2m$, let $\gamma_{g'_{X_j}}$ denote a geodesic segment joining $\mathfrak{g}'(X_1\cdots X_{{j-1}})$ to $\mathfrak{g}'(X_1\cdots X_{{j}})$. Then, there exist decompositions for $\{\gamma_{g'_{X_j}}\}_{1\leq j\leq 2m}$:
    \begin{equation}\label{eq: geodesic encoding another}
        \gamma_{g'_{X_j}}=\tilde{\beta}^{(0)}_j\tilde{\gamma}_{X_j}\tilde{\beta}^{(1)}_j \;(1\leq j\leq 2m)
    \end{equation}
    such that:
    \begin{enumerate}[font=\normalfont, resume]
        \item[(\textbf{P7})]\linkdest{encoding: extra} For all $1\leq j\leq 2m$, $l(\tilde{\beta}^{(0)}_j)\leq \tilde{D}_3$,  $l(\tilde{\beta}^{(1)}_j)\leq \tilde{D}_3$, $l(\tilde{\gamma}_{X_j})=l(\gamma_{X_j})$, $\gamma_{X_j}$ and $\tilde{\gamma}_{X_j}$ $4(K+30\delta)$-fellow travel, and $\tilde{\gamma}_{X_j}$ and $\pi_{X_j}$ $4(K+40\delta)$-fellow travel.
    \end{enumerate}
    \begin{figure}[h]
        \centering
        \includegraphics[width=0.7\columnwidth]{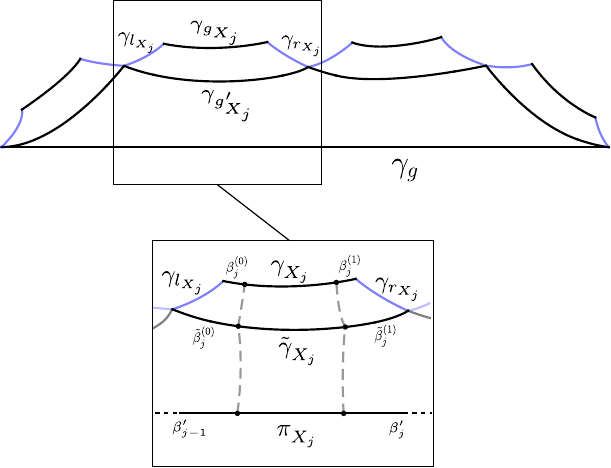}
        \caption{The notations introduced in \Cref{construction: geodesicencoding}}
        \label{fig: double tracking notations}
    \end{figure}
    The data in \Cref{eq: geodesic encoding} and \Cref{eq: geodesic encoding another} will be called the \emph{geodesic encoding} of $\mathfrak{g}$ into $\gamma_g$, and it will be denoted by $\mathcal{E}=\mathcal{E}(\mathfrak{g},\gamma_g)$.

    Throughout this section we will only be using the decompositions in \Cref{eq: geodesic encoding} and properties \linkref{encoding: matching segments 1}{(P1)}-\linkref{encoding: fellow travel bound}{(P6)}. The decomposition in \Cref{eq: geodesic encoding another} and property \linkref{encoding: extra}{(P7)} will only be used in \Cref{subsec: Conjugacy class count: Lower bounds: Hyperbolic groups}.
    
\end{construction}

\noindent\textit{Canonical partition induced by geodesic encoding:}\linkdest{canonical partition} We continue with the set-up and notation used in \Cref{construction: geodesicencoding}. Property \linkref{encoding: pairing}{(P3)} implies that $\gvec{l}_{\pi}\coloneq(l(\pi_{X_j}))_{1\leq j\leq 2m}$, is a length type for $W$. Let $\gvec{b}_{\pi}\coloneq(l(\beta'_j))_{0\leq j\leq 2m}$. There is a tautological correspondence between the canonical partition induced by the data $(W,\gvec{l}_{\pi},\gvec{b}_{\pi})$ on the set of positions $\{1,\dots, L\}$ (see \Cref{def: matching position}) and the decomposition $\gamma_g=\beta'_0\pi_{X_1}\beta'_1\dots\pi_{X_{2m}}\beta'_{2m}$: for every position $i$, we denote by $x_i\in S$ the $i^{\text{th}}$ letter of the word representing $\gamma_g$.\\We will use the correspondence between the positions and the letters of the word representing $\gamma_g$, to relate the combinatorial reduction in \Cref{subsec: combinatorialreduction} with the degeneracy in substitutions seen in \Cref{prop: hyperbolic degenerate substitutions}.

Consider a sequence of $k$ consecutive positions $\mathbf{i}=(a,\dots, a+(k-1))$. We need to set up some notation:
    \begin{itemize}
        \item $\pi(\mathbf{i})$ will denote the subsegment of $\gamma_g$ starting at the $a^{\text{th}}$ vertex and ending at the $(a+k)^{\text{th}}$ vertex of $\gamma_g$. $x(\mathbf{i})$ will denote the group element corresponding to $\pi(\mathbf{i})$, that is, $x(\mathbf{i})\coloneq x_{a}\cdots x_{a+(k-1)}$.
        \item Suppose $\mathbf{i}$ is contained in the letter $X_j$ with respect to $\gvec{l}_{\pi}$. This means that $\pi(\mathbf{i})\subseteq\pi_{X_j}$. Let $\gamma_{X_j}(\mathbf{i})$ denote the subsegment of $\gamma_{X_j}$ fellow-traveling with $\pi(\mathbf{i})$ (see properties \linkref{encoding: matching segments 2}{(P2)} and \linkref{encoding: fellow travel bound}{(P6)}). $g(\mathbf{i})$ will denote the group element corresponding to $\gamma_{X_j}(\mathbf{i})$.
    \end{itemize}

    Suppose $\mathbf{i}=(a,\dots, a+(k-1))$ is contained in the letter $X_j$ with respect to $\gvec{l}_{\pi}$ and let $\varepsilon_{X_j}$ be the sign of $X_j$. Then the sequence of matching positions $\underline{\mathbf{i}}=(\underline{a},\dots,\underline{a+(k-1)})$.
    It follows from properties \linkref{encoding: matching segments 1}{(P1)}, \linkref{encoding: matching segments 2}{(P2)}, \linkref{encoding: pairing}{(P3)}, and \linkref{encoding: mismatch bound}{(P5)} that there are group elements $\eta_1,\eta_2\in B(D_3)$ such that:
    \begin{equation}\label{eq: matching position equation hyperbolic}
    g(\underline{\mathbf{i}})=\eta_1g(\mathbf{i})^{\varepsilon_{X_j}}\eta_2
    \end{equation}

\begin{proposition}
    Suppose we have $(W,\mathfrak{g}_1,\vec{l})$ and $(W',\mathfrak{g}_2,\vec{l}')$ such that $\mathfrak{g}_1'(W)=\mathfrak{g}'_2(W)=g$. Fix a geodesic $\gamma_g$ joining the identity element to $g$ and let $\mathcal{E}_1$ and $\mathcal{E}_2$ denote the geodesic encoding of $\mathfrak{g}_1$ and $\mathfrak{g}_2$ into $\gamma_g$ respectively. Suppose there is a sequence $\mathbf{i}$ of consecutive positions such that:
    \begin{itemize}
        \item $\mathbf{i}$ is contained in the letter $X$ (of $W$) with respect to $\gvec{l}_{\pi}$ and it is contained in the letter $X'$ (of $W'$) with respect to $\gvec{l}^{'}_{\pi}$.
        \item Both $\underline{\mathbf{i}}$ and $\underline{\mathbf{i}}'$ are admissible.
    \end{itemize}
   Let $\varepsilon_X$ and $\varepsilon_{X'}$ denote the signs of $X$ and $X'$ respectively. Then:
   \begin{enumerate}[font=\normalfont]
       \item If $\underline{\mathbf{i}}$ and $\underline{\mathbf{i}}'$ are disjoint, then the substitutions $\mathfrak{g}_1$ and $\mathfrak{g}_2$ satisfy the following non-trivial equations:
   \begin{equation}\label{eq: reduction equations - symmetric hyperbolic}
       g_1(\underline{\mathbf{i}})^{\varepsilon_X}=\xi'_1g_1(\underline{\mathbf{i}}')^{\varepsilon_{X'}}\xi'_2,\quad g_2(\underline{\mathbf{i}})^{\varepsilon_X}=\xi'_3g_2(\underline{\mathbf{i}}')^{\varepsilon_{X'}}\xi'_4
   \end{equation}
   where $\xi'_1,\xi'_2,\xi'_3,\xi'_4\in B(16(K+40\delta)+2D_3)$. 
    \item If $\mathbf{i}$ and $\underline{\mathbf{i}}'$ are disjoint, then the substitution $\mathfrak{g}_1$ satisfies the following non-trivial equation:
   \begin{equation}\label{eq: reduction equation 1 - asymmetric hyperbolic}
       g_1(\mathbf{i})=\xi''_1 g_1(\underline{\mathbf{i}}')^{\varepsilon_{X'}}\xi''_2
   \end{equation}
   where $\xi''_1,\xi''_2\in B(16(K+40\delta)+D_3)$. 
   \item If $\mathbf{i}$ and $\underline{\mathbf{i}}$ are disjoint, then the substitution $\mathfrak{g}_2$ satisfies the following non-trivial equation:
   \begin{equation}\label{eq: reduction equation 2 - asymmetric hyperbolic}
       g_2(\mathbf{i})=\xi''_3 g_2(\underline{\mathbf{i}})^{\varepsilon_{X}}\xi''_4
   \end{equation}
   where $\xi''_3,\xi''_4\in B(16(K+40\delta)+D_3)$.
   \end{enumerate}
\end{proposition}

\begin{proof}
    Firstly, using \Cref{eq: matching position equation hyperbolic}, we get the following:
   \begin{equation}\label{eq: matching position equations - hyperbolic}
       g_1(\underline{\mathbf{i}})=\eta_1g_1(\mathbf{i})^{\varepsilon_{X}}\eta_2,\quad g_2(\underline{\mathbf{i}}')=\eta_3g_2(\mathbf{i})^{\varepsilon_{X'}}\eta_4
   \end{equation}
   where $\eta_1,\eta_2,\eta_3,\eta_4\in B(D_3)$. Secondly, using the admissibility conditions, we get:
   \begin{equation}\label{eq: admissibility equations - hyperbolic}
       g_1(\mathbf{i})=\xi_1g_2(\mathbf{i})\xi_2,\quad g_1(\underline{\mathbf{i}})=\xi_3g_2(\underline{\mathbf{i}})\xi_4,\quad g_1(\underline{\mathbf{i}}')=\xi_5g_2(\underline{\mathbf{i}}')\xi_6 
   \end{equation}
   where $\xi_1,\xi_2,\xi_3,\xi_4,\xi_5,\xi_6\in B(8(K+40\delta))$.  \Cref{eq: reduction equations - symmetric hyperbolic}, \Cref{eq: reduction equation 1 - asymmetric hyperbolic}, and \Cref{eq: reduction equation 2 - asymmetric hyperbolic} follow directly from \Cref{eq: matching position equations - hyperbolic} and \Cref{eq: admissibility equations - hyperbolic}.
\end{proof}

\subsubsection{Proof of \Cref{lem: hyperbolic element intersection count for reduction in degrees of freedom}}

\begin{proposition}\label{prop: hyperbolic degenerate substitutions}
    Let $N$ and $D_0$ be positive integers. Let $W$ and $W'$ be maximal pre-Wicks forms of length $2m$ and $2m'$ respectively. Let $\vec{l}$ and $\vec{l}'$ be two length types for $W$ and $W'$ respectively, such that:
    \begin{enumerate}[font=\normalfont]
        \item $||\vec{l}||=||\vec{l}'||(\coloneq L)$
        \item For all $1\leq i\leq m$, $1\leq j\leq m'$, $l_i,l'_j>\max\left\{6N+5D_2+2,\frac{R_1}{D_0}\right\}$
        \item $(W,\vec{l})\neq (W',\vec{l}')$
    \end{enumerate}
    For all $D_0> D_2$ and for every $g\in G_W(10D_0\cdot\vec{l})\cap G_{W'}(10D_0\cdot \vec{l}')$, there is a substitution $\mathfrak{g}$ for $W$, of length type $10D_0\cdot\vec{l}$ with $\mathfrak{g}'(W)=g$, having one of the following \emph{degenerate} algebraic forms:
    \begin{enumerate}[font=\normalfont]
        \item \linkdest{algebraic form 1} There are two distinct unknowns $X$ and $X'$ of $W$ such that:
        \begin{equation*}
            g_X=uzv,\quad g_{X'}=u'\xi_1z^{\varepsilon}\xi_2v'
        \end{equation*}
        where $\varepsilon\in\{+1,-1\}$, $\xi_1,\xi_2\in B(16(K+40\delta)+2D_3)$, and:
        \begin{equation}\label{eq: hyp degen len eqn 1}
            l(z)=ND_0,\quad l(u)+l(v)+ND_0=10l_XD_0,\quad l(u')+l(v')+ND_0=10l_{X'}D_0 
        \end{equation}
        \item \linkdest{algebraic form 2} There is an unknown $X$ of $W$ such that:
        \begin{equation*}
            g_X=uzw\xi_1z^{\varepsilon}\xi_2v
        \end{equation*}
        where $\varepsilon\in\{+1,-1\}$, $\xi_1,\xi_2\in B(16(K+40\delta)+2D_3)$, and:
        \begin{equation}\label{eq: hyp degen len eqn 2}
            l(z)=ND_0,\quad l(u)+l(w)+l(v)+2ND_0=10l_XD_0
        \end{equation}
        \item \linkdest{algebraic form 3} There is an unknown $X$ of $W$, and a positive integer $T\geq D_0$ such that:
        \begin{equation*}
            g_X=uz_0z_1\cdots z_kv
        \end{equation*}
        where $k=\left\lfloor\frac{ND_0}{T}\right\rfloor+1$, $l(z_0)=l(z_1)=\cdots=l(z_{k-1})=T$, $l(z_k)<T$. Moreover:
        \begin{equation*}
             z_i=\xi_iz_{i-1}\xi'_i \;\; (1\leq i\leq k-1),\quad z_k=\xi_kz'_{k-1}\xi'_k
        \end{equation*}
        where $\xi_1,\xi'_1,\dots,\xi_k,\xi'_k\in B(16(K+40
        \delta)+2D_3)$, and $z'_{k-1}$ is a group element of length $ND_0-(k-1)T$ which appears as a prefix in some geodesic word representing $z_{k-1}$. Furthermore:
        \begin{equation}\label{eq: hyp degen len eqn 3}
            l(u)+l(v)+l(z_k)+kT=10l_XD_0,\quad (k-1)T+l(z_k)=ND_0
        \end{equation}
        
        \item \linkdest{algebraic form 4} There is an unknown $X$ of $W$ such that:
        \begin{equation*}
            g_X=uzw\xi_1z^{-1}\xi_2v
        \end{equation*}
        where $\xi_1,\xi_2\in B(16(K+40\delta)+2D_3)$, $w\in B(3D_0)$, and:
        \begin{equation}\label{eq: hyp degen len eqn 4}
            l(u)+l(w)+l(v)+2l(z)=10l_XD_0,\quad l(z)+l(w)=ND_0
        \end{equation}
        \item \linkdest{algebraic form 5} There is an unknown $X$ of $W$ such that:
        \begin{equation*}
            g_X=uz_0z_1w\xi_1z_1^{-1}\xi'_1\xi_2z_0^{-1}\xi'_2v
        \end{equation*}
        where $\xi_1,\xi_2,\xi'_1,\xi'_2\in B(16(K+40\delta)+2D_3)$, $w\in B(D_0)$, and:
        \begin{equation}\label{eq: hyp degen len eqn 5}
            l(z_1)\geq D_0,\quad l(u)+l(w)+l(v)+2l(z_0)+2l(z_1)=10l_XD_0,\quad l(z_0)+2l(z_1)+l(w)=ND_0
        \end{equation}
    \end{enumerate}
\end{proposition}

\begin{proof}\label{proof: proof of degenerate hyperbolic substitutions}
Let $D_0$ be any positive integer greater than $D_2$.
Suppose there are substitutions $\mathfrak{g}_1$ and $\mathfrak{g}_2$ for $W$ and $W'$ of length types $\gvec{l}_{10D_0}\coloneq 10D_0\vec{l}$ and $\gvec{l'}_{10D_0}\coloneq 10D_0\vec{l}'$ respectively, such that $\mathfrak{g}_1'(W)=\mathfrak{g}_2'(W')=g$. Let $\gamma_g$ be a geodesic joining the identity element to $g$. Let $\mathcal{E}_1$ and $\mathcal{E}_2$ denote geodesic encodings of $\mathfrak{g}_1$ and $\mathfrak{g}_2$ into $\gamma_g$ respectively. The encodings, in particular, the data $(W,\overrightarrow{(l_{10D_0})}_{\pi},\vec{b}_{\pi})$ and $(W',\overrightarrow{(l'_{10D_0})}_{\pi},\vec{b}'_{\pi})$ induce two canonical partitions on the set of positions $\{1,\dots, l(g)\}$.

\begin{claim}
    There is a $D_0$-reduction in degrees of freedom by $N$ between $(W,\overrightarrow{(l_{10D_0})}_{\pi},\vec{b}_{\pi})$ and $(W',\overrightarrow{(l'_{10D_0})}_{\pi},\vec{b}'_{\pi})$. That is, there exists an admissible sequence of $ND_0$ consecutive positions $\mathbf{i}=(i_1,\dots,i_{ND_0})$ such that both $\underline{\mathbf{i}}$ and $\underline{\mathbf{i}}'$ are admissible, and $|\underline{i_1}-\underline{i_1}'|\geq D_0$.
\end{claim}
\begin{proof}
    The data $(W,\overrightarrow{(l_{10D_0})}_{\pi},\vec{b}_{\pi})$, $(W',\overrightarrow{(l'_{10D_0})}_{\pi},\vec{b}'_{\pi})$ satisfies conditions \ref{item: L1} and \ref{item: L2} of \Cref{thm: the main reduction theorem}. Note that each entry of $\vec{b}_{\pi}$ and $\vec{b'}_{\pi}$ lies in the interval $[0,D_2]$ (see \linkref{encoding: pairing}{(P4)}). Therefore, it follows from \Cref{thm: the main reduction theorem} that either there is $D_0$-reduction in degrees of freedom by $N$ or $W=W'$, $m=m'$, and for all $1\leq i\leq 2m$, $|10D_0l_i-10D_0l'_i|\leq 3D_0+4D_2$. Since $D_0>D_2$, this implies that $\vec{l}=\vec{l}'$. But this contradicts our hypothesis that $(W,\vec{l})\neq(W',\vec{l}')$. Hence, there is $D_0$-reduction in degrees of freedom by $N$.
\end{proof}

The rest of the proof is purely algebraic. Suppose $\mathbf{i}$ is contained in the letter $X_{j_1}$ with respect to $\gvec{l}_{D_0}$ and is contained in the letter $X'_{j_2}$ with respect to $\gvec{l'}_{D_0}$. Set $\varepsilon_{j_1}=\text{sign}(X_{j_1})$ and $\varepsilon_{j_2}=\text{sign}(X_{j_2})$. There are three possibilities for $\underline{\mathbf{i}}$ and $\underline{\mathbf{i}}'$:
\begin{enumerate}[label=(\Roman*)]
    \item $\underline{\mathbf{i}}$ and $\underline{\mathbf{i}}'$ are disjoint and they are either contained in the same letter or in different, unpaired letters with respect to $\overrightarrow{(l_{D_0})}_{\pi}$.
    \item $\underline{\mathbf{i}}$ and $\underline{\mathbf{i}}'$ are not disjoint and they are contained in the same letter with respect to $\overrightarrow{(l_{D_0})}_{\pi}$.
    
    \item $\underline{\mathbf{i}}$ and $\underline{\mathbf{i}}'$ are not disjoint and they are contained in different but paired letters with respect to $\overrightarrow{(l_{D_0})}_{\pi}$.
\end{enumerate}
\noindent\textbf{Case - \RomanCaps{1}}: See \Cref{fig: degenerate subs 1}. In this case, by \Cref{eq: reduction equations - symmetric hyperbolic} for $\mathfrak{g}$, there exist element $\xi_1,\xi_2\in B(16(K+40\delta)+2D_3)$ such that:
\begin{align*}
g_1(\underline{\mathbf{i}})=\xi_1g_1(\underline{\mathbf{i}}')^{(\varepsilon_{j_1}\varepsilon_{j_2})}\xi_2
\end{align*}
If $\underline{\mathbf{i}}$ and $\underline{\mathbf{i}}'$ belong to the same letter, we get the algebraic form given in \linkref{algebraic form 2}{$(1)$}, otherwise we get the algebraic form given in \linkref{algebraic form 1}{$(2)$}.
    \begin{figure}[h]
        \centering
        \includegraphics[width=0.8\columnwidth]{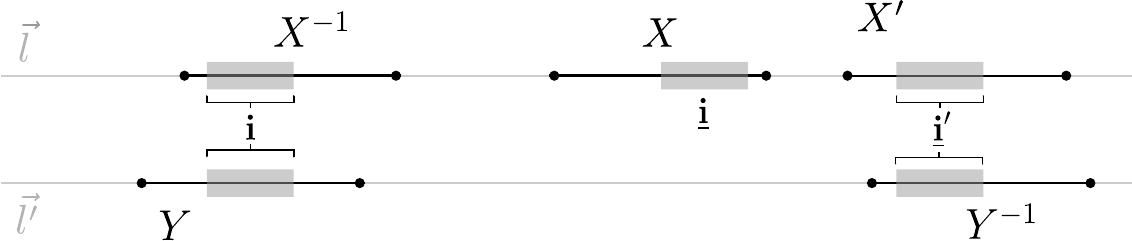}
        \caption{When $\underline{\mathbf{i}}$ and $\underline{\mathbf{i}}'$ are disjoint.}
        \label{fig: degenerate subs 1}
    \end{figure}
\\\textbf{Case - \RomanCaps{2}}: See \Cref{fig: degenerate subs 2}. Let $\underline{\mathbf{i}}\cup\underline{\mathbf{i}}'$ denote the sequence of consecutive positions comprising $\underline{\mathbf{i}}$ and $\underline{\mathbf{i}}'$.\\Suppose $\varepsilon_{j_1}=\varepsilon_{j_2}=\varepsilon$. There are four subcases here, depending on whether $\varepsilon=+1$ or $-1$ and whether the letter is an unknown, say $X$ or its inverse $X^{-1}$ (in case $\text{sign}(X)=-1$). We will only argue for the case when $\varepsilon=+1$, and $\underline{\mathbf{i}}$ and $\underline{\mathbf{i}}'$ are contained in the unknown $X$. The other cases can be dealt with entirely similar arguments. Set $T=|\underline{i_1}-\underline{i_1}'|$ and write $\mathbf{i}=(\mathbf{i_1},\dots,\mathbf{i_k})$, where:
\begin{align*}
    k=\left\lfloor\frac{ND_0}{T}\right\rfloor+1,\quad|\mathbf{i_1}|=\cdots=|\mathbf{i_{k-1}}|=T,\; |\mathbf{i_{k}}|<T
\end{align*}
Note that $T\geq D_0$. Define $\mathbf{a_j}=\underline{\mathbf{i_{j}}}$ and $\mathbf{b_j}=\underline{\mathbf{i_{j}}}'$ for all $1\leq j\leq k$. Suppose $\underline{i_1}<\underline{i_1}'$. Because $X_{j_1}$ and $X_{j_2}$ have the same sign, it follows from \Cref{eq: matching position 1} that: $\mathbf{a_{j+1}}=\mathbf{b_j}$ for all $1\leq j\leq k-2$. Therefore:
\begin{equation}\label{eq: degenerate substitutions proof case II eq 1}
    \text{for all }1\leq j\leq k-2,\quad g(\mathbf{b_j})=g(\mathbf{a_{j+1}})
\end{equation}
Also, by definition, $\mathbf{a_j}=\underline{\mathbf{i_j}}$ and $\mathbf{b_j}=\underline{\mathbf{i_j}}'$ are disjoint for all $1\leq j\leq k$. Thus, using \Cref{eq: reduction equations - symmetric hyperbolic}, for each $1\leq j\leq k$, there exist $\xi_j,\xi'_j\in B(16(K+40\delta)+2D_3)$  such that:
\begin{equation}\label{eq: degenerate substitutions proof case II eq 2}
    g(\mathbf{b_j}) = \xi_jg(\mathbf{a_j})\xi_j'
\end{equation}
Also note that $g(\mathbf{b_{k-1}})$ determines $g(\mathbf{a_k})$. Looking at the geodesic segment $\gamma_{g_X}$ used in the encoding $\mathcal{E}_1$, we see that $g(\mathbf{a_k})$ can be read off from a prefix of the (sub)word representing $g(\mathbf{b_{k-1}})$ in $\gamma_{g_X}$. Setting $z_0=g(\mathbf{a_1})$, $z'_{k-1}=g(\mathbf{a_k})$, and $z_j=g(\mathbf{b_j})$ for all $1\leq j\leq k$, and using \Cref{eq: degenerate substitutions proof case II eq 1} and \Cref{eq: degenerate substitutions proof case II eq 2}, we get:
\begin{align*}
    g(\underline{\mathbf{i}}\cup\underline{\mathbf{i}}')=z_0z_1\cdots z_k,\quad z_j=\xi_jz_{j-1}\xi'_j \;\; (1\leq j\leq k-1),\quad z_k=\xi_kz'_{k-1}\xi'_k
\end{align*}
Evidently, $\mathfrak{g}$ has the algebraic form given in \linkref{algebraic form 3}{$(3)$}.
    \begin{figure}[h]
        \centering
        \includegraphics[width=0.8\columnwidth]{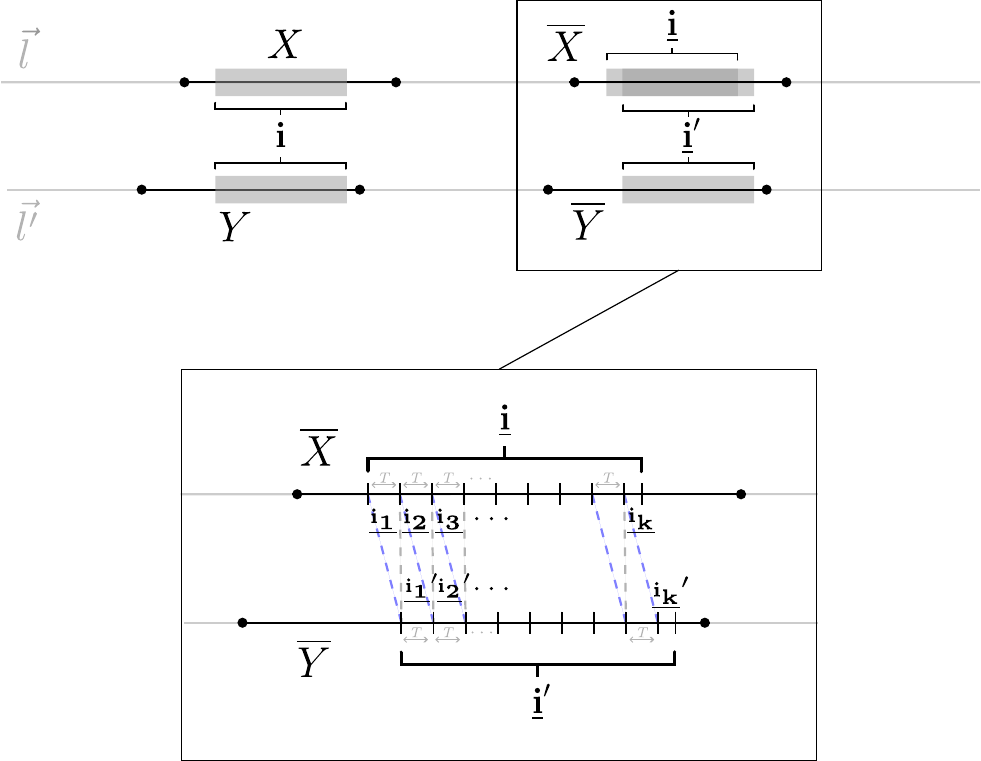}
        \caption{When $\underline{\mathbf{i}}$ and $\underline{\mathbf{i}}'$ overlap and are contained in the same letter, assuming $\text{sign}(X) = \text{sign}(Y) = 1$. The gray vertical lines indicate equality between positions, and the blue slanted ones indicate near-equality between the respective $\mathfrak{g}$ and $\mathfrak{g}'$ substitutions.}
        \label{fig: degenerate subs 2}
    \end{figure}
\\Next, suppose $\varepsilon_{j_1}\neq\varepsilon_{j_2}$. We will argue for the case when $\underline{\mathbf{i}}$ and $\underline{\mathbf{i}}'$ are contained in the unknown $X$, $\varepsilon_{j_1}=-1$, $\varepsilon_{j_2}=+1$, and $\underline{i_{ND_0}}<\underline{i_1}'<\underline{i_1}<\underline{i_{ND_0}}
'$. The other cases can be argued similarly.  In this case, either we have ({\romannumeral 1}) $\underline{i_1}-\underline{i_1}'< 3D_0$, or ({\romannumeral 2}) there exists a positive integer $T\geq D_0$ such that $\underline{i_1}-\underline{i_1}'=2T+D_0$.
\begin{enumerate}[label=(\roman*)]
    \item $\underline{i_1}-\underline{i_1}'\leq 3D_0$. We write $\mathbf{i}=(\mathbf{i_1},\mathbf{i_2})$, where $|\mathbf{i_2}|=\underline{i_1}-\underline{i_1}'+1$. Set $\mathbf{a_j}=\underline{\mathbf{i_j}}$ and $\mathbf{b_j}=\underline{\mathbf{i_j}}'$ for all $1\leq j\leq 2$. Then, $\mathbf{a_1}=\mathbf{b_1}$. Also, $\mathbf{a_2}=\underline{\mathbf{i_2}}$ and $\mathbf{b_2}=\underline{\mathbf{i_2}}'$ are disjoint. So by \Cref{eq: reduction equations - symmetric hyperbolic}, there exist elements $\xi_1,\xi_2\in B(16(K+40\delta)+2D_3)$ such that:
    \begin{align*}
        g(\mathbf{b_2})=\xi_1g(\mathbf{a_2})^{-1}\xi_2
    \end{align*}
    Thus, setting $z=g(\mathbf{a_2})$ and $w=g(\mathbf{a_1})$, we have:
    \begin{equation*}
        g(\underline{\mathbf{i}}\cup\underline{\mathbf{i}}')=zw\xi_1z^{-1}\xi_2
    \end{equation*}
    Therefore, $\mathfrak{g}$ has the algebraic form given in \linkref{algebraic form 4}{$(4)$}.
    \item There exists a positive integer $T\geq D_0$ such that $\underline{i_1}-\underline{i_1}'=2T+D_0$. we write $\mathbf{i}=(\mathbf{i_1},\mathbf{i_2},\mathbf{i_3},\mathbf{i_4})$ where $|\mathbf{i_1}|=T$, $|\mathbf{i_2}|=D_0$, and $|\mathbf{i_3}|=T$. Set $\mathbf{a_j}=\underline{\mathbf{i_j}}$ and $\mathbf{b_j}=\underline{\mathbf{i_j}}'$ for all $1\leq j\leq 4$. Then, $\mathbf{b_1}=\mathbf{a_3}$, $\mathbf{b_2}=\mathbf{a_2}$, and $\mathbf{b_3}=\mathbf{a_1}$. Now note that $\mathbf{a_1}=\underline{\mathbf{i_1}}$ and $\mathbf{b_1}=\underline{\mathbf{i_1}}'$ are disjoint, and $\mathbf{a_4}=\underline{\mathbf{i_4}}$ and $\mathbf{b_4}=\underline{\mathbf{i_4}}'$ are disjoint. So, it follows from \Cref{eq: reduction equations - symmetric hyperbolic} that there exist elements $\xi_1,\xi_2,\xi'_1,\xi'_2\in B(16(K+40\delta)+2D_3)$ such that:
    \begin{equation*}
        g(\mathbf{a_1})=\xi_1g(\mathbf{b_1})^{-1}\xi'_1,\quad g(\mathbf{b_4})=\xi_2g(\mathbf{a_4})^{-1}\xi'_2
    \end{equation*}
    Therefore, setting $z_0=g(\mathbf{a_4})$, $z_1=g(\mathbf{a_3})$, and $w=g(\mathbf{a_2})$, we have:
    \begin{equation*}
        g(\underline{\mathbf{i}}\cup\underline{\mathbf{i}}')=z_0z_1w\xi_1z_1^{-1}\xi'_1\xi_2z_0^{-1}\xi'_2
    \end{equation*}
    Therefore, $\mathfrak{g}$ has the algebraic form given in \linkref{algebraic form 5}{$(5)$}.
\end{enumerate}
\textbf{Case - \RomanCaps{3}}: $\underline{\mathbf{i}}$ and $\underline{\mathbf{i}}'$ are contained in different but paired letters with respect to $\overrightarrow{(l_{D_0})}_{\pi}$. Therefore $\mathbf{i}$ and $\underline{\mathbf{i}}'$ are contained in the same letter with respect to $\overrightarrow{(l_{D_0})}_{\pi}$. There are two subcases here, depending on the sign of $X'_{j_2}$:
\begin{enumerate}[label=(\roman*)]
    \item Suppose $\varepsilon_{j_2}=-1$. In this case, an analysis identical to subcase ({\romannumeral 1}) in Case - \RomanCaps{2}, using \Cref{eq: reduction equation 1 - asymmetric hyperbolic} instead of \Cref{eq: reduction equations - symmetric hyperbolic} implies that $\mathfrak{g}$ has the algebraic form given in \linkref{algebraic form 4}{$(4)$}.
    \item Suppose $\varepsilon_{j_2}=+1$. We observe that:
    \begin{equation*}
        |i_1 - \underline{i_1}'|=|i-(s^{\mathcal{E}_2}_{X'_{\bar{j_2}}}+i-s^{\mathcal{E}_2}_{X'_{j_2}})|=|s^{\mathcal{E}_2}_{X'_{j_2}}-s^{\mathcal{E}_2}_{X'_{\bar{j_2}}}|>D_0
    \end{equation*}
    where $s^{\mathcal{E}_2}_{X'_{\bar{j_2}}}$ and $s^{\mathcal{E}_2}_{X'_{j_2}}$ denote the starting positions of $X'_{\bar{j_2}}$ and $X'_{j_2}$ with respect to $\overrightarrow{(l'_{D_0})}_{\pi},\gvec{b'}_{\pi}$. Thus, an analysis identical to that made in the same sign subcase in Case - \RomanCaps{2}, using \Cref{eq: reduction equation 1 - asymmetric hyperbolic} instead of \Cref{eq: reduction equations - symmetric hyperbolic} implies that $\mathfrak{g}$ has the algebraic form given in \linkref{algebraic form 3}{$(3)$}.
\end{enumerate}
\end{proof}

\begin{proof}[Proof of \Cref{lem: hyperbolic element intersection count for reduction in degrees of freedom}]
    Throughout this proof, we will be using the notation introduced in \Cref{prop: hyperbolic degenerate substitutions}. Let $D_4\coloneq 16(K+40\delta)+2D_3$. $(W,\vec{l})$ and $(W,\vec{l}')$ satisfy the hypotheses of \Cref{prop: hyperbolic degenerate substitutions}, so for all $D_0> D_2$ and for any $g\in G_W(D_0\vec{l})\cap G_{W'}(D_0\vec{l}')$, there exists a substitution $\mathfrak{g}$ for $(W,D_0\vec{l})$ such that $\mathfrak{g}'(W)=g$ and it has an algebraic form given by \linkref{algebraic form 1}{$(1)$}-\linkref{algebraic form 5}{$(5)$}. For all $1\leq j\leq 5$, let $\text{Subs}_j(W,10D_0\cdot \vec{l})\subset\prod_{i=1}^{m}B(10D_0l_i)$ denote the set of all substitutions $\mathfrak{g}$ for $(W,10D_0\cdot\vec{l})$ having an algebraic form of type $j$. Then:
    \begin{equation*}
        |G_W(10D_0\cdot\vec{l})\cap G_{W'}(10D_0\cdot\vec{l}')|\leq \left|\bigcup_{j=1}^{5}\text{Subs}_j(W,10D_0\cdot\vec{l})\right |\leq \sum_{j=1}^{5}\left|\text{Subs}_j(W,10D_0\cdot\vec{l})\right|
    \end{equation*}
    In order to estimate $|\text{Subs}_j(W,10D_0\cdot\vec{l})|$, we will be using the upper bound in \Cref{thm: ballgrowth} repeatedly, that is, for all $L>0$ we have $|B(L)|\leq C_2e^{hL}$. We also fix some $D_0>\max\{D_2,4(D_4+\frac{\log C_2}{h})\}$.
    \\\\\noindent \underline{$\mathbf{j=1}$}: Let $\text{Subs}_1(X,X',l(u)=r_1,l(u')=r_2)$ be the set of all substitutions of type \linkref{algebraic form 1}{$(1)$} where the distinct unknowns $X$ and $X'$ are fixed and $l(u)=r_1,l(u')=r_2$. We have:
    \begin{align*}
        &|\text{Subs}_1(X,X',l(u)=r_1,l(u')=r_2)|\\
        &\leq \left(C_2^{m-2} e^{h\sum_{i\neq X,X'}10D_0l_i}\right)\times\left( C_2^3e^{h(l(u)+ND_0+l(v))}\right)\times\left( C_2^2e^{h(l(u')+l(v'))}\right)\times\left(C_2e^{hD_3}\right)^2\\
        &=\left(C_2^{m+5}e^{2hD_4}\right)\times e^{h\cdot 10D_0\sum_{i=1}^ml_i} \times e^{-hND_0}\quad[\text{Using \Cref{eq: hyp degen len eqn 1}}]
    \end{align*}
    Therefore,
    \begin{align*}
        |\text{Subs}_1(W,10D_0\cdot\vec{l})|&=\sum_{X,X',r_1,r_2}|\text{Subs}_1(X,X',l(u)=r_1,l(u')=r_2)|\\
        &\leq {m\choose 2}(10D_0L)^2\times \left(C_2^{m+5}e^{2hD_4}\right)\times e^{\frac{h\cdot 10D_0L}{2}}\times e^{-hND_0}\\
        &\leq \left(100D_0^2m^2C_2^{m+5}e^{2hD_4}\right)\times e^{\frac{h\cdot 10D_0L}{2}}\times L^2e^{-hND_0}
    \end{align*}
    \noindent\underline{$\mathbf{j=2}$}: This case is essentially the same as the one above. Let $\text{Subs}_2(X,l(u)=r_1,l(w)=r_2)$ be the set of all substitutions of type \linkref{algebraic form 2}{$(2)$} where the unknown $X$ is fixed and $l(u)=r_1,l(w)=r_2$. Using \Cref{eq: hyp degen len eqn 2}, we get:
    \begin{align*}
        |\text{Subs}_2(X,l(u)=r_1,l(w)=r_2)|&\leq \left(C_2^{m+3}e^{2hD_4}\right)\times e^{h\cdot 10D_0\sum_{i=1}^ml_i}\times e^{-hND_0}\\
        |\text{Subs}_2(W,10D_0\cdot\vec{l})|&=\sum_{X,r_1,r_2}|\text{Subs}_2(X,l(u)=r_1,l(w)=r_2)|\\
        &\leq \left(100D_0^2mC_2^{m+3}e^{2hD_4}\right)\times e^{\frac{h\cdot 10D_0L}{2}}\times L^2e^{-hND_0}
    \end{align*}
    \noindent\underline{$\mathbf{j=3}$}: Let $\text{Subs}_3(X,l(u)=r_1,l(z_0)=T)$ be the set of all substitutions of type \linkref{algebraic form 3}{$(3)$} where the unknown $X$ is fixed, $l(u)=r_1$, and $l(z_0)=T$, for $T\geq D_0$. By \Cref{eq: hyp degen len eqn 3}, $r_1$ and $T$ determine $l(v)$. Therefore, we get:
    \begin{align*}
        &|\text{Subs}_3(X,l(u)=r_1,l(z_0)=T)|\\
        &\leq \left(C_2^{m-1}e^{h\cdot 10D_0\sum_{i\neq X}l_i}\right)\times\left(C_2^3e^{h(l(u)+l(z_0)+l(v))}\right)\times\left(C_2e^{hD_4}\right)^{2k}\times\left(C_2e^{2h\delta}\right)\\
        &\leq \left(C_2^{m+3}e^{2h\delta}\right)\times e^{h\cdot 10D_0\sum_{i=1}^ml_i}\times e^{-hND_0}\times e^{2kh\left(D_4+\frac{\log C_2}{h}\right)} \quad [\text{Using \Cref{eq: hyp degen len eqn 3}}]\\
        &\leq \left(C_2^{m+5}e^{2h(D_4+\delta)}\right)\times e^{\frac{h\cdot 10D_0L}{2}}\times e^{-hN\left(D_0-2\left(D_4+\frac{\log C_2}{h}\right)\right)} \quad\left[\text{Since }k\leq N+1 \right]\\
        &\leq\left(C_2^{m+5}e^{2h(D_4+\delta)}\right)\times e^{\frac{h\cdot 10D_0L}{2}}\times e^{-\frac{hND_0}{2}} 
    \end{align*}
    where in the last step we have used the fact that $D_0\geq4(D_4+\frac{\log C_2}{h})$. 
    \begin{align*}
        |\text{Subs}_3(W,10D_0\cdot\vec{l})|&=\sum_{X,r_1,T}|\text{Subs}_3(X,l(u)=r_1,l(z_0)=T)|\\
        &\leq \left(100D_0^2mC_2^{m+5}e^{2h(D_4+\delta)}\right)\times e^{\frac{h\cdot 10D_0L}{2}}\times L^2e^{-\frac{hND_0}{2}} 
    \end{align*}
    \noindent\underline{$\mathbf{j=4}$}: Let $\text{Subs}_4(X,l(u)=r_1,l(w)=r_2)$ be the set of all substitutions of type \linkref{algebraic form 4}{$(4)$} where the unknown $X$ is fixed, $l(u)=r_1$ and $l(w)=r_2$. By \Cref{eq: hyp degen len eqn 4}, $r_1$ and $r_2$ determine $l(z)$ and $l(v)$. Also, note that $r_2\leq 3D_0$. Therefore:
    \begin{align*}
        |\text{Subs}_4(X,l(u)=r_1,l(w)=r_2)|
        &\leq\left(C_2^{m-1}e^{h\cdot 10D_0\sum_{i\neq X}l_i}\right)\times \left(C_2^4e^{h(l(u)+l(w)+l(v)+l(z))}\right)\times\left(C_2e^{hD_4}\right)^2\\
        &\leq \left(C_2^{m+5}e^{h(2D_4+3D_0)}\right)\times e^{\frac{h\cdot 10D_0L}{2}}\times e^{-hND_0}\quad[\text{Using \Cref{eq: hyp degen len eqn 4}}]\\
        |\text{Subs}_4(W,10D_0\cdot\vec{l})|&=\sum_{X,r_1,r_2}|\text{Subs}_4(X,l(u)=r_1,l(w)=r_2)|\\
        &\leq \left(30mD_0^2C_2^{m+5}e^{h(2D_4+3D_0)}\right)\times e^{\frac{h\cdot 10D_0L}{2}}\times Le^{-hND_0}
    \end{align*}
    \noindent\underline{$\mathbf{j=5}$}: Let $\text{Subs}_5(X,l(u)=r_1,l(z_1)=r_2,l(w)=r_3)$ denote the set of all substitutions of type \linkref{algebraic form 5}{$(5)$} such that the unknown $X$ is fixed, $l(u)=r_1$, $l(z_1)=r_2$, and $l(w)=r_3$. By \Cref{eq: hyp degen len eqn 5}, $r_1, r_2, r_3$ determine $l(z_0)$ and $l(v)$. Also, note that $r_3\leq D_0$. Therefore,
    \begin{align*}
        &|\text{Subs}_5(X,l(u)=r_1,l(z_1)=r_2,l(w)=r_3)|\\
        &\leq \left(C_2^{m-1}e^{h\cdot 10D_0\sum_{i\neq X}l_i}\right)\times\left(C_2^5e^{h(l(u)+l(w)+l(v)+l(z_0)+l(z_1))}\right)\times\left(C_2e^{hD_4}\right)^4\\
        &\leq \left(C_2^{m+8}e^{h(4D_4+D_0)}\right)\times e^{\frac{h\cdot 10D_0L}{2}}\times e^{-\frac{h(N-1)D_0}{2}}
    \end{align*}
    \begin{align*}
        |\text{Subs}_5(W,10D_0\cdot\vec{l})|&=\sum_{X,r_1,r_2,r_3}|\text{Subs}_5(X,l(u)=r_1,l(z_1)=r_2,l(w)=r_3)|\\
        &\leq \left(100mD_0^3C_2^{m+8}e^{h(4D_4+2D_0)}\right)\times e^{\frac{h\cdot 10D_0L}{2}}\times L^2e^{-\frac{hND_0}{2}}
    \end{align*}
    Inspecting the constants in the five estimates, for all $1\leq j\leq 5$ we get:
    \begin{equation*}
        |\text{Subs}_j(W,10D_0\cdot\vec{l})|\leq \left(100m^2D^3_0C_2^{m+8}e^{h(4D_4+3D_0+\delta)}\right)\times e^{\frac{h\cdot 10D_0L}{2}}\times L^2e^{-\frac{hND_0}{2}}
    \end{equation*}
    Using the estimate in \Cref{lem: lowerboundgoodsubs} we get:
    \begin{equation*}
        \frac{|G_W(10D_0\cdot\vec{l})\cap G_{W'}(10D_0\cdot\vec{l}')|}{|G_W(10D_0\cdot\vec{l})|}\leq \left(500m^2D_0^3(C'_1)^{-1}C_2^{m+8}e^{h(4D_4+3D_0+\delta)}\right)\times L^2e^{-\frac{hND_0}{2}}
    \end{equation*}
\end{proof}

\section{Conjugacy class count: Upper bounds}\label{sec: Conjugacy class count: Upper bounds}

Suppose $\mathcal{X}$ is a set of unknowns and $W$ is a Wicks form with unknowns in $\mathcal{X}$. Throughout this section and the next, a Wicks form $W$ will be considered as a cyclic word. Moreover, whenever $W$ is non-orientable, it is assumed that $\text{genus}^-(W)>1$. We will often need to look at a linear word representing a Wicks form. In this case we will either denote this linear word by a different symbol, or write $W=X_1\cdots X_{2m}$. The notion of a length type needs some clarification. A length type $\vec{l}$ for $W$ is an assignment of non-negative integers to the unknowns appearing in $W$ (see \Cref{rmk: length type as an m-tuple}). For the purpose of specifying length types, we will write $\vec{l}=(l_i)_{1\leq i\leq m}$ and it will be understood that some total order has been chosen on the set of unknowns. Suppose $V=X_1\cdots X_{2m}$ is a linear word representing $W$. Then $\vec{l}_V$ denotes the tuple $(l_{X_i})_{1\leq i\leq 2m}$.

\subsection{Counting in Free groups}\label{subsec: Conjugacy class count: Upper bounds: Free groups}
Fix a positive integer $r\geq 2$ and let $F_r$ denote the free group of rank $r$, equipped with the standard symmetric generating set of $2r$ elements. Note that the volume entropy of $F_r$ equals $h_r\coloneq \log(2r-1)$. 

\begin{lemma}[\cite{culler1981using},\cite{ComerfordEdmunds1994}]\label{lem: canonical forms for conjugacy classes}
    Let $n$ be a positive integer.
    \begin{enumerate}[font=\normalfont]
        \item If $[g]$ has orientable genus $n$, then there exists an orientable Wicks form $W$ of genus $n$ and a cyclically reduced non-canceling substitution $\mathfrak{g}$ such that $[\mathfrak{g}(W)]=[g]$ and $l(\mathfrak{g}(W))=l([g])$.
        \item If $[g]$ has non-orientable genus $n$, then either of the following holds true:
        \begin{enumerate}[font=\normalfont]
            \item There exists a non-orientable Wicks form $W$ of non-orientable genus $n$ and a cyclically non-canceling substitution $\mathfrak{g}$ such that $[\mathfrak{g}(W)]=[g]$ and $l(\mathfrak{g}(W))=l([g])$.
            \item There exists an orientable Wicks form $W'$ of orientable genus $(n-1)/2$ and a cyclically non-canceling substitution $\mathfrak{g}'$ such that $[\mathfrak{g}'(W')]=[g]$ and $l(\mathfrak{g}(W))=l([g])$.
        \end{enumerate}
        The latter is possible only when $n$ is odd and $\textnormal{genus}^+([g])=(n-1)/2$.
    \end{enumerate}
\end{lemma}
\begin{proof}
    In either case, consider any representative $g$ of $[g]$ and let $(W_t,\mathfrak{g}_t)$ be the tracking pair associated with $g$ (cf. \Cref{corollary: terminalwordexistence}). Note that $W_t$ is apriori only a pre-Wicks form. We observe that the substitution $\mathfrak{g}_t$ is non-canceling: Consider two consecutive letters $X_i$ and $X_{i+1}$ in $W_t$. The fact that \linkref{Move 1}{Move 1} and \linkref{Move 3}{Move 3} cannot be applied to $(W_t,\mathfrak{g}_t)$ implies that there cannot be any cancellation between $g_{X_{i}}$ and $g_{X_{i+1}}$. As described in \Cref{prop: what does a preWicks form look like}, a pre-Wicks form comes in four forms. If $W_t$ is a Wicks form, then by applying our moves to cyclic permutations of $W_t$ (if necessary), we get a tuple $(W,\mathfrak{g})$ such that $\mathfrak{g}$ is cyclically non-canceling, $W$ is a Wicks form, and $[\mathfrak{g}(W)]=[g]$. Since $\mathfrak{g}(W)$ is cyclically reduced, $l(\mathfrak{g}(W))=l([g])$. If $W_t$ is of the form \ref{item: preWicks form 2 - tail reduction}, or \ref{item: preWicks form 3 - tail redundancy}, or \ref{item: preWicks form 4 - tail redundancy nested in tail reduction}, then the tail reduction or redundancy in $W_t$ can be removed by a conjugation.   
\end{proof}

\noindent It follows from \Cref{lem: canonical forms for conjugacy classes} that the set $\partial \mathcal{C}^{\pm}_n(L)$ is empty when $L$ is odd. Therefore in the rest of this subsection, the parameter $L$ is always assumed to be even.

\begin{definition}\label{def: set of cyclically cancellation free subs - conjugacy classes}
    Let $W$ be a Wicks form and let $\vec{l}$ be a length type for $W$. Let $V$ be a linear word representing $W$. Define $G_V(\vec{l}\,)$ to be the set of group elements coming from cyclically non-canceling substitutions for $V$ of length type $\vec{l}_V$. Define $\mathcal{G}_W(\vec{l}\,)$ to be the set of conjugacy classes of group elements in $G_V(\vec{l})$. 
\end{definition}

Let $V=X_1\cdots X_{2m}$ be a linear Wicks form and let $\mathfrak{g}$ be a substitution. For any cyclic permutation $\sigma\in\text{Aut}(V)$, we define a new substitution $\sigma\cdot\mathfrak{g}$ as follows: for all $1\leq i\leq 2m$, $\sigma\cdot\mathfrak{g}(X_i)\coloneq g_{X_{\sigma(i)}}$. We will also say that a length type $\vec{l}$ is \emph{non-degenerate} if no two unknowns get assigned the same value.

\begin{definition}\label{def: the numbers kappa}
    Given any Wicks form $W=X_1\cdots X_{2m}$, define $\kappa_W=\kappa_W(r)$ to be the number of distinct cyclically non-canceling choices of starting and ending letters of $X_1,\dots, X_{2m}$. That is, $\kappa_W$ is the number of distinct choices of $\{(s_j, t_j)\}_{j=1}^{2m}$, with each $s_j, t_j \in S$ such that:
         \begin{enumerate}
         \item For all $1\leq j\leq 2m$: 
         \begin{align*}
             (s_{\overline{j}},t_{\overline{j}})=\begin{cases}
             (t_j^{-1},s_j^{-1}) &\text{if sign}(X_j)=-1\\
             (s_j,t_j) &\text{if sign}(X_j)=+1
             \end{cases}
         \end{align*}
        \item For all $1\leq j\leq 2m-1$, $t_j \cdot s_{j+1}\neq 1$, and $t_{2m} \cdot s_1\neq 1$.
    \end{enumerate}
    Given any positive integers $n\geq 1$ and $r\geq 2$, define:
    \begin{align*}
        \kappa^+_n&=\kappa^+_n(r)\coloneq (2r(2r-1)(2r-2))^{4n-2}\\ \kappa^-_n&=\kappa^-_n(r)\coloneq (2r(2r-1)(2r-2))^{2n-2} 
    \end{align*}
\end{definition}

\begin{lemma}\label{lem: kappa for maximal Wicks forms}
    For all maximal orientable pre-Wicks forms $W$ of genus $n$, $\kappa_W=\kappa^+_n$. For all maximal non-orientable pre-Wicks forms $W$ of genus $n$, $\kappa_W=\kappa^-_n$. 
\end{lemma}
% \begin{proof}
    
% \end{proof}

\begin{definition}\label{def: automorphism group and mass}
    Let $W$ be a Wicks form of length $2m$ and let $V$ be a linear word representing it. $\text{Aut}(V)$ is defined to be the set of cyclic permutations $\sigma\in\mathbb{Z}/2m\mathbb{Z}$ such that $\sigma(V)$ is linearly equivalent to $V$. For any cyclic permutation $\sigma$ of $V$, $\text{Aut}(\sigma(V))=\text{Aut}(V)$, so $\text{Aut}(W)$ is well-defined. Let $n$ be a positive integer. Define the \emph{mass} of orientable (and non-orientable respectively) Wicks forms of genus $n$ as:
    \[\mathfrak{M}_n^+\coloneq \sum_{\substack{\text{genus}^+(W)=n\\\text{maximal}\\\text{orientable}}}\frac{1}{|\text{Aut}(W)|},\quad \mathfrak{M}_n^-\coloneq \sum_{\substack{\text{genus}^-(W)=n\\\text{maximal}\\\text{non-orientable}}}\frac{1}{|\text{Aut}(W)|}\]
    
\end{definition}

Continuing with the notation in \Cref{def: automorphism group and mass}, given a Wicks form $W$ in the set of unknowns $\mathcal{X}$, and $V$ a linear word representing it, $\text{Aut}(V)$ acts on the set of length types for $W$. To see this, note that for any cyclic permutation $\sigma$ of $V$, $\sigma\cdot \vec{l}_V$ is a length type for $V$ if and only if for all $i, \sigma(\overline{i})=\overline{\sigma(i)}$, that is $\sigma\in\text{Aut}(V)$. Moreover, note that any $\sigma\in\text{Aut}(V)$ determines a of permutation of $\mathcal{X}\cup\mathcal{X}^{-1}$.  We describe some examples below:

\begin{enumerate}[label=(\roman*)]
    \item $V=ABCA^{-1}B^{-1}C^{-1}$. This is the orientable genus $1$ maximal Wicks form. The cyclic permutations $\sigma=(123456)$ generates $\text{Aut}(V)$, and therefore $|\text{Aut}(V)|=6$. The induced action on the set of unknowns is given by:
    \[A\mapsto B,\quad B\mapsto C,\quad C\mapsto A^{-1}\]
    \item $V= ABCA^{-1}CB$. This is a maximal non-orientable Wicks form of genus $2$. In this case, $\sigma=(123456)^3=(14)(23)(56)$ generates $\text{Aut}(V)$, and therefore $|\text{Aut}(V)|=2$. The induced action on the set of unknowns is given by:
    \[A\mapsto A^{-1},\quad B\mapsto C,\quad  C\mapsto B\]
\end{enumerate}

    \begin{lemma}\label{lemma: action of automorphisms on a length type}
        Let $W$ be a maximal Wicks form and let $\vec{l}$ be a non-degenerate length type for $W$. Let $\sigma$ be the cyclic permutation generating $\text{Aut}(W)$, and suppose it has order $n>1$. Consider the collection of length types $\vec{l},\sigma\cdot\vec{l},\dots,\sigma^{n-1}\cdot\vec{l}$.
        \begin{enumerate}[font=\normalfont]
            \item\label{item: When V is genus 1 maximal} If $W=ABCA^{-1}B^{-1}C^{-1}$, then $\sigma^3=1$, therefore there are three distinct length types in the collection: $\vec{l}$, $\sigma\cdot\vec{l}$, and $\sigma^2\cdot\vec{l}$.
            \item\label{item: When V is not genus 1 maximal} If $W\neq ABCA^{-1}B^{-1}C^{-1}$, then for all $0\leq i\neq j\leq n-1$, $\sigma^i\cdot\vec{l}\neq \sigma^j\cdot\vec{l}$.
        \end{enumerate}
    \end{lemma}
\begin{proof}
    Let $V=X_1\cdots X_{2m}$ be a linear maximal Wicks form representing $W$. The case when $V=ABCA^{-1}B^{-1}C^{-1}$ is easily verified. So we assume that $V\neq ABCA^{-1}B^{-1}C^{-1}$.
    
    Let $0\leq i_1<i_2\leq n-1$ be such that $\sigma^{i_1}\cdot\vec{l}_V=\sigma^{i_2}\cdot\vec{l}_V$, that is $\vec{l}_V=\sigma^{i_2-i_1}\cdot\vec{l}_V$. Setting $\alpha=\sigma^{i_2-i_1}$, this means that for all $1\leq i\leq 2m$, $l_{X_i}=l_{X_{\alpha(i)}}$. Since $\vec{l}$ is non-degenerate and $\alpha$ is non-trivial, this implies that for all $1\leq i\leq 2m$, $\alpha(i)=\overline{i}$. But $\alpha$ is a cyclic permutation, so $V=X_1X_2\cdots  X_mX_{\overline{1}}X_{\overline{2}}\cdots X_{\overline{m}}$. We argue below that such a \textit{palindromic} $V$ cannot be a maximal Wicks form. Thus, $\sigma^{i_1}\cdot\vec{l}_V\neq \sigma^{i_2}\cdot\vec{l}_V$.
    
    Note that since $V\neq ABCA^{-1}B^{-1}C^{-1}$ and $V\neq A^2$, so if $V$ is orientable, $\text{genus}^+(V)\geq 2$, so $m\geq 9$. If $V$ is non-orientable, $\text{genus}^-(V)\geq 2$, so $m\geq 3$.
    
    Recall the graph $\Gamma_V$ associated with $V$ (\Cref{def: wicks form genus}). Before we begin, let us note down a simple observation:
    \begin{quote}
        In a maximal Wicks form $V$, if a letter $X$ has sign $-1$, then its image in $\Gamma_V$ cannot be a loop.
    \end{quote}
    
    \begin{enumerate}[label=(\roman*)]
        \item Suppose $\text{sign}(X_1)=-1$. We write $X_1=X$, $X_2=Y$, $X_3=Z$, $X_4=U$.
        \begin{itemize}
            \item Suppose $\text{sign}(Y)=+1$. Then $V$ looks like $XY\cdots X^{-1}Y$. But this means that the image of $X$ in $\Gamma_V$ is a loop, which is not possible.
            \item Suppose $\text{sign}(Y)=-1$. Note that $Z\neq X^{-1}$, otherwise $V=XYX^{-1}Y^{-1}$ which is not a maximal Wicks form. Now if $\text{sign}(Z)=+1$, then $V$ looks like $XYZ\cdots X^{-1}Y^{-1}Z$ and the image of $Y$ in $\Gamma_V$ is a loop which is not possible. Therefore $\text{sign}(Z)=-1$, and $V$ is of the form $XYZ\cdots X^{-1}Y^{-1}Z^{-1}$. Therefore, in $\Gamma_V$, we see a trivalent vertex, say $v_0$, with $X$, $Y$, and $Z$ sharing it as a vertex. Since $V\neq XYZX^{-1}Y^{-1}Z^{-1}$, there has to be a distinct letter $U$, so that $V$ is of the form: $XYZU\cdots X^{-1}Y^{-1}Z^{-1}U^{\varepsilon}$. But this is not possible as it forces $v_0$ to have degree more than $3$, contradicting maximality of $V$. 
        \end{itemize}
        Thus, $\text{sign}(X_1)\neq-1$.
        \item Suppose $\text{sign}(X_1)=+1$. We write $X_1=X$ and $X_2=Y$. We must have $\text{sign}(Y)=-1$, otherwise $V$ has a redundancy. This means $V$ is of the form $XY\cdots XY^{-1}$. But this means that the image of $Y$ in $\Gamma_V$ is a loop, which is not possible.  
    \end{enumerate}
\end{proof}

\begin{theorem}{\cite[Theorem 1.1]{BACHER200213}}\label{thm: orientable mass formula}
    Let $W$ be a maximal orientable Wicks form. Then, $\text{Aut}(W)$ is cyclic of order $1$, $2$, $3$, or $6$, and for all $n\geq 1$:
    \[\mathfrak{M}_n^+=\frac{2}{12^n}\frac{(6n-5)!}{n!(3n-3)!}\]
\end{theorem}

\begin{lemma}\label{lem: substitution count asymptote conj class}
    Fix \(k > 0\). Let $V$ be a linear Wicks form of genus $n$ and length $2m$. Let $\vec{l}$ be a length type for $V$ such that each entry of \(\vec{l}\) is at least \(k \log L \). Let $C_1$, $C_2$ be the constants in \Cref{lem: countingpinnedreducedwords} and let $L=||\vec{l}||$. Then:
    \begin{equation*}
        \left (1+ \frac{C_1}{r^{k\log L-2}} \right)^m\leq \frac{|G_V(\vec{l})|}{\displaystyle{\kappa_W \cdot\frac{(2r - 1)^{\frac{||\vec{l}||}{2} - m}}{\displaystyle (2r)^m}}}\leq \left(1+\frac{C_2}{r^{k\log L-2}}\right)^{m}
    \end{equation*}
\end{lemma}
\begin{proof}
    The proof of this lemma is identical to that of \Cref{lem: substitution count asymptote}, the only difference lies in the fact that we're counting cyclically non-canceling substitutions this time.
\end{proof}

\begin{lemma}\label{lem: conj class substitution count asymptote}
    Fix \(k > 0\). Let $W$ be a Wicks form of length $2m$ and let $\vec{l}$ be a non-degenerate length type for $W$ such that each entry of \(\vec{l}\) is at least \(k \log L \), where $L=||\vec{l}||$. Let $C_2$ be the constant appearing in \Cref{lem: countingpinnedreducedwords}. Then:
    \begin{align*}
        \text{If }W\neq ABCA^{-1}B^{-1}C^{-1},\quad |\mathcal{G}_W(\vec{l})|&\leq \left(1+\frac{C_2}{r^{k\log L-2}}\right)^{m} \left(\displaystyle{\kappa_W \cdot\frac{(2r - 1)^{\frac{||\vec{l}||}{2} - m}}{\displaystyle (2r)^m}}\right)\\
        \text{If }W= ABCA^{-1}B^{-1}C^{-1},\quad |\mathcal{G}_W(\vec{l})|&\leq \frac{1}{2}\left(1+\frac{C_2}{r^{k\log L-2}}\right)^{3} \left(\displaystyle{\kappa_W \cdot\frac{(2r - 1)^{\frac{||\vec{l}||}{2} - 3}}{\displaystyle (2r)^3}}\right)
    \end{align*}
\end{lemma}
\begin{proof}
    Let $V$ be any linear Wicks form representing $W$. Then, by definition, $|\mathcal{G}_W(\vec{l})|\leq |G_V(\vec{l}_V)|$. Since $(V,\vec{l}_V)$ satisfies the hypotheses of \Cref{lem: substitution count asymptote conj class}, the estimates follows directly when $W\neq ABCA^{-1}B^{-1}C^{-1}$. When $W=ABCA^{-1}B^{-1}C^{-1}$, for any cyclically non-canceling substitution $\mathfrak{g}$ for $ABCA^{-1}B^{-1}C^{-1}$ of length type $\vec{l}$, consider the substitution $\mathfrak{g}'$ defined by: $g'_A\coloneq g_A^{-1}, g'_B\coloneq g_B^{-1}, g'_C\coloneq g_C^{-1}$. Then $[\mathfrak{g}(ABCA^{-1}B^{-1}C^{-1})]=[\mathfrak{g}'(ABCA^{-1}B^{-1}C^{-1})]$ but $\mathfrak{g}(ABCA^{-1}B^{-1}C^{-1})\neq \mathfrak{g}'(ABCA^{-1}B^{-1}C^{-1})$. Therefore, in this case, $|\mathcal{G}_W(\vec{l})|\leq \frac{1}{2}|G_V(\vec{l}_V)|$, and by \Cref{lem: substitution count asymptote conj class}, we are done.
\end{proof}

\begin{proposition}\label{thm: free grp shell conj class upper bound asymptote}
    Let $n$ and $r$ be positive integers with $r\geq 2$.
    \begin{equation*}
        \limsup_{\substack{L\to\infty\\L\textnormal{ even}}}\frac{|\partial \mathcal{C}^{+}_n(L)|}{L^{6n-4}e^{\frac{h_rL}{2}}}\leq \frac{\kappa^+_n(r)\mathfrak{M}^+_n}{2^{6n-4}(6n-4)!(2r)^{6n-3}(2r-1)^{6n-3}}
    \end{equation*}
\end{proposition}
\begin{proof}
We begin with a simple observation: given a Wicks form $W$ of length $2m$, a length type $\vec{l}$, and a $\sigma\in\text{Aut}(W)$, $\mathcal{G}_W(\vec{l})=\mathcal{G}_W(\sigma\cdot \vec{l})$. 

Further, given $L>0$ we define two families of length types for $W$ as follows:
    \begin{align*}
        &\mathcal{L}_0(W,m,r,L)\coloneq\left\{\vec{l}\in\mathbb{Z}_{\geq 0}^{m} \middle |\; \sum_{i=1}^{m}l_i= \frac{L}{2} \right\}\\
        &\mathcal{L}_1(W,m,r,L)\coloneq \left\{\vec{l}\in\mathbb{Z}_{\geq 0}^{m} \middle |\; \sum_{i=1}^{m}l_i= \frac{L}{2}, l_i>2\;(1\leq i\leq m), \vec{l}\text{ is non-degenerate}\right\}\\
        &\mathcal{L}_2(W,m,r,L)\coloneq \left\{\vec{l}\in\mathbb{Z}_{\geq 0}^{m} \middle |\; \sum_{i=1}^{m}l_i= \frac{L}{2},\, l_i\geq \log L \; (1\leq i\leq m), \vec{l}\text{ is non-degenerate} \right\}\\
        &\tilde{\mathcal{L}}_i(W,m,r,L)\coloneq \mathcal{L}_i(W,m,r,L)/\text{Aut(W)}, \quad (0\leq i\leq 2)
    \end{align*}
    
    We will suppress the arguments of $\mathcal{L}_1,\mathcal{L}_2$, $\tilde{\mathcal{L}}_1,\tilde{\mathcal{L}}_2$ below whenever they are clearly understood.
    \begin{claim}\label{claim: claim in free grp conj class upper bound asymptote proof}
    \leavevmode
        \begin{enumerate}
            \item For any non-maximal orientable Wicks form $W$, we have:
            \begin{equation*}
                \lim_{\substack{L\to\infty\\L\textnormal{ even}}} \frac{\sum_{\vec{l}\in \tilde{\mathcal{L}}_1} |\mathcal{G}_W(\vec{l}\,)|}{L^{6n-4}e^{\frac{h_rL}{2}}} =0 
            \end{equation*}
            \item For any maximal orientable Wicks form $W$:
            \begin{equation*}
                \lim_{\substack{L\to\infty\\L\textnormal{ even}}}\frac{\sum_{\vec{l}\in\tilde{\mathcal{L}_1}-\tilde{\mathcal{L}}_2}|\mathcal{G}_W(\vec{l})|}{L^{6n-4}e^{\frac{h_rL}{2}}}=0
            \end{equation*}
            where $\tilde{\mathcal{L}}_1-\tilde{\mathcal{L}}_2\coloneq (\mathcal{L}_1-\mathcal{L}_2)/\text{Aut}(W)$.
            \item For any maximal orientable Wicks form $W$:
            \begin{equation*}
                \limsup_{\substack{L\to\infty\\L\textnormal{ even}}}\frac{\sum_{\vec{l}\in \tilde{\mathcal{L}}_2}|\mathcal{G}_{W}(\vec{l})|}{L^{6n-4}e^{\frac{h_rL}{2}}}\leq \frac{1}{\text{Aut}(W)}\frac{\kappa_n^+(r)}{2^{6n-4}(6n-4)!(2r)^{6n-3}(2r-1)^{6n-3}}
            \end{equation*}
            where $\kappa^+_n(r)$ is the constant defined in \Cref{def: the numbers kappa}.
        \end{enumerate}
    \end{claim}
     The proof of \Cref{claim: claim in free grp conj class upper bound asymptote proof} is completely analogous to that of \Cref{claim: claim in free grp upper bound asymptote proof}, where \Cref{lem: conj class substitution count asymptote} serves the role of the upper bound in \Cref{lem: substitution count asymptote}. A point to be noted here is that when $W=ABCA^{-1}B^{-1}C^{-1}$, we do get the desired $\frac{1}{\text{Aut}(W)}$ factor: $\frac{1}{3}$ comes from counting length types, while $\frac{1}{2}$ comes from the $\mathcal{G}_W(\vec{l})$ estimate.
        
    For part $(3)$, note that $\lim_{L\to\infty}\frac{|\mathcal{L}_2(W,n,r,L)|}{L^{6n-4}}=\frac{1}{2^{6n-4}(6n-4)!}$, because asymptotically, most length types satisfying $\sum l_i=L/2$ are non-degenerate. This is easily verified by a computation, but the heuristic here is that all degenerate length types lie on a finite number of hyperplanes. Therefore, for a maximal Wicks form $W$, by \Cref{lemma: action of automorphisms on a length type}:
    \begin{align*}
        \text{If }W\neq ABCA^{-1}B^{-1}C^{-1},\quad\lim_{L\to\infty}\frac{|\tilde{\mathcal{L}}_2(W,n,r,L)|}{L^{6n-4}}&=\frac{1}{\text{Aut}(W)}\frac{1}{2^{6n-4}(6n-4)!}\\
        \text{If }W= ABCA^{-1}B^{-1}C^{-1},\quad\lim_{L\to\infty}\frac{|\tilde{\mathcal{L}}_2(W,n,r,L)|}{L^{6n-4}}&=\frac{1}{3}\frac{1}{2^{6n-4}(6n-4)!}\quad (n=1)
    \end{align*}
   
    The proposition follows from the \Cref{claim: claim in free grp conj class upper bound asymptote proof}, and the fact that (\Cref{lem: canonical forms for conjugacy classes}):
    \[\partial \mathcal{C}^+_n(L)\subset\bigcup_{W}\bigcup_{\vec{l}\in \tilde{\mathcal{L}}_0(W,L)}\mathcal{G}_W(\vec{l}\,)\]
    and the fact that for any two maximal orientable Wicks forms of the same genus, $|\mathcal{L}_2(W,n,r,L)|=|\mathcal{L}_2(W',n,r,L)|$. 
\end{proof}
Similarly, we have:
\begin{proposition}\label{thm: free grp shell conj class upper bound asymptote non-orientable}
    Let $n>1$ and $r\geq 2$ be positive integers.
    \begin{equation*}
        \limsup_{\substack{L\to\infty\\L\textnormal{ even}}}\frac{|\partial \mathcal{C}^{-}_n(L)|}{L^{3n-4}e^{\frac{h_rL}{2}}}\leq \frac{\kappa^-_n(r)\mathfrak{M}^-_n}{2^{3n-4}(3n-4)!(2r)^{3n-3}(2r-1)^{3n-3}}
    \end{equation*}
\end{proposition}
\begin{proof}
    The proof follows using the same arguments used in the proof of \Cref{thm: free grp shell conj class upper bound asymptote}, with the estimates in \Cref{claim: claim in free grp conj class upper bound asymptote proof} replaced by their non-orientable analogues. However there is one subtlety that needs to be accounted for. By \Cref{lem: canonical forms for conjugacy classes}, there are two types of conjugacy classes in $\partial\mathcal{C}^-_n(L)$: those that arise from cyclically non-canceling substitutions of non-orientable Wicks forms $W$ with $\textnormal{genus}^-(W)=n$, and those that arise from cyclically non-canceling substitutions of orientable Wicks forms $W'$ with $\textnormal{genus}^+(W')=(n-1)/2$. However, the number of conjugacy classes of the second type is $O(L^{3n-7}e^{\frac{h_rL}{2}})$ and therefore do not contribute to the estimate in the statement of the proposition.
\end{proof}

\subsection{Counting in Hyperbolic groups}\label{subsec: Conjugacy class count: Upper bounds: Hyperbolic groups}
Throughout this section, $G$ will denote a hyperbolic group. We will equip $G$ with the word metric induced by some symmetric generating set $S$.

Let $\mathcal{C}^+_n(L)$ denote the number of conjugacy classes in $G$ of orientable genus $n$ and length at most $L$. Let $\mathcal{C}^-_n(L)$ denote the number of conjugacy classes in $G$ of non-orientable genus $n$ and length at most $L$.

\begin{theorem}\label{thm: conjugacy count upper bound orientable hyperbolic}

Let $n\geq 1$ be a positive integer. Suppose $G$ is a hyperbolic group and $S$ is a finite generating set so that $G$ is $\delta$-hyperbolic. Then there are constants $C, L_0>0$ depending only on $\delta, n, G$ and the choice of the generating set, so that 

\begin{align*}
    \forall L\geq L_0,\quad  |\mathcal{C}^{+}_n(L)| \leq  C L^{6n-4} e^{hL/2}
\end{align*}

where $h$ is the volume entropy of $G$.
\end{theorem}

\begin{theorem}\label{thm: conjugacy count upper bound non-orientable hyperbolic}
    Let $n>1$ be a positive integer. Suppose $G$ is a hyperbolic group and $S$ is a finite generating set so that $G$ is $\delta$-hyperbolic. Then there are constants $C', L'_0 > 0$ depending only on $\delta,n,G$ and the choice of generating set, such that:
    \[
        \forall L\geq L'_0,\quad  |\mathcal{C}^{-}_n(L)| \leq  C' L^{3n-4} e^{hL/2}
    \]
\end{theorem}

The proofs of \Cref{thm: conjugacy count upper bound orientable hyperbolic} and \Cref{thm: conjugacy count upper bound non-orientable hyperbolic} are analogous to those of \Cref{thm: element count upper bound orientable hyperbolic} and \Cref{thm: element count upper bound non-orientable hyperbolic} respectively. It suffices to prove the following:

\begin{proposition}\label{prop: long short decomposition in conjugacy}
    Let $n\geq1$ be a positive integer. Let $[g]$ be a conjugacy class in $G$ of orientable genus $n$ and length $L$. Then there exists a pair $(W,\mathfrak{g})$ where $W$ is an orientable Wicks form of genus $n$ and $\mathfrak{g}$ is a substitution for $W$ such that $[\mathfrak{g}(W)]=[g]$ and $W$ admits a decomposition as $W=S_1L_1S_2\dots S_mL_mS_{m+1}$ where $L_1,\dots,L_m$ are the long edges of $W$ and $S_1,\dots, S_{m+1}$ are products of short edges (allowed to be the empty word). Suppose $S_j$ is a product of $s_j$ many short edges, for $1\leq j\leq m+1$. Then we have:
    \begin{align*}
        m+\sum_{j=1}^{m+1}s_j&\leq 12n-6\\
        \sum_{j=1}^{m}l(g_{L_j})&\leq L + 12nD^+
    \end{align*}
    where, $D^+=D^{+}(n)=4l_0^{+}(n)+M+2$, $l_0^+(n)=2\delta\lceil\log_2(12n-1)\rceil$ and $M$ is the constant appearing in \Cref{lemma: conjugationbound}.
\end{proposition}
\begin{proof}
    Let $g$ be an element of the conjugacy class $[g]$ which realizes its length, that is, $l(g)=l([g])=L$. By \Cref{corollary: terminalwordexistence}, there exists a tracking pair (\Cref{def: tracking pair}) $(W_t,\mathfrak{g}_t)$ such that $\mathfrak{g}_t(W_t)=g$. $W_t$ can be written as $W_t=S_1L_1\dots S_mL_mS_{m+1}$ where $L_i$ are the long edges of $W_t$ and each $S_j$ is a product of short edges (allowed to be the empty word), say $s_j$ many, and we have $m+\sum_{j=1}^{m+1}s_j=l(W_t)$. By \Cref{lem: length of a preWicks form}, $l(W_t)\leq 12n-2$ and by \Cref{lem: longedgeparallelfellow-travelling}, we have:
    \begin{align}\label{eqn: sum of lengths}
        \sum_{j=1}^{m}l(g_{L_j})\leq L + (12n-2)D^+
    \end{align}
    If $W_t$ is a genuine Wicks form, then we set $W\coloneq W_t$, $\mathfrak{g}\coloneq \mathfrak{g}_t$, and in this case $l(W_t)\leq 12n-6$, so there's nothing extra to be done. But $W_t$ is a pre-Wicks form, so it could either have a tail reduction, a tail redundancy, or a tail redundancy nested in a reduction (cf. \Cref{def: pre-Wicks form}). The last case encapsulates the others. So suppose $W_t=DAU_1A^{-1}B^{-1}U_2BD^{-1}$. We set $W\coloneq U_1X^{-1}U_2X$ and we set $\mathfrak{g}$ the same as $\mathfrak{g}_t$ for every unknown other than $X, A, B$ and $D$, with $g_X\coloneq (g_t)_B(g_t)_A$. By definition, $W$ is a Wicks form, and clearly $\mathfrak{g}(W)\in c$. $W$ also admits a decomposition $W=S'_1L'_1\dots S'_{m'}L'_{m'}S'_{m'+1}$ where $L'_i$ are the long edges of $W$ and each $S'_j$ is a product of short edges, say $s'_j$ many. Note that $m'+\sum_{j=1}^{m'+1}s'_j=l(W_t)-4\leq 12n-6$. In the process of obtaining $W$ from $W_t$, we have dropped the letter $D$ and therefore potentially lost a long letter, and resolved a redundancy $BA$ and therefore potentially gained a long letter (in case $A$ and $B$ are short letters but $l((g_t)_B(g_t)_A)\geq D^+$). Thus \Cref{eqn: sum of lengths} holds true for the long letters $L'_j$ $(1\leq j\leq m')$, but with an additive factor of $2D^+$ to the upper bound.
\end{proof}

Reasoning similarly, we have a non-orientable version of \Cref{prop: long short decomposition in conjugacy}
\begin{proposition}\label{prop: long short decomposition in conjugacy non-orientable}
    Let $n\geq1$ be a positive integer. Let $[g]$ be a conjugacy class in $G$ of non-orientable genus $n$ and length $L$. Also assume that $\textnormal{genus}^+([g])\geq n/2$. Then there exists a pair $(W,\mathfrak{g})$ where $W$ is an non-orientable Wicks form of genus $n$ and $\mathfrak{g}$ is a substitution for $W$ such that $[\mathfrak{g}(W)]=[g]$ and $W$ admits a decomposition as $W=S_1L_1S_2\dots S_mL_mS_{m+1}$ where $L_1,\dots,L_m$ are the long edges of $W$ and $S_1,\dots, S_{m+1}$ are products of short edges (allowed to be the empty word). Suppose $S_j$ is a product of $s_j$ many short edges, for $1\leq j\leq m+1$. Then we have:
    \begin{align*}
        m+\sum_{j=1}^{m+1}s_j&\leq 6n-6\\
        \sum_{j=1}^{m}l(g_{L_j})&\leq L + 6nD^- 
    \end{align*}
    where $D^-=D^{-}(n)=4l_0^{-}(n)+M+2$, $l_0^-(n)=2\delta\lceil\log_2(6n-1)\rceil$ and $M$ is the constant appearing in \Cref{lemma: conjugationbound}.
\end{proposition}

\section{Conjugacy class count: Lower bounds}\label{sec: Conjugacy class count: Lower bounds}

\subsection{Combinatorial Reduction: Revisited}\label{subsec: Combreduction revisited}

In order to obtain the lower bounds for counting conjugacy classes, we need to refine our analysis of combinatorial reduction. We will be using the same notation and definitions introduced in \Cref{subsec: combinatorialreduction}.

\begin{definition}[Split Wicks forms]\label{def: split wicks forms}
A \emph{split Wicks form} is a pre-Wicks form of the form $AU(BA)^{\varepsilon}VB$ such that $XUX^{\varepsilon}V$ is a Wicks form, where $\varepsilon\in \{\pm 1\}$, $A,B$, and $X$ are letters that do not appear in $U$ and $V$. These were introduced in \Cref{prop: what does a preWicks form look like} part \ref{item: preWicks form 3 - tail redundancy}. 
\end{definition}

\begin{remark}\label{rmk: split Wicks forms}
    In what follows, we will write a Wicks form of length $2m$ as $W=X_1\cdots X_{2m}$, and a split Wicks form of length $2m+2$ as $W'_{\text{split}}=X'_1\cdots (BX'_1)^{\varepsilon}\cdots X'_{2m}B$. Associated with $W'_{\text{split}}$ is the Wicks form $W'$ obtained by cycling $B$ up to the front and removing the redundancy $BX'_1$ from $W'_{\text{split}}$. Moreover, if $W'$ is a non-orientable Wicks form then $\text{genus}^-(W')>1$. This is because there is exactly one non-orientable genus $1$ Wicks form, namely $X^2$, and while it can be `split' as $XBXB$, it is not a pre-Wicks form. 
\end{remark}

\begin{proposition}\label{prop: another main reduction theorem}
    Let $c, A_0, H, D_0$ be non-negative integers with $c, D_0\geq 1$. Let $W$ be a maximal Wicks form of length $2m$ and let $W'_{\text{split}}$ be a split Wicks form of length $2m+2$. Let $\vec{l}$ and $\vec{l}'$ be length types for $W$ and $W'_{\text{split}}$ respectively. Furthermore, let $\vec{b}\in [A_0,A_0+H]^{2m+1}$ and $\vec{b}'\in [A_0,A_0+H]^{2m+3}$. Now suppose the following holds:
    \begin{enumerate}[font=\normalfont]
        \item[$(\mathbf{L'1})$]\label{item: L'1} $L(W,\vec{l},\vec{b})=L(W'_{\text{split}},\vec{l}',\vec{b}')$.
        \item[$(\mathbf{L'2})$]\label{item: L'2} For all $1\leq i\leq 2m$, $1\leq j\leq 2m$, $j\neq 1,\overline{1}$, $l_i,l'_j> 8c+10(A_0+H)+4D_0$, $l'_1=l'_{\overline{1}'}> 4c+4(A_0+H)+2D_0$, and $l'_B\leq 4c+4(A_0+H)+2D_0$.  
    \end{enumerate}
    Then, either:
    \begin{enumerate}[font=\normalfont]
        \item There is $D_0$-reduction in degrees of freedom by $c$, or
        \item $W=W'$, and $l'_B\leq 2(D_0+H-1)$. 
    \end{enumerate}
\end{proposition}

\begin{proposition}\label{prop: reduction or no splitting}
    Let $c, A_0, H,D_0$ be non-negative integers with $c,D_0\geq 1$. Let $W$ be a maximal Wicks form of length $2m$ and let $W'_{\text{split}}$ be a split Wicks form of length $2m+2$. Let $\vec{l}$ and $\vec{l}'$ be length types for $W$ and $W'_{\text{split}}$ respectively. Furthermore, let $\vec{b}\in [A_0,A_0+H]^{2m+1}$ and $\vec{b}'\in [A_0,A_0+H]^{2m+3}$. Now suppose the following holds:
    \begin{enumerate}[font=\normalfont]
        \item[$(\mathbf{L'1})$] $L(W,\vec{l},\vec{b})=L(W'_{\text{split}},\vec{l}',\vec{b}')$.
        \item[$(\mathbf{L'3})$] For all $1\leq i\leq 2m$, $1\leq j\leq 2m$, $j\neq 1,\overline{1}'$, $l_i,l'_j> 10c+25(A_0+H+D_0)$, and $l'_{X'_1}+l'_B> 10c+25(A_0+H+D_0)$.  
    \end{enumerate}
    Then, either:
    \begin{enumerate}[font=\normalfont]
        \item There is $D_0-$reduction in degrees of freedom by $c$, or
        \item Either $W=W'$, for all $1\leq i\leq 2m$, $|l_i-l'_i|\leq 11(D_0-1) + 16(A_0+H)$, and $l'_B\leq 2(D_0-1+H)$, or $W=\sigma(W')$, for all $1\leq i\leq 2m$, $|l_i-l'_{\sigma(i)}|\leq 11(D_0-1) + 16(A_0+H)$, and $l'_{X'_1}\leq 2(D_0-1+H)$, where $\sigma=(1 2 \dots 2m)$.
    \end{enumerate}
\end{proposition}
\begin{proof}
    Suppose there is no $D_0$-reduction in degrees of freedom by $c$. Since $l'_{X'_1}+l'_B> 10c+20(A_0+H+D_0)$, we cannot have both $l'_{X'_1}\leq 4c+4(A_0+H)+2D_0$ and $l'_B\leq 4c+4(A_0+H)+2D_0$. If both $l'_{X'_1}> 4c+4(A_0+H)+2D_0$ and $l'_B> 4c+4(A_0+H)+2D_0$, then $(W,\vec{l},\vec{b})$ and $(W'_{\text{split}},\vec{l}',\vec{b}')$ satisfy conditions \ref{item: L1} and \ref{item: L3} in \Cref{subsec: combinatorialreduction}. Therefore, by \Cref{prop: reduction or shifted matching for different wicks forms} and \Cref{rmk: reduction or equality strengthened}, $W$ and $W'_{\text{split}}$ are equivalent as pre-Wicks forms, and in particular have equal length, which is not possible. Therefore, exactly one among $l'_{X'_1}$ and $l'_B$ is strictly greater than $4c+4(A_0+H)+2D_0$. It follows from \Cref{prop: another main reduction theorem}, applied to either $W$, $W'_{\text{split}}$ or their inverses, that we must have: $W=W'$ and $l'_B\leq 2(D_0-1+H)$ or $W=\sigma(W')$ and $l'_{X'_1}\leq 2(D-0-1+H)$, where $\sigma$ is the cyclic permutation $(12\dots 2m)$. We argue the case when $l'_B\leq 2(D_0-1+H)$ as the other case is entirely similar.
    We formally promote $B$ to a junction. That is, we define $\vec{b}''$ as follows, for $0\leq i\leq 2m$:
    \[
    \text{For }\varepsilon=-1, b''_i\coloneq \begin{cases}
        b'_i&; i< \overline{1}\\
        b'_{\overline{1}}+l'_B+b'_{\overline{1}+1}&; i= \overline{1}\\
        b'_{i+1}&; \overline{1}<i<2m\\
        b'_{2m+1}+l'_B+b'_{2m+2}&;  i=2m\\
    \end{cases} 
    \]

    \[
    \text{For }\varepsilon=+1, b''_i\coloneq \begin{cases}
        b'_i&; i< \overline{1}-1\\
        b'_{\overline{1}-1}+l'_B+b'_{\overline{1}}&; i= \overline{1}-1\\
        b'_{i+1}&;  \overline{1}-1<i<2m\\
        b'_{2m+1}+l'_B+b'_{2m+2}&;  i=2m
    \end{cases} 
    \]
    Furthermore, we define $A'_0\coloneq A_0$, $H'=A_0+2(D_0-1)+4H$, and note that $\vec{b}''\in [A'_0,A'_0+H']^{2m+1}$.
    Next, we observe that $(W,\vec{l},\vec{b})$ and $(W',\vec{l}',\vec{b}'')$ satisfy conditions \ref{item: L1} and \ref{item: L2} of \Cref{thm: the main reduction theorem}. For \ref{item: L1}, it suffices to observe that $||\vec{b''}||=||\vec{b}'||$. For \ref{item: L2}, firstly note that because $l'_{X'_1}+l'_B>10c+ 25(A_0+H+D_0)$ and $l'_B\leq 2(D_0-1+H)$, we have $l'_{X'_1}>10c + 25A_0 +23(D_0+H) > 6c +5(A'_0+H') +2D_0 = 6c+ 10A_0+20H+12D_0$ and similarly, $l'_{X'_j}>6c +5(A'_0+H') +2D_0$ for all $2\leq j\leq 2m$. Moreover, there is no $D_0$-reduction by $c$ between $(W,\vec{l},\vec{b})$ and $(W',\vec{l}',\vec{b}'')$. So it follows from \Cref{thm: the main reduction theorem}, that for all $1\leq i\leq 2m$ $|l_i-l'_i|\leq 3(D_0-1)+4H'=11(D_0-1)+4(A_0+4H)$. For future reference, we note down the following estimate using \Cref{eq: shift bound in element count}:
    \begin{equation}\label{eq: shift bound conj class count}
        \text{for all }1\leq j\leq 2m, \; |\text{shift}(j)|\leq \frac{3}{2}(D_0-1+H')\leq 5(D_0-1)+ 2(A_0+4H)
    \end{equation}
\end{proof}
The case when $A_0=H=0$ and $D_0=1$ will be used in \Cref{subsec: Conjugacy class count: Lower bounds: Free groups}, and we record it separately.
\begin{proposition}\label{cor: reduction or no splitting - free groups}
    Let $c\geq 1$ be a positive integer. Let $W$ be a maximal Wicks form of length $2m$ and let $W'_{\text{split}}$ be a split Wicks form of length $2m+2$. Let $\vec{l}$ and $\vec{l}'$ be length types for $W$ and $W'_{\text{split}}$ respectively. Suppose the following holds:
    \begin{enumerate}[font=\normalfont, label=(\alph*)]
        \item\label{item: L'1 free} $||\vec{l}||=||\vec{l}'||$.
        \item\label{item: L'3 free} For all $1\leq i\leq 2m$, $1\leq j\leq 2m$, $j\neq 1,\overline{1}'$, $l_i,l'_j> 10c+25$, and $l'_{X'_1}+l'_B> 10c+25$.  
    \end{enumerate}
    Then, either:
    \begin{enumerate}[font=\normalfont]
        \item There is $1$-reduction in degrees of freedom by $c$, or
        \item Either $W=W'$, $\vec{l}=\vec{l}'$, and $l'_B=0$ or $W=\sigma(W')$, for all $1\leq i\leq 2m$, $l_i=l_{\sigma(i)}$, and $l'_{X'_1}=0$, where $\sigma = (1 2 \dots 2m)$.
    \end{enumerate}
\end{proposition}

\subsubsection{Tools for the proof of \Cref{prop: another main reduction theorem}}
Let $c,A_0,H$, $W,m,\vec{b},\vec{l}$, $W'_{\text{split}},\vec{b}',\vec{l}'$ satisfy $(\mathbf{L'1})$, $(\mathbf{L'2})$ as in \Cref{prop: another main reduction theorem}. These will be the standing assumptions for the rest of this subsection.
\begin{lemma}\label{lem: matching position extended conj class}
    Suppose there is a position $i$, contained in a letter $X_j$ of $W$ with respect to $\vec{l}$, while it is contained in a letter $X'_k$ of $W'_{\text{split}}$ with respect to $\vec{l}'$. Moreover, suppose that $\text{sign}(X_j)=\text{sign}(X'_k)$. Then:
    \begin{enumerate}[font=\normalfont]
        \item If $\underline{i}\simeq_{D_0-1}\underline{i}'$, then for every position $p$ common to $X_j$ and $X'_k$, we have $\underline{p}\simeq_{D_0-1}\underline{p}'$.
        \item The condition $\underline{i}\simeq_{D_0-1}\underline{i}'$ is equivalent to:
        \begin{itemize}
            \item $t_{X_{\overline{j}}} - t'_{X'_{\overline{k}'}} \simeq_{D_0-1} s'_{X'_{k}} - s_{X_j}$, when $\text{sign}(X_j)=\text{sign}(X'_k)=-1$.
            \item $s_{X_{\overline{j}}} - s'_{X'_{\overline{k}'}} \simeq_{D_0-1} s_{X_j} - s'_{X'_{k}}$, when $\text{sign}(X_j)=\text{sign}(X'_k)=+1$.
        \end{itemize}
    \end{enumerate}
\end{lemma}
\begin{proof}
    Identical to the proof of \Cref{lem: matching position extended}.
\end{proof}

\begin{lemma}\label{lem: common positions conj class}
    Suppose there is no $D_0$-reduction in degrees of freedom by $c$ between $(W,\vec{l},\vec{b})$ and $(W'_{\text{split}},\vec{l}',\vec{b}')$. Moreover, suppose there are $c_0\geq 4c+3(A_0+H)+D_0$ many consecutive positions $p+1,\dots, p+c_0$ common between $X_j$, a letter of $W$, and $X'_k$, a letter of $W'$. Then:
    \begin{enumerate}[font=\normalfont]
        \item $\text{sign}(X_j)=\text{sign}(X'_k)$.
        \item For all $1\leq i\leq c_0$, $\underline{p+i}\simeq_{D_0-1}\underline{(p+i)}'$.
    \end{enumerate}
\end{lemma}
\begin{proof}
    The proof of this lemma follows along the same lines as that of \Cref{lem: common positions}, with the nuance that we don't have a lower bound on $l'_B$. In particular, when $l'_B<c$, the arguments in the proof \Cref{claim: consecutiveadmissiblepositions} run unchanged, except that $B$ serves as a dummy junction of size at most $c$. This accounts for the additive factor of $c+(A_0+H)$ in the lower bound for $c_0$ when compared with \Cref{lem: common positions}.
    
\end{proof}

\begin{lemma}\label{lem: reduction or redundancy conj class}
    Suppose there is no $D_0$-reduction in degrees of freedom by $c$ between $(W,\vec{l},\vec{b})$ and $(W'_{\text{split}},\vec{l}',\vec{b}')$. Then none of the following are possible:
    \begin{enumerate}[font=\normalfont]
        \item A letter $X'_k$ of $W'_{\text{split}}$ has at least $4c+3(A_0+H)+D_0$ positions in common with each of two consecutive letters $X_j$, and $X_{j+1}$ of $W$.
        \item A letter $X_k$ of $W$ has at least $4c+3(A_0+H)+D_0$ positions in common with each of two consecutive letters $X'_j$ and $X'_{j+1}$ of $W_{\text{split}}$.
    \end{enumerate}
\end{lemma}
\begin{proof}
    $(1)$ and $(2)$ follow from \Cref{lem: common positions conj class} in the same way in which \Cref{lem: reduction or redundancy} follows from \Cref{lem: common positions}, using the fact that $W$ and $W'_{\text{split}}$ have no redundancies. 
\end{proof}

\subsubsection{Proof of \Cref{prop: another main reduction theorem}}
Let $c,D_0,A_0,H$, $W,m,\vec{b},\vec{l}$, $W'_{\text{split}},\vec{b}',\vec{l}'$ satisfy $(\mathbf{L'1})$ and $(\mathbf{L'2})$ as in \Cref{prop: another main reduction theorem}. Let $W'$ denote the Wicks form associated to $W'_{\text{split}}$ (see \Cref{rmk: split Wicks forms}). These will be the standing assumptions for the rest of this subsection.

\begin{proposition}\label{prop: reduction or equality of Wicks forms conj class}
    Suppose that there is no reduction in degrees of freedom by $c$ between $(W,\vec{l},\vec{b})$ and $(W'_{\text{split}},\vec{l}',\vec{b}')$. Then $W$ is equivalent to $W'$, and for every $1\leq i\leq 2m$, $X_i$ and $X'_i$ have at least $4c+3(A_0+H)+D_0$ positions in common. 
\end{proposition}
\begin{proof}
    Before we begin, note that $W=X_1\cdots X_{2m}$ and $W'_{\text{split}}=X'_1\cdots(BX'_1)^{\varepsilon}\cdots X'_{2m}B$. In order to prove that $W$ is equivalent to $W'$, it suffices to show that for every $1\leq i\leq 2m$, $X_i$ and $X'_i$ have at least $4c+3(A_0+H)+D_0$ positions in common. We prove this by induction on $i$. Let $X'_j$ be the letter in $W'_{\text{split}}$ preceding the sub-word $(BX'_1)^{\varepsilon}$, that is $j=\overline{1}'-1$ and:
    \[W'_{\text{split}}=\begin{cases}
        X'_1\cdots X'_jBX'_{\overline{1}'}X'_{j+2}\cdots X'_{2m}B&\text{; when }\varepsilon=+1\\
        X'_1\cdots X'_jX'_{\overline{1}'}B^{-1}X'_{j+2}\cdots X'_{2m}B&\text{; when }\varepsilon=-1
    \end{cases}\]
    \textbf{Base case: $\mathbf{i=1}$} $|s_1-s'_1|=|b_0-b'_0|\leq H$. Since $l_{X_1},l'_{X'_1}>(4c+3(A_0+H)+2D_0)+H$, $X_1$ and $X'_1$ have at least $4c+3(A_0+H)+2D_0$ positions in common.\\\\\textbf{Induction step:} Let $i\geq 1$. Suppose that for all $1\leq k\leq i$, $X_k$ and $X'_k$ have at least $4c+3(A_0+H)+D_0$ positions in common. We will show that $X_{i+1}$ and $X'_{i+1}$ have at least $4c+3(A_0+H)+D_0$ positions in common.\\\emph{Case - 1:} Suppose $i<j$ or $j+1<i$. By induction hypothesis, $X'_i$ and $X_i$ have at least $4c+3(A_0+H)+D_0$ positions in common. If $X'_i$ and $X_{i+1}$ have any positions in common, they are at most $4c+3(A_0+H)+D_0-1$ in number. But $l_{X_{i+1}}>8c+7(A_0+H)+2D_0$, so $X_{i+1}$ and $X'_{i+1}$ must have at least $4c+3(A_0+H)+D_0$ positions in common. Similarly, if $X_i$ and $X'_{i+1}$ share any positions in common, we conclude that $X_{i+1}$ and $X'_{i+1}$ must have at least $4c+3(A_0+H)+D_0$ positions in common, since $l'_{X'_{i+1}}>8c+7(A_0+H)+2D_0$. Otherwise, the junction between $X_i$ and $X_{i+1}$ must contain the junction between $X'_i$ and $X'_{i+1}$, and we infer the same.\\\emph{Case - 2:} Suppose $i=j$, that is, $i=\overline{1}'-1$. By \Cref{lem: common positions conj class} and the base case, $X'_{\overline{1}'}$ and $X_{\overline{1}}$ have at least $4c+3(A_0+H)+D_0$ positions in common. Therefore it suffices to show that $\overline{1}=j+1=\overline{1}'$.  We need to argue separately for different values $\varepsilon$.
    \begin{itemize}
        \item \underline{$\varepsilon=-1$}: By \Cref{lem: reduction or redundancy conj class}, the letter $X_{\overline{1}-1}$ has at most $4c+3(A_0+H)+D_0-1$ positions in common with $X'_{\overline{1}'}$. Since the size of $B'_{j}$ is at most $A_0+H$ and $l_{X_{\overline{1}-1}}>8c+7(A_0+H)+2D_0$, it follows that $X_{\overline{1}-1}$ has at least $4c+3(A_0+H)+D_0$ positions in common with $X'_j$. But $X_j$ also shares at least $4c+3(A_0+H)+D_0$ positions in common with $X'_j$. So it follows from \Cref{lem: reduction or redundancy conj class} that $j=\overline{1}-1$.

        \item \underline{$\varepsilon=+1$}: $|s_{X_1}-s_{X'_1}|\leq H$, and $X_1$ and $X'_1$ have at least $4c+3(A_0+H)+D_0$ positions in common, so it follows from \Cref{lem: matching position extended conj class} that $|s_{X_{\overline{1}}}-s_{X_{\overline{1}'}}|\leq H+D_0-1$.
        This means that the letter $X_{\overline{1}-1}$ has at most $H+D_0-1$ many positions in common with $X'_{\overline{1}'}$. Because the sizes of $B'_j$ and $B'_{j+1}$ are at most $A_0+H$, $l'_B\leq 4c+4(A_0+H)+2D_0$, and $l_{X_{\overline{1}-1}}>8c+10(A_0+H)+4D_0$, it follows that $X_{\overline{1}-1}$ has at least $4c+3(A_0+H)+D_0$ positions in common with $X'_j$. But $X_j$ also shares at least $4c+3(A_0+H)+D_0$ positions in common with $X'_j$. So it follows from \Cref{lem: reduction or redundancy conj class} that $j=\overline{1}-1$.
    \end{itemize}
    \emph{Case - 3:} Suppose $i=j+1$. 
    \begin{itemize}
        \item \underline{$\varepsilon=-1$}: $|s_{X_1}-s_{X'_1}|\leq H$, and $X_1$ and $X'_1$ have at least $4c+3(A_0+H)+D_0$ positions in common, so it follows from \Cref{lem: matching position extended conj class} that $|t_{X_{\overline{1}}}-t_{X_{\overline{1}'}}|\leq H+D_0-1$.
        This means that the letter $X_{j+2}$ has at most $H+D_0-1$ positions in common with $X'_{j+1}(=X'_{\overline{1}'})$. Because the sizes of $B'_{j+1}$ and $B'_{j+2}$ are at most $A_0+H$, $l'_B\leq 4c+4(A_0+H)+2D_0$, and $l_{X_{j+2}}>8c+10(A_0+H)+4D_0$, it follows that $X_{j+2}$ has at least $4c+3(A_0+H)+D_0$ positions in common with $X'_{j+2}$. 
        
        \item \underline{$\varepsilon=+1$}: The same arguments made in Case - $1$ work.
        \end{itemize}
\end{proof}

\begin{remark}
    The proof of \Cref{prop: reduction or equality of Wicks forms conj class} in fact tells us that if there is no $D_0$-reduction in degrees of freedom by $c$ between $(W,\vec{l},\vec{b})$ and $(W'_{\text{split}},\vec{l}',\vec{b}')$, then $B$ is at the head of $X_{\overline{1}+1}$ when $\varepsilon=-1$, and $B$ is at the tail of $X_{\overline{1}-1}$ when $\varepsilon=+1$. 
\end{remark}

The proof of \Cref{prop: another main reduction theorem} will proceed by using shift equations (see \Cref{prop: reduction or shifted matching}) just like we did in \Cref{subsubsec: Proof of the reduction theorem}. Because there is no lower bound on $l'_B$, it essentially behaves like a dummy junction. We do this formally in \Cref{def: shift conj class} below. The shift equations will imply that if there is no $D_0$-reduction by $c$, the size of such a dummy junction has to be small.

\begin{definition}\label{def: shift conj class}
    Given $(W,\vec{l},\vec{b})$ and $(W'_{\text{split}},\vec{l}',\vec{b}')$, for every letter $X_i$ of $W$, define $\text{shift}(i)\coloneq t_{X_i}-t'_{X'_i}$. Recall that $W'_{\text{split}}$ has the following form:
    \[W'_{\text{split}}=\begin{cases}
        X'_1\cdots X'_jBX'_{\overline{1}}X'_{j+2}\cdots X'_{2m}B&\text{; when }\varepsilon=+1\\
        X'_1\cdots X'_jX'_{\overline{1}}B^{-1}X'_{j+2}\cdots X'_{2m}B&\text{; when }\varepsilon=-1
    \end{cases}\]
    In order to treat the letter $B^\varepsilon$ as a junction, we make the following definition. When $\varepsilon=+1$, define $b''_{\overline{1}-1}\coloneq b'_{\overline{1}-1}+l'_B+b'_{\overline{1}}$. When $\varepsilon=-1$, define $b''_{\overline{1}}\coloneq b'_{\overline{1}}+l'_B+b'_{\overline{1}+1}$.\\ Next, define $\Delta(i)$ as the difference in sizes of the junctions immediately following the letters $X_i$ and $X'_i$. More precisely, for all $0\leq i\leq 2m$:
    \[ \text{For }\varepsilon=-1, \Delta(i)\coloneq \begin{cases}
        b_i-b'_i&\text{; if } i< \overline{1}\\
        b_{\overline{1}}-b''_{\overline{1}}&\text{; if } i= \overline{1}\\
        b_i-b'_{i+1}&\text{; if } i> \overline{1}
    \end{cases}\;\; \text{For }\varepsilon=+1, \Delta(i)\coloneq \begin{cases}
        b_i-b'_i&\text{; if } i< \overline{1}-1\\
        b_{\overline{1}-1}-b''_{\overline{1}-1}&\text{; if } i= \overline{1}-1\\
        b_i-b'_{i+1}&\text{; if } i> \overline{1}-1
    \end{cases} \]
\end{definition}
 We make a simplification in notation in view of \Cref{prop: reduction or equality of Wicks forms conj class}. For $(W,\vec{l},\vec{b})$ and $(W'_{\text{split}},\vec{l}',\vec{b}')$, $t_j\coloneq t_{X_j}$ and $t'_j\coloneq t'_{X'_j}$. Also note that $t_0=t'_0\coloneq 0$, so $\text{shift}(0)=0$.

\begin{lemma}\label{lem: shift equations conj class}
    Suppose there is no reduction in degrees of freedom by $c$ between $(W,\vec{l},\vec{b})$ and $(W'_{\text{split}},\vec{l}',\vec{b}')$. Then the following equations hold for all $1\leq j\leq 2m$:
    \begin{enumerate}[font=\normalfont, label=(\arabic*)]
            \item If \(\text{sign}(X_j) = -1\):
            \begin{enumerate}[font=\normalfont,label=(1\alph*)]
                \item\label{item: shift equation 1a conj}  If $j-1\neq\overline{1}$, then $\text{shift}(j-1)\simeq_{D_0+H-1} - \text{shift}(\overline{j})$.
                \item\label{item: shift equation 1b conj}  If $\overline{j}-1\neq\overline{1}$, then  $\text{shift}(j)\simeq_{D_0+H-1} -\text{shift}(\overline{j} - 1)$.
                \item\label{item: shift equation 1c conj} If $j-1=\overline{1}$, then $\text{shift}(j-1)\simeq_{D_0+H-1}-\text{shift}(\overline{j})+l'_B+b'_{\overline{1}+1}$, and if $\overline{j}-1=\overline{1}$, then  $\text{shift}(\overline{j}-1)\simeq_{D_0+H-1}-\text{shift}(j)+l'_B+b'_{\overline{1}+1}$.
            \end{enumerate}
            
            \item If \(\text{sign}(X_j) = 1\):
            \begin{enumerate}[font=\normalfont,label=(2\alph*)]
                \item\label{item: shift equation 2a conj} If $j\neq 1,\overline{1}$, then $\text{shift}(j-1) \simeq_{D_0+2H-1} \text{shift}(\overline{j}-1)$.
                \item\label{item: shift equation 2b conj} $\text{shift}(j) \simeq_{D_0-1} \text{shift}(\overline{j})$.
                \item\label{item: shift equation 2c conj} If $j=\overline{1}$, then $\text{shift}(j-1)\simeq_{D_0+2H-1}l'_B+b'_{\overline{1}}$, and if $j=1$, then $\text{shift}(\overline{j}-1)\simeq_{D_0+2H-1}l'_B+b'_{\overline{1}}$.
            \end{enumerate}

            \item\label{item: shift equation 3 conj} $\text{shift}(2m)\simeq_{H}l'_B+b'_{2m+2}$
    \end{enumerate}
\end{lemma}
\begin{proof}
    $(W,\vec{l},\vec{b})$ and $(W'_{\text{split}},\vec{l}',\vec{b}')$ satisfy $(\mathbf{L'1})$. So, we have:
    \begin{align*}
        t_{2m}+b_{2m}=t'_{2m}+b'_{2m+1}+l'_B+b'_{2m+2} 
    \end{align*}
    Therefore, $\text{shift}(2m)=l'_B-\Delta(2m)+b'_{2m+2}$. Recalling that $W'_{\text{split}}=X'_1\cdots (BX'_1)^{\varepsilon}\cdots X'_{2m}B$, we note that when $\varepsilon=+1$, $\overline{1}\leq 2m$, and when $\varepsilon=-1$, $\overline{1}<2m$. In any case, we have $|\Delta(2m)|=|b_{2m}-b'_{2m+1}|\leq H$. Thus, \ref{item: shift equation 3 conj} follows. 
    
    By \Cref{prop: reduction or equality of Wicks forms conj class}, for all $1\leq j\leq 2m$, $X_j$ and $X'_j$ have at least $4c+3(A_0+H)+D_0$ positions in common. Consider a position $i$, which is contained in $X_j$ with respect to $\vec{l}$ and also contained in $X'_j$ with respect to $\vec{l}'$. By \Cref{lem: common positions conj class}, we have $\underline{i}\simeq_{D_0-1}\underline{i}'$. Moreover, using \Cref{def: shift conj class}, we have the following equation for all $1\leq j\leq 2m$:
    \[\quad  s_j-s'_j=\text{shift}(j-1)+\Delta(j-1)\]
    Together with \Cref{lem: matching position extended}, the same arguments used in \Cref{cor: shift equations}, give:
    \begin{itemize}
        \item If $\text{sign}(X_j)=-1$, then $\text{shift}(\overline{j})\simeq_{D_0-1}-(\text{shift}(j-1)+\Delta(j-1)$).
        \item If $\text{sign}(X_j)=+1$, then $\text{shift}(j-1)+\Delta(j-1)\simeq_{D_0-1}\text{shift}(\overline{j}-1)+\Delta(\overline{j}-1)$, and $\text{shift}(j)\simeq_{D_0-1}\text{shift}(\overline{j})$.
    \end{itemize}
    So, the rest of the parts follow from the above in conjunction with the following observation: Because $\vec{b}\in [A_0,A_0+H]^{2m}$ and $\vec{b}'\in[A_0,A_0+H]^{2m+3}$, we have:
    \begin{itemize}
        \item When $\varepsilon=-1$, and $i\neq \overline{1}$, then $|\Delta(i)|\leq H$. $|\Delta(\overline{1})+l'_B+b'_{\overline{1}+1}|=|b_{\overline{1}}-b'_{\overline{1}}|\leq H$.
        \item  When $\varepsilon=+1$, and $i\neq\overline{1}-1$, then $|\Delta(i)|\leq H$. $|\Delta(\overline{1}-1)+l'_B+b'_{\overline{1}}|=|b_{\overline{1}-1}-b'_{\overline{1}-1}|\leq H$.
    \end{itemize}
\end{proof}

\begin{proof}[Proof of \Cref{prop: another main reduction theorem}] 
    By \Cref{prop: reduction or equality of Wicks forms conj class}, we already know that $W=W'$. It remains to prove that $l'_B\leq 2(D_0-1+H)$.
    Recall the definition of the graph $\Gamma_W$ associated with $W$ (see \Cref{def: wicks form genus}). As $W$ is a maximal Wicks form, $\Gamma_W$ is a trivalent graph. The cycle $\partial P_W$ starting at the edge labeled $X_1$ (corresponding to reading the letters of $W$ from left to right), projects to a cycle, $\alpha$ in $\Gamma_W$. We think of $\alpha$ as a vertex path $(v_0,\dots,v_{l(W)})$.
    The proof now splits into several cases depending on the  $1$-neighborhood of $v_0$, which can be exactly one among five types, as depicted in \Cref{fig: circuit in 1-neighbourhood (conjugacy)}.

    \begin{figure}[h]
            \centering
            \includegraphics[width=0.85\columnwidth]{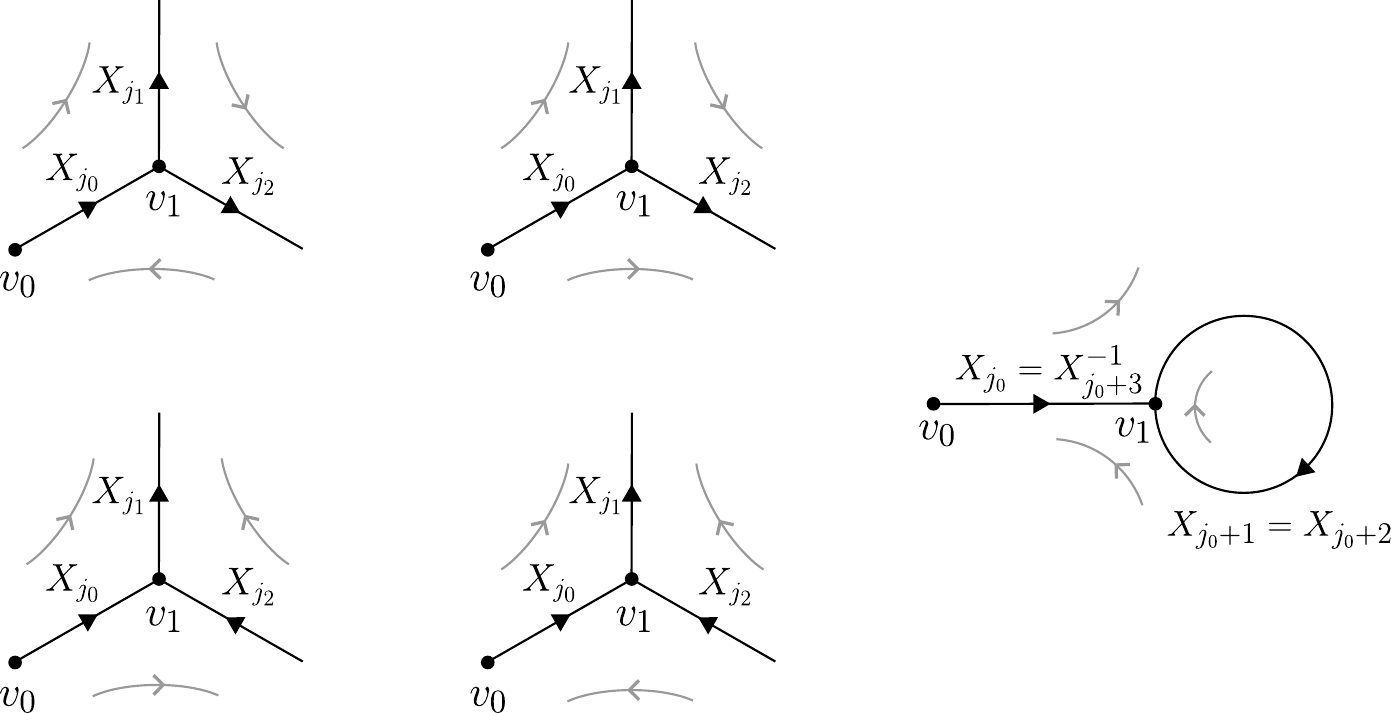}
            \caption{The possible ways of the circuit corresponding to $W$ intersecting the 1-neighbourhood of $v_1$. This figure is nearly identical to \Cref{fig: circuit in 1-neighbourhood}, and is reproduced here for reference.}
            \label{fig: circuit in 1-neighbourhood (conjugacy)}
        \end{figure}
      We define indices \(j_0, j_1, j_2\) as follows:  Set \(j_0 = 1\). When the $1$-neighborhood of $v_0$ in $\Gamma_W$ is a tripod, we define $j_1$ and $j_2$ as follows: \(X_{j_1}\) is the unique letter adjacent to \(X_{\overline{j_0}} = X_{\overline{1}}\) such that between the corresponding edges in \(\Gamma_W\), \(v_0\) is the common vertex. Similarly, \(X_{j_2}\) and \(X_{\overline{j_1}}\) are adjacent letters such that the corresponding edges share the vertex \(v_0\) in \(\Gamma_W\). When the $1$-neighborhood of $v_0$ has a loop, $X_{j_1}$ is the unique letter adjacent to $X_{\overline{j_0}}$ (other than possibly $X_{j_0}$, which could also be adjacent to $X_{\overline{j_0}}$) such that for the corresponding edges in $\Gamma_W$, $v_0$ is the common vertex. Using this notation, the following is an algebraic description of the figures in \Cref{fig: circuit in 1-neighbourhood (conjugacy)}:
    \begin{enumerate}[label=(\alph*)]
        \item \underline{Top left, without loop:} $\text{sign}(X_{j_0})=\text{sign}(X_{j_1})=\text{sign}(X_{j_2})=-1$ and $W$ is of the form: 
        \[X_{j_0}\cdots X_{\overline{j_0}}X_{j_1}\cdots X_{\overline{j_1}}X_{j_2}\cdots X_{\overline{j_2}}\]
        \item \underline{Top right, without loop:} $\text{sign}(X_{j_0})=\text{sign}(X_{j_2})=+1$ $\text{sign}(X_{j_1})=-1$, and $W$ is of the form: 
        \[X_{j_0}\cdots X_{j_1}X_{\overline{j_0}}\cdots X_{j_2}X_{\overline{j_1}}\cdots X_{\overline{j_2}}\]
        \item \underline{Bottom left, without loop:} $\text{sign}(X_{j_0})=\text{sign}(X_{j_1})=+1$, $\text{sign}(X_{j_2})=-1$ and $W$ is of the form: 
        \[X_{j_0}\cdots X_{j_1}X_{\overline{j_0}}\cdots X_{\overline{j_1}}X_{j_2}\cdots X_{\overline{j_2}}\]
        \item \underline{Bottom right, without loop:} $\text{sign}(X_{j_0})=-1$, $\text{sign}(X_{j_1})=\text{sign}(X_{j_2})=+1$ and $W$ is of the form: 
        \[X_{j_0}\cdots X_{\overline{j_0}}X_{j_1}\cdots X_{j_2}X_{\overline{j_1}}\cdots X_{\overline{j_2}}\]
        \item \underline{with loop:} This splits into three sub-cases:
        \begin{enumerate}[label=(\roman*)]
            \item $\text{sign}(X_{j_0})=+1$, $\text{sign}(X_{j_1})=-1$, and $W$ is of the form:
            \[X_{j_0}X_{\overline{j_0}}X_{j_1}\cdots X_{\overline{j_1}}\]
            \item $\text{sign}(X_{j_0})=+1$, $\text{sign}(X_{j_1})=-1$, and $W$ is of the form:
            \[X_{j_0}X_{\overline{j_1}}\cdots X_{j_1} X_{\overline{j_0}}\]
            \item $\text{sign}(X_{j_0})=-1$, $\text{sign}(X_{j_1})=+1$, and $W$ is of the form:
            \[X_{j_0}\cdots X_{\overline{j_0}} X_{j_1} X_{\overline{j_1}} \]
        \end{enumerate}
    \end{enumerate}
    Note that $(a)$ is the only relevant neighborhood type when $W$ is orientable. But when $W$ is non-orientable, all five types need to be taken account for. The arguments in all the cases follow the same template to obtain the desired upper bound for $l'_B$: We get two expressions $\text{shift}(2m)$, one from part \ref{item: shift equation 3 conj} of \Cref{lem: shift equations conj class}, and another using appropriate shift equations. $l'_B$ appears in both expressions, but with different signs, thus leading to the upper bound. The arguments are essentially the same as those used in the proof of \Cref{prop: no reduction implies length types are equal}, so we only demonstrate two cases:
    \begin{enumerate}[label=(\alph*)]
        \item \underline{Top left, without loop:} In this case, we have $\text{sign}(X_{j_0})=\text{sign}(X_{j_1})=\text{sign}(X_{j_2})=-1$ and $j_0=1$, $j_1=\overline{j_0}+1$, $j_2=\overline{j_1}+1$, and $\overline{j_2}=2m$.
        \begin{align*}
            0=\text{shift}(j_0-1)&\simeq_{D_0+H-1}-\text{shift}(\overline{j_0}) &(\text{using }\ref{item: shift equation 1a conj})\\
            &= -\text{shift}(j_1-1)  & \\
            &\simeq_{D_0+H-1} -(-\text{shift}(\overline{j_1})+l'_B+b'_{\overline{1}+1})  &(\text{using }\ref{item: shift equation 1c conj})\\
            &= \text{shift}(j_2-1)-l'_B-b'_{\overline{1}+1}  &\\
            &\simeq_{D_0+H-1}-\text{shift}(\overline{j_2})-l'_B-b'_{\overline{1}+1}  &(\text{using }\ref{item: shift equation 1a conj})\\
            \therefore \text{shift}(\overline{j_2})&\simeq_{3(D_0+H-1)}-l'_B-b'_{\overline{1}+1} &
        \end{align*}
        Using \ref{item: shift equation 3 conj}, we have: $l'_B+b'_{2m+2}\simeq_{3(D_0-1)+4H}-l'_B-b'_{\overline{1}+1}$. Therefore, since $b'_{\overline{1}+1}$, $b'_{2m+2}$ are non-negative integers, it follows that $l'_B\leq 2(D_0+H-1)$.
        \end{enumerate}
        \begin{enumerate}[label=(\alph*), start=3]
        \item \underline{Bottom left, without loop:} In this case, we have $\text{sign}(X_{j_0})=\text{sign}(X_{j_1})=+1$, and $\text{sign}(X_{j_2})=-1$ and $j_0=1$, $j_1=\overline{j_0}-1$, $j_2=\overline{j_1}+1$, and $\overline{j_2}=2m$.
        \begin{align*}
            l'_B+b'_{\overline{1}}&\simeq_{D_0+2H-1}\text{shift}(\overline{j_0}-1) &(\text{using }\ref{item: shift equation 2c conj})\\
            &= \text{shift}(j_1) & \\
            &\simeq_{D_0-1} \text{shift}(\overline{j_1}) &(\text{using }\ref{item: shift equation 2b conj})\\
            &= \text{shift}(j_2-1) & \\
            &\simeq_{D_0+H-1} -\text{shift}(\overline{j_2})&(\text{using }\ref{item: shift equation 1a conj})\\
            \therefore \text{shift}(\overline{j_2})&\simeq_{3(D_0+H-1)} -l'_B-b'_{\overline{1}} & 
        \end{align*}
        Comparing with \ref{item: shift equation 3 conj} we again get $l'_B\leq 2(D_0+H-1)$.
    \end{enumerate} 
    Furthermore, using \ref{item: shift equation 3 conj} and the bound on $l'_B$, we get:
    \begin{equation}\label{eq: shift 2m bound}
        |\text{shift}(2m)|\leq 2(D_0-1 +H)+A_0+2H
    \end{equation}
    \end{proof}

\subsection{Counting in Free groups}\label{subsec: Conjugacy class count: Lower bounds: Free groups}

We continue to use the notation and definitions made in \Cref{subsec: Conjugacy class count: Upper bounds: Free groups}. Moreover, recall that due to \Cref{lem: canonical forms for conjugacy classes}, the set $\partial \mathcal{C}^{\pm}_n(L)$ is empty when $L$ is odd. Therefore in the rest of this subsection, the parameter $L$ is always assumed to be even.

\begin{definition}\label{def: set of degenerate group elements free groups}
    Let $V$ be a linear Wicks form and let $\vec{l}$ be a length type for $V$. Let $N\geq 1$ be a positive integer. Define $G^{\text{deg}}_{V}(N,\vec{l})$ to be the set of group elements coming from cyclically non-canceling substitutions $\mathfrak{g}$ of length type $\vec{l}$ which satisfy one of the following:
\begin{enumerate}[font=\normalfont]
        \item  There are two distinct unknowns $X$ and $Y$ of $W$ such that:
        \begin{equation*}
            g_X=uzv,\quad g_{Y}=u'z^{\varepsilon}v'
        \end{equation*}
        where $\varepsilon\in\{+1,-1\}$, and:
        \begin{equation*}
            l(z)=N,\quad l(u)+l(v)+N=l_X,\quad l(u')+l(v')+N=l_{Y}
        \end{equation*}
        \item  There is an unknown $X$ of $W$ such that:
        \begin{equation*}
            g_X=uzwz^{\varepsilon}v
        \end{equation*}
        where $\varepsilon\in\{+1,-1\}$, and:
        \begin{equation*}
            l(z)=N,\quad l(u)+l(w)+l(v)+2N=l_X
        \end{equation*}
        \item There is an unknown $X$ of $W$, and a positive integer $T\geq 1$ such that:
        \begin{equation*}
            g_X=uz_0^kz_1v
        \end{equation*}
        where $l(z_0)=T$, $k=\left\lfloor\frac{N}{T}\right\rfloor+1$, $l(z_1)<T$, and $z_{1}$ is a group element of length $N-(k-1)T$ which appears as the prefix in the reduced word representing $z_0$. Furthermore:
        \begin{equation*}
            l(u)+l(v)+l(z_1)+kT=l_X,\quad (k-1)T+l(z_1)=N
        \end{equation*}
    \end{enumerate}
\end{definition}

\begin{lemma}\label{lem: degenerate count free group upper bound}
    Let $V$ be a linear Wicks form of length $2m$ and let $\vec{l}$ be a length type for $V$. Then there is a constant $c'_1$ depending only on $r$ and $m$ such that:
    \begin{align*}
        |G^{\text{deg}}_{V}(N,\vec{l})|\leq c'_1L^2e^{h_r\left(\frac{||\vec{l}||}{2}-N \right)}
    \end{align*}
\end{lemma}
\begin{proof}
    This follows readily from \Cref{lem: countingpinnedreducedwords}, and the fact that in each of the algebraic forms in \Cref{def: set of degenerate group elements free groups}, there is a reduction in choosing substitutions by a length of $N$. 
\end{proof}

\begin{lemma}\label{lem: key lemma for conj class lower bounds}
    Let $N\geq 1$ be a positive integer. Let $V$ and $V'$ be linear maximal Wicks forms. Let $\vec{l}$ and $\vec{l'}$ be length types for $V$ and $V'$ respectively, satisfying:
    \begin{itemize}
        \item $||\vec{l}||=||\vec{l}'||$
        \item For all $1\leq i\leq l(V)$, $1\leq j\leq l(V')$, $l_i,l'_j> 10N+25$.
    \end{itemize}
    
    Let $\mathfrak{g}$ and $\mathfrak{g}'$ be cyclically non-canceling substitutions of length type $\vec{l}$ and $\vec{l}'$ for $V$ and $V'$ respectively, such that $[\mathfrak{g}(V))]=[\mathfrak{g}'(V')]$. Then either:
    \begin{enumerate}[font=\normalfont]
        \item $\mathfrak{g}(V)\in G^{\text{deg}}_{V}(N,\vec{l})$, or
        \item There is a cyclic permutation $\sigma$ such that $V$ is linearly equivalent to $\sigma(V')$, $\vec{l}=\sigma\cdot \vec{l}'$, and $\mathfrak{g}=\sigma\cdot\mathfrak{g}'$.
    \end{enumerate}
\end{lemma}

\begin{proof}
    Since $\mathfrak{g}$ and $\mathfrak{g}'$ are cyclically non-canceling substitutions, $\mathfrak{g}(V)$ and $\mathfrak{g}'(V')$ are cyclically reduced words. Therefore, the fact that $[\mathfrak{g}(V))]=[\mathfrak{g}'(V')]$ implies that $\mathfrak{g}'(V')$ and $\mathfrak{g}(V)$ are cyclic permutations of each other. This can happen in the following possible ways:
    \begin{enumerate}
        \item There is a cyclic permutation $\sigma$ of $V'$, so that $\mathfrak{g}(V)=\mathfrak{g}'(\sigma(V'))$. 

        In this case, consider the data $(V,\vec{l})$ and $(\sigma(V'),\sigma(\vec{l}'))$. They satisfy conditions \ref{item: L1 free} and \ref{item: L2 free} of \Cref{cor: the main reduction theorem for free groups}. So, either there is reduction in degrees of freedom by $N$ or $V=\sigma(V')$, $\vec{l}=\sigma\cdot \vec{l}'$. Recalling \Cref{rmk: Relating combinatorial reduction to degeneracy in substitutions}, let us note that in the former case, $\mathfrak{g}$ (and $\mathfrak{g}'$) has one of the degenerate algebraic forms listed in \Cref{def: set of degenerate group elements free groups} (see proof of \Cref{lem: element intersection count for reduction in degrees of freedom}), and so $\mathfrak{g}(V)\in G^{\text{deg}}_{V}(N,\vec{l}) $. While in the latter case, $\mathfrak{g}(V)=\mathfrak{g}'(\sigma(V'))$ together with $V=\sigma(V')$, $\vec{l}=\sigma(\vec{l}')$ implies that $\mathfrak{g}=\sigma\cdot\mathfrak{g}'$.

        \item There is a cyclic permutation $\sigma$ of $V'$, and a splitting of its first letter, so that $\mathfrak{g}(V)=\mathfrak{g}'_{\text{split}}(V'_{\text{split}})$. More precisely, we write $\sigma(V')=X'U'_1(X')^{\varepsilon}U'_2$, split $X'=B'X'_1$, and set $V'_{\text{split}}=X'_1U'_1(BX'_1)^{\varepsilon}U'_2B$. This is a split Wicks form (see \Cref{def: split wicks forms}) and it comes with a substitution which is the same as $\mathfrak{g}'$ for all letters except $X'$, which has been split. A splitting of the reduced word $g'_{X'}=g'_Bg'_{X'_1}$ (where $g'_B$ and $g'_{X'_1}$ are non-trivial reduced words), determines the new substitution. We denote the new substitution by $\mathfrak{g}'_{\text{split}}$. 

        Now consider the data $(V,\vec{l})$ and $(V'_{\text{split}},\mathfrak{g}'_{\text{split}})$. They satisfy conditions \ref{item: L'1 free} and \ref{item: L'3 free} of \Cref{cor: reduction or no splitting - free groups}. So, either there is reduction in degrees of freedom by $N$ or $V=\sigma(V')$, $\vec{l}=\sigma\cdot \vec{l}'$ and either of $l'_B$ or $l'_{X'_1}$ equals zero. The latter is not possible since $g'_B$ and $g'_{X'_1}$ are non-trivial reduced words by assumption. Therefore, there is reduction in degrees of freedom by $N$, and by \Cref{rmk: Relating combinatorial reduction to degeneracy in substitutions}, $\mathfrak{g}$ (and $\mathfrak{g}'$) has one of the degenerate algebraic forms listed in \Cref{def: set of degenerate group elements free groups}, and so $\mathfrak{g}(V)\in G^{\text{deg}}_{V}(N,\vec{l})$.
        
    \end{enumerate}
\end{proof}

\Cref{lem: key lemma for conj class lower bounds} will allow us to derive the key estimates required to prove the conjugacy class lower bound estimate. 

\begin{lemma}\label{lem: conj class substitution count asymptote part 2}
    Let $N\geq 1$ be a positive integer. Let $W$ be a maximal Wicks form of length $2m$ and let $\vec{l}$ be a non-degenerate length type for $W$ such that each entry of $\vec{l}$ is greater than $10N+25$, where $L=||\vec{l}||$. Let $c'_1$ be the constant appearing in \Cref{lem: degenerate count free group upper bound}. Then:
    \begin{align*}
        \frac{|\mathcal{G}_W(\vec{l})|}{\displaystyle{\kappa_W \cdot\frac{(2r - 1)^{\frac{||\vec{l}||}{2} - m}}{\displaystyle (2r)^m}}}\geq \begin{cases}
            \left(1-\frac{c'_1L^2}{e^{h_rN}}\right) &; W\neq ABCA^{-1}B^{-1}C^{-1}\\
            \frac{1}{2}\left(1-\frac{c'_1L^2}{e^{h_rN}}\right) &; W= ABCA^{-1}B^{-1}C^{-1}
        \end{cases}
    \end{align*}
\end{lemma}

\begin{proof}
    Let $V=X_1\cdots X_{2m}$ be a linear word representing $W$. Consider the map $\psi:G_V(\vec{l}_V)-G^{\text{deg}}_V(\vec{l}_V)\rightarrow \mathcal{G}_W(\vec{l})$ which maps group elements to their conjugacy classes. 
    \begin{claim}\label{claim: finite-to-one-ness of psi}
        If $V\neq ABCA^{-1}B^{-1}C^{-1}$, $\psi$ is injective. If $V= ABCA^{-1}B^{-1}C^{-1}$, $\psi$ is $2$-to-$1$.
    \end{claim}
    \begin{proof}
        Let $\sigma$ be a generator of $\text{Aut}(V)$ and denote its order by $n$. Let $\mathfrak{g}$ and $\mathfrak{g}'$ be two cyclically non-canceling substitutions such that $\mathfrak{g}(V),\mathfrak{g}'(V)\in G_V(\vec{l}_V)-G^{\text{deg}}_V(\vec{l}_V)$, and $[\mathfrak{g}(V)]=[\mathfrak{g}'(V)]$. The data $(V,\mathfrak{g},\vec{l}_V)$ and $(V,\mathfrak{g}',\vec{l}_V)$ satisfies conditions \ref{item: L'1 free} and \ref{item: L'3 free} of \Cref{cor: reduction or no splitting - free groups}, and $\mathfrak{g}(V)\notin G^{\text{deg}}_V(\vec{l}_V)$ , so  it follows that there is a cyclic permutation $\alpha$ of $V$ such that $V$ is linearly equivalent to $\alpha(V)$, that is, $\alpha=\sigma^j$ for some $0\leq j\leq n-1$. Moreover, $\vec{l}_V=\sigma^j\cdot\vec{l}_V$, and $\mathfrak{g}=\sigma^j\cdot\mathfrak{g}'$. Since $\vec{l}$ is non-degenerate, it follows from \Cref{lemma: action of automorphisms on a length type} that when $V\neq ABCA^{-1}B^{-1}C^{-1}$, $\sigma^j$ has to be the identity, which proves that $\psi$ is injective. When $V= ABCA^{-1}B^{-1}C^{-1}$, either $j=0$ or $j=3$. Note that $\mathfrak{g}'(V)\neq \sigma^{3}\cdot\mathfrak{g}'(V)$. Therefore, in this case $\psi$ is $2$-to-$1$.  
    \end{proof}
    By \Cref{lem: degenerate count free group upper bound} and \Cref{lem: substitution count asymptote conj class}, we have:
    \begin{align*}
        G_V(\vec{l}_V)-G^{\text{deg}}_V(\vec{l}_V)\geq \left(1-\frac{c'_1L^2}{e^{h_rN}}\right)G_V(\vec{l}_V)
    \end{align*}
    Together with \Cref{claim: finite-to-one-ness of psi}, the lemma follows.
\end{proof}
\begin{lemma}\label{lem: key lemma for free group conj class lower bound counting argument}
    Let $W$ be a maximal orientable Wicks form of genus $n$. Then, for all $L\geq 1$, there exists a collection of length types $\mathcal{L}=\mathcal{L}(W,n,r,L)$ and there are constants $L_0,c'_2\geq 1$ depending only on $r$ and $n$ satisfying, for all $L\geq L_0$:
    \begin{enumerate}[font=\normalfont]
        \item For all $\vec{l}\in \mathcal{L}(W,n,r,L)$, $||\vec{l}||=L$, $\mathcal{L}(W,n,r,L)$ is $\text{Aut}(W)$-invariant,  and \[\lim_{L\to\infty}\frac{|\mathcal{L}(W,n,r,L)|}{L^{6n-4}}=\frac{1}{2^{6n-4}(6n-4)!}
        \]
    \end{enumerate}
    Let $W'$ be another maximal orientable Wicks form of genus $n$, and let $[\vec{l}]\in \mathcal{L}(W,n,r,L)/\text{Aut}(W)$, $[\vec{l}']\in \mathcal{L}(W',n,r,L)/\text{Aut}(W')$.
    \begin{enumerate}[font=\normalfont, resume]
        \item If $(W,[\vec{l}])\neq(W',[\vec{l}'])$, then:
        \begin{align*}
            |\mathcal{G}_W(\vec{l})\cap \mathcal{G}_{W'}(\vec{l}')|\leq c'_2\frac{e^{\frac{h_rL}{2}}}{L^{6n}}
        \end{align*}
    \end{enumerate}
\end{lemma}
\begin{proof}
    Suppose $W$ is a maximal orientable Wicks form. Then $W$ has length equal to $12n-6$. Consider the following collection of length types for $W$:
    \begin{align*}
         &\mathcal{L}(W,n,r,L)\\
         &\coloneq\left\{\vec{l}\in\mathbb{Z}^{6n-3}_+\middle|\quad  \sum_{i=1}^{6n-3}l_i = \frac{L}{2};1\leq i\leq 6n-3, l_i> 10\left \lfloor \frac{(6n+2)\log L}{h_r} \right \rfloor+25,\vec{l}\text{ non-degenerate}\right\}
     \end{align*}
    $\mathcal{L}$ evidently satisfies $(1)$, because asymptotically, most length types satisfying $\sum l_i=\frac{L}{2}$ are non-degenerate.
    Note that every pair of length types in $\mathcal{L}$ satisfies the hypotheses of \Cref{lem: key lemma for conj class lower bounds}, with $N=\left \lfloor \frac{(6n+2)\log L}{h_r}\right \rfloor$.
    
    Next, let $W'$ is another Wicks form, and $\vec{l}\in \mathcal{L}(W,n,r,L)$, $\vec{l}'\in \mathcal{L}(W',n,r,L)$. Let $V$ and $V'$ be linear Wicks forms representing $W$ and $W'$ respectively. Let $\mathfrak{g}$ and $\mathfrak{g}'$ be two cyclically non-canceling substitutions such that $\mathfrak{g}(V)\in G_V(\vec{l}_V)$, $\mathfrak{g}'(V')\in G_V(\vec{l'}_{V'})$, and $[\mathfrak{g}(V)]=[\mathfrak{g}'(V')]$. The data $(V,\mathfrak{g},\vec{l})$ and $(V',\mathfrak{g}',\vec{l}')$ satisfies the hypotheses of \Cref{lem: key lemma for conj class lower bounds}. If $(W,[\vec{l}])\neq(W',[\vec{l}'])$, $\mathfrak{g}(V)\in G^{\text{deg}}_{V}(N,\vec{l}_V)$, where $N=\left \lfloor \frac{(6n+2)\log L}{h_r}\right \rfloor$. Therefore, by \Cref{lem: degenerate count free group upper bound}, we get:
    \begin{align*}
        |\mathcal{G}_W(\vec{l})\cap \mathcal{G}_{W'}(\vec{l}')|\leq \left|G^{\text{deg}}_{V}\left(\left \lfloor \frac{(6n+2)\log L}{h_r}\right \rfloor,\vec{l}_V\right)\right|\leq \frac{(2r-1)c'_1e^{\frac{h_rL}{2}}}{L^{6n}}
    \end{align*}
\end{proof}

\begin{proposition}\label{thm: free grp shell conj class lower bound asymptote}
    Let $n$ and $r$ be positive integers with $r\geq 2$. Then:
    \begin{equation*}
        \liminf_{\substack{L\to\infty\\L\textnormal{ even}}}\frac{|\partial \mathcal{C}^{+}_n(L)|}{L^{6n-4}e^{\frac{h_rL}{2}}}\geq \frac{\kappa^+_n(r)\mathfrak{M}^+_n}{2^{6n-4}(6n-4)!(2r)^{6n-3}(2r-1)^{6n-3}}
    \end{equation*}
\end{proposition}

\begin{proof}
    The proof follows from \Cref{lem: conj class substitution count asymptote part 2}, \Cref{lem: key lemma for free group conj class lower bound counting argument}, and \Cref{thm: free grp shell conj class upper bound asymptote}, in the same way in which \Cref{thm: free grp shell lower bound asymptote} follows from \Cref{lem: substitution count asymptote}, \Cref{lem: non-cancellingcount2}, and \Cref{thm: free grp shell upper bound asymptote}. We provide a summary below for the sake of completeness.
    
    Let $W$ be a maximal orientable Wicks form of genus $n$, thus having length \(12n - 6\) and \(6n-3\) unknowns. Let $\mathcal{L}(W,n,r,L), L_0$ and $c'_2$ be as in \Cref{lem: key lemma for free group conj class lower bound counting argument}. For any $\vec{l}\in\mathcal{L}(W)$, a conjugacy class in $\mathcal{G}_W(\vec{l})$ has genus at most $n$. That is we have:
    \[\bigcup_{\substack{W\\\text{maximal}\\\text{orientable}\\\text{genus }n}}\bigcup_{\vec{l} \in \mathcal{L}}\mathcal{G}_W(\vec{l}\,)\subset \partial \mathcal{C}_{\leq n}^+(L)\]
    
    where \(\partial \mathcal{C}_{\leq n}^+(L)\coloneq \bigcup_{i=1}^n\partial \mathcal{C}^{+}_i(L)\) is the set of all elements in \(F_r\) having orientable genus lesser than or equal to \(n\) and length equal to $L$. 

    Recall that for any $\sigma\in\text{Aut}(W)$, $\mathcal{G}_W(\vec{l})=\mathcal{G}_W(\sigma\cdot \vec{l})$. Therefore, we define: $\tilde{\mathcal{L}}(W,n,r,L)\coloneq \mathcal{L}(W,n,r,L)/\text{Aut}(W)$.
    So we have:
    \begin{align*}
        |\mathcal{C}_{\leq n}^+(L)|\geq \sum_W\sum_{[\vec{l}]\in \tilde{\mathcal{L}}(W)}|\mathcal{G}_W(\vec{l})| - \frac{1}{2}\sum_{\substack{W,W'\\\vec{l}\in\tilde{\mathcal{L}}(W),\vec{l}'\in \tilde{\mathcal{L}}(W')\\(W,[\vec{l}])\neq(W',[\vec{l}'])}}|\mathcal{G}_W(\vec{l})\cap \mathcal{G}_W(\vec{l}')|
    \end{align*}
    
    Using \Cref{lem: non-cancellingcount2}, \Cref{lem: key lemma for free group conj class lower bound counting argument}, and the fact that for any two maximal orientable Wicks forms of the same genus, $|\mathcal{L}(W,n,r,L)|=|\mathcal{L}(W',n,r,L)|$ :
    \begin{align*}
        \frac{1}{L^{6n-4}e^{\frac{h_rL}{2}}}\sum_{\substack{W,W'\\\vec{l}\in\tilde{\mathcal{L}}(W),\vec{l}'\in \tilde{\mathcal{L}}(W')\\(W,[\vec{l}])\neq(W',[\vec{l}'])}}|\mathcal{G}_W(\vec{l})\cap \mathcal{G}_{W'}(\vec{l}')|&\leq \frac{c'_2(\mathcal{N}^+(n))^2e^{-O(n\log n)}}{L^4}
    \end{align*}
    \begin{equation*}
        \text{Therefore, }\lim_{L\to\infty}\frac{1}{L^{6n-4}e^{\frac{h_rL}{2}}}\sum_{\substack{W,W'\\\vec{l}\in\tilde{\mathcal{L}}(W),\vec{l}'\in \tilde{\mathcal{L}}(W')\\(W,[\vec{l}])\neq(W',[\vec{l}'])}}|\mathcal{G}_W(\vec{l})\cap \mathcal{G}_{W'}(\vec{l}')|=0
    \end{equation*}
    Using \Cref{lem: conj class substitution count asymptote part 2}, \Cref{lem: kappa for maximal Wicks forms}, and the fact that:
    \[\lim_{L\to\infty}\frac{|\tilde{\mathcal{L}}(W,n,r,L)|}{L^{6n-4}}=\begin{cases}
            \frac{1}{|\text{Aut}(W)|}\frac{1}{2^{6n-4}(6n-4)!} &\text{; }W\neq ABCA^{-1}B^{-1}C^{-1}\\
            \frac{1}{3}\frac{1}{2^{2}(2)!} &\text{; }W= ABCA^{-1}B^{-1}C^{-1}
        \end{cases}\]
    we get:
    \begin{equation*}
        \lim_{L\to\infty}\frac{1}{L^{6n-4}e^{\frac{h_rL}{2}}} \sum_{[\vec{l}]\in \tilde{\mathcal{L}}(W)}|\mathcal{G}_W(\vec{l})|=\frac{1}{|\text{Aut}(W)|} \frac{\kappa^+_n(r)}{2^{6n-4}(6n-4)!(2r)^{6n-3}(2r-1)^{6n-3}}
    \end{equation*}
    \begin{align*}
        \text{Therefore, } \liminf_{\substack{L\to\infty\\L\textnormal{ even}}}\frac{|\partial \mathcal{C}^{+}_{\leq n}(L)|}{L^{6n-4}e^{\frac{h_rL}{2}}}\geq \frac{\kappa^+_n(r)\mathfrak{M}^+_n}{2^{6n-4}(6n-4)!(2r)^{6n-3}(2r-1)^{6n-3}}
    \end{align*}
    Now the proposition follows from \Cref{thm: free grp shell conj class upper bound asymptote}.
\end{proof}

And entirely analogously, we have:
\begin{lemma}\label{lem: key lemma for free group conj class lower bound counting argument non-orientable}
    Let $W$ be a maximal non-orientable Wicks form of genus $n$. Then, for all $L\geq 1$, there exists a collection of length types $\mathcal{L}=\mathcal{L}(W,n,r,L)$ and there are constants $L_0,c'_3\geq 1$ depending only on $r$ and $n$ satisfying, for all $L\geq L_0$:
    \begin{enumerate}[font=\normalfont]
        \item For all $\vec{l}\in \mathcal{L}(W,n,r,L)$, $||\vec{l}||=L$, $\mathcal{L}(W,n,r,L)$ is $\text{Aut}(W)$-invariant,  and \[\lim_{L\to\infty}\frac{|\mathcal{L}(W,n,r,L)|}{L^{3n-4}}=\frac{1}{2^{3n-4}(3n-4)!}
        \]
    \end{enumerate}
    Let $W'$ be another maximal non-orientable Wicks form of genus $n$, and let $[\vec{l}]\in \mathcal{L}(W,n,r,L)/\text{Aut}(W)$, $[\vec{l}']\in \mathcal{L}(W',n,r,L)/\text{Aut}(W')$.
    \begin{enumerate}[font=\normalfont, resume]
        \item If $(W,[\vec{l}])\neq(W',[\vec{l}'])$, then:
        \begin{align*}
            |\mathcal{G}_W(\vec{l})\cap \mathcal{G}_{W'}(\vec{l}')|\leq c'_3\frac{e^{\frac{h_rL}{2}}}{L^{3n}}
        \end{align*}
    \end{enumerate}
\end{lemma}
\begin{proposition}\label{thm: free grp shell conj class lower bound asymptote non-or}
    Let $n$ and $r$ be positive integers with $r\geq 2$. Then:
    \begin{equation*}
        \liminf_{\substack{L\to\infty\\L\textnormal{ even}}}\frac{|\partial \mathcal{C}^{-}_n(L)|}{L^{3n-4}e^{\frac{h_rL}{2}}}\geq \frac{\kappa^-_n(r)\mathfrak{M}^-_n}{2^{3n-4}(3n-4)!(2r)^{3n-3}(2r-1)^{3n-3}}
    \end{equation*}
\end{proposition}

\subsection{Counting in Hyperbolic groups}\label{subsec: Conjugacy class count: Lower bounds: Hyperbolic groups}
The aim of this section is to provide lower bounds for conjugacy class counts in hyperbolic groups. The proofs adapt the arguments in \Cref{subsec: Conjugacy class count: Lower bounds: Free groups} to a general hyperbolic setting in the same way in which \Cref{subsec: Element count: Lower bounds: Hyperbolic groups} adapts the proofs in \Cref{subsec: Element count: Lower bounds: Free groups}, using the geodesic encoding construction. The main work done in this section is to provide analogues for \Cref{lem: key lemma for conj class lower bounds}, \Cref{lem: conj class substitution count asymptote part 2}, and \Cref{lem: key lemma for free group conj class lower bound counting argument}.

We will use definitions and notation introduced in \Cref{subsec: Element count: Lower bounds: Hyperbolic groups}. $G$ will denote a hyperbolic group equipped with the word metric induced by some symmetric generating set $S$. $h$ will denote the volume entropy of $G$. Let $K, R_0, K'(R_0), R_1$, $D_1,D_2,D_3,\tilde{D}_3$ be the constants appearing in \Cref{def: goodsubs} and \Cref{construction: geodesicencoding}. In particular, $R_0, R_1$ are constants chosen to satisfy the same conditions as in \Cref{construction: geodesicencoding}. Let $W$ be a Wicks form and let $\vec{l}$ be a length type for $W$. Let $V$ be a linear word representing $W$. $G_V(\vec{l})$ will denote the set of group elements coming from good substitutions of the length type $\vec{l}_V$ (see \Cref{def: goodsubs}). Define $\mathcal{G}_W(\vec{l})$ to be the set of conjugacy classes of group elements in $G_V(\vec{l})$.

Furthermore, let $\vec{l},\vec{l}'$ be length types for $W$, and $\mathfrak{g}_1,\mathfrak{g}_1$ be substitutions for $W$ of length types $\vec{l},\vec{l}'$ respectively. Let $c>0$. We write $\vec{l}\simeq_{c}\vec{l}'$ if for every unknown $X$ of $W$, $|l_X-l'_X|\leq c$. We write $\mathfrak{g}_1\simeq_c\mathfrak{g}_2$ if for every unknown $X$ of $W$, there exist group elements $\xi_X,\xi'_X\in B(c)$ such that $(g_1)_{X}=\xi_X (g_2)_{X}\xi'_X$.

\begin{definition}\label{def: set of degenerate group elements hyperbolic groups}
    Let $V$ be a linear Wicks form and let $\vec{l}$ be a length type for $V$. Let $N,D$, and $R$ be positive integers. Define $G^{\text{deg}}_V(N,D,R,\vec{l})$ to be the set of group elements $\mathfrak{g}'(V)$, where $\mathfrak{g}'$ is the good substitution associated to a substitution $\mathfrak{g}$ of length type $\vec{l}$, which has either of the following algebraic forms:
    \begin{enumerate}[font=\normalfont]
        \item \linkdest{algebraic form 1 conj} There are two distinct unknowns $X$ and $X'$ of $V$ such that:
        \begin{equation*}
            g_X=uzv,\quad g_{X'}=u'\xi_1z^{\varepsilon}\xi_2v'
        \end{equation*}
        where $\varepsilon\in\{+1,-1\}$, $\xi_1,\xi_2\in B(R)$, and:
        \begin{equation*}\label{eq: hyp degen len eqn 1 conj}
            l(z)=ND,\quad l(u)+l(v)+ND=l_X,\quad l(u')+l(v')+ND=l_{X'}
        \end{equation*}
        \item \linkdest{algebraic form 2 conj} There is an unknown $X$ of $W$ such that:
        \begin{equation*}
            g_X=uzw\xi_1z^{\varepsilon}\xi_2v
        \end{equation*}
        where $\varepsilon\in\{+1,-1\}$, $\xi_1,\xi_2\in B(R)$, and:
        \begin{equation*}\label{eq: hyp degen len eqn 2 conj}
            l(z)=ND,\quad l(u)+l(w)+l(v)+2ND=l_X
        \end{equation*}
        \item \linkdest{algebraic form 3 conj} There is an unknown $X$ of $W$, and a positive integer $T\geq D$ such that:
        \begin{equation*}
            g_X=uz_0z_1\cdots z_kv
        \end{equation*}
        where $k=\left\lfloor\frac{ND}{T}\right\rfloor+1$, $l(z_0)=l(z_1)=\cdots=l(z_{k-1})=T$, $l(z_k)<T$. Moreover:
        \begin{equation*}
             z_i=\xi_iz_{i-1}\xi'_i \;\; (1\leq i\leq k-1),\quad z_k=\xi_kz'_{k-1}\xi'_k
        \end{equation*}
        where $\xi_1,\xi'_1,\dots,\xi_k,\xi'_k\in B(R)$, and $z'_{k-1}$ is a group element of length $ND-(k-1)T$ which appears as a prefix in some geodesic word representing $z_{k-1}$. Furthermore:
        \begin{equation*}\label{eq: hyp degen len eqn 3 conj}
            l(u)+l(v)+l(z_k)+kT=l_X,\quad (k-1)T+l(z_k)=ND
        \end{equation*}
        \item \linkdest{algebraic form 4 conj} There is an unknown $X$ of $W$ such that:
        \begin{equation*}
            g_X=uzw\xi_1z^{-1}\xi_2v
        \end{equation*}
        where $\xi_1,\xi_2\in B(R)$, $w\in B(3D)$, and:
        \begin{equation}\label{eq: hyp degen len eqn 4 conj}
            l(u)+l(w)+l(v)+2l(z)=l_X,\quad l(z)+l(w)=ND
        \end{equation}
        \item \linkdest{algebraic form 5 conj} There is an unknown $X$ of $W$ such that:
        \begin{equation*}
            g_X=uz_0z_1w\xi_1z_1^{-1}\xi'_1\xi_2z_0^{-1}\xi'_2v
        \end{equation*}
        where $\xi_1,\xi_2,\xi'_1,\xi'_2\in B(R)$, $w\in B(D)$, and:
        \begin{equation*}\label{eq: hyp degen len eqn 5 conj}
            l(z_1)\geq D,\quad l(u)+l(w)+l(v)+2l(z_0)+2l(z_1)=l_X,\quad l(z_0)+2l(z_1)+l(w)=ND
        \end{equation*}
    \end{enumerate}

\end{definition}
\begin{lemma}\label{lem: degenerate count hyp group upper bound}
    Let $N,D,R$ be positive integers. Let $V$ be a linear Wicks form of length $2m$ and let $\vec{l}$ be a length type for $V$. Then there is a constant $c_1>0$ depending only on $m,\delta,h,D,R$ and a positive integer $\tilde{D}_4$ depending only on $h,R,\delta$, such that for all $D>\tilde{D}_4$, we have:
    \begin{align*}
        |G^{\text{deg}}_V(N,D,R,\vec{l})|\leq c_1L^2e^{\frac{h||\vec{l}||}{2} }e^{-\frac{hND}{2}}
    \end{align*}
\end{lemma}
\begin{proof}
    The proof is identical to that of \Cref{lem: hyperbolic element intersection count for reduction in degrees of freedom}. Formally, let $C_2$ be the constant appearing in \Cref{thm: ballgrowth} and set $\tilde{D}_4=\lceil4(R+\frac{\log C_2}{h})\rceil$. Then for all $D>\tilde{D}_4$, the analysis done in the proof of \Cref{lem: hyperbolic element intersection count for reduction in degrees of freedom} works out, with $R$ and $D$ assuming the roles of $D_4$ and $D_0$ respectively. 
\end{proof}

\begin{lemma}\label{lem: key lemma for conj class lower bounds hyp}
    There exist positive numbers $D_{\text{conj}}$, $\tilde{D}_0$, $R_{\text{conj}}$ depending only on $\delta, K,R_0,R_1$ such that the following holds true: Let $N,D$ be positive integers. Let $V$ and $V'$ be linear maximal Wicks forms of length 2m. Let $\vec{l}$ and $\vec{l}'$ be length types for $V$ and $V'$ respectively, satisfying:
    \begin{enumerate}[font=\normalfont, label=(\alph*)]
        \item\label{item: key lemma conj class hyp 1} $||\vec{l}||=||\vec{l}'||$
        \item\label{item: key lemma conj class hyp 2} For all $1\leq i\leq 2m$, $1\leq j\leq 2m$, $l_i,l'_j> 10ND + 30(D+(D_2+2D_{\text{conj}})) + 2R_1 $
    \end{enumerate}
    Let $\mathfrak{g}_1$ and $\mathfrak{g}_2$ be substitutions of length type $\vec{l}$ and $\vec{l}'$ for $V$ and $V'$ respectively, such that $[\mathfrak{g}'_1(V)]=[\mathfrak{g}'_2(V')]$. Then, for all $D>\tilde{D}_0$ either:
    \begin{enumerate}[font=\normalfont]
        \item $\mathfrak{g}'_1(V)\in G^{\text{deg}}_{V}(N,D,R_{\text{conj}},\vec{l})$, or
        \item There is a cyclic permutation $\sigma$ such that $V$ is linearly equivalent to $\sigma(V')$, $\vec{l}\simeq_{25D} \sigma\cdot\vec{l}'$, and $\mathfrak{g}_1\simeq_{25D}\sigma\cdot\mathfrak{g}_2$.
    \end{enumerate}
\end{lemma}
\begin{proof}
    \begin{claim}
        There exists a positive integer $D_{\text{conj}}$ depending only on $K,R_0,R_1,\delta$ such that either:
        \begin{enumerate}[label=(\Roman*)]
            \item There are cyclic permutations $\sigma$ and $\sigma'$ such that $\mathfrak{g}_1(\sigma(V))=\alpha\mathfrak{g}_2(\sigma'(V'))\alpha^{-1}$, where $\alpha\in B(D_{\text{conj}})$ or
            \item There are cyclic permutations $\sigma$ and $\sigma'$, and a split Wicks form $V'_{\text{split}}$ such that $\mathfrak{g}'_1(\sigma(V))=\alpha(\mathfrak{g}'_2)_{\text{split}}(V'_{\text{split}})\alpha^{-1}$, where $V'_{\text{split}}$ and $(\mathfrak{g}'_2)_{\text{split}}$ are obtained from a splitting of the first letter of $\sigma'(V')$, and $\alpha\in B(D_{\text{conj}})$. More precisely, we write $\sigma'(V')=X'_{\sigma'(1)}U'_1(X')^{\varepsilon}U'_2$, split $X'_{\sigma'(1)}=B'X''_{\sigma'(1)}$ and set $V'_{\text{split}}=X''_{\sigma'(1)}U'_1(BX''_{\sigma'(1)})^{\varepsilon}U'_2B$. $V'_{\text{split}}$ comes with a substitution $(\mathfrak{g}_2)_{\text{split}}$, which is the same as $\mathfrak{g}_2$ for all letters except $X'_{\sigma'(1)}$, which has been split.  A geodesic word representing $(g_2)_{X'_{\sigma'(1)}}$ is chosen, and there is a splitting $(g_2)_{X'_{\sigma'(1)}}=(g_2)_{B}(g_2)_{X''_{\sigma'(1)}}$, where $(g_2)_{B}$ and $(g_2)_{X''_{\sigma'(1)}}$ are group elements corresponding to a splitting of the geodesic word chosen for $(g_2)_{X'_{\sigma'(1)}}$, such that $l((g_2)_{B})>R_1+2D_{\text{conj}}$ and $l((g_2)_{X''_{\sigma'(1)}})>R_1+2D_{\text{conj}}$. $((g_2)_{\text{split}})_{X''_{\sigma'(1)}}\coloneq(g_2)_{X''_{\sigma'(1)}}$ and $((g_2)_{\text{split}})_{B}\coloneq(g_2)_{B}$. The splitting is also recorded in a length type for $V'_{\text{split}}$, which we continue to denote as $\vec{l}'$. The substitution $(\mathfrak{g}'_2)_{\text{split}}$ is constructed from the good substitution associated to $\mathfrak{g}_2$, that is, from $\mathfrak{g}'_2$ in the natural way: $(\mathfrak{g}'_2)_{\text{split}}$ is the same as $\mathfrak{g}_2$ for all letters except $X'_{\sigma'(1)}$, and $((g'_2)_{\text{split}})_{B}\coloneq l_{X'_{\sigma'(1)}}(g_2)_{B}$, $((g'_2)_{\text{split}})_{X''_{\sigma'(1)}}\coloneq (g_2)_{X''_{\sigma'(1)}}r_{X'_{\sigma'(1)}}$, where $l_{X'_{\sigma'(1)}}$ and $r_{X'_{\sigma'(1)}}$ are Schottky elements corresponding to $((g_2)_{\text{split}})_{X'_{\sigma'(1)}}$. 
        \end{enumerate}
    \end{claim}
    \begin{proof}
        This can be easily proved using the ideas in \cite[Chapter {\RomanCaps{3}}.$\Gamma$, Lemma 2.9]{bridson2013metric} adapted for chain substitutions. It is worth noting that the length lower bounds for the split letter are chosen to be greater than $R_1$ so that the geodesic encoding construction can be done.
    \end{proof}
    We take up each of the cases listed in the claim individually:\\\textbf{Case - {\RomanCaps{1}}}: Let $g_1=\mathfrak{g}'_1(\sigma(V))$ and $g_{2}=\mathfrak{g}'_2(\sigma(V'))$. Let  $\gamma_{g_1}$ be a geodesic joining the identity to $g_1$, and let $\gamma_{g_2}$ be a geodesic joining $\alpha$ to $\alpha g_{2}$. Comparing the lengths of $\gamma_{g_1}$ and $\gamma_{g_2}$, we chop off an initial subsegment of length $\Delta(g_1,g_2)\coloneq |l(\gamma_{g_1})-l(\gamma_{g_2})|\leq 2D_{conj}$ from the longer segment, and we denote the new segments thus obtained as $\tilde{\gamma}_{g_1}$ and $\tilde{\gamma}_{g_2}$. It follows from \Cref{prop: ftbound} that $\tilde{\gamma}_{g_1}$ and $\tilde{\gamma}_{g_2}$ $D_{FT}$-fellow travel, where $D_{FT}=6D_{\text{conj}}+4\delta$.\\Let $\mathcal{E}_1$ and $\mathcal{E}_2$ denote geodesic encodings of $\mathfrak{g}_1$ and $\mathfrak{g}_2$ into $\gamma_{g_1}$ and $\gamma_{g_2}$ respectively (encodings can be constructed because good substitutions are cyclically good, as explained in \Cref{rmk: good subs for Wicks forms are cyclically good}). The encodings induce two canonical partitions on the set of positions $\{1,\dots, L\}$ where $L=\min\{l(\gamma_{g_1}),l(\gamma_{g_2})\}$, as follows: Without any loss of generality, suppose $l(\gamma_{g_1})\leq l(\gamma_{g_2})$. $\mathcal{E}_1$ comes with the data $(\vec{l_{\pi}},\vec{b_{\pi}})$ which induces one canonical partition on $\{1,\dots, L\}$. Similarly, $\mathcal{E}_2$ comes with the data $(\vec{l'_{\pi}},\vec{b'_{\pi}})$. We restrict it to the subsegment $\tilde{\gamma}_{g_2}$, that is, we discard the initial subsegment of $\gamma_{g_2}$ of length $\Delta(g_1,g_2)$ and modify the data comprising the encoding $\mathcal{E}_2$ in the obvious way: if $(b'_{\pi})_0\geq \Delta(g_1,g_2)$, re-define it to be equal to $(b'_{\pi})_0-\Delta(g_1,g_2)$. Otherwise, we set $(b'_{\pi})_0=0$ and we re-define $(l'_{\pi})_{X'_{\sigma'(1)}}$ to equal $(l'_{\pi})_{X'_{\sigma'(1)}}-\Delta(g_1,g_2)$. We do the same for the paired letter $X_{\overline{\sigma'(1)}}$, and depending on the sign of $X_{\sigma'(1)}$, we increase the size of the junction preceding or following $X_{\overline{\sigma'(1)}}$ by $\Delta(g_1,g_2)$. Thus, the segments $\tilde{\gamma}_{g_1}$ and $\tilde{\gamma}_{g_2}$ correspond to the set of positions $\{1,\dots, L\}$ and there are two canonical partitions on them, coming from $(\sigma(V),\vec{l_{\pi}},\vec{b_{\pi}})$ and $(\sigma'(V'),\vec{l'_{\pi}},\vec{b'_{\pi}})$. Moreover, $\vec{b}_{\pi},\vec{b'_{\pi}}\in [0, H]^{2m+1}$, where $H\coloneq D_2+2D_{\text{conj}}$ (see property \linkref{encoding: junctions}{(P4)}). Therefore, together with \linkref{encoding: mismatch bound}{(P5)}, we have:
    \begin{equation}\label{eq: first length discrepancy bound}
        \text{for all }1\leq j\leq 2m, |l_{X_{\sigma(j)}}-(l_{\pi})_{X_{\sigma(j)}}|\leq 2D_3+2D_{\text{conj}}, |l'_{X'_{\sigma'(j)}}-(l'_{\pi})_{X'_{\sigma'(j)}}|\leq 2D_3+2D_{\text{conj}}
    \end{equation}
    The data $(\sigma(V),\vec{l_{\pi}},\vec{b_{\pi}})$ and $(\sigma'(V'),\vec{l'_{\pi}},\vec{b'_{\pi}})$ satisfy the hypotheses of \Cref{thm: the main reduction theorem}, using which we infer that either of the following occurs:
    \begin{enumerate}
        \item There is $D$-reduction by $ND$. In this case, $g_1=\mathfrak{g}'_1(\sigma(V))\in G^{\text{deg}}_{\sigma(V)}(N,D,R_{\text{conj}},\vec{l})$, where $R_{\text{conj}}=2D_3+16(K+40\delta)+2D_{FT}+4D_{\text{conj}}$. This follows from arguments identical to those in the algebraic part of the proof of \Cref{prop: hyperbolic degenerate substitutions}, with the only differences being that there is an additive factor of $2D_{FT}$ due to the fact that $\tilde{\gamma}_{g_1}$ and $\tilde{\gamma}_{g_2}$ $D_{FT}$-fellow travel, and another additive factor of $4D_{\text{conj}}$ which comes from the $\Delta(g_1,g_2)$ adjustment. 
        \item There is no $D$-reduction by $ND$, therefore $\sigma(V)$ is linearly equivalent to $\sigma'(V')$. and:
        \[\text{for all }1\leq j\leq 2m, |(l_\pi)_{X_{\sigma(i)}}-(l'_\pi)_{X'_{\sigma'(i)}}|\leq 3(D-1)+4(D_2+2D_{\text{conj}})\]
        So, using \Cref{eq: first length discrepancy bound}, we get, for all $1\leq j\leq 2m$:
        \begin{align*}
            |l_{X_{\sigma(i)}}-l'_{X'_{\sigma'(i)}}|&\leq 3(D-1)+4(D_2+2D_{\text{conj}})+2(2D_3+2D_{\text{conj}})\\
            &<3D+12(D_2+D_3+D_{\text{conj}})
        \end{align*}
        Next, consider any $1\leq j\leq 2m$, and the segments $\gamma_{X_{\sigma(j)}}$ and $\gamma_{X'_{\sigma'(j)}}$ which $8(K+40\delta)$-fellow travel with $\pi_{X_{\sigma(j)}}$ and $\pi_{X'_{\sigma'(j)}}$ respectively (see \Cref{construction: geodesicencoding}, and property \linkref{encoding: fellow travel bound}{(P6)}). Recalling the definition of shifts, \Cref{def: shifts}, we note that:
        \begin{align*}
            d(i(\pi_{X_{\sigma(j)}}),i(\pi_{X'_{\sigma'(j)}}))&\leq D_{FT}+ |\text{shift}(j-1)|+ |(b_{\pi})_{j}-(b'_{\pi})_j|\\
            d(\tau(\pi_{X_{\sigma(j)}}),\tau(\pi_{X'_{\sigma'(j)}}))&\leq D_{FT}+ |\text{shift}(j)|
        \end{align*}
        Moreover, by \Cref{eq: shift bound in element count}, for all $1\leq j\leq 2m $, $|\text{shift}(j)|\leq \frac{3}{2}(D-1+ (D_2+2D_{\text{conj}}))$. Thus, we compute:
        \begin{align*}
            &\max\{d(i(\pi_{X_{\sigma(j)}}),i(\pi_{X'_{\sigma'(j)}})),d(\tau(\pi_{X_{\sigma(j)}}),\tau(\pi_{X'_{\sigma'(j)}}))\}\\&\leq \frac{3}{2}(D-1+ (D_2+2D_{\text{conj}})) + (D_2+2D_{\text{conj}})+D_{FT}\\
            &<2D+(3D_2+6D_{\text{conj}}+D_{FT})
        \end{align*}
        Now it follows that $(g_1)_{X_{\sigma(j)}}=\xi_j (g_2)_{X'_{\sigma'(j)}}\xi'_j $, where $\xi_j,\xi'_j$ are elements of length at most $2D+(3D_2+6D_{\text{conj}}+D_{FT})+(2D_3+2D_{\text{conj}})$ where the extra additive term $2D_3$ is due to property \linkref{encoding: mismatch bound}{(P5)} and $2D_{\text{conj}}$ is due to the modification done to the encodings at the start. Therefore, for $D>100(D_2+D_3+D_{FT}+D_{\text{conj}})$: $ \vec{l}\simeq_{25D}\tilde{\sigma}\cdot \vec{l}'$, and $\mathfrak{g}_1\simeq_{25D}\tilde{\sigma}\cdot\mathfrak{g}_2$, where $\tilde{\sigma}=\sigma^{-1}\circ\sigma'$.
    \end{enumerate}
    \noindent\textbf{Case - {\RomanCaps{2}}}: This case is analyzed just as in the previous case, and the role of \Cref{thm: the main reduction theorem} is now assumed by \Cref{prop: reduction or no splitting}. To begin with, let $g_1=\mathfrak{g}'_1(\sigma(V))$ and $g_2=(\mathfrak{g}_2')_{\text{split}}(V'_{\text{split}})$. Let $\gamma_{g_1}$, $\gamma_{g_2}$, $\mathcal{E}_1$,$\mathcal{E}_2$, $\tilde{\gamma}_{g_1}$,$\tilde{\gamma}_{g_2}$, $\Delta(g_1,g_2)$ be defined as above. The data $(\vec{l_{\pi}},\vec{b_{\pi}})$ and $(\vec{l'_{\pi}},\vec{b'_{\pi}})$ and the encodings $\mathcal{E}_1$ and $\mathcal{E}_2$ are adjusted to correct for $\Delta(g_1,g_2)$. \Cref{eq: first length discrepancy bound} continues to hold true, mutatis mutandis.
    The data $(\sigma(V),\vec{l_{\pi}},\vec{b_{\pi}})$ and $(V'_{\text{split}},\vec{l'_{\pi}},\vec{b'_{\pi}})$ satisfy the hypotheses of \Cref{prop: reduction or no splitting}, using which we infer that either of the following occurs:
    \begin{enumerate}
        \item There is $D$-reduction by $ND$. In this case, $g_1=\mathfrak{g}'_1(\sigma(V))\in G^{\text{deg}}_{\sigma(V)}(N,D,R_{\text{conj}},\vec{l})$, where $R_{\text{conj}}=2D_3+16(K+40\delta)+2D_{FT}$, as before.
        \item There is no $D$-reduction by $ND$, therefore either $\sigma(V)$ is linearly equivalent to $\sigma'(V')$, or to $\sigma''\circ\sigma'(V')$, where $\sigma''$ is the cyclic permutation $(1 2 \dots 2m)$. Without any loss of generality, we assume the former. Therefore, we also have:
        \begin{align*}
            &(l'_{\pi})_B\leq 2(D-1)+(D_2+2D_{\text{conj}}),\\
            & |(l_\pi)_{X_{\sigma(i)}}-(l'_\pi)_{X'_{\sigma'(i)}}|\leq  11(D-1)+ 16(D_2+2D_{\text{conj}})\;(1<j\leq 2m, j\neq (\sigma')^{-1}(\overline{\sigma'(1)}))\\
            &|(l_\pi)_{X_{\sigma(1)}}-(l'_\pi)_{X''_{\sigma'(1)}}|\leq  11(D-1)+ 16(D_2+2D_{\text{conj}})
        \end{align*}
        Using the analogous version of \Cref{eq: first length discrepancy bound}, for $1<j\leq 2m, j\neq (\sigma')^{-1}(\overline{\sigma'(1)})$:
        \begin{align*}
            |l_{X_{\sigma(i)}}-l'_{X'_{\sigma'(i)}}|&\leq  11(D-1)+ 16(D_2+2D_{\text{conj}})+2(2D_3+2D_{\text{conj}})\\
            &<11D+18(D_2+2D_{\text{conj}})+4D_3
        \end{align*}
        Using the fact that $l'_{X'_{\sigma'(1)}}=l'_{B}+l'_{X''_{\sigma'(1)}}$, we have:
        \begin{align*}
            |l_{X_{\sigma(1)}}-l'_{X'_{\sigma'(1)}}|&\leq |l_{X_{\sigma(1)}}-l'_{X''_{\sigma'(1)}}|+(l'_{\pi})_B+|l'_B-(l'_{\pi})_B|\\
            &<13D+20(D_2+2D_{\text{conj}})+6D_3
        \end{align*}
        As before, consider any $1\leq j\leq 2m$, and the segments $\gamma_{X_{\sigma(j)}}$ and $\gamma_{X'_{\sigma'(j)}}$ which $8(K+40\delta)$-fellow travel with $\pi_{X_{\sigma(j)}}$ and $\pi_{X'_{\sigma'(j)}}$ respectively. By \Cref{eq: shift bound conj class count}, we have: $\text{shift}(j)<5D+8(D_2+2D_{\text{conj}})$, and therefore:
        \begin{align*}
            &\max\{d(i(\pi_{X_{\sigma(j)}}),i(\pi_{X'_{\sigma'(j)}})),d(\tau(\pi_{X_{\sigma(j)}}),\tau(\pi_{X'_{\sigma'(j)}}))\}\\
            &< (5D+8(D_2+2D_{\text{conj}}))+ (D_2+2D_{\text{conj}})+D_{FT}
        \end{align*}
         As before, it follows that for all $j\neq1, (\sigma')^{-1}(\overline{\sigma'(1)})$  , $(g_1)_{X_{\sigma(j)}}=\xi_j (g_2)_{X'_{\sigma'(j)}}\xi'_j $, $(g_1)_{X_{\sigma(1)}}=\xi_1 (g_2)_{X''_{\sigma'(1)}}\xi'_1$ where all $\xi_j,\xi'_j$ are elements of length at most $5D+9(D_2+2D_{\text{conj}})+D_{FT}+(2D_3+2D_{\text{conj}})$. Because $l_{B}\leq (l'_{\pi})_B+|l'_B-(l'_{\pi})_B|<2D+2(D_2+2D_{\text{conj}})+2D_3$, it follows that: $(g_1)_{X_{\sigma(1)}}=\xi''_1 (g_2)_{X'_{\sigma'(1)}}\xi'_1$, $\xi''_1=\xi'_1l_{B}^{-1}$ has length at most $7D+12(D_2+2D_{\text{conj}})+4D_3+D_{FT}$.

         Therefore, for $D>100(D_2+D_3+D_{FT}+D_{\text{conj}})$: $ \vec{l}\simeq_{25D}\tilde{\sigma}\cdot \vec{l}'$, and $\mathfrak{g}_1\simeq_{25D}\tilde{\sigma}\cdot\mathfrak{g}_2$, where $\tilde{\sigma}=\sigma^{-1}\circ\sigma'$.
    \end{enumerate}
    
\end{proof}

In view of \Cref{lem: degenerate count hyp group upper bound} and \Cref{lem: key lemma for conj class lower bounds hyp}, for the rest of this section, we fix $D$ to be some integer greater than both $\tilde{D}_0$ and $\tilde{D}_4$. 

\begin{lemma}\label{lem: conj class substitution count asymptote part 2 hyp}
    Let $N\geq 1$ be a positive integer. Let $W$ be a maximal Wicks form of length $2m$ and let $\vec{l}$ be a non-degenerate length type for $W$ satisfying condition \ref{item: key lemma conj class hyp 2} of \Cref{lem: key lemma for conj class lower bounds hyp}. Let $L=||\vec{l}||$. Then, there are constants $c_2,c_3>0$ depending only on $m,\delta,h,D,K,R_0,R_1$ such that:
    \begin{align*}
        \frac{\mathcal{G}_W(100D\cdot \vec{l})}{c_2e^{\frac{h\cdot 100DL}{2}}}\geq \begin{cases}
            \left ( 1-\frac{c_3(100DL)^2}{e^{\frac{hND}{2}}}\right) &; W\neq ABCA^{-1}B^{-1}C^{-1}\\
            \frac{1}{2}\left ( 1-\frac{c_3(100DL)^2}{e^{\frac{hND}{2}}}\right) &; W= ABCA^{-1}B^{-1}C^{-1}
        \end{cases}
    \end{align*}
\end{lemma}
\begin{proof}
    We write $\overrightarrow{l_D}=100D\cdot\vec{l}$. Let $V=X_1\cdots X_{2m}$ be a linear word representing $W$. Consider the map $\psi:G_V((\overrightarrow{l_D})_V)-G^{\text{deg}}_V(N,D,R_{\text{conj}},(\overrightarrow{l_D})_V)\rightarrow \mathcal{G}_W(\overrightarrow{l_D})$ which maps group elements to their conjugacy classes. The lemma follows from \Cref{claim: finite-to-one-ness of psi hyp}, \Cref{lem: degenerate count hyp group upper bound}, and \Cref{lem: lowerboundgoodsubs}.
    \begin{claim}\label{claim: finite-to-one-ness of psi hyp}
        There is a constant $c_2>0$ depending only on $m,\delta,h,K,R_0,R_1$ such that: if $V\neq ABCA^{-1}B^{-1}C^{-1}$, then $\psi$ is at most $c_2$-to-$1$, and if $V= ABCA^{-1}B^{-1}C^{-1}$, $\psi$ is atmost $2c_2$-to-$1$.
    \end{claim}
    \begin{proof}
        Let $\sigma$ be a generator of $\text{Aut}(V)$ and denote its order by $n$. Let $\mathfrak{g}_1$ and $\mathfrak{g}_2$ be two substitutions of length type $(\overrightarrow{l_D})_V$, such that $\mathfrak{g}'_1(V)$ and $\mathfrak{g}'_2(V)$ are members of the set $G_V((\overrightarrow{l_D})_V)-G^{\text{deg}}_V(N,D,R_{\text{conj}},(\overrightarrow{l_D})_V)$, and moreover, suppose that $[\mathfrak{g}'_1(V)]=[\mathfrak{g}'_2(V)]$. The data $(V,\mathfrak{g}_1,(\overrightarrow{l_D})_V)$ and $(V,\mathfrak{g}_2,(\overrightarrow{l_D})_V)$ satisfies conditions \ref{item: key lemma conj class hyp 1} and \ref{item: key lemma conj class hyp 2} of \Cref{lem: key lemma for conj class lower bounds hyp}, and $\mathfrak{g}'_1(V)\notin G^{\text{deg}}_V(N,D,R_{\text{conj}},(\overrightarrow{l_D})_V)$, so it follows that there is a cyclic permutation $\sigma^{j}$ such that $V$ is linearly equivalent to $\sigma^{j}(V)$, for some $0\leq j\leq n-1$. Moreover, $(\overrightarrow{l_D})_V\simeq_{25D} \sigma^j\cdot(\overrightarrow{l_D})_V$, and $\mathfrak{g}_1\simeq_{25D}\sigma^j\cdot\mathfrak{g}_2$. But $(\overrightarrow{l_D})_V\simeq_{25D} \sigma^j\cdot(\overrightarrow{l_D})_V$ implies that $\vec{l}_V\simeq_{1/4}\sigma^{j}\vec{l}_V$, and therefore, $\vec{l}_V=\sigma^{j}\vec{l}_V$. Since $\vec{l}_V$ is non-degenerate, it follows from \Cref{lemma: action of automorphisms on a length type} that when $V\neq ABCA^{-1}B^{-1}C^{-1}$, $\sigma^j$ has to be the identity and $\mathfrak{g}_1\simeq_{25D}\mathfrak{g}_2$. On the other hand, when $V= ABCA^{-1}B^{-1}C^{-1}$, either $j=0$ or $j=3$. That is, either $\mathfrak{g}_1\simeq_{25D}\mathfrak{g}_2$ or $\mathfrak{g}_1\simeq_{25D}\sigma^3\cdot\mathfrak{g}_2$. Thus, setting $c_2=(2|B(25D)|)^m$, the claim follows.  
    \end{proof}
\end{proof}

\begin{lemma}\label{lem: key lemma for hyp group conj class lower bound counting argument}
    Let $W$ be a maximal orientable Wicks form of genus $n$. There are positive constants $L_0$ and $c_4$, depending only on $n,\delta, K,R_0,R_1$ such that for all $L\geq L_0$, there is a collection of length types $\mathcal{L}(W,n,L)$ satisfying:
    \begin{enumerate}[font=\normalfont]
        \item For all $\vec{l}\in\mathcal{L}(W,n,L)$, $||\vec{l}||=2\left\lceil\frac{L}{2}\right\rceil$, $\mathcal{L}(W,n,L)$ is $\text{Aut}(W)$-invariant, and:
        \[\lim_{L\to\infty}\frac{\mathcal{L}(W,n,L)}{L^{6n-4}}=\frac{1}{2^{6n-4}(6n-4)!}\]
    \end{enumerate}
    Let $W'$ be another maximal orientable Wicks form of genus $n$, and let $[\vec{l}]\in\mathcal{L}(W,n,L)/\text{Aut}(W)$, $[\vec{l}']\in\mathcal{L}(W',n,L)/\text{Aut}(W')$.
    \begin{enumerate}[font=\normalfont, resume]
        \item If $(W,[\vec{l}])\neq (W',[\vec{l}'])$, then:
        \[\mathcal{G}_W(100D\cdot\vec{l})\cap \mathcal{G}_{W'}(100D\cdot\vec{l}')\leq c_4\frac{e^{\frac{h\cdot 100DL}{2}}}{L^{6n}}\]
    \end{enumerate}
\end{lemma}
\begin{proof}
    Let $W$ be a maximal orientable Wicks form of genus $n$, thus having a length of $12n-6$. Consider the following collection of length types for $W$:
    \begin{align*}
         \mathcal{L}(W,n,r,L)=
         \coloneq\left\{\vec{l}\in\mathbb{Z}^{6n-3}_+\middle|\quad \begin{aligned}
             &\sum_{i=1}^{6n-3}l_i = \left\lceil\frac{L}{2}\right\rceil, \vec{l}\text{ non-degenerate}, 1\leq i\leq 6n-3,\\
             & l_i> 10\left \lceil \frac{2(6n+2)\log L}{hD} \right \rceil+30(D+(D_2+2D_{\text{conj}})) + 2R_1
         \end{aligned}\right\}
     \end{align*}
    $\mathcal{L}(W,n,L)$  satisfies $(1)$, because asymptotically, most length types satisfying $\sum l_i=\left\lceil\frac{L}{2}\right\rceil$ are non-degenerate.
    Note that every pair of length types in $\mathcal{L}$ satisfies the hypotheses of \Cref{lem: key lemma for conj class lower bounds hyp}, with $N=\left \lceil \frac{2(6n+2)\log L}{hD}\right \rceil$. We write $\overrightarrow{l_D}=100D\cdot\vec{l}$ for any length type $\vec{l}$.
    
    Next, let $W'$ be another Wicks form, and $\vec{l}\in \mathcal{L}(W,n,L)$, $\vec{l}'\in \mathcal{L}(W',n,L)$. Let $V$ and $V'$ be linear Wicks forms representing $W$ and $W'$ respectively. Let $\mathfrak{g}_1$ and $\mathfrak{g}_2$ be two substitutions for $V$ and $V'$ of length types $\overrightarrow{l_D}$ and $\overrightarrow{l'_D}$ respectively, such that $[\mathfrak{g}'_1(V)]=[\mathfrak{g}'_2(V')]$. The data $(V,\mathfrak{g}_1,(\overrightarrow{l_D})_V)$ and $(V',\mathfrak{g}_2,(\overrightarrow{l'_D})_V)$ satisfy the hypotheses of \Cref{lem: key lemma for conj class lower bounds hyp}. If $\mathfrak{g}(V)\notin G^{\text{deg}}_{V}(N,D,R_{\text{conj}},(\overrightarrow{l_D})_V)$, where $N=\left \lceil \frac{2(6n+2)\log L}{hD}\right \rceil$, then it follows from \Cref{lem: key lemma for conj class lower bounds hyp} that there exists a cyclic permutation $\sigma$ such that $V$ is linearly equivalent to $\sigma(V')$, that is $W=W'$, and $(\overrightarrow{l_D})_V=\sigma\cdot(\overrightarrow{l'_D})_V$, that is $(W,[\vec{l}])= (W',[\vec{l}'])$, which is not possible. 

    Therefore, $\mathfrak{g}(V)\in G^{\text{deg}}_{V}(N,D,R_{\text{conj}},(\overrightarrow{l_D})_V)$, and using \Cref{lem: degenerate count hyp group upper bound}, we get:
    \begin{align*}
        |\mathcal{G}_W(100D\cdot\vec{l})\cap \mathcal{G}_{W'}(100D\cdot\vec{l}')|&\leq \left|G^{\text{deg}}_{V}\left(\left \lceil \frac{2(6n+2)\log L}{hD}\right \rceil,D,R_{\text{conj}},(\overrightarrow{l_D})_V\right)\right|\\
        &\leq c_1L^2\times e^{\frac{h\cdot 100D(L+2)}{2}}\times e^{-\frac{hD}{2}\left \lceil \frac{2(6n+2)\log L}{hD}\right \rceil}\\
        &\leq (c_1e^{100hD})\times\frac{e^{\frac{h\cdot 100DL}{2}}}{L^{6n}}
    \end{align*}
\end{proof}

\begin{theorem}\label{thm: lowerboundhyp orientable conj class}
    Let $G$ be a hyperbolic group equipped with the word metric induced by some finite symmetric generating set $S$. Let $h$ denote the volume entropy of $G$. There exist positive constants $c_5,L_0$ depending only on $n$, $G$, and the choice of the generating set, such that:
    \begin{equation*}
        \forall L\geq L_0, |\mathcal{C}^{+}_n(L)|\geq c_5 L^{6n-4}e^{hL/2}
    \end{equation*}
\end{theorem}
\begin{proof}
    This is easily proved using an inclusion-exclusion argument, using \Cref{lem: conj class substitution count asymptote part 2 hyp} and \Cref{lem: key lemma for hyp group conj class lower bound counting argument}.
\end{proof}

Similarly, we have:

\begin{lemma}
    Let $W$ be a maximal non-orientable Wicks form of genus $n$. There are positive constants $L'_0$ and $c'_4$, depending only on $n,\delta, K,R_0,R_1$ such that for all $L\geq L_0$, there is a collection of length types $\mathcal{L}(W,n,L)$ satisfying:
    \begin{enumerate}[font=\normalfont]
        \item For all $\vec{l}\in\mathcal{L}(W,n,L)$, $||\vec{l}||=2\left\lceil\frac{L}{2}\right\rceil$, $\mathcal{L}(W,n,L)$ is $\text{Aut}(W)$-invariant, and:
        \[\lim_{L\to\infty}\frac{\mathcal{L}(W,n,L)}{L^{3n-4}}=\frac{1}{2^{3n-4}(3n-4)!}\]
    \end{enumerate}
    Let $W'$ be another maximal non-orientable Wicks form of genus $n$, and let $[\vec{l}]\in\mathcal{L}(W,n,L)/\text{Aut}(W)$, $[\vec{l}']\in\mathcal{L}(W',n,L)/\text{Aut}(W')$.
    \begin{enumerate}[font=\normalfont, resume]
        \item If $(W,[\vec{l}])\neq (W',[\vec{l}'])$, then:
        \[\mathcal{G}_W(100D\cdot\vec{l})\cap \mathcal{G}_{W'}(100D\cdot\vec{l}')\leq c_4\frac{e^{h\cdot 100DL}}{L^{3n}}\]
    \end{enumerate}
\end{lemma}

\begin{theorem}\label{thm: lowerboundhyp nonorientable conj class}
    Let $G$ be a hyperbolic group equipped with the word metric induced by some finite symmetric generating set $S$. Let $h$ denote the volume entropy of $G$. There exist positive constants $c'_5,L'_0$ depending only on $n$, $G$, and the choice of the generating set, such that:
    \begin{equation*}
        \forall L\geq L'_0, |\mathcal{C}^{-}_n(L)|\geq c'_5 L^{3n-4}e^{hL/2}
    \end{equation*}
\end{theorem}

  \bibliographystyle{alpha}
  \bibliography{bibliography}
\end{document}